\documentclass[12pt,twoside]{article}
\usepackage{amsmath}
\usepackage{latexsym}
\usepackage{amssymb}
\usepackage{a4wide}
\newtheorem{thm}{Theorem}[section]
\newtheorem{la}[thm]{Lemma}
\newtheorem{Defn}[thm]{Definition}
\newtheorem{Cnv}[thm]{Convention}
\newtheorem{Remark}[thm]{Remark}
\newtheorem{Note}[thm]{Note}
\newtheorem{prop}[thm]{Proposition}

\newtheorem{cor}[thm]{Corollary}
\newtheorem{Example}[thm]{Example}
\newtheorem{Examples}[thm]{Examples}
\newtheorem{Problems}[thm]{Problems}

\newtheorem{Problem}[thm]{Problem}
\newtheorem{Number}[thm]{\!\!}
\newenvironment{defn}{\begin{Defn}\rm}{\end{Defn}}
\newenvironment{cnv}{\begin{Cnv}\rm}{\end{Cnv}}

\newenvironment{example}{\begin{Example}\rm}{\end{Example}}
\newenvironment{examples}{\begin{Examples}\rm}{\end{Examples}}

\newenvironment{problem}{\begin{Problem}\rm}{\end{Problem}}
\newenvironment{rem}{\begin{Remark}\rm}{\end{Remark}}
\newenvironment{numba}{\begin{Number}\rm}{\end{Number}}
\newenvironment{proof}{{\noindent\bf Proof.}}%
                  {\nopagebreak\hspace*{\fill}$\Box$\medskip\medskip\par}   
\newcommand{\Punkt}{\nopagebreak\hspace*{\fill}$\Box$}
\newcommand{\wb}{\overline}

\newcommand{\wt}{\widetilde}

\newcommand{\tensor}{\otimes}

\newcommand{\n}{\rm}

\newcommand{\aeq}{\Leftrightarrow}

\newcommand{\mto}{\mapsto}
\newcommand{\emb}{\hookrightarrow}

\newcommand{\ve}{\varepsilon}

\newcommand{\isom}{\cong}

\newcommand{\N}{{\mathbb N}}

\newcommand{\R}{{\mathbb R}}
\newcommand{\F}{{\mathbb F}}
\newcommand{\bO}{{\mathbb O}}
\newcommand{\Q}{{\mathbb Q}}
\newcommand{\sph}{{\mathbb S}}
\newcommand{\Z}{{\mathbb Z}}
\newcommand{\C}{{\mathbb C}}

\newcommand{\K}{{\mathbb K}}

\newcommand{\cU}{{\cal U}}
\newcommand{\cO}{{\cal O}}

\newcommand{\cV}{{\cal V}}
\newcommand{\cF}{{\cal F}}

\newcommand{\dl}{{\displaystyle \lim_{\longrightarrow}}}

\newcommand{\End}{\mbox{\n End}}

\newcommand{\Iso}{{\mbox{\n Iso}}}

\newcommand{\wh}{\widehat}

\newcommand{\one}{\mbox{\rm \bf 1}}

\newcommand{\dsemi}{\mbox{$\times\!$\rule{.15 mm}{2.1 mm}$\,$}}

\newcommand{\take}{\backslash}
\newcommand{\sub}{\subseteq}

\newcommand{\GL}{\mbox{\rm GL}}

\newcommand{\im}{\mbox{\n im}}

\newcommand{\pr}{\mbox{\rm pr}}

\newcommand{\Aff}{\mbox{{\n Aff}}}
\newcommand{\aff}{\mbox{{\n aff}}}

\newcommand{\id}{\mbox{\n id}}

\newcommand{\cB}{{\cal B}}
\newcommand{\cA}{{\cal A}}
\newcommand{\cD}{{\cal D}}

\newcommand{\cK}{{\cal K}}
\newcommand{\cL}{{\cal L}}
\newcommand{\cE}{{\cal E}}

\newcommand{\bL}{{\mathbb L}}

\newcommand{\sSup}{\mbox{\n\footnotesize supp}}
\newcommand{\Supp}{\mbox{\n supp}}

\newcommand{\sbull}{{\scriptscriptstyle \bullet}}
\newcommand{\Diff}{{\mbox{{\rm Diff}}}}

\newcommand{\pl}{{\displaystyle\lim_{\longleftarrow}}}
\newcommand{\bx}{{\text{{\rm box}}}}
\newcommand{\tvs}{{\text{{\rm tvs}}}}
\newcommand{\lcx}{{\text{{\rm lcx}}}}
\newcommand{\spa}{\makebox[4 mm]{}}
\begin{document}
\thispagestyle{plain}
\pagenumbering{roman}
%
%
%
\renewcommand{\thefootnote}{\fnsymbol{footnote}}
$\;$\vspace*{-17mm}\\
\begin{center}
{\Large \bf Lie groups over non-discrete topological fields}\vspace{5mm}\\
{\bf Helge Gl\"{o}ckner\footnote{Darmstadt University of Technology,
Department of Mathematics~AG~5,
Schlossgartenstr.\,7,
\hspace*{6.3mm}64289~Darmstadt, Germany.
\,\,E-Mail: gloeckner@mathematik.tu-darmstadt.de}}\vspace{15mm}
\end{center}
\renewcommand{\thefootnote}{\arabic{footnote}}
\setcounter{footnote}{0}
%
%
\begin{abstract}\noindent
We generalize the classical
construction principles of infinite-dimensional
real (or complex) Lie groups
to the case of Lie groups over non-discrete
topological fields.
In particular,
we discuss linear Lie groups, mapping groups,
test function groups,
diffeomorphism groups,
and weak direct products
of Lie groups.
The specific tools of differential calculus
required for the Lie group constructions
are developed.
Notably, we establish
differentiability
properties of composition and evaluation,
as well as exponential laws
for function spaces.
We also present
techniques to deal with the subtle
differentiability and continuity properties of
non-linear mappings between spaces of test functions.
Most of the results are independent
of any specific properties
of the topological vector spaces involved;
in particular, we can deal with real
and complex Lie groups modeled
on non-locally convex spaces.\vspace{27mm}
\end{abstract}
{\small
{\bf Classification:}
22E65, 
22E67, 
58D05, 
26E30 (main); 
26E15, 
26E20, 
46A16, 
46S10, 
58C20\\[5mm] 
}
%
%
%
\vfill\pagebreak

\thispagestyle{plain}
\noindent
{\bf\large Contents:}\vspace{6mm}

\noindent
Introduction\dotfill\pageref{introdu}\vspace{3.7mm}

\noindent
Preliminaries:\vspace{1mm}

1. Differential calculus over topological
fields\dotfill\pageref{seccalc}\vspace{3mm}

\noindent
Linear Lie groups:\vspace{1mm}

2. Continuous inverse algebras and their unit
groups\dotfill\pageref{secALG}\vspace{3mm}

\noindent
Mapping groups and related constructions:\vspace{1mm}

3. Spaces of continuous mappings and mappings between
them\dotfill\pageref{secveccts}

4. Spaces of $C^r$-maps and mappings between
them\dotfill\pageref{secvecCr}

5. Mapping groups and mapping algebras\dotfill\pageref{secmapgps}

6. Mappings between direct sums\dotfill\pageref{secdirsum}

7. Weak direct products of Lie groups\dotfill\pageref{secweakprod}

8. Spaces of test functions and mappings between
them\dotfill\pageref{secmaptf}

9. Test function groups and algebras of test
functions\dotfill\pageref{sectfgps}

10. Differentiability of almost
local mappings\dotfill\pageref{secalmloc}\vspace{3mm}

\noindent
Diffeomorphism groups:\vspace{1mm}

11. Smoothness of evaluation and composition\dotfill\pageref{seccompo}

12. Basic exponential law for smooth mappings\dotfill\pageref{secexplaw}

13. Diffeomorphism groups of finite-dimensional,
paracompact\\
\hspace*{1.4cm}smooth manifolds
over local fields\dotfill\pageref{secdiffeos}

14. The diffeomorphism groups $\text{Diff}^r(M)$ and
$\text{Diff}^\infty(M)\wt{\,}$\dotfill\pageref{secdiff2}\vspace{3mm}

\noindent
Appendices:\vspace{1mm}

A. Proof of Proposition~\ref{propprop}\dotfill\pageref{appfun}

B. Proof of Proposition~\ref{crucial}\dotfill\pageref{appcruc}\vspace{-.1mm}

C. Proof of Proposition~\ref{comparbit}\dotfill\pageref{appcomparbit}

D. Smoothness vs.\ weak smoothness over local
fields\dotfill\pageref{grothend}

E. Towards a $p$-adic analogue of Boman's theorem\dotfill\pageref{appboman}

F. Spaces of sections in vector bundles and mappings between
them\dotfill\pageref{appsections}\vspace{4mm}

\noindent
Bibliography\dotfill\pageref{bibliog}\\[2mm]
Index\dotfill\pageref{nowindex}\\[2mm]
List of Symbols\dotfill\pageref{nowsymb}\vfill

\pagebreak
\thispagestyle{plain}
\pagenumbering{arabic}
\setcounter{page}{1}
\begin{center}\label{introdu}
{\large\bf Introduction}
\end{center}
Most of the known examples of
infinite-dimensional real or complex Lie groups
can be subsumed under (at least) one of the following main classes
of Lie groups:
1. linear Lie groups; 2. mapping groups;
3. diffeomorphism groups.
In the present article, we show that
the general construction principles
underlying
these classes of Lie groups
work just as well
beyond the real and complex cases.
Thus, we are able to discuss linear Lie groups
and groups of continuous Lie group-valued mappings
over arbitrary non-discrete topological fields;
groups of smooth Lie group-valued mappings on finite-dimensional
smooth manifolds over locally compact fields;
and diffeomorphism groups of paracompact
finite-dimensional smooth manifolds
over local fields (of arbitrary characteristic).
In the real and complex cases, it becomes possible
to construct Lie groups modeled on arbitrary
(not necessarily locally convex)
topological vector spaces.\\[3mm]
A fourth main class of infinite-dimensional
Lie groups are Lie groups
obtained from direct limit constructions,
in particular direct limits of finite-dimensional
Lie groups (see \cite{NRW}, \cite{NRW2}, \cite{DIR}, \cite{FUN}).
Direct limits of finite-dimensional Lie groups
over local fields have been discussed in~\cite{FUN}
(cf.\ also \cite{DIR}). Here,
we construct weak direct products of
(finite- or infinite-dimensional)
Lie groups over valued fields,
generalizing the discussion of weak direct products
of Lie groups
modeled on real or complex locally convex spaces from~\cite{MEA}.\\[3mm]
Our studies are based on a
concept of $C^k$-maps (and smooth maps)
between open subsets of topological vector spaces over
non-discrete topological fields
introduced in~\cite{Ber}, where more generally
an axiomatic approach to differential
calculus over arbitrary infinite fields
(and suitable commutative rings) is developed.
A map between
open subsets of {\em locally convex\/} real topological vector
spaces is of class~$C^k$ in the sense considered here
if and only if it is a $C^k$-map
in the sense of Michal-Bastiani,
{\em i.e.}, a $C^k_c$-map in the terminology of
Keller's monograph~\cite{Kel}
(see \cite{Ber}).
$C^k$-maps in the latter sense
have been used as the basis of differential
calculus and infinite-dimensional Lie theory
by many authors (see, {\em e.g.},
\cite{BCR}, \cite{RES}--\cite{HOL},
\cite{DRn}, \cite{INF}, \cite{FUN}--\cite{COM},
\cite{Ham}, \cite{KaR},
\cite{Mic}--
\cite{Mil}, \cite{Nee}--\cite{CUR}, \cite{PaS}, and
\cite{Wur}).
Furthermore, a map between open subsets
of complex locally convex spaces is complex $C^\infty$
if and only if it is complex analytic in the usual sense
(as in \cite{BaS}).\\[3mm]
Taken together,
the papers \cite{Ber}, \cite{IMP}
and the present work
develop, from first principles,
a comprehensive theory of differential calculus
and Lie groups over arbitrary non-discrete topological fields.
In \cite{Ber} already mentioned,
an exposition of differential
calculus over topological fields
and the corresponding basic theory of manifolds
and Lie groups is given;
this article is directed to a broad audience
including readers without prior
knowledge of differential calculus
over topological fields.
In \cite{IMP} (needed here only for
the discussion of diffeomorphism groups),
implicit (and inverse) function theorems
for $C^k$-maps over complete valued fields
are established.
The present article, then,
provides concrete examples
of Lie groups over topological fields, and develops
the specific aspects
of differential calculus required for this purpose.
Further papers related to the ``General
Differential Calculus'' over topological fields
are \cite{Be2}, \cite{BaN}, \cite{BaN2},
\cite{NCX}, \cite{ANA}, \cite{FUN} and \cite{HOE}. 
In \cite{ANA}, it was shown that every finite-dimensional
smooth $p$-adic Lie group is a $p$-adic analytic
Lie group in the usual sense. Important aspects
of differential geometry
over topological fields have been worked out
in~\cite{Be2}.
Jordan theoretic applications
are described in \cite{BaN} and \cite{BaN2}.\\[4mm]
The article commences
with a brief introduction to differential calculus
over non-discrete topological fields (Section~\ref{seccalc}),
and ends with appendices covering material
which is best taken on faith on a first reading,
and whose presentation within the text
would have distracted from the main line of thought.
Apart from these sections, the main body of the
text is divided into three parts,
devoted to the three main classes and construction
principles of infinite-dimensional
Lie groups described above:\vspace{2mm}
\begin{center}
{\bf I. Linear Lie groups}
\end{center}
Paradigms of real or complex Lie groups are {\em linear\/} Lie groups,
{\em i.e.}, unit groups of unital Banach algebras
(or other well-behaved topological algebras)
and their Lie subgroups
(see \cite{ALG}, \cite{dlH}, \cite[Ch.\,5]{HaM},
\cite{MaR}).
We begin our studies with
a discussion of linear Lie groups
over topological fields (Section~\ref{secALG}),
as this only requires
a minimum of technical machinery.
If $\K$ is a non-discrete
topological field,
a good class of topological $\K$-algebras
to look at are the {\em continuous inverse algebras\/}
(or CIAs),
{\em viz.}\ unital associative topological
$\K$-algebras~$A$
such that the group of units~$A^\times$
is open in~$A$ and the inversion map
$\iota\!: A^\times \to A$, $a\mto a^{-1}$
is continuous.
We describe examples of CIAs and construction
principles for CIAs over
arbitrary non-discrete topological fields.
Since the unit group~$A^\times$
is a $\K$-Lie group for
any continuous inverse algebra $A$ (Proposition~\ref{ciasmooth}),
we thus always have a certain supply of
$\K$-Lie groups, for any~$\K$ (beyond the trivial
examples, the additive groups of topological $\K$-vector spaces).
Algebras of continuous or differentiable
maps on compact
topological spaces or compact
manifolds, with values in a CIA,
are again CIAs
(Proposition~\ref{mapcia}).
For further typical examples of CIAs in the real or complex cases,
see \cite{ALG},
\cite{HOL}, \cite[1.15]{Gra}, \cite{YMO}.\vspace{2mm}
\begin{center}
{\bf II. Mapping groups and related constructions}
\end{center}
The second widely studied class
of infinite-dimensional real (or complex)
Lie groups are the mapping groups.
Typical examples
are the ``loop groups''
$C(\sph^1,G)$ and $C^\infty(\sph^1,G)$,
where $G$ is a finite-dimensional real
(or complex)
Lie group~\cite{PaS}.
More generally, let $G$ be a real or complex Lie group
modeled on a locally convex space,
$r\in \N_0\cup\{\infty\}$, and
$M$ be a finite-dimensional
smooth manifold (or topological space
if $r=0$).
Among the types of mapping groups
encountered in the literature,
we mention: the groups $C^r(M,G)$
of $G$-valued $C^r$-maps,
for compact~$M$;
the groups
$C^r_K(M,G)$ of $G$-valued
$C^r$-maps supported in a compact set
$K\sub M$;
and, for $\sigma$-compact manifolds~$M$,
the ``test function groups''
$C^r_c(M,G):=\bigcup_K\, C^r_K(M,G)$
of compactly supported $G$-valued $C^r$-maps
(see \cite{Alb}, \cite{GCX}, \cite{INF},
\cite{KaM}, \cite{Mil}, \cite{NRW2},
\cite{CEN}, \cite{CUR}).
\\[3mm]
In the second main part
of this article, we
construct Lie group
structures on analogous
mapping groups in the case
of Lie groups over topological fields.
The results include:\\[3mm]
{\bf Groups of continuous mappings.\/}
{\em If $\K$ is a non-discrete topological field,
$G$ a $\K$-Lie group, $X$ a topological space,
and $K\sub X$ a compact subset, then the
group
\[
C_K(X,G)\, :=\, \{\gamma\in C(X,G)\!: \gamma|_{X\setminus K}=1\}
\]
of continuous $G$-valued mappings supported in~$K$
can be made a $\K$-Lie group modeled on
$C_K(X,L(G))$, in a natural way.
In particular, $C(K,G)=C_K(K,G)$ is a $\K$-Lie group, for every
compact
topological space~$K$ and any $\K$-Lie group~$G$.}\\[3mm]
{\bf Groups of differentiable mappings.\/}
{\em Let $\F$ be a locally compact, non-discrete topological field,
$\K$ be a topological extension field of~$\F$,
$G$ be a $\K$-Lie group,
$r\in \N_0\cup\{\infty\}$,
$M$ be a finite-dimensional $C^r_\F$-manifold,
$K\sub M$ a compact subset, and
\[
C_K^r(M,G)\, :=\, \{\gamma\in C^r(M,G)\!: \gamma|_{M\setminus K}=1\}
\]
be the group of $G$-valued $C^r_\F$-functions on~$M$
supported in~$K$.
Then $C^r_K(M,G)$
is a $\K$-Lie group modeled on $C^r_K(M,L(G))$, in a natural way.
If $M$ is paracompact and
the topology on~$\K$ arises from an absolute value,
then also the group
$C^r_c(M,G):=\bigcup_K\, C^r_K(M,G)$
of $G$-valued test functions of class~$C^r_\F$
is a $\K$-Lie group, modeled on $C^r_c(M,L(G))$.}\\[3mm]
(See Sections~\ref{secmapgps}
and~\ref{sectfgps}).
Typically, we might choose
$\K:=\F$ here, or $\F:=\R$, $\K:=\C$.
As the basis for the construction of Lie group
structures on mapping groups, we study
continuity and differentiability properties
of mappings between function spaces.
The three cases of interest
(mappings between spaces of $C_K$-functions,
$C_K^r$-functions, and $C^r_c$-functions,
respectively)
are discussed in turn in
Sections~\ref{secveccts}, \ref{secvecCr}, resp.,
\ref{secmaptf} and \ref{secalmloc}.
For simplicity, let
us assume that $\K:=\F$ now, for the remainder of
the introduction.
The results obtained
subsume, for example, that
the mappings
\[
\begin{array}{rcl}
C^\infty_K(M,g) : & C^\infty_K(M,U)\to C^\infty_K(M,F),&
\gamma\mto g\circ \gamma\quad \mbox{\hspace*{1.4cm} and}\\
f_* : & C^\infty_K(M,U)\to C^\infty_K(M,F),&
\gamma\mto f\circ (\id_M,\gamma)
\end{array}
\]
are of class $C^\infty_\K$,
for any $C^\infty_\K$-maps
$g\!: U\to F$ and $f\!: M\times U\to F$
such that $g(0)=0$ and $f|_{(M\setminus K)\times\{0\}}=0$,
where $\K$ is a locally compact,
non-discrete topological field,
$E$ and $F$ are
topological $\K$-vector spaces,
$U\sub E$ an open zero-neighbourhood,
$M$ a finite-dimensional
$C^\infty_\K$-manifold, and $K\sub M$ a compact subset.
For paracompact~$M$, analogous conclusions are
valid for
$C^\infty_c(M,g)$ and $f_*\!:C^\infty_c(M,U)\to
C^\infty_c(M,F)$.
More generally,
results of the preceding type
are established
for mappings between spaces of sections
in vector bundles,
whose fibres are arbitrary topological vector spaces
(Appendix~\ref{appsections}).
%
For the real locally convex case,
the reader may compare
Michor~\cite{Mic},
in particular his ``$\Omega$-Lemma'' \cite[Thm.\,8.7]{Mic}
for finite-dimensional real vector bundles;
\cite{GCX} (for maps between spaces of
test functions),
and~\cite{SEC}.
It is a peculiarity
of differential calculus over general topological fields
that, when we are trying to prove
differentiability properties of $f_*$ (or related results),
{\em parameter-dependent\/}
variants invariably pop up in the natural inductive
arguments, even when we are only interested
in the case of $f_*$ (not involving parameters).
On the one hand, this makes
the proofs more complicated;
but, on the other hand,
we are rewarded with stronger,
parameter-dependent versions of the basic
results
(like our ``$\Omega$-Lemma with Parameters'', Theorem~\ref{OmegaP}),
which are novel
even in the real locally convex case.\\[3mm]
{\em Topologies on spaces of test functions.}
If $\K\not=\C$ is locally compact, $E$ a topological
$\K$-vector space and $M$ a $\sigma$-compact,
finite-dimensional
$C^r_\K$-manifold,
we equip $C^r_c(M,E)=\dl\,C^r_K(M,E)$\vspace{-.8mm}
with the topology making it the direct limit of its subspaces
$C^r_K(M,E)$ in the category
of topological $\K$-vector spaces.
Although little is known on
direct limits of general topological vector spaces
(in contrast to the
real or complex locally convex case,
which has been studied extensively),
we get a perfectly firm grip
on the topology of $C^r_c(M,E)$
by showing that the linear map
\begin{equation}\label{getgrip}
\rho\!: C^r_c(M,E)\to \bigoplus_{i\in I} C^r(U_i,E),\quad
\rho(\gamma):=(\gamma|_{U_i})_{i\in I}
\end{equation}
is a topological embedding onto
a closed vector subspace,
for any locally finite cover
$(U_i)_{i\in I}$ of~$M$
by relatively compact, open subsets~$U_i$.
Here, the direct sum is equipped with the
box topology,
which is extremely simple to work with.\\[3mm]
If $\K=\C$, or if $M$ is merely paracompact,
then the topology making
$C^r_c(M,E)$ the direct limit topological vector space
$\dl\, C^r_K(M,E)$\vspace{-.8mm} is too elusive
to be useful for us.
Instead of excluding these cases altogether
from our considerations,
we simply replace the direct limit topology
on $C^r_c(M,E)$ with the topology induced by~$\rho$,
which turns out to be independent of the choice
of open cover $(U_i)_{i\in I}$.
This enables us to carry out most of our constructions
also in the complex case, and also for paracompact manifolds
(see Proposition~\ref{comparetop}
and Remark~\ref{justffy}
for further discussions of the box topology
and explanations why we prefer to use it).
Note that, for a non-locally convex complex
topological vector space~$E$,
the space $C^r_c(M,E)$ of compactly supported $E$-valued $C^r_\C$-maps
need not be reduced to the locally constant functions.
For example,
there are non-zero compactly supported $C^\infty_\C$-maps
$\C\to E$, for
suitable~$E$ (see \cite{NCX}).
As we do not have cut-off functions
(nor partitions of unity)
available in the complex case,
we have to proceed with particular care.
In Section~\ref{secalmloc},
their use cannot be avoided any longer,
and the complex case has to be excluded then.\\[3mm]
{\em Mappings between direct sums.}
The embedding $\rho$ from (\ref{getgrip})
allows us to reduce the study
of continuity and differentiability properties
of mappings between spaces of test functions (or
compactly supported sections in vector bundles)
almost entirely to the
study of differentiability properties
of mappings of the form
\[
\oplus_{i \in I}\,f_i : \;\bigoplus_{i\in I} U_i\to \bigoplus_{i\in I} F_i,
\quad
(x_i)_{i\in I}\mto (f_i(x_i))_{i\in I}
\]
on open boxes $\bigoplus_{i\in I}U_i$ in
direct sums
$\bigoplus_{i\in I}E_i$ of topological $\K$-vector spaces.
In Section~\ref{secdirsum}, for arbitrary valued
fields~$\K$, we show
that a mapping $\oplus_{i\in I}f_i$ as before
is $C^k_\K$ provided each $f_i$ is $C^k_\K$.
Although the proof of the special case
where $I$ is countable and we are dealing
with real or complex locally convex topological
vector spaces is almost trivial (see \cite{MEA}),
the proof of the general case requires a
substantial amount of work.
In order to control simultaneously, for each $i\in I$,
the dependence on
the parameter $t$ of the extended difference
quotient maps
\[
f_i^{[1]}\!: U_i\times E_i\times \K\supseteq U_i^{[1]}\to F_i\, ,
\quad
(x_i,y_i,t)\mto f_i^{[1]}(x_i,y_i,t)
\]
and, more generally,
analogous parameters in the mappings $f^{[j]}_i$, where $j\in \N$, $j\leq k$,
which are encountered in the inductive arguments,
we are forced to investigate in some depth the symmetry
properties
of the higher difference quotient maps
\[
f^{[j]}_i\! : U^{[j]}_i\to F_i\,,\]
which are rather complicated
functions depending on $\,2^{j+1}\!-\!1$ variables
(Proposition~\ref{mapsdirsums}).
In the case of locally compact fields,
the dependence on parameters can
be controlled more easily
by means of compactness
arguments, which in fact permit us
to formulate stronger results,
involving additional parameters (Proposition~\ref{sumspara}).
As a first straightforward application,
mappings between direct sums are used to
define a Lie group structure on (countable
or uncountable) weak direct products
$\prod_{i\in I}^*G_i$ of Lie groups $G_i$
over a valued field~$\K$, based on the
box topology on direct sums (Section~\ref{secweakprod}).
As we shall see later,
such groups are encountered quite frequently in the case
of ultrametric fields; for example,
they shall play an important role in our discussion
of diffeomorphism groups over local fields.
Weak direct products
of Lie groups
modeled on real or complex locally convex spaces,
based on locally convex direct sums,
have first been considered in \cite{MEA}.\\[3mm]
In the real finite-dimensional case,
embeddings in locally convex direct sums
are implicit in \cite[\S4.7]{Mic}
in connection with descriptions
of the ``$\cD$-topology'' on mapping spaces,
which are used there for similar purposes.
The usefulness of
embeddings into real and complex locally convex direct sums
for the study of mappings between spaces of test
functions (and compactly supported sections)
has been pointed out explicitly in~\cite{SEC},
\cite{DIF}.
The arguments in~\cite{Mic}
(based on jet bundles) are restricted
to vector bundles over finite-dimensional bases.
We approach function spaces
and mappings between them
in a more direct way.
This allows us, for instance,
to prove smoothness
of the pushforward
\[
f_*\!: C^\infty(M,E)\to C^\infty(M,F)
\]
in the case of a globally defined
$C^\infty_\K$-map
$f\!: M\times E\to F$
in utmost generality,
namely, for $\K$ an arbitrary topological field,
$M$ a $C^\infty_\K$-manifold
modeled on an arbitrary
topological $\K$-vector space, and $E$, $F$ arbitrary topological
$\K$-vector spaces (Proposition~\ref{globcruc}).\\[3mm]
{\em Differentiability properties of
almost local mappings.}
In order to motivate the most general
results we have to offer which exploit embeddings into direct
sums (presented in Section~\ref{secalmloc}), 
we recall from \cite{DIS}
that the self-map
\[
f\!: C^\infty_c(\R,\R)\to C^\infty_c(\R,\R),\quad 
f(\gamma):=\gamma\circ \gamma\,-\,\gamma(0)
\]
of the space $C^\infty_c(\R,\R)$ of real-valued test functions
on the real line
is discontinuous (and hence not smooth),
although the restriction of $f$ to
$C^\infty_K(\R,\R)$ is smooth, for every
compact subset $K\sub \R$.
In the real and complex case,
and, more generally, in the case of locally compact
topological fields~$\K$,
this poses the question whether
it is possible to specify
simple and easily verified
{\em additional properties\/} ensuring that
a mapping
\[
f\!: C^r_c(M,E)\to C^s_c(N,F)
\]
between spaces of vector-valued test functions
is indeed $C^k_\K$, provided
the restriction of $f$ to $C^r_K(M,E)$
is $C^k_\K$ for each
compact subset $K\sub M$
(and likewise for mappings between
spaces of compactly
supported sections, or open subsets thereof).\\[2mm]
Generalizing the real locally
convex case (see \cite{INF}, \cite{SEC}, \cite{DIF}),
we show that also in the case of general locally compact fields
$\K\not=\C$,
the requirement that $f$ be {\em almost local\/}
is a suitable additional property on~$f$
(Theorem~\ref{smoothy}),
meaning that there exist locally finite, relatively compact open covers
$(U_i)_{i\in I}$ of~$M$ and $(V_i)_{i\in I}$ of~$N$
such that  $f(\gamma)|_{V_i}$ only depends on
$\gamma|_{U_i}$, for any~$i$.
The class of almost local maps
includes most mappings of interest.
For example, in the case $M=N$, every
pushforward of sections
associated with a fibre-preserving bundle map
is almost local.
Furthermore,
all mappings encountered
in the construction of Lie group
structures on diffeomorphism groups
of $\sigma$-compact finite-dimensional
real $C^\infty$-manifolds
are (locally) almost local
(see \cite{DRn} and \cite{DIF},
where diffeomorphism groups are discussed
along lines
independent of the earlier work \cite{Mic}).\\[3mm]
For a highly developed theory of
spaces (and manifolds) of mappings and mappings
between these in the ``convenient setting of analysis''
(based on Mackey complete
real or complex locally convex spaces),
which is inequivalent to the setting
of analysis we
are working in here, see
\cite{FaK}, \cite{KaM}
and further works by the same authors.
\begin{center}
{\bf III. Diffeomorphism groups}
\end{center}
The third main part of the article
is devoted to diffeomorphism groups
of finite-dimensional manifolds over local fields,
and related material.
We begin with a discussion of
continuity and differentiability
properties of evaluation
and composition of maps in the context of locally
compact fields~$\K$ (Section~\ref{seccompo}).
Among variants and related results,
we show in particular
that the evaluation map
\[
\ve\!: C^r(M,E)\times M\to E, \quad \ve(\gamma,x):=\gamma(x)
\]
is of class $C^r_\K$, for every
finite-dimensional $C^r_\K$-manifold~$M$
and topological $\K$-vector space~$E$
(Proposition~\ref{evalCk}),
and
that the composition map
\[
\Gamma\!: C^{r+k}(U,E)\times C^r_K(M,U)\to C^r(M,E),\quad
\Gamma(\gamma,\eta):=\gamma\circ \eta
\]
is of class $C^k_\K$, for any $r,k\in \N_0\cup\{\infty\}$,
finite-dimensional $C^r_\K$-manifold~$M$,
topological $\K$-vector space~$E$,
compact subset $K\sub M$, and
open subset~$U$ of a finite-dimensional
$\K$-vector space~$F$ (Proposition~\ref{compcomp}).
We then turn to the exponential law
for smooth mappings (Section~\ref{secexplaw}).
Given any topological field~$\K$,
topological $\K$-vector space~$E$,
$r,k\in \N_0\cup\{\infty\}$
and arbitrary $C^{r+k}_\K$-manifolds $M$ and $N$, we show that
\[
f^\vee\!: M\to C^r(N,E),\quad f^\vee(x)(y):=f(x,y)
\]
is of class $C^k_\K$ for all $f\in C^{r+k}(M\times N,E)$,
and that the mapping
\begin{equation}\label{firstPhi}
\Phi\!: C^{r+k}(M\times N,E)\to C^k(M,C^r(N,E)),\quad \Phi(f):=
f^\vee
\end{equation}
is a continuous linear injection
(Lemma~\ref{halfcartesian}).
If $\K$ is locally compact and $N$ is finite-dimensional
(but $M$ arbitrary),
we show that $\Phi$ is an isomorphism of
topological vector spaces, in the case
$r=k=\infty$
(Proposition~\ref{nownow}).\footnote{Analogous results for the case
where $k=r=\infty$, $E$ is locally convex and both
$M$ and $N$ are open subsets of real locally convex spaces
(instead of manifolds)
have been obtained earlier in~\cite{Bil},
along with interesting additional information.}
Similar (slightly weaker) conclusions hold
if both $M$ and $N$ are modeled on metrizable
topological vector spaces and
$\K$ is $\R$ or an ultrametric field
(Proposition~\ref{expmetriz}).
To deduce the surjectivity of~$\Phi$
in the metrizable case,
we make use of techniques of convenient differential
calculus (already mentioned),
suitably adapted to non-locally convex or
ultrametric analysis
by means of preparatory results
provided in~\cite{Ber}.
We can only broach on the subject of
``ultrametric convenient differential calculus''
here, and have to confine ourselves
to what is actually needed for the concrete purpose.
A further application of the exponential laws
is given in Appendix~\ref{appboman}.
Combining the latter and an ultrametric analogue
of Grothendieck's Theorem (relating smoothness
and weak smoothness of suitable maps)
provided in Appendix~\ref{grothend},
it is shown there that Boman's theorem
will hold for mappings $f\!: E\supseteq U\to F$ from an open subset
of a metrizable topological vector space~$E$
over a local field~$\K$
to a Mackey complete
locally convex topological $\K$-vector space~$F$
provided Boman's theorem holds for all mappings $f\!: \K^2\to \K$.
Recall that Boman's classical theorem~\cite[Thm.\,1]{Bom}
asserts that a map $f\!: \R^n\to\R$
is smooth if and only if $f\circ \gamma\!: \R\to\R$
is smooth for each smooth curve $\gamma\!: \R\to\R^n$.
Whether Boman's theorem transfers to maps
$f\!: \K^2\to \K$ is still unknown.\\[3mm]
In the final Sections~\ref{secdiffeos} and \ref{secdiff2},
which can be read independently,
we describe two\linebreak
approaches to diffeomorphism groups
of finite-dimensional
smooth manifolds over local fields.
The first approach (Section~\ref{secdiffeos})
produces a Lie group structure on $\Diff^\infty(M)$,
for every paracompact, finite-dimensional smooth
manifold~$M$ over a local field~$\K$.
The second approach (Section~\ref{secdiff2})
is restricted to $\sigma$-compact~$M$.
It produces {\em two\/}
Lie group structures on $\Diff^\infty(M)$
(one of which coincides with the one
constructed in Section~\ref{secdiffeos}).
Both approaches make use of many of the results and techniques
prepared before (and hence also
illustrate the usefulness and typical applications
of these results).\\[3mm]
{\bf First approach\/} (Section~\ref{secdiffeos}).
Let $M$ be a finite-dimensional, paracompact $C^\infty_\K$-manifold
over a local field~$\K$, and
$\Diff^\infty(M)$ be the group of all $C^\infty_\K$-diffeomorphisms
of~$M$.
Our first construction
of a Lie group structure on $\Diff^\infty(M)$
relies on the fact $M$ is
a disjoint union of a family $(B_i)_{i\in I}$
of open and compact balls $B_i\sub M$
(i.e., subsets $B_i\sub M$ which are
$C^\infty_\K$-diffeomorphic
to balls in $\K^d$ with respect to the supremum-norm).
Motivated by this decomposition,
we first turn the diffeomorphism group
$\Diff^\infty(B)$ of a ball $B\sub \K^d$ (where $d\in \N_0$)
into a Lie group; this is quite easy,
because $\Diff^\infty(B)$ is a mere open subset
of $C^\infty(B,\K^d)$.
Next, we form the weak direct product
of Lie groups
\[
{\textstyle {\prod^*}_{\!\!\! i\in I} \;\Diff^\infty(B_i)}
\]
modeled on
$\bigoplus_{i\in I} C^\infty(B_i,TB_i)=
C^\infty_c(M,TM)$.
This weak direct product
can be identified with a subgroup
of the group $\Diff^\infty_c(M)$ of ``compactly
supported''
diffeomorphisms.
It then only remains to equip
$\Diff^\infty(M)$
with a smooth $\K$-manifold
structure making it a Lie group
with $\prod^*_{i\in I} \,\Diff^\infty(B_i)$
as an open subgroup.
We remark that smoothness
of the inversion map
$\Diff^\infty(B)\to \Diff^\infty(B)$, $\gamma\mto \gamma^{-1}$
is a simple consequence of the exponential
laws established here and
a version of the implicit function theorem
(the ``Inverse Function Theorem with Parameters'')
for mappings from metrizable topological vector spaces
to Banach spaces~\cite{IMP}.\\[3mm]
{\bf Second approach\/} (Section~\ref{secdiff2}).
For our second approach to diffeomorphism
groups, we assume that $M$ is $\sigma$-compact
and of positive dimension over~$\K$.
In this case, $M$ is $C^\infty_\K$-diffeomorphic
to an open subset $U$ of its modeling space~$\K^d$
(Lemma~\ref{onlyopen}\,(a); cf.\ \cite{Lu4}),
making it quite easy to deal with~$M$.
Given $r\in \N\cup\{\infty\}$,
we consider the monoid $\End_c^r(U)$
of all $C^r_\K$-maps $U\to U$ which coincide with $\id_U$
outside some compact set. We show that $\gamma\mto\gamma-\id_U$
identifies $\End^r_c(U)$
with an open subset
of $C^r_c(U,\K^d)$. Thus $\End^r_c(U)$ is a $C^\infty_\K$-manifold
with a global chart.
We show that $\Diff^r_c(U)=\End_c^r(U)^\times$
is open in $\End^r_c(U)$,
and we show that, for each $k\in \N_0\cup\{\infty\}$,
the composition map
\[
\Diff^{r+k}_c(U)\times \Diff^r_c(U)\to \Diff^r_c(U)
\]
and the inversion map
$\Diff^{r+k}_c(U)\to \Diff^r_c(U)$
are $C^k_\K$.
This enables us to
turn $\Diff^\infty(M)$ into a Lie group with $\Diff_c^\infty(M)\isom
\Diff^\infty_c(U)$ as an open subgroup,
modeled on the space
$C^\infty_c(M,TM)$ of compactly supported
smooth vector fields on~$M$,
equipped with its natural LF vector topology.\footnote{The
Lie group structure so obtained coincides with the one
from Section~\ref{secdiffeos}.}
But it also enables us to define a second
Lie group structure on $\Diff^\infty(M)$
(which we then denote by $\Diff^\infty(M)\wt{\,}$\,).
It is modeled on the same vector space $C^\infty_c(M,TM)$,
equipped however with the (in general properly) coarser
vector topology making this space the projective limit
\[
\bigcap_{k\in \N_0}C^k_c(M,TM)=\pl_{k\in \N_0}
\, C^k_c(M,TM)\,.
\]
Apparently, the definition of this second Lie group structure
is close in spirit to the ILB-approach to diffeomorphism groups
in the works of Omori \cite{Omm}, \cite{Omo}. An analogous construction
for diffeomorphism groups of $\sigma$-compact
real manifolds
had been proposed in \cite{MiP} (and was fully worked
out in \cite{DIF}).
As in the real case,
we can also give $\Diff^r(M)$ a smooth manifold
structure for each finite~$r$,
with $\Diff^r_c(M)\isom \Diff^r_c(U)$
as an open subgroup, such that $\Diff^r(M)$
becomes a topological group and
all right translations are smooth.\\[3mm]
We remark that
certain groups
$\Diff(t,M)$, $G(t,M)$,
and $GC(t,M)$ of diffeomorphisms of class of smoothness
$C(t)$ for manifolds
over local fields of
characteristic zero have
already been discussed
in~\cite{Lu1}, \cite{Lu3}, \cite{Lu10}
and further works of
S.\,V. Ludkovsky,
where they are considered mainly
as manifolds and topological groups
(rather than Lie groups).
He discusses irreducible
representations of these groups (\cite{Lu3})
and measures on such groups (or $M$)
which are quasi-invariant with respect to
dense subgroups~\cite{Lu1}.
Ludkovsky's approach to differential calculus
(which is necessarily restricted to
local fields~$\K$ of characteristic zero,
and differs from the one we use),
is described in \cite[I, \S2.3]{Lu8}
(for maps between ultrametric Banach spaces)
and extended to the case of
locally convex topological $\K$-vector spaces
in \cite[II, Rem.\,4.4]{Lu8}
and \cite{Lu10},
where diffeomorphism groups are discussed
in further detail and generality.
Ludkovsky
also turns groups
$\Diff^{\,t}_{\beta,\gamma}(M)$
of certain diffeomorphisms of real Banach manifolds
(subject to H\"{o}lder-type conditions)
into topological groups \cite[Thm.\,3.1]{Lu6},
and discusses
irreducible representations
and quasi-invariant measures
for such groups~\cite{Lu6}, \cite{Lu9}.
Note that
the ``non-archimedian loop groups''
discussed in \cite{Lu8}
are not mapping groups
in the sense considered in the present article,
but something
different.\\[3mm]
Finally, let us mention that
weak direct products of Lie groups are also useful
to obtain information on diffeomorphism groups
of real manifolds (although they are not simply 
open subgroups here, as in the ultrametric case).
In \cite{COM}, weak direct products
are used to show that the Lie group
$\Diff_c^\infty(M)$
of compactly supported diffeomorphisms
of a $\sigma$-compact, finite-dimensional smooth manifold~$M$
is the direct limit $\dl\,\Diff_K^\infty(M)$\vspace{-.8mm}
both in the category of Lie groups modeled
on real locally convex spaces and
in the category of topological groups
(where $K$ ranges through the set of compact subsets
of~$M$, and $\Diff^\infty_K(M):=\{\gamma\in \Diff^\infty(M)\!:
(\forall x\in M\setminus K) \; \gamma(x)=x\}$).
This is remarkable, because,
as a consequence of results from~\cite{DIS},
in general the Lie group $\Diff_c^\infty(M)$
neither is the direct limit
$\dl\,\Diff_K^\infty(M)$\vspace{-.8mm}
in the category of topological spaces,
nor in the category of smooth manifolds
(\cite{COM}; cf.\ also \cite{TSH}).
Analogous results can be obtained for the test functions
groups $C^\infty_c(M,G)$, for $G$ a finite-dimensional
real Lie group~\cite{COM}.

\begin{center}
{\bf Concluding remarks and guidance for the reader}
\end{center}
Readers who wish to get quickly to the main results
can skip part of the material.
For example, since all of the vector bundles
required for the discussion of diffeomorphism groups over local fields
are trivial bundles,
only very few of our results on
spaces of sections in vector bundles
(discussed in Appendix~\ref{appsections})
will actually be used,
and these are easy to take on faith
(cf.\ Remark~\ref{remnosecs}).
Proposition~\ref{globcruc}
(concerning pushforwards $f_*$ for globally defined~$f$)
is only needed for the discussion of spaces
of sections in vector bundles, while
its more complicated variants (Propositions~\ref{pushforw2}
and \ref{crucial}) are vital for the Lie group
constructions.
Nonetheless, the author recommends
to study the proof of
the simpler Proposition~\ref{globcruc}
first, and to leave the proofs
of Propositions~\ref{pushforw2}
and \ref{crucial} for a second reading.
Section~\ref{secalmloc}
is only needed for our second approach
to diffeomorphism groups (Section~\ref{secdiff2}),
but not for the first approach (Section~\ref{secdiffeos}).
The general case of Proposition~\ref{comparbit}
(proved in Appendix~\ref{appcomparbit})
is not used elsewhere,
and only part of
Section~\ref{secexplaw} (concerning the exponential law)
is needed for the discussion of diffeomorphism
groups: Lemma~\ref{halfcartesian}\,(a),
and Lemma~\ref{halfcartesian}\,(b) for $k=0$ suffice.\\[3mm]
It is clear, however, that
it would be inefficient not to include
such closely related results, when this can be done
without much additional effort.
Besides their obvious potential
for applications,
the additional results also serve
to put the Lie theoretic
constructions in a larger perspective,
and thus foster their understanding.\\[3mm]
Let us remark in closing
that it was necessary to refrain
from developing the surrounding theory
up to the possible limits of generality,
in order not to be carried away too far
from the subject matter of Lie group
constructions, to increase the readability,
and to avoid the discussions from getting even
more technical.\\[3mm]
We mainly think of
two possible generalizations.
In the real locally convex case,
a more refined discussion
of the maps $C^r_c(M,f)\!: C^r_c(M,U)\to C^r_c(M,F)$
between open subsets of spaces of compactly supported
sections is possible~\cite{SEC};
in this case, $C^r_c(M,f)$ is $C^k$ provided,
in local coordinates, $f$ is $C^k$ along the fibre,
with fibre derivatives jointly~$C^r$.
An analogous condition, based on iterated
partial difference quotient maps,
should be sufficient to ensure
that $C^r_c(M,f)$ be $C^k$, in the general case
of bundles
of topological vector spaces over
finite-dimensional paracompact $C^r$-manifolds
over locally compact fields.
This would substantially generalize
our version of the $\Omega$-Lemma (which, however,
already incorporates the cases
of main relevance),
but would inflict complicated
technical arguments
on us, which are irrelevant for the Lie group
constructions.\\[3mm]
The second possible generalization
concerns the exponential law.
If $k=r=\infty$, the map $\Phi$ from (\ref{firstPhi})
should always be a topological embedding
(cf.\ \cite{Bil}).
Furthermore, for general $r$ and $k$,
a more detailed
analysis of the problem should reveal that
$\Phi$
can be written
as a composition
\[
C^{r+k}(M\times N,E)\to C^{k,r}(M\times N,E)\to
C^k(M,C^r(N,E))
\]
of continuous linear injections
for a suitably defined space $C^{k,r}(M\times N,E)$
of $E$-valued $C^{k,r}$-maps on $M\times N$.
Here, the first mapping
is the inclusion map.
The second map, $C^{k,r}(M\times N,E)\ni f\mto f^\vee \in
C^k(M,C^r(N,E))$, should always be
a topological embedding.
Again, the author feels
that the immense additional technical effort
would not be justified in the present context.
The problems may be analyzed elsewhere.
\section{Differential calculus over topological fields}\label{seccalc}
It is possible to define $C^k$-mappings
and smooth mappings
once a topologized ring and a so-called $C^0$-concept is given,
satisfying suitable axioms
(see \cite{Ber}).
In the present paper, we exclusively
consider the special case where the given
topologized ring is a non-discrete
topological field~$\K$
(Hausdorff, as all topological spaces we consider),
where $C^0$-maps are defined as
continuous maps between open subsets of topological
vector spaces over~$\K$, and where the product
topology is used on products of topological vector spaces.
In this section, we briefly describe the resulting setting
of differential calculus.
\begin{numba}\label{convents} {\bf Conventions.}
All topological spaces occurring in this paper
are assumed Hausdorff.
All topological fields are supposed to be non-discrete.
A field $\K$, together with an absolute value
$|.|\!:\K\to [0,\infty[\,$
giving rise to a non-discrete topology on~$\K$
will be called a {\em valued field}.
An {\em ultrametric field\/}
is a valued field $(\K,|.|)$ whose
absolute value is ultrametric,
i.e.,
\[
|x+y|\,\leq\, \max\{|x|,|y|\}\qquad\mbox{for all $x,y\in \K$.}
\]
If $(\K,|.|)$ is an ultrametric field,
then $\bO:=\{x\in \K\!: |x|\leq 1\}$ is
a subring of~$\K$, called the {\em valuation ring}.
The valuation ring is an open and closed subset of~$\K$.
Totally disconnected, locally compact
topological fields will be referred to as {\em local fields}.
It is well known that every
local field admits an ultrametric absolute value
defining its topology, and can therefore be considered
as a complete ultrametric field.
The valuation ring of a local field
is open and compact.
It is also known that
every locally compact topological field
is either isomorphic to $\R$, $\C$, or a local field~\cite{Wei}.
A complete\footnote{The requirement is
that $(\K,d)$ be a complete metric space,
where $d\!: \K\times \K\to [0,\infty[$,
$d(x,y):=|x-y|_\K$.}
valued field $(\K,|.|_\K)$ is either ultrametric
or isomorphic as a valued field
to $(\R,|.|^\alpha)$ or $(\C,|.|^\alpha)$
for some $\alpha\in \;]0,1]$,
where $|.|$ is the usual absolute value
\cite[VI, \S6]{BCA}.  
\end{numba}
\begin{numba}\label{deflcx}
A topological vector space~$E$
over an ultrametric field $(\K,|.|)$
is called {\em locally convex\/}
if every zero-neighbourhood
of~$E$ contains an open $\bO$-submodule
of~$E$, where
$\bO$
is the valuation ring of~$\K$
(see \cite{Mon}, Ch.\,III,
\S2, Prop.\,4 and
\S3, D\'{e}f.\,1 when~$\K$ is complete).
It is well known that a topological vector space~$E$
over an ultrametric field is locally
convex if and only if its topology arises
from a family $(\|.\|_\gamma)_\gamma$ of {\em ultrametric\/}
continuous seminorms $\|.\|_\gamma\!: E\to [0,\infty[$,
satisfying $\|x+y\|_\gamma\leq\max\{\|x\|_\gamma,\|y\|_\gamma\}$
for all $x,y\in E$.
\end{numba}
\begin{numba}\label{dfblls}
If $(E,\|.\|)$
is a normed space over a valued field~$\K$,
given $\ve>0$ and $x\in E$ we write $B^E_\ve(x):=\{y\in E\!:
\|y-x\|<\ve\}$, or simply $B_\ve(x):=B^E_\ve(x)$
if $E$ is understood. Note that,
since the image of a norm $\|.\|$ need not be
contained in the image $|\K|$ of the absolute value,
it is not possible in general to normalize
elements: Given $0\not=x\in E$ we need not find
$t\in \K^\times$ such that $\|tx\|=1$.
To ensure that $\|Ax\|\leq \|A\|\,\|x\|$,
the operator norm of a linear operator $A\!: E\to F$
between normed spaces $(E,\|.\|_E)$
and $(F,\|.\|_F)$ therefore has to be defined as
$\|A\|:=\min\{C\geq 0\!: (\forall x\in E)\; \|Ax\|_F\leq C\|x\|_E\}$.
If $E=F=\K^d$ for some $d\in \N_0$ and $\|.\|_E=\|.\|_F$
is the maximum norm $\|.\|_\infty\!: \K^d\to[0,\infty[$,
$\|(x_1,\ldots, x_d)\|_\infty:=\max\{|x_1|,\ldots, |x_d|\}$,
then every non-zero vector can be normalized
and thus
$\|A\|=\max\{\|Ax\|_\infty\!: x\in \K^d,\, \|x\|_\infty\leq 1\}$.
As usual,
given topological vector spaces $E$ and $F$
over a topological field~$\K$,
we let $\cL(E,F)$ denote the set of all
continuous linear maps from~$E$ to~$F$;
we abbreviate $\cL(E):=\cL(E,E)$.
\end{numba}
Throughout the remainder of this section,
$\K$ denotes a
(non-discrete) topological field.\\[3mm]
Before we define $C^k$-maps,
we need an efficient notation for
the domains of certain mappings associated
with~$C^k$-maps.
\begin{defn}\label{domains}
If $E$ is a topological $\K$-vector
space and $U\sub E$ an open subset,
we define $U^{[0]}:=U$
and
\[
U^{[1]}:=\{(x,y,t)\in U\times E\times \K\!:
x+ty\in U\},
\]
which is an open subset of the topological
$\K$-vector space $E\times E\times \K$.
Having defined $U^{[j]}$ inductively
for a natural number $j\geq 1$,
we set $U^{[j+1]}:=(U^{[j]})^{[1]}$.
\end{defn}
In particular, $E^{[1]}=E\times E\times \K$,
$E^{[2]}=E\times E\times \K\times E\times E\times \K\times \K$,
etc.\
\begin{defn}\label{defnCr}
Let $E$ and $F$ be topological $\K$-vector
spaces and $f\!: U\to F$ be a mapping,
defined on an open subset $U\sub E$.
We say that $f$ is {\em of class $C^0_\K$}
if $f$ is continuous, in which case
we set $f^{[0]}:=f$
and call $f^{[0]}$ the $0$th extended difference
quotient map of~$f$.
If $f$ is continuous and there exists
a continuous mapping $f^{[1]}\!:
U^{[1]}\to F$ such that
\begin{equation}\label{diffquot}
{\textstyle 
\frac{1}{t}(f(x+ty)-f(x))=f^{[1]}(x,y,t)}\;\;
\mbox{for all $(x,y,t)\in U^{[1]}$ such that $t\not=0$,}
\end{equation}
we say that $f$ is {\em of class~$C^1_\K$},
and call $f^{[1]}$ the (first) {\em extended difference quotient
map\/} of~$f$ (note that, as $\K$ is non-discrete,
the continuous map $f^{[1]}$ is uniquely
determined by (\ref{diffquot})). 
Recursively, having defined $C^j_\K$-maps
and $j$th extended difference quotient maps
for $j=0,\ldots, k-1$ for some natural number $k\geq 2$,
we call $f$ a mapping {\em of class $C^k_\K$\/}
if $f$ is of class $C^{k-1}_\K$ and $f^{[k-1]}$
is of class~$C^1_\K$.
In this case, we define the $k$th extended difference
quotient map of~$f$ via
\[
f^{[k]}:= (f^{[k-1]})^{[1]}\!: U^{[k]}\to F\,.
\]
The mapping~$f$ is {\em of class $C^\infty_\K$\/}
(or {\em $\K$-smooth\/})
if it is of class $C^k_\K$ for all $k\in \N_0$.
If $\K$ is understood, we simply write $C^k$ instead
of $C^k_\K$, and call $f$ smooth or of class $C^\infty$
if it is $\K$-smooth.
If, conversely, we want to stress the fact that the
field~$\K$ is used, we shall write
$U^{[k]}_\K$ for $U^{[k]}$ and $f^{[k]}_\K$
for $f^{[k]}$.
\end{defn}
\begin{examples}\label{linbilin}
Every continuous $\K$-linear
mapping $\lambda\!: E\to F$ between topological $\K$-vector
spaces is smooth,
with $\lambda^{[1]}(x,y,t)=\lambda(y)$
for all $(x,y,t)\in E\times E\times \K$.
If $V, W$ and $F$ are topological
$\K$-vector spaces and $\beta\!: V\times W\to F$
is a continuous bilinear map,
then $\beta$ is smooth,
with
\[
\beta^{[1]}((v,w),(v',w'),t)=
\beta(v,w')+\beta(v',w)+t\beta(v',w')
\]
for all $v,v'\in V$, $w,w'\in W$, and $t\in \K$
(cf.\ \cite{Ber}).
\end{examples}
\begin{numba}\label{usefulsimp}
Note that, for $k\geq 2$,
a mapping
$f$ as above is of class $C^k_\K$
if and only if $f$ is of class $C^1_\K$
and $f^{[1]}$ is of class $C^{k-1}_\K$;
in this case, $f^{[k]}=(f^{[1]})^{[k-1]}$
(these claims are proved by a trivial induction).
\end{numba}
\begin{numba}\label{gargel}
Given a map $f\!: U\to F$ as before,
let $f^{]1[}\!: U^{]1[}\to F$,
$f^{]1[}(x,y,t):=\frac{1}{t}(f(x+ty)-f(x))$
be the associated difference quotient map, defined
on $U^{]1[}:=\{(x,y,t)\in U^{[1]}\!: t\not=0\}$.
Then $f$ is $C^1_\K$
if and only if it is continuous
and $f^{]1[}$ extends to a continuous function,
$f^{[1]}$, on $U^{[1]}$.
The set $U^{]1[}$ is open and dense
in $U^{[1]}$,
and we have $U^{[1]}=U^{]1[}\,\cup \,(U\!\times\!E\times\!\{0\})$,
as a disjoint union. If $f$ is $C^k_\K$,
then so is $f^{]1[}$ (cf.\ {\bf \ref{chainr}} below).
\end{numba}
\begin{numba}\label{differentials}
Given a $C^1_\K$-map $f\!: U\to F$ as before,
we define the {\em directional derivative\/}
of~$f$ at $x\in U$ in the direction $v\in E$
via
\[
df(x,v)\,:=\, (D_vf)(x)\, :=\, \lim_{0\not=t\to 0} {\textstyle
\frac{1}{t}(f(x+tv)-f(x))}
\, =\, f^{[1]}(x,v,0)\,.
\]
Then $df\!:U\times E\to F$ is continuous, being a partial map
of~$f^{[1]}$, and it can be shown that
the ``differential'' $df(x,\sbull)\!: E\to F$
of $f$ at~$x$ is a continuous $\K$-linear map,
for each $x\in U$ \cite[Prop.\,2.2]{Ber}.
If
$f$ is $C^2$,
we define a continuous map $d^2f\!: U\times E^2\to F$
via
\begin{eqnarray*}
d^2f(x,v_1,v_2) &:=& (D_{v_2}(D_{v_1}f))(x)\, =\, \lim_{0\not=t\to 0}
{\textstyle
\frac{1}{t}(df(x+tv_2,v_1)-df(x,v_1))}\\
&=& f^{[2]}(x,v_1,0,v_2,0,0,0).
\end{eqnarray*}
Similarly, if $f$ is of class $C^k_\K$,
we obtain continuous
mappings (the ``higher differentials'')
$d^jf\!: U\times E^j\to F$,
$d^jf(x,v_1,\ldots, v_j):=(D_{v_j}\cdots D_{v_1}f)(x)$
for all $j\in \N_0$
such that $j\leq k$
(where $d^0f:=f$).
Here $d^jf(x,\sbull)\!:E^j\to F$
is symmetric and $j$-multilinear \cite[La.\,4.8]{Ber}.
\end{numba}
\begin{numba}\label{MichalB}
If $\K$ is $\R$ or $\C$
and the range~$F$ is locally convex,
the considerations in {\bf \ref{differentials}}
show that every $C^k_\K$-map
in the preceding sense
is a $C^k$-map in the sense
of Michal-Bastiani (a $C^k_{MB}$-map for short),
{\em i.e.}, the iterated directional derivatives needed to define
$d^jf$ exist for all $j\in \N$
such that $j\leq k$,
and the mappings $d^jf\!: U\times E^j\to F$
so obtained (as well as~$f$)
are continuous
(such mappings are also called
``Keller's $C^k_c$-maps'' in the literature,
following~\cite{Kel}).
Exploiting the Fundamental Theorem
of Calculus,
it can be shown that,
conversely, every $C^k_{MB}$-map with locally convex range
is of class~$C^k_\K$ \cite[Prop.\,7.4]{Ber}.
Thus, when dealing with maps
into real or complex locally convex spaces,
it is fully sufficient (and much more convenient)
to work
with the differentials $d^jf$,
no use has to be made
of the mappings $f^{[j]}$.
However, as soon as we turn to
mappings into non-locally convex real or complex
topological vector spaces,
and also in the case of base fields
other than $\R$ and $\C$,
the differentials alone do
not encode enough information
on~$f$, and it is necessary to work
with the mappings~$f^{[j]}$.
For example,
consider the function $f\!: \;]0,1[\, \to L^0[0,1]$,
$f(t):=\one_{[0,t]}$ taking $t$ to the characteristic
function of the interval $[0,t]$;
here $[0,1]$ is equipped with Lebesgue measure,
and $L^0[0,1]$ denotes the space of equivalence classes
of measurable
real-valued functions on $[0,1]$ (modulo equality a.e.),
equipped with the topology of convergence
in measure (see \cite{KPR}). Then $f$ is injective.
It can be shown that $f$
is of class $C^\infty_\R$, with $d^jf$ vanishing
identically for all~$j\in \N$ (cf.\ \cite{NCX}).
\end{numba}
\begin{numba}\label{chainr}
(Chain Rule).
If $E$, $F$ and $H$ are topological
$\K$-vector spaces, $U\sub E$ and $V\sub F$ are
open subsets, and $f\!: U\to V\sub F$,
$g\!: V\to H$ are mappings of class~$C^k$,
then also the composition
$g\circ f\!: U\to H$ is of class~$C^k$.
If $k\geq 1$,
we have
$(f(x),f^{[1]}(x,y,t),t)\in V^{[1]}$
for all $(x,y,t)\in U^{[1]}$, and
\begin{equation}\label{formchain}
(g\circ f)^{[1]}(x,y,t)=g^{[1]}(f(x),f^{[1]}(x,y,t),t)\,.
\end{equation}
In particular, $d(g\circ f)(x,y)=dg(f(x),df(x,y))$
for all $(x,y)\in U\times E$ (see \cite{Ber},
Prop.\,3.1 and Prop.\,4.5).
\end{numba}
We recall that being of class $C^k$ is a local
property~\cite[La.\,4.9]{Ber}:
\begin{la}\label{Crlocal}
Let $E$ and $F$ be topological $\K$-vector spaces
and $f\!: U\to F$ be a mapping, defined
on an open subset
$U$ of~$E$. Let $k\in \N_0\cup\{\infty\}$. 
If there is an open cover $(U_i)_{i\in I}$
of~$U$ such that
$f|_{U_i}\!: U_i\to F$ is of class
$C^k$ for each $i\in I$,
then $f$ is of class~$C^k$.\Punkt
\end{la}
\begin{numba}\label{verybas}
Compositions of composable $C^k$-mappings
being of class~$C^k$, we can define $C^k$-manifolds
modeled on topological $\K$-vector
spaces in the usual way,
namely as Hausdorff
topological spaces~$M$, together with a
set (``atlas'') of
homeomorphisms (``charts'') from
open subsets of~$M$ onto open subsets
of the modeling topological $\K$-vector space~$Z$,
whose domains cover~$M$, and such that
the transition maps are of class~$C^k$.
A {\em Lie group over~$\K$}
is a group~$G$, equipped with
a smooth manifold structure
modeled on a topological $\K$-vector
space~$Z$, with respect to which
inversion and
the group multiplication are smooth mappings.
For every $\K$-Lie group~$G$,
the geometric tangent space $T_1(G)$
can be turned into a topological $\K$-Lie
algebra $L(G)$ in a natural way
(see \cite{Ber} for more information).
\end{numba}
\begin{rem}\label{findimgp}
It can be shown that
$\K$-analytic maps
from open subsets of (ultrametric) normed
spaces to locally convex
topological $\K$-vector spaces
(as in \cite{BoF}, where locally convex spaces are called ``polynormed'')
are $C^\infty_\K$,
for every ultrametric field
$(\K,|.|)$ \cite[Prop.\,7.20]{Ber}.
As a consequence,
every finite-dimensional $\K$-analytic 
Lie group~$G$
in the usual sense
(as defined in \cite[p.\,102]{Ser})
can be considered as a $\K$-Lie group
in our sense, and likewise
for the analytic Lie groups
modeled on ultrametric Banach spaces
considered in~\cite{Bou}.
\end{rem}
We recall three simple,
but very useful facts (\cite{Ber}, Lemmas~10.1--10.3):
\begin{la}\label{corestr}
Let $E$ and $F$ be topological $\K$-vector
spaces, $U\sub E$ be open,
and $f\!: U\to F$ be a mapping of class~$C^k$, where $k\in \N_0
\cup\{\infty\}$.
Let $F_0$ be a vector subspace of~$F$ containing
the image of~$f$.
If $F_0$ is closed or if $F_0$ is sequentially closed
and~$\K$ is metrizable,
then the co-restriction
$f|^{F_0}\!: U\to F_0$
is $C^k$ as a map into~$F_0$.\Punkt
\end{la}
\begin{la}\label{inprod}
Suppose that $E$ is a topological $\K$-vector
space, $(F_i)_{i\in I}$ a family of
topological $\K$-vector spaces,
$U\sub E$ an open subset,
and $f\!: U\to P$ a mapping,
where $P:=\prod_{i\in I}\, F_i$,
with canonical projections
$\pr_i\!: P\to F_i$.
Let $k\in \N_0\cup\{\infty\}$.
Then $f$ is of class~$C^k$ if and only if
$\pr_i\circ f$ is of class~$C^k$
for each $i\in I$.\Punkt
\end{la}
\begin{la}\label{inpl}
Let $E$ be a topological $\K$-vector
space, $(F_i)_{i\in I}$ be a family of
topological $\K$-vector spaces,
$U\sub E$ be open,
and $f\!: U\to F$ be a map,
where $F=\pl_{i\in I}\, F_i$,\vspace{-.8mm}
with limit maps $\pi_i\!: F\to F_i$.
Let $k\in \N_0\cup\{\infty\}$.
Then $f$ is~$C^k\,$ iff
$\,\pi_i\circ f$ is~$C^k$
for each $i\in I$.\Punkt
\end{la}
As in the case of Banach-Lie groups,
general Lie groups can be described locally:
\begin{prop}[Local description
of Lie group structures]\label{locchar}
Suppose that a subset $\,U$ of a group $\,G$ is
equipped with a smooth manifold structure
modeled on a topological $\K$-vector space~$E$,
and suppose that $\,V$ is an open subset of~$\,U$
such that $\,1\in V$, $V=V^{-1}$, $VV\sub U$,
and such that the multiplication map
$V\times V \to U$, $(g,h)\mto gh$ is smooth
as well as inversion $V\to V$, $g\mto g^{-1}$;
here $V$ is considered as an open submanifold of~$\,U$.
Suppose that for every element $x$ in a symmetric
generating set
of~$\,G$, there is an open identity-neighbourhood
$W\sub U$ such that $xWx^{-1}\sub U$,
and such that the mapping $W\to U$, $w\mto xwx^{-1}$
is smooth.\footnote{This condition is automatically
satisfied if~$V$ generates~$G$.}
Then there is a unique
$\K$-Lie group structure on~$G$
which makes~$V$,
equipped with the above
manifold structure,
an open submanifold of~$\,G$.
\end{prop}
\begin{proof}
The proof of \cite{Bou}, Chapter~3, \S1.9, Proposition~18
can easily be adapted.
\end{proof}
\section{Unit groups of continuous inverse algebras}\label{secALG}
Let $\K$ be an arbitrary topological field.
In this section,
we show that the groups~$A^\times$ of invertible
elements in suitable topological
$\K$-algebras $A$ (the continuous
inverse $\K$-algebras) are $\K$-Lie groups.
We describe constructions
producing new continuous inverse algebras
from given ones.
In this way, we obtain
a supply of continuous inverse $\K$-algebras
and thus also of $\K$-Lie groups.
For much more information
concerning locally convex
continuous inverse algebras
over the real or complex field,
and their unit groups,
the reader is referred to~\cite{ALG}.
Further examples can be found in \cite{HOL},
\cite[1.15]{Gra}, and \cite{YMO}.
\begin{defn}\label{defnci}
A {\em continuous inverse algebra\/} (over $\K$) is
a unital associative topological $\K$-algebra~$A$
whose group of units $A^\times$ is open in~$A$
and such that inversion $\iota\!: A^\times\to A$,
$a\mto a^{-1}$
is continuous.
\end{defn}
Continuous inverse algebras are
of interest in the present context,
for the following reason:
\begin{prop}\label{ciasmooth}
For every continuous inverse $\K$-algebra $A$,
inversion
$\iota\!: A^\times\to A$ is of class $C^\infty_\K$,
and thus $A^\times$ is a $\K$-Lie group
when considered as an open
submanifold of~$A$.
\end{prop}
\begin{proof}
The algebra multiplication $A\times A\to A$ is continuous bilinear
and hence smooth. Hence so is the group multiplication
$A^\times\times A^\times \to A^\times$.
We now show by induction that $\iota$ is $C^k_\K$ for each $k\in \N_0$.
By hypothesis, $\iota$ is $C^0_\K$.
Suppose that $\iota$ is $C^k_\K$.
Using that $b^{-1}-a^{-1}=b^{-1}(a-b)a^{-1}$
for $a,b\in A^\times$, we obtain for any
$(x,v,t)\in (A^\times)^{[1]}\sub A^\times\times A\times \K$:
\begin{equation}\label{showsC1}
\iota(x+tv)-\iota(x)=
(x+tv)^{-1}-x^{-1}=-t((x+tv)^{-1}vx^{-1}
=tF(x,v,t)\,,
\end{equation}
where $F\!: (A^\times)^{[1]}\to A$,
$F(x,v,t):=-\iota(x+tv)v\iota(x)$.
Since $\iota$ is $C^k_\K$ by the induction hypotheses,
$F$ is of class $C^k_\K$, in particular of class~$C^0_\K$.
Thus (\ref{showsC1})
shows that $\iota$ is of class $C^1_\K$,
with $\iota^{[1]}=F$
a mapping of class~$C^k_\K$.
Therefore $\iota$ is of class $C^{k+1}_\K$
(see {\bf \ref{usefulsimp}}).
\end{proof}
For example, $\K$ is a continuous inverse
algebra over~$\K$, and thus $\K^\times$
is a $\K$-Lie group.
\begin{prop}\label{matrices}
If $A$ is a continuous inverse $\K$-algebra,
then so is the algebra $M_n(A)$ of $n\!\times\! n$-matrices with entries
in~$A$,
when equipped with the natural vector topology
$\isom A^{n\times n}$.
\end{prop}
\begin{proof}
Apparently, $M_n(A)$ is a topological $\K$-algebra.
Its unit group is open and inversion is
continuous by \cite[Cor.\,1.2]{Swa}.
\end{proof}
\begin{la}\label{subfindim}
If $A$ is a finite-dimensional unital associative
$\K$-algebra and $B\sub A$ a unital subalgebra,
then $B^\times=A^\times\cap B$.
\end{la}
\begin{proof}
This is a well-known fact
(cf.\ \cite[La.\,9.4]{ALG} if $\K=\R$
or $\C$).
\end{proof}
\begin{numba}\label{canonic}
Recall that, among the vector topologies
on a finite-dimensional
$\K$-vector space~$V$ of dimension~$d$, there is a uniquely determined
vector topology making~$V$ isomorphic
to the direct product $\K^d$ as a topological vector space.
It is called the {\em canonical $\K$-vector space topology\/}
(see \cite{BTV}, Ch.\,I, \S1, no.\,1, Example~5).
\end{numba}
We remark that it
is possible to characterize those
topological fields~$\K$ having the special property
that any finite-dimensional
$\K$-vector space admits only one
(Hausdorff) vector topology
(\cite{Nac},
also \cite{Wie}, Section~5.4, Theorem~10).
They are necessarily complete.
\begin{prop}\label{findimalg}
Let $A$ be a finite-dimensional
unital associative $\K$-algebra.
Then the canonical $\K$-vector space topology turns~$A$
into a continuous inverse algebra
over~$\K$.
\end{prop}
\begin{proof}
Let $n:=\dim_\K(A)$.
It is clear that the canonical $\K$-vector space
topology
turns~$A$ into a topological $\K$-algebra,
and it is clear that
the left
regular representation
\[
\lambda\!: A\to \cL(A)\isom M_n(\K),\;\;\;
\lambda(a)(b):=ab
\]
is a topological embedding
(where $\cL(A)$ denotes the $\K$-algebra
of $\K$-linear self-maps of~$A$,
equipped with the canonical $\K$-vector space topology).
It therefore suffices to assume
that $A$ is a subalgebra of $M_n(\K)$.
Now, $M_n(\K)$ being a finite-dimensional
$\K$-algebra, we have
\begin{equation}\label{intsect}
A^\times=M_n(\K)^\times\cap A
\end{equation}
by Lemma~\ref{subfindim}.
Since $M_n(\K)$ is a continuous
inverse algebra by Proposition~\ref{matrices},
$M_n(\K)^\times$ is open in $M_n(\K)$
and thus $A^\times$ is open in~$A$ by (\ref{intsect}).
The inversion map $\iota\!: A^\times\to A$
being a restriction
of the continuous inversion map
$M_n(\K)^\times\to M_n(\K)$,
we deduce that~$\iota$ is continuous.
\end{proof} 
Tensor products of finite-dimensional
algebras and continuous inverse algebras
are again continuous inverse algebras.
\begin{prop}\label{tensorprod}
Given a continuous inverse $\K$-algebra~$A$
and finite-dimensional unital associative
$\K$-algebra~$F$,
consider the associative unital
$\K$-algebra $F\tensor_\K A$.
Pick any $\K$-basis $e_1,\ldots,e_n$
of~$F$, and
equip $F\tensor_\K A$ with the topology making
$\phi\!: A^n\to F\tensor_\K A$, $(a_i)_{i=1}^n
\mto \sum_{i=1}^n
e_i\tensor a_i$ an isomorphism of topological
$\K$-vector spaces. Then this topology does not depend on the choice of basis,
and it turns $F\tensor_\K A$ into a continuous inverse algebra
over~$\K$.
\end{prop}
\begin{proof}
The natural map $M_n(\K)\times A^n\to A^n$
being continuous, we readily deduce that the topology on
$F\tensor_\K A$ is independent of the choice
of $\K$-basis
for~$F$.
Given $i,j\in\{1,\ldots, n\}$,
we have $e_ie_j=\sum_{k=1}^n t_{i,j,k}e_k$
for uniquely determined elements (``structure
constants'') $t_{i,j,k}\in \K$.
Given $z=(z_i)_{i=1}^n$,
$v=(v_i)_{i=1}^n$ in $A^n$,
we have
\[
\phi(z)\cdot \phi(v)=\sum_{k=1}^n e_k\tensor
\left(\sum_{i,j=1}^n t_{i,j,k} z_iv_j\right)=\phi\left(\Big(\sum_{i,j=1}^n
t_{i,j,k}z_iv_j\Big)_{k=1}^n\right)\, .
\]
As~$A$ is a topological $\K$-algebra,
we readily deduce from the preceding formula
that multiplication in $F\tensor_\K A$ is
continuous. Thus $F\tensor_\K A$ is a topological
$\K$-algebra. Given $z$ and~$v$ as before,
we calculate
\begin{eqnarray}
(1+\phi(z))\cdot(1+\phi(v)) & = &
(1+\sum_{i=1}^n e_i \tensor z_i)\cdot
(1+\sum_{j=1}^n e_j\tensor v_j) \nonumber\\
& = &
1+\sum_{i,j=1}^n \underbrace{e_ie_j}_{=\sum_{k=1}^n t_{i,j,k}e_k}
\!\!\!\!\tensor \;z_iv_j
+\sum_{k=1}^n e_k\tensor z_k+\sum_{k=1}^n e_k\tensor v_k\nonumber\\
& = &
1+\sum_{k=1}^n e_k\tensor (z_k + (S(z).v)_k)\, ,\label{lastl}
\end{eqnarray}
where
$S(z):=(a_{kj}(z))_{k,j=1}^n\in M_n(A)$
with
$a_{kj}(z):=\delta_{k,j}+\sum_{i=1}^n z_it_{i,j,k}$,
and where
$(S(z).v)_k$ denotes the $k$th coordinate of the vector
$S(z).v\in A^n$.
We strive to show that for
$z\in A^n$ in some zero-neighbourhood,
we can choose~$v$ such that $S(z).v=-z$.
Then, by Equation~(\ref{lastl}),
the element $1+\phi(v)$ will be a right inverse for $1+\phi(z)$.\\[3mm]
Since $M_n(A)$ is a continuous inverse $\K$-algebra
(Proposition~\ref{matrices}),
and $S\!: A^n\to M_n(A)$ is a continuous
mapping such that $S(0)=\one\in M_n(A)^\times$,
there is a zero-neighbourhood~$U$ in~$A^n$
such that $S(U)\sub M_n(A)^\times$.
Apparently, the mapping
\[
\rho\!: U \to A^n,\;\;\; \rho(z):=  -S(z)^{-1}.z
\]
is continuous.
For each $z\in U$, we have
$S(z)\rho(z)=-z$,
and thus $(1+\phi(z))\cdot (1+\phi(\rho(z)))=1$
by (\ref{lastl}).
Thus, for each
$a$ in the open identity neighbourhood
$V:=1+\phi(U)\sub F\tensor_\K A$, the element
\[r(a):=1+\phi(\rho(\phi^{-1}(a-1)))\in F\tensor_\K A\]
is a right inverse for~$a$ in $F\tensor_\K A$,
and the mapping $r\!:V\to F\tensor_\K A$ is continuous.
Very similar arguments show that there
is an identity neighbourhood $W$ in $F\tensor_\K A$
such that every $a\in W$ has a left inverse in $F\tensor_\K A$.
Then $P:=V\cap W$ is an identity neighbourhood
in $F\tensor_\K A$ such that $P\sub (F\tensor_\K A)^\times$,
and the inversion map $\iota\!:
(F\tensor_\K A)^\times \to F\tensor_\K A$
satisfies $\iota|_P=r|_P$ an thus is continuous
on~$P$. As a consequence, $(F\tensor_\K A)^\times$
is open in $F\tensor_\K A$ and $\iota$
is continuous (cf.\ \cite[La.\,2.8]{ALG}).
\end{proof}
Concerning extension of scalars,
we readily deduce:
\begin{cor}\label{extendscale}
For every continuous inverse algebra~$A$ over~$\K$
and finite extension field~$\bL$ of~$\K$,
$A_\bL:=\bL\tensor_\K A$ is a continuous inverse
algebra over~$\bL$ $($where $\bL$ is
equipped with the canonical $\K$-vector space topology$)$.
\end{cor}
\begin{proof}
We equip $A_\bL=\bL\tensor_\K A$
with the topological $\K$-algebra
structure defined in Proposition~\ref{tensorprod},
which makes it a continuous inverse
algebra over~$\K$. It is easy to see that
the mapping
$\bL\to \bL\tensor_\K A$, $z\mto z\tensor 1$
is a continuous
$\K$-algebra homomorphism.
The continuity of scalar multiplication
$\bL\times A_\bL\to A_\bL$
therefore follows from the continuity of
the multiplication map $A_\bL\times A_\bL\to A_\bL$.
\end{proof}
For the final result
of this section,
we specialize to the case where
$\K$ is a locally compact
topological field.
\begin{numba}\label{refback}
We consider a unital associative
$\K$-algebra~$A$
which is {\em locally finite\/}
in the sense that every finite subset
is contained in a finite-dimensional
subalgebra of~$A$. We also assume
that~$A$ is of countable dimension
as a $\K$-vector space.
As a consequence,
there exists an ascending sequence
$A_1\sub A_2\sub \cdots$ of finite-dimensional unital
subalgebras $A_n$ of~$A$ such that $A=\bigcup_{n\in \N}A_n$.
\end{numba}
\begin{numba}\label{defnftop}
We equip~$A$ with the so-called
``finite topology,'' {\em i.e.},
the final topology
with respect to the inclusion maps
$\lambda_F\!: F\to A$, where $F$ runs
through the set $\cF$ of finite-dimensional
vector subspaces of~$A$.
The set $\{A_n\!: n\in \N\}$
being co-final in $\cF$ (directed with respect to
inclusion), the finite topology
on~$A$ is also the final topology with respect to
the family $(\lambda_{A_n})_{n\in \N}$.
Then $A=\dl\,A_n$\vspace{-.8mm} as a topological space,
furthermore
$A\times A=\dl\, (A_n\times A_n)$
and $\K\times A=\dl\, (\K\times A_n)$,\vspace{-.8mm}
each $A_n$ and ~$\K$ being locally compact
(\cite{HSTH} or \cite[Prop.\,3.3]{DIR}).
As a consequence, $A$ is a topological $\K$-algebra
(cf.\ \cite{DIR}).
\end{numba}
\begin{prop}\label{dirlimalgs}
Every countable-dimensional,
locally finite associative unital
algebra over
a locally compact
topological field~$\K$
is a continuous inverse $\K$-algebra
when equipped with the finite
topology.
\end{prop}
\begin{proof}
We have already shown that~$A$ is a topological
$\K$-algebra. The openness of $A^\times$ in~$A$
as well as continuity of inversion can be
shown as in the real and complex
special cases (see \cite[Prop.\,9.5]{ALG}).
\end{proof}
\begin{rem}\label{remmentBCH}
If $A$ is a real or complex locally convex
CIA, then $A^\times$ is in fact an analytic Lie group.
If, furthermore, $A$ is complete
(or, at least, Mackey complete),
then $A^\times$ is a Baker-Campbell-Hausdorff
(BCH) Lie group, {\em viz.}
it possesses a locally analytically diffeomorphic
exponential function, and its multiplication is
given locally by the BCH-series
(see \cite{ALG}).
In this case, the results of \cite{Rob}
and \cite{GCX} facilitate to integrate
closed Lie subalgebras of~$A$ to analytic subgroups
of~$A^\times$, providing us with a much richer
supply of ``linear Lie groups''
then the mere full unit groups~$A^\times$.
In the case where $\K$ is a complete valued field,
in some cases subgroups of unit groups of continuous
inverse $\K$-algebras
may be turned into Lie groups using the inverse
function and
implicit function theorems from~\cite{IMP}.
\end{rem}
\section{Spaces of continuous
mappings and mappings\protect\\
between them}\label{secveccts}
As a preliminary for our studies
in Section~\ref{secmapgps}, where we
shall turn
the group $C(K,G)$ of continuous
mappings from a compact topological
space~$K$ to a $\K$-Lie group~$G$
into a $\K$-Lie group,
in the present section
we study differentiability properties
of certain types of mappings between
spaces of continuous
vector-valued functions on compact
topological spaces.
More generally, here (and in Section~\ref{secmapgps})
we can consider mappings on non-compact spaces
supported in given compact sets.\\[3mm]
Throughout
this section,
$\K$ denotes a
topological field,
$X$ a topological space,
and $K$ a compact subset of~$X$.
\begin{numba}\label{defncsp}
If $E$ is a topological
$\K$-vector space, we let
$C_K(X,E)\sub E^X$ denote the $\K$-vector
space of all
continuous mappings $\gamma\!: X\to E$
such that $\Supp(\gamma)\sub K$.
We equip $C_K(X,E)$ with the topology of uniform
convergence, which apparently
makes $C_K(X,E)$ a topological $\K$-vector
space. A basis of open
zero-neighbourhoods
is given by the sets $C_K(X,U):=$\linebreak
$\{\gamma\in
C_K(X,E)\!: \im\,\gamma\sub U\}$,
where $U$ ranges through the open zero-neighbourhoods
in~$E$.
\end{numba}
\begin{numba}\label{lcx1}
Note that if $\K\in\{\R,\C\}$ and $E$ is locally convex,
then $C_K(X,U)$ is convex for each convex,
open $0$-neighbourhood $U\sub E$
and thus $C_K(X,E)$ is locally convex.
If $\K$ is an ultrametric field
with valuation ring~$\bO$
and $E$ is locally convex (see {\bf \ref{deflcx}}), then
$C_K(X,U)$ is an open $\bO$-submodule of
$C_K(X,E)$ for each open $\bO$-submodule
$U\sub E$, and hence $C_K(X,E)$ is locally convex.
\end{numba}
The following proposition (and a $C^r$-analogue
to be proved later)
is the technical
backbone of our discussion of mapping groups.
\begin{prop}\label{backbone}
Let $E$, $F$, and $Z$ be topological
$\K$-vector spaces, $U\sub E$ and $P\sub Z$
be open subsets, $k\in \N_0\cup\{\infty\}$, and
$f\!: X\times U\times P\to F$
be a mapping.
Suppose that
\begin{itemize}
\item[\n (a)]
$f(x,\sbull)=0$
for all $x\in X\,\take\, K$;
\item[\n (b)]
$f(x,\sbull)\!: U\times P\to F$
is of class~$C^k$ for each $x\in X$, and
\item[\n (c)]
$X\times (U\times P)^{[j]}\to F$,
$(x,y)\mto f(x,\sbull)^{[j]}(y)$
is a continuous map, for each $j\in \N_0$
such that $j\leq k$.
\end{itemize}
Then $C_K(X,U):=C_K(X,E)\cap U^X$ is a $($possibly empty$)$
open subset of $C_K(X,E)$, and
\[
\phi\!:
C_K(X,U)\times P \to C_K(X,F),\;\;\;\;
\phi(\gamma,p):=f(\sbull,p)_*(\gamma)
\]
is a mapping of class $C^k$
$($where $f(\sbull,p)_*(\gamma)(x):=f(x,\gamma(x),p)$
for $x\in X)$.
\end{prop}
\begin{proof}
It is clear that $C_K(X,U)$ is open, and that
$\phi(\gamma,p)\in C_K(X,F)$
indeed. To show that $\phi$ is of class~$C^k$,
we clearly may assume that $k<\infty$.
The proof is by induction.\\[2mm]
{\em The case $k=0$.}
Let $\xi \in C_K(X,U)$,
$p\in P$, and $V\sub F$ be an open zero-neighbourhood.
Let $W\sub F$ be an open zero-neighbourhood
such that $W-W\sub V$.
For each $x\in K$,
we find an open neighbourhood $A_x\sub K$
of~$x$ in~$K$ and open zero-neighbourhoods
$B_x\sub E$
and $C_x\sub Z$
such that $\xi(A_x)+B_x\sub U$,
$p+C_x\sub P$,
and
\[
f(y,u,q)-f(x,\xi(x),p)\in W
\]
for all $y\in A_x$, $u\in \xi(A_x)+B_x$,
and $q\in p+C_x$.
By compactness, $K\sub \bigcup_{x\in I}A_x$
for some finite subset $I\sub K$.
Then $B:=\bigcap_{x\in I}B_x\sub E$
and $C:=\bigcap_{x\in I}C_x\sub Z$
are open zero-neighbourhoods.
Let $\eta\in \xi+C_K(X,B)$
and $q\in p+C\sub P$.
Given $y\in K$,
there is $x\in I$ such that
$y\in A_x$.
Hence
\begin{eqnarray*}
f(y,\eta(y),q)-f(y,\xi(y),p)
& = & f(y,\eta(y),q)-f(x,\xi(x),p)
-(f(y,\xi(y),p)-f(x,\xi(x),p))\\
& \in & W-W\sub V.
\end{eqnarray*}
For $y\in X\,\take\,K$ on the other hand, we have
$f(y,\eta(y),q)=f(y,\xi(y),p)=0$
and thus
$f(y,\eta(y),q)-f(y,\xi(y),p)=0\in V$ trivially.
We have shown that
$\phi(\eta,q)-\phi(\xi,p)\in C_K(X,V)$
for all $(\eta,q)$
in the open neighbourhood
$(\xi+C_K(X,B))\times (p+C)$
of $(\xi,p)$.
Thus $\phi$ is continuous, as required.\\[3mm]
{\em Induction step.}
Suppose that $k\geq 1$, and suppose
that the proposition holds for~$k$ replaced with $k-1$.
Abbreviate $Q:=(C_K(X,U)\times P)^{[1]}\sub
C_K(X,U)\times P\times C_K(X,E)\times Z \times \K$
and $Q^\times:=\{
(\xi,p,\eta,q,t)\in Q\!: t\not=0\}$.
For all $(\xi,p,\eta,q,t)\in Q^\times$,
we have
\begin{eqnarray}
{\textstyle \frac{1}{t}(\phi(\xi+t \eta,p+tq)-\phi(\xi,p))(x)}
& = &
{\textstyle
\frac{1}{t}
(f(x,\xi(x)+t\eta(x),p+tq)-f(x,\xi(x),p))}\nonumber\\
& = &
f(x,\sbull)^{[1]}((\xi(x),p),\;
(\eta(x),q),\; t)\label{suggesti}
\end{eqnarray}
for all $x\in X$, which
suggests to define
\[
\psi\!: Q\to C_K(X,F),\;\;\;
\psi(\xi,p,\eta,q,t)(x):=
f(x,\sbull)^{[1]}((\xi(x),p),\;
(\eta(x),q),\; t)\;\;\;\;
\mbox{for $x\in X$}.
\]
If we can show that $\psi$
is continuous,
then $\phi$ will be $C^1$ with
$\phi^{[1]}=\psi$, by~(\ref{suggesti}).\\[3mm]
{\em Claim~$1$.
$\psi$ is of class $C^{k-1}$ on $Q^\times$.}
In fact, inversion $\K^\times\to \K^\times$
being smooth,
addition and scalar multiplication
in $C_K(X,E)$ and $C_K(X,F)$
being continuous linear (resp., bilinear)
and thus smooth,
and $\phi$ being of class $C^{k-1}$ by induction,
the claim readily follows from
the formula $\psi(\xi,p,\eta,q,t)=
\frac{1}{t}(\phi(\xi+t\eta,p+tq)-\phi(\xi,p))$
for $(\xi,p,\eta,q,t)\in Q^\times$.\\[3mm]
{\em
Claim~$2$. Every $(\xi,p,\eta,q,0)\in Q$
has an open neighbourhood on which
$\psi$ is of class~$C^{k-1}$.}
In fact, since
$\im\,\xi\sub \xi(K)\cup\{0\}$ and $\im\,\eta\sub \eta(K)\cup\{0\}$
are compact
subsets of~$E$, and $\im\,\xi \sub U$,
there exist open neighbourhoods
$A\sub U$ of $\im\,\xi$,
$B\sub E$ of $\im\,\eta$
and an open zero-neighbourhood
$C\sub \K$
such that $A+C\cdot B\sub U$.
Shrinking~$C$ if necessary,
we may furthermore assume
that there exist open neighbourhoods
$D\sub P$ of~$p$ and
$G\sub Z$ of~$q$
such that $D+C\cdot G\sub P$.
Then $U_1:=A\times B$ is an open subset
of $E\times E$ containing
$\im (\xi,\eta)$,
and $P_1:= D\times G\times C$
is an open neighbourhood
of $(p,q,0)$ in $P\times Z\times \K$.
The definition
\[
f_1\!: X\times U_1\times P_1\to F,\;\;\;\;
f_1(x,(a,b),(p',q',t)):=
f(x,\sbull)^{[1]}((a,p'),(b,q'),t)
\]
makes sense by choice of $U_1$ and $P_1$
({\em i.e.}, the expression on
the right
hand side is defined). As an immediate
consequence of hypothesis\,(c),
the mapping~$f_1(x,\sbull)$ is of class $C^{k-1}$,
for each $x\in X$,
and $X\times (U_1\times P_1)^{[j]}\to F$,
$(x,y)\mto f_1(x,\sbull)^{[j]}(y)$
is continuous
for all $j\in \N_0$ such that $j\leq k-1$.
By induction,
\[
\phi_1\!:
C_K(X,U_1)\times P_1\to C_K(X,F),\;\;\;\;
\phi_1(\gamma,p_1)(x):=f_1(x,\gamma(x),p_1)
\]
is a mapping of class $C^{k-1}$.
Since $\phi_1((\sigma,\tau),(p',q',t))=
\psi(\sigma,p',\tau,q',t)$
for all $(\sigma,p',\tau, q',t)\in
C_K(X,A)\times D\times C_K(X,B)\times G\times C
\sub Q$,
Claim~2 is established.\\[3mm]
In view of Claims~1 and~2,
Lemma~\ref{Crlocal}
shows that $\psi$ is a mapping of
class~$C^{k-1}$.
In particular, $\psi$ is continuous
and thus, in view of (\ref{suggesti}),
the mapping~$\phi$ is of class~$C^1$
with $\phi^{[1]}=\psi$
of class~$C^{k-1}$.
Thus $\phi$ is of class~$C^k$.
\end{proof}
We readily deduce:
\begin{cor}\label{pushforw}
Let $E$ and $F$ be topological
$\K$-vector spaces, $U\sub E$
an open subset, $k\in \N_0\cup\{\infty\}$, and
$f\!: X\times U \to F$
be a mapping.
Suppose that
\begin{itemize}
\item[\n (a)]
$f(x,\sbull)=0$
for all $x\in X\,\take\, K$;
\item[\n (b)]
$f(x,\sbull)\!: U\to F$
is of class~$C^k$ for each $x\in X$, and
\item[\n (c)]
$X\times U^{[j]}\to F$, $(x,y)\mto f(x,\sbull)^{[j]}(y)$
is continuous, for each $j\in \N_0$
with $j\leq k$.
\end{itemize}
Then
\[
f_*\!:
C_K(X,U) \to C_K(X,F),\;\;\;\;
f_*(\gamma)(x):=f(x,\gamma(x))
\]
is a mapping of class $C^k$.\Punkt
\end{cor}
\begin{cor}\label{Cf}
Let $E$ and $F$ be topological $\K$-vector
spaces
and $f\!:U\to F$ be a mapping of class~$C^k$,
defined on an open subset~$U$ of~$E$.
If $K\not=X$, assume that $0\in U$ and $f(0)=0$.
Then
\[
C_K(X,f)\!: C_K(X,U)\to C_K(X,F),\;\;\;\;
\gamma\mto f\circ \gamma
\]
is a mapping of class~$C^k$.
\end{cor}
\begin{proof}
We have $C_K(X,f)=g_*$,
where $g\!:X\times U\to F$,
$g(x,y):=f(y)$ is easily seen to satisfy
Conditions (a), (b), (c)
of Corollary~\ref{pushforw}.
\end{proof}
Before working through the analogues
of the preceding facts for spaces of $C^r$-maps
stated in the next section---which are
considerably harder to prove---the
reader may wish to pass directly to
the construction of continuous mapping groups
in Section~\ref{secmapgps} (assuming $r=0$
there), to see what
the results just proved are good for.
\section{Spaces of {\boldmath $C^r$}-maps
and mappings between them}\label{secvecCr}
In this section, we discuss
spaces of vector-valued $C^r$-maps,
and mappings between such spaces,
to facilitate the construction of
a manifold structure on groups of $C^r$-maps
in Section~\ref{secmapgps}.
We begin with the special case of
vector-valued mappings on
open subsets of topological vector spaces.\\[3mm]
In this section,
$\F$ denotes a
topological field, and $\K$ a topological field
extending~$\F$, meaning that $\K$ contains $\F$
as a subfield, and
that the inclusion map $\F\to \K$
is continuous.\footnote{Typical examples
are: 1. $\K=\F$; 2. $\F=\R$, $\K=\C$.}
Starting with Proposition~\ref{pushforw2},
we shall assume that~$\F$ is locally compact.
We remark that,
if $\F$ is a valued field (for instance, if $\F$ is locally compact),
then the inclusion map
$\F\to \K$ is a topological embedding
automatically,
as every $1$-dimensional (Hausdorff) topological $\F$-vector
space is topologically isomorphic to~$\F$
(\cite{Wie}, \S5.1, Example~1
and \S5.4, Theorem~9).
%
%
%
%
\begin{center}
{\bf The spaces {\boldmath $C^r(U,E)$},
when {\boldmath $U$} is an open subset of the modeling space}
\end{center}
The preparatory results concerning mappings on open subsets
$U$ of topological vector spaces provided in this
subsection are essential
for our later discussion of the general case,
where $U$ is replaced with a manifold.
\begin{numba}\label{top1}
Given a
topological $\K$-vector space~$E$
and open subset $U$ of a
topological $\F$-vector space~$Z$,
we let $C^r(U,E)$ be the set of all
mappings $\gamma\!: U\to E$ of class~$C^r_\F$
(where $r\in \N_0\cup\{\infty\}$).
It is clear that pointwise
operations turn
$C^r(U,E)$ into a $\K$-vector space.
We give $C^r(U,E)$ the initial topology with respect
to the family of mappings
\[
C^r(U,E)\to C(U^{[j]},E),\;\;\;\;
\gamma\mto \gamma^{[j]},
\]
where $j\in \N_0$ such that $j\leq r$,
and where $C(U^{[j]},E)$ is equipped with the
topology of uniform convergence on compact sets
(which coincides with the compact-open topology).
It is clear that
$C^r(U,E)$ becomes a topological $\K$-vector space
in this way.
\end{numba}
The sets\label{thesets}
\[
\lfloor K,W\rfloor:=\{ \gamma\in C(U^{[j]},E)\!:
\gamma(K)\sub W\}
\]
form a basis of open zero-neighbourhoods
for the topology on $C(U^{[j]},E)$
when $K$ ranges through the compact subsets of~$U^{[j]}$
and $W$ through the open zero-neighbourhoods of~$E$.
\begin{rem}\label{simplobs}
The following assertions readily
follow from the definitions:
\begin{itemize}
\item[\n (a)]
For every $r\geq s$, the inclusion
map $C^r(U,E)\to C^s(U,E)$ is a continuous
linear map.
The topology on
$C^\infty(U,E)$
is initial with respect
to the family of inclusion maps
$C^\infty(U,E)\to C^k(U,E)$, where $k\in \N_0$.
Furthermore,
$C^\infty(U,E)=\pl\, C^k(U,E)$.
\item[\n (b)]
For every $k\in \N_0$,
the topology on $C^{k+1}(U,E)$ is initial
with respect to the inclusion
map $C^{k+1}(U,E)\to C(U,E)$
together with the mapping
\[
C^{k+1}(U,E)\to C^k(U^{[1]},E),\;\;\;
\gamma\mto \gamma^{[1]}\,.
\]
\end{itemize}
\end{rem}
Strengthening (b), we have:
\begin{la}\label{closdinprod}
In the preceding situation, the map
\[
\Lambda : \, C^{k+1}(U,E)\, \to \; C(U,E)\times C^k(U^{[1]},E)\,,
\quad \Lambda(\gamma)\, :=\, (\gamma,\gamma^{[1]})
\]
is a topological embedding onto a closed
vector subspace of $C(U,E)\times C^k(U^{[1]},E)$.
\end{la}
\begin{proof}
By Remark~\ref{simplobs}\,(b), the map $\Lambda$ is a topological embedding.
To see that $\, \im(\Lambda)$ is closed,
let $(\gamma_\alpha)$ be a net in $C^{k+1}(U,E)$
such that $\Lambda(\gamma_\alpha)$ converges
in $C(U,E)\times C^k(U^{[1]},E)$, say to $(\gamma,\eta)$
with $\gamma\in C(U,E)$ and $\eta\in C^k(U^{[1]},E)$.
Let $(x,y,t)\in U^{]1[}$ (see {\bf \ref{gargel}} for the notation).
Then $\gamma_\alpha(x)\to \gamma(x)$
and $\gamma_\alpha(x+ty)\to\gamma(x+ty)$,
entailing that $(\gamma_\alpha)^{[1]}(x,y,t)=
\frac{1}{t}(\gamma_\alpha(x+ty) -\gamma_\alpha(x))
\to \frac{1}{t}(\gamma(x+ty) -\gamma(x))$.
Since also $(\gamma_\alpha)^{[1]}(x,y,t)\to \eta(x,y,t)$,
we deduce that $\eta(x,y,t)=\frac{1}{t}(\gamma(x+ty) -\gamma(x))$.
The map $\eta$ being continuous,
this means that $\gamma$ is $C^1_\F$,
with $\gamma^{[1]}=\eta$. The map $\gamma^{[1]}=\eta$
being $C^k_\F$, we deduce that $\gamma$ is $C^{k+1}_\F$
and thus
$\gamma\in C^{k+1}(U,E)$.
Then $\lim \Lambda(\gamma_\alpha)=
(\gamma,\eta)=(\gamma,\gamma^{[1]})=\Lambda(\gamma)$.
Thus $\,\im(\Lambda)$ is closed.
\end{proof}
\begin{la}\label{pullback}
Suppose that $Z$ and $Y$ are topological
$\F$-vector spaces, $U\sub Z$ and $V\sub Y$
open subsets, and $f\!: U\to V$ a $C^r_\F$-map.
Then the ``pullback''
\[
C^r(f,E)\!: C^r(V,E)\to C^r(U,E),\;\;\;
\gamma\mto \gamma\circ f
\]
is a continuous $\K$-linear map.
\end{la}
\begin{proof}
Given elements $r\geq s\in \N_0\cup\{\infty\}$, let $i_{s,r}\!: C^r(U,E)\to
C^s(U,E)$ and\linebreak
$j_{s,r}\!: C^r(V,E)\to
C^s(V,E)$ be the respective inclusion maps.
Since $i_{k,\infty}\circ C^\infty(f,E)$\linebreak
$=C^k(f,E)\circ j_{k,\infty}$
if $f$ is of class $C^\infty_\F$,
in view of Remark~\ref{simplobs}\,(a)
the continuity of $C^\infty(f,E)$
follows if we can show that $C^k(f,E)$
is continuous for each $k\in \N_0$.
Thus, we may assume that $r\in \N_0$,
and prove the assertion by induction on~$r$.
The case $r=0$ is a standard fact, see \cite{Eng},
p.\,157, Assertion\,(2).\\[3mm]
{\em Induction step.} Suppose the lemma is correct
for some $r\in \N_0$, and suppose that
$f\!: U\to V$ is of class $C^{r+1}_\F$.
The mapping $i_{0,r+1}\circ C^{r+1}(f,E)
=C^0(f,E)\circ j_{0,r+1}$
being continuous, in view of
Remark~\ref{simplobs}\,(b)
it only remains to show that
\[
\phi\!: C^{r+1}(V,E)\to C^r(U^{[1]},E),\;\;\;
\phi(\gamma):=\left(C^{r+1}(f,E)(\gamma)\right)^{[1]}=
(\gamma\circ f)^{[1]}
\]
is continuous.
By the Chain Rule, we have
\[ \phi(\gamma)=(\gamma\circ f)^{[1]}=C^r(\Phi,E)(\gamma^{[1]}),
\]
where $\Phi\!: U^{[1]}\to V^{[1]}$,
$\Phi(u,y,t):=(f(u),f^{[1]}(u,y,t),t)$
is of class~$C^r$.
By induction, $C^r(\Phi,E)\!: C^r(V^{[1]},E)\to
C^r(U^{[1]},E)$ is continuous,
and also $\psi\!:C^{r+1}(V,E)\to
C^r(V^{[1]},E)$,
$\gamma\mto \gamma^{[1]}$
is continuous (Remark~\ref{simplobs}\,(b)).
Thus $\phi=C^r(\Phi,E)\circ \psi$
is continuous.
\end{proof}
\begin{la}\label{multop}
Let $E$ be a topological $\K$-vector space,
$Z$ a topological $\F$-vector space,
$U\sub Z$ an open subset, and
$f\!: U\to \K$ be a mapping of class~$C^r_\F$.
Then the ``multiplication operator''
\[
m_f\!: C^r(U,E)\to C^r(U,E),\;\;\;
(m_f(\gamma))(x):=f(x)\cdot \gamma(x)
\]
is a continuous $\K$-linear map.
\end{la}
\begin{proof}
Arguing as before,
we find that
it suffices to discuss the case where $r\in \N_0$.
The proof is by induction.
In the following, let us write $m_{f,r}$
for $m_f$, to stress its dependence on~$r$.
Given $r\geq s\in \N_0$, we let
$i_{s,r}\!: C^r(U,E)\to C^s(U,E)$ be the inclusion map.\\[3mm]
{\em Case $r=0$}: Let $K\sub U$ be a compact
subset and $V\sub E$ be an open zero-neighbourhood.
Since $f(K)$ is compact, $f(K)\cdot 0\sub \{0\}\sub V$,
and scalar multiplication is continuous,
there is an open zero-neighbourhood
$W\sub E$ such that $f(K)\cdot W\sub V$.
As a consequence, $m_f\big(\lfloor K,W\rfloor\big)\sub
\lfloor K,V\rfloor$. Being continuous at~$0$
by the preceding, the linear map~$m_{f,0}$ is continuous.\\[3mm]
{\em Induction step.}
Suppose that the assertion of the lemma is
correct for some $r\in \N_0$,
and let $f\!: U\to \K$ be a mapping of class
$C^{r+1}_\F$.
Then $i_{0,r+1}\circ m_{f,r+1}=m_{f,0}\circ i_{0,r+1}$
shows that $i_{0,r+1}\circ m_{f,r+1}$ is a continuous
linear map.
Using that scalar multiplication
$\beta\!: \K\times E\to E$
is a continuous $\K$-bilinear (and thus
$\F$-bilinear) map,
the formula for $\beta^{[1]}$ (see Examples~\ref{linbilin})
combined with the Chain Rule shows that
\begin{eqnarray*}
\lefteqn{(m_{f,r+1}(\gamma))^{[1]}(x,y,t)}\\
&=&
(\beta\circ (f,\gamma))^{[1]}(x,y,t)\\
& = & \beta(f(x),\gamma^{[1]}(x,y,t))+
\beta(f^{[1]}(x,y,t),\gamma(x))
+t\beta(f^{[1]}(x,y,t),\gamma^{[1]}(x,y,t))\\
& = &
f(x)\cdot \gamma^{[1]}(x,y,t)+
f^{[1]}(x,y,t)\cdot \gamma(x)
+t f^{[1]}(x,y,t)\cdot \gamma^{[1]}(x,y,t),
\end{eqnarray*}
whence
\begin{equation}\label{givesmult}
(m_{f,r+1}(\gamma))^{[1]}=
(m_{f\circ \pi+\tau\cdot f^{[1]},r}\circ \phi)(\gamma)
+
(m_{f^{[1]},r}\circ C^r(\pi,E)\circ i_{r,r+1})(\gamma)\, ,
\end{equation}
where $\pi\!: U^{[1]}\to U$, $\pi(x,y,t):=x$
and $\tau\!: U^{[1]}\to\K$,
$\tau(x,y,t):=t$ are smooth and thus $C^r_\F$,
multiplication operators
are denoted in the apparent way,
and $\phi\!: C^{r+1}(U,E)\to C^r(U^{[1]},E)$
denotes the continuous linear map
$\gamma\mto \gamma^{[1]}$.
In view of the induction
hypothesis and Lemma~\ref{pullback}, Equation\,(\ref{givesmult})
shows that $C^{r+1}(U,E)\to C^r(U^{[1]},E)$,
$\gamma\mto (m_{f,r+1}(\gamma))^{[1]}$
is a continuous $\K$-linear map.
By Remark~\ref{simplobs}\,(b),
$m_{f,r+1}$ is continuous. 
\end{proof}
\begin{la}\label{cover}
Let $Z$ be a topological $\F$-vector space,
$U\sub Z$ be an open subset, and
$(U_i)_{i\in I}$ be an open cover of~$U$.
For $i\in I$, let
$\lambda_i\!: U_i\emb U$ be the inclusion map,
and
\[
\rho_i:=C^r(\lambda_i,E)\!:C^r(U,E)\to C^r(U_i,E),
\;\;\;
\rho_i(\gamma):=\gamma|_{U_i}
\]
be the corresponding restriction map.
Then the topology on $C^r(U,E)$ is initial
with respect to the family
$(\rho_i)_{i\in I}$.
\end{la}
\begin{proof}
Arguing as usual, we may assume that~$r$ is finite.
The proof is by induction.
By Lemma~\ref{pullback},
each map $\rho_i$ is continuous linear
and thus the initial topology $\cO_r$ on $C^r(U,E)$
with respect to $(\rho_i)_{i\in I}$
is a (Hausdorff) vector topology on
$C^r(U,E)$ which is coarser than the given topology.
We shall write $\rho_{i,r}$ for $\rho_i$,
to stress its dependence on~$r$.\\[3mm]
{\em The case $r=0$.}
Suppose that $K$ is a compact subset of~$U$,
and $W\sub E$ an open zero-neighbourhood.
Given $x\in K$,
there exists $i\in I$ such that $x\sub U_i$.
Since $K$ is compact Hausdorff,
there exists
a compact neighbourhood $V_x$ of~$x$ in~$K$
such that $V_x\sub K\cap U_i$.
As a consequence, using the compactness
of~$K$ we find finitely many compact subsets $A_1,\ldots, A_n$
of~$K$ covering~$K$, such that,
for each $j=1,\ldots, n$, there
is $i_j\in I$ with $A_j\sub U_{i_j}$.
Then $\lfloor K,W\rfloor\sub C^0(U,E)$
coincides with
$\bigcap_{j=1}^n \rho_j^{-1}(\lfloor A_j, W\rfloor)$,
where $\lfloor A_j,W\rfloor\sub C^0(U_j,E)$.
As a consequence, the vector topology $\cO_0$
on $C^0(U,E)$ is finer than the given topology and thus
coincides with it.\\[3mm]
{\em Induction step.}
Suppose that the assertion of
the lemma is correct for some $r\in \N_0$.
In view of Remark~\ref{simplobs}\,(b),
we have to show that the mappings
\[
\;\;\phi\!: (C^{r+1}(U,E),\cO_{r+1})\to C^0(U,E),\;\;\;
\phi(x)=x\;\;\;\;\;\mbox{and}
\]
\[
\!\!\!\psi\!: (C^{r+1}(U,E),\cO_{r+1})\to C^r(U^{[1]},E),\;\;\;
\psi(\gamma):=\gamma^{[1]}
\]
are continuous, using the usual topology
on the spaces on the right hand side.
Let $j_i\!: C^{r+1}(U_i,E)\to C^0(U_i,E)$
be the inclusion map, which is continuous linear.
As $\rho_{i,0}\circ \phi
=j_i \circ \rho_{i,r+1}$
is continuous, we deduce from the $C^0$-case
of the lemma already proved that $\phi$
is continuous.\\[3mm]
To see that also $\psi$ is continuous, let
$(x,y,t)\in U^{]1[}$
and $\gamma\in C^{r+1}(U,E)$ (see {\bf\ref{gargel}}).
Then
\[
\psi(\gamma)(x,y,t)={\textstyle \frac{1}{t}(\gamma(x+ty)-\gamma(x))}
\]
and thus
\begin{equation}\label{depcts}
\psi(\gamma)|_{U^{]1[}}= 
(m_f\circ (C^r(s,E)-C^r(\pi,E)) \circ \mu_{r,r+1})(\gamma),
\end{equation}
where the inclusion map
$\mu_{r,r+1}\!: (C^{r+1}(U,E),\cO_{r+1})\to (C^r(U,E),\cO_r)=C^r(U,E)$
(induction hypothesis\,!)
is apparently continuous linear, and
$f\!: U^{]1[}\to \K$, $(x,y,t)\mto t^{-1}$
is of class~$C^r_\F$
as well as the mappings
$s\!: U^{]1[}\to U$, $s(x,y,t):=x+ty$ and
$\pi\!: U^{]1[}\to U$, $\pi(x,y,t):=x$.
By Lemma~\ref{multop},
the multiplication operator
$m_f\!: C^r(U^{]1[},E)\to C^r(U^{]1[},E)$
is continuous, and
by Lemma~\ref{pullback}, the mappings $C^r(s,E)$
and $C^r(\pi,E)$ are continuous.
Thus Equation\,(\ref{depcts})
shows that
\begin{equation}\label{eq1}
(C^{r+1}(U,E),\cO_{r+1})\to
C^r(U^{]1[},E), \;\;\;\;\gamma\mto \psi(\gamma)|_{U^{]1[}}
\end{equation}
is a continuous mapping.\\[3mm]
Next, suppose that $p=(x_0,y_0,t_0)\in U^{[1]}$
is given
such that $t_0=0$.
There exists $i\in I$ such that $x_0\in U_i$.
Then $(x_0,y_0,0)\in (U_i)^{[1]}$,
which is an open subset of $U^{[1]}$.
As $\rho_{i,r+1}$
is continuous on $(C^{r+1}(U,E),\cO_{r+1})$
and also $C^{r+1}(U_i,E)\to C^r(U_i^{[1]},E)$,
$\gamma\mto\gamma^{[1]}$ is continuous,
we deduce that the mapping
$(C^{r+1}(U,E),\cO_{r+1}) \to C^r(U_i^{[1]},E)$,
\begin{equation}\label{eq2}
\gamma \mto (\gamma|_{U_i})^{[1]}=\gamma^{[1]}|_{U_i^{[1]}}
=\psi(\gamma)|_{U_i^{[1]}}
\end{equation}
is continuous.
Now $\{U^{]1[}\}\cup\{ U_i^{[1]}\!: i\in I\}$
being an open cover of $U^{[1]}$,
using the induction hypothesis
we deduce from the continuity
of the mappings described in (\ref{eq1})
and (\ref{eq2}) that~$\psi$ is continuous.
\end{proof}
\begin{center}
{\bf The spaces {\boldmath $C^r(M,E)$\/} and mappings between them}
\end{center}
\begin{numba}\label{top2}
Given a
topological $\K$-vector space~$E$
and $\F$-manifold~$M$ of class $C^r_\F$,
modeled on a topological $\F$-vector space~$Z$,
we let $C^r(M,E)$ be the set of all
mappings $\gamma\!: M\to E$ of class~$C^r_\F$.
It is clear that pointwise
operations turn
$C^r(M,E)$ into a $\K$-vector space.
We give $C^r(M,E)$ the initial topology with respect
to the mappings
\begin{equation}\label{dag}
\theta_\kappa\!:
C^r(M,E)\to C^r(V,E),\;\;\;\;
\gamma\mto \gamma|_U\circ \kappa^{-1}\, ,
\end{equation}
where $\kappa\!: U\to V\sub Z$ ranges through the charts
of~$M$.
It is clear that this topology makes
$C^r(M,E)$ a topological $\K$-vector space.
\end{numba}
\begin{numba}\label{newold}
Since an open subset
$U\sub Z$ may be considered as an $\F$-manifold,
we now have two definitions of a
topology on $C^r(U,E)$,
described in {\bf \ref{top1}} and {\bf \ref{top2}}.
As a consequence of Lemma~\ref{cover},
both topologies coincide:
\end{numba}
\begin{la}\label{atlas}
If $\cA$ is an atlas of charts
for~$M$, then the topology
on $C^r(M,E)$ is initial
with respect to the family
$(\theta_\kappa)_{\kappa\in \cA}$.
\end{la}
\begin{proof}
Apparently,
the initial topology $\cO$ with respect to
$(\theta_\kappa)_{\kappa\in \cA}$
is coarser than the given topology on $C^r(M,E)$.
To see that it is also finer, we have to show
that $\cO$ makes $\theta_\eta$ continuous,
for every chart $\eta\!: U\to V$ of~$M$.
For $\kappa\in \cA$, say $\kappa\!: U_\kappa\to
W_\kappa$, define
$V_\kappa:=\eta(U_\kappa\cap U)$.
Then
$(V_\kappa)_{\kappa\in \cA}$
is an open cover of~$V$,
and as $(C^r(M,E),\cO) \to C^r(V_\kappa,E)$,
\[
\gamma \mto \theta_\eta(\gamma)|_{V_\kappa}
=\gamma \circ \eta^{-1}|_{V_\kappa}
=\theta_\kappa(\gamma)\circ \kappa\circ \eta^{-1}|_{V_\kappa}
= (C^r(\kappa\circ \eta^{-1}|_{V_\kappa},E)\circ \theta_\kappa)(\gamma)
\]
is a continuous function of~$\gamma$
by Lemma~\ref{pullback}
and definition of~$\cO$,
for each $\kappa\in \cA$,
we deduce from Lemma~\ref{cover}
that $\theta_\eta$ is continuous,
which completes the proof.
\end{proof}
\begin{rem}\label{compareco}
The topology on $C^0(M,E)=C(M,E)$
just defined
coincides with the compact-open topology.
Indeed, the new topology obviously
is coarser than the compact open topology,
but it is also finer, by the argument used in
the proof of Lemma~\ref{cover}, case $r=0$
(see also Lemma~\ref{flooropen}).
\end{rem}
\begin{la}\label{pb2}
Let $M$ and $N$ be $C^r_\F$-manifolds
modeled on topological $\F$-vector spaces,
$E$ be a topological $\K$-vector space,
and $f\!: M\to N$ be a $C^r_\F$-map.
Then the ``pullback''
\[
C^r(f,E)\!: C^r(N,E)\to C^r(M,E),\;\;\;
\gamma\mto \gamma\circ f
\]
is a continuous $\K$-linear map.
\end{la}
\begin{proof}
It is clear that $C^r(f,E)$ is $\K$-linear.
There exists an atlas $\{\kappa_i\!: i\in I\}$
of charts $\kappa_i\!: U_i\to V_i$
of $M$ such that,
for each $i\in I$, $\,f(U_i)\sub A_i$
for some chart $\phi_i\!: A_i\to B_i$
of~$N$.
Given $i\in I$, consider
$\theta_i\!: C^r(M,E)\to C^r(V_i,E)$,
$\theta_i(\gamma):=\gamma\circ \kappa_i^{-1}$
and $\Theta_i\!:
C^r(N,E)\to C^r(B_i,E)$, $\Theta_i(\gamma):=
\gamma\circ \phi_i^{-1}$.
In view of Lemma~\ref{atlas},
the mapping $C^r(f,E)$ is continuous if and only
if $\theta_i\circ C^r(f,E)$
is continuous for each~$i$.
But $\theta_i\circ C^r(f,E)=
C^r( \phi_i \circ
f|_{U_i}^{A_i}\circ \kappa_i^{-1},E)   \circ \Theta_i$
is a composition of continuous
mappings (see Lemma~\ref{pullback})
\end{proof}
\begin{la}\label{amend}
Let $M$ be a $C^r_\F$-manifold, modeled
on a topological $\F$-vector space, $E$ be a topological $\K$-vector
space, and $(U_i)_{i\in I}$
be an open cover of~$M$.
Then
\[
\rho:=(\rho_i)_{i\in I}\!: C^r(M,E)\to \prod_{i\in I}C^r(U_i,E),\;\;\;
\gamma\mto (\gamma|_{U_i})_{i\in I}
\]
is a topological embedding, whose image is
a closed vector subspace of $\prod_{i\in I}C^r(U_i,E)$.
\end{la}
\begin{proof}
Let $\lambda_i\!: U_i\to M$ be the inclusion maps.
The coordinate functions of $\rho$
are the restriction maps $\rho_i=C^r(\lambda_i,E)\!:
C^r(M,E)\to C^r(U_i,E)$, $\gamma\mto \gamma\circ \lambda_i=\gamma|_{U_i}$,
which are continuous linear by Lemma~\ref{pb2}.
Hence $\rho$ is continuous linear, and apparently it is
injective. Let $\cA$ be the set of all charts of~$M$
whose domain is contained in some~$U_i$.
Then $\cA$ is an atlas for~$M$.
Given $\kappa\in \cA$, say $\kappa\!: U\to V$ with
$U\sub U_i$, we can write
$\theta_\kappa\!: C^r(M,E)\to C^r(V,E)$,
$\gamma\mto \gamma\circ \kappa^{-1}$
as $\theta_\kappa= \Theta_\kappa\circ
\rho_i$, where $\Theta_\kappa\!: C^r(U_i,E)\to C^r(V,E)$,
$\eta\mto \eta\circ \kappa^{-1}$.
As a consequence,
$\theta_\kappa$ is continuous with respect to the topology~$\cO$
induced by $\rho$ on $C^r(M,E)$.
Hence, by Lemma~\ref{atlas}, $\cO$ has to be finer
than the given topology on $C^r(M,E)$.
Being also coarser (since $\rho$ is continuous),
it coincides with the given topology.
Thus $\rho$ is a topological embedding.

Let $F:=\im\,\rho$, and $\wb{F}$ be its closure.
Given $j,k\in I$, and $x\in U_j\cap U_k$,
define $f_{j,k,x}\!: \prod_{i\in I}C^r(U_i,E)\to E$
via $(\gamma_i)_{i\in I}\mto \gamma_j(x)-\gamma_k(x)$.
Then $f_{j,k,x}$ is a continuous linear
map. From $f_{j,k,x}(F)=\{0\}$
we deduce $f_{j,k,x}(\wb{F})\sub \wb{\{0\}}=\{0\}$.
Thus $\gamma_j|_{U_j\cap U_k}=
\gamma_k|_{U_j\cap U_k}$
for all $(\gamma_i)_{i\in I}\in \wb{F}$ and
$j,k\in I$.
As a consequence,
given $(\gamma_i)_{i\in I}\in \wb{F}$,
we can unambiguously define
a mapping $\gamma\!: U\to E$
via $\gamma(x):=\gamma_i(x)$
if $x\in U_i$.
Since $\gamma|_{U_i}=\gamma_i$
is of class~$C^r$ for each $i\in I$,
Lemma~\ref{Crlocal} shows that
$\gamma$ is a mapping of class~$C^r$.
It remains to note that
$(\gamma_i)_{i\in I}=\rho(\gamma)\in F$.
\end{proof}
Various simple observations will be useful.
\begin{la}\label{linearcase}
Suppose that $\lambda\!: E\to F$
is a continuous $\K$-linear map between topological $\K$-vector spaces.
Then
\[
C^r(M,\lambda)\!: C^r(M,E)\to C^r(M,F),\quad
\gamma\mto \lambda\circ \gamma
\]
is a continuous linear map.
\end{la}
\begin{proof}
Given a chart $\kappa\!: U\to V$ of~$M$,
we have $\theta_\kappa^F\circ C^r(M,\lambda)=C^r(V,\lambda)
\circ \theta^E_\kappa$,
where $\theta_\kappa^E\!: C^r(M,E)\to C^r(V,E)$,
$\gamma\mto \gamma\circ \kappa^{-1}$,
and $\theta_\kappa^F\!: C^r(M,F)\to C^r(V,F)$
is defined analogously.
The topology on $C^r(M,F)$ being initial with respect
to the mappings $\theta_\kappa^F$,
it therefore suffices to show that $C^r(V,\lambda)$
is continuous, for any open subset
$V$ of the modeling space of~$M$.
Let $j\in \N_0$ such that
$j\leq r$. For each $\gamma\in C^r(V,E)$,
we have
$(C^r(V,\lambda)(\gamma))^{[j]}=(\lambda\circ\gamma)^{[j]}=
\lambda\circ (\gamma^{[j]})=C(V^{[j]},\lambda)(\gamma^{[j]})$
since $\lambda$ is continuous linear.
Here $C^r(V,E)\to C(V^{[j]},E)$, $\gamma\mto \gamma^{[j]}$
is continuous linear by definition of the
$C^r$-topology,
and $C(V^{[j]},\lambda)\!: C(V^{[j]},E)\to C(V^{[j]},F)$,
$\eta\mto \lambda\circ \eta$
is continuous with respect to the compact-open topologies
by \cite[\S3.4, Assertion\,(1)]{Eng}. The topology on
$C^r(V,F)$ being initial with respect to the maps
$(\sbull)^{[j]}$,
we deduce that $C^r(V,\lambda)$ is continuous.
\end{proof}
If the topology on a topological space~$X$ is initial with respect to a family
of maps into topological spaces whose topology is
again initial with respect to certain families
of maps, then the topology on~$X$ is initial with respect
to the family of composed maps.
This well-known fact will be referred to as
``transitivity of initial topologies''
in the following.
\begin{la}\label{Ctopinit}
Suppose that the topology on $E$ is initial with respect to
a family $(\lambda_i)_{i\in I}$ of $\K$-linear
maps $\lambda_i\!: E\to E_i$
into topological $\K$-vector spaces~$E_i$.
Then the topology on $C^r(M,E)$ is
initial with respect to the family
$(C^r(M,\lambda_i))_{i\in I}$
of the linear mappings $C^r(M,\lambda_i)\!: C^r(M,E)\to C^r(M,E_i)$.
\end{la}
\begin{proof}
The topologies on $C^r(M,E)$ and $C^r(M,E_i)$ are
initial with respect to the mappings
$\theta_\kappa\!: C^r(M,E)\to C^r(V_\kappa,E)$, resp.,
$\theta_{\kappa,i}\!: C^r(M,E_i)\to C^r(V_\kappa,E_i)$
(as in (\ref{dag})),
where $\kappa\!: U_\kappa\to V_\kappa$ ranges through
the set of charts of~$M$.
Hence,
we deduce from $C^r(V_\kappa,\lambda_i)\circ \theta_\kappa
=\theta_{\kappa,i}\circ C^r(M,\lambda_i)$
and the transitivity of initial topologies that
the assertion will hold if we can prove it
when $M$ is an open subset of a topological $\F$-vector space~$Z$
(like the sets $V_\kappa$).
Using Remark~\ref{simplobs}\,(a)
in a similar way, we may furthermore assume that $r\in \N_0$
is finite. Now, the proof is by induction.

For $r=0$, in view of Remark~\ref{compareco}
the assertion is immediate from \cite[La.\,3.4.6]{Eng}.
If $r\in \N$ and the assertion holds when $r$ is replaced
with $r-1$, we recall that, for $M=U\sub Z$, the topology on
$C^r(U,E)$ is initial
with respect to the inclusion map $f\!: C^r(U,E)\to
C(U,E)$ and the map $(\sbull)^{[1]}\!:
C^r(U,E)\to C^{r-1}(U^{[1]},E)$.
Let $f_i\!: C^r(U,E_i)\to C(U,E_i)$
be inclusion. 
Since $C^{r-1}(U^{[1]},\lambda_i)\circ
(\sbull)^{[1]}=(\sbull)^{[1]}\circ C^r(U,\lambda_i)$
and $C(U,\lambda_i)\circ f=f_i\circ C^r(U,\lambda_i)$,
we deduce from the induction hypothesis,
the case $r=0$ and the transitivity of initial
topologies that the topology on $C^r(M,E)$ is indeed initial
with respect to the maps $C^r(U,\lambda_i)$.
\end{proof}
As an immediate consequence, we have:
\begin{la}\label{corprod}
Let $E_1$ and $E_2$ be topological $\K$-vector spaces,
and $\pr_1\!: E_1\times E_2\to E_1$,
$\pr_2\!: E_1\times E_2\to E_2$ be the coordinate projections.
Then
\[
\big(C^r(M,\pr_1),\, C^r(M,\pr_2)\big)
:\,
C^r(M,E_1\times E_2)\to C^r(M,E_1)\times C^r(M,E_2)
\]
is an isomorphism of topological $\K$-vector spaces.\Punkt
\end{la}
Using the latter isomorphism,
we shall frequently identify a function
$\gamma\in C^r(M,E_1\times E_2)$
with its pair of coordinate functions
$(\gamma_1,\gamma_2)$, $\gamma_i:=\pr_i\circ \gamma$.
\begin{prop}\label{globcruc}
Let $E$, $F$, $H$ and $\wt{Z}$ be topological
$\K$-vector spaces, $P\sub H$
be an open subset, $r,k\in \N_0\cup\{\infty\}$,
$\wt{M}$ be a $\K$-manifold
of class $C^{r+k}_\K$ modeled on~$\wt{Z}$,
and $\tilde{f}\!: \wt{M}\times E \times P\to F$
be a mapping of class $C^{r+k}_\K$.
Let $M$ be an $\F$-manifold
of class $C^r_\F$, modeled on a topological
$\F$-vector space~$Z$.
Given a $C^r_\F$-map $\sigma\!: M\to \wt{M}$,
define
\[
f:=\tilde{f}\circ (\sigma\times \id_E\times \id_P)\!: M\times E\times P\to F\,.
\]
Then
\[
\phi\!:
C^r(M,E)\times P \to C^r(M ,F),\;\;\;\;
\phi(\gamma,p):=f(\sbull,p)_*(\gamma)
\]
is a mapping of class $C^k_\K$,
where $f(\sbull,p)_*(\gamma)(x):=f(x,\gamma(x),p)$
for $x\in M$.
\end{prop}
\begin{proof}
Since $C^\infty(M,F)={\pl}_{\ell\in \N_0}\,C^\ell(M,F)$\vspace{-.8mm}
apparently and the inclusion map\linebreak
$C^\infty(M,E)\to
C^\ell(M,E)$ is continuous linear and thus
of class $C^\infty_\K$,
we easily deduce with Lemma~\ref{inpl}
that $\phi$ is of class $C^\infty_\K$ in the case $r=\infty$,
provided the proposition holds for all $r\in \N_0$.
It therefore suffices to consider
the case $r\in \N_0$.
\begin{center}
{\bf Reduction to open subsets of~{\boldmath $Z$ and $\wt{Z}$}}
\end{center}
There is an atlas $\cA=\{\kappa_i\!: i\in I\}$
of charts
$\kappa_i\!: U_i\to V_i\sub Z$
of~$M$ such that $\sigma(U_i)$ is
contained in the domain
$S_i$ of a chart $\tau_i\!: S_i\to R_i\sub \wt{Z}$
of $\wt{M}$.
In view of Lemma~\ref{amend}, Lemma~\ref{pb2},
Lemma~\ref{corestr} and Lemma~\ref{inprod},
the map $\phi$ will be $C^k_\K$ if we can show that
\[
h_i\!: C^r(M,E)\times P\to C^r(V_i,F), \quad (\gamma,p)\mto \phi(\gamma,p)\circ
\kappa_i^{-1}
\]
is of class $C^k_\K$,
for every $i\in I$.
Then $\tilde{f}_i:=\tilde{f}\circ (\tau_i^{-1}\times \id_E\times
\id_P)\!: R_i\times E\times P\to F$
is a $C^{r+k}_\K$-map, and
$\sigma_i:=\tau_i\circ \sigma|_{U_i}^{S_i}\circ \kappa_i^{-1}\!:
V_i\to R_i$ is of class $C^r_\F$.
We set $f_i:=\tilde{f}_i\circ (\sigma_i\times \id_E\times \id_P)\!:
V_i\times E\times P\to F$ and
define
\[
\phi_i\!: C^r(V_i,E)\times P\to C^r(V_i,F),\quad
\phi_i(\gamma,p):=f(\sbull,p)\circ (\id_{V_i},\gamma)\,.
\]
In view of Lemma~\ref{pb2},
the formula $h_i(\gamma,p)=\phi_i(C^r(\kappa_i^{-1},E)(\gamma),p)$
shows that $h_i$ will be $C^k_\K$
if so is $\phi_i$. Replacing $M$ with $V_i$ and $\wt{M}$ with
$R_i$, we may therefore assume that $M$ and $\wt{M}$
are open subsets of $Z$, resp.,
$\wt{Z}$, for the rest of the proof.\\[3mm]
Apparently, it suffices to consider the case
where $k\in \N_0$; the proof is by induction
on~$k$.
\begin{center}
{\bf The case {\boldmath $k=0$}.}
\end{center}
The proof is by induction on~$r$.
If $r=0$, then the topology on $C^0(M,E)$ and $C^0(M,F)$
is the topology of uniform convergence on compact
sets (see {\bf \ref{newold}}).
Let $\gamma \in C(M,E)$,
$p\in P$, $L$ be a compact subset of~$M$,
and $V\sub F$ be an open zero-neighbourhood.
Let $W\sub F$ be an open zero-neighbourhood
such that $W-W\sub V$.
For each $x\in L$,
we find an open neighbourhood $A_x\sub M$
of~$x$ and open zero-neighbourhoods
$B_x\sub E$
and $C_x\sub H$
such that
$p+C_x\sub P$
and
\[
f(y,u,q)-f(x,\gamma(x),p)\in W
\]
for all $y\in A_x$, $u\in \gamma(A_x)+B_x$,
and $q\in p+C_x$.
By compactness, $L\sub \bigcup_{x\in I}A_x$
for some finite subset $I\sub L$.
Then $B:=\bigcap_{x\in I}B_x\sub E$
and $C:=\bigcap_{x\in I}C_x\sub H$
are open zero-neighbourhoods.
Let $\xi\in \gamma+\lfloor L,B\rfloor$
and $q\in p+C\sub P$.
Given $y\in L$,
there is $x\in I$ such that
$y\in A_x$.
Then
\begin{eqnarray*}
f(y,\xi(y),q)-f(y,\gamma(y),p)
\!& = & \!f(y,\xi(y),q)-f(x,\gamma(x),p)
-(f(y,\gamma(y),p)-f(x,\gamma(x),p))\\
& \in & \! W-W\sub V.
\end{eqnarray*}
We have shown that
$\phi(\xi,q)-\phi(\gamma,p)\in \lfloor L,V\rfloor\sub
C(M,F)$
for all $(\xi,q)$
in the open neighbourhood
$(\gamma+\lfloor L,B\rfloor)\times (p+C)$
of $(\gamma,p)$.
Thus $\phi$ is continuous, as required.\\[3mm]
{\em Induction step on~$r$.}
We write $\phi_r$ for~$\phi$,
to emphasize its dependence on~$r$.
Suppose the assertion of the lemma is correct
for $k=0$ and some $r\in \N_0$.
Suppose that the hypotheses of the lemma
are satisfied by $\tilde{f}$ and $\sigma$,
with $r$ replaced by $r+1$.
Let $i\!: C^{r+1}(M,E)\to C(M,E)$ and
$j\!:C^{r+1}(M,F)\to C(M,F)$ be the inclusion maps.
The mapping
\[
C^{r+1}(M,F)\to C(M,F)\times C^r(M^{[1]},F),\;\;\;\;
\gamma\mto (\gamma,\gamma^{[1]})
\]
is an embedding of topological $\K$-vector
spaces by Remark~\ref{simplobs}\,(b).
Thus $\phi_{r+1}$ will be continuous
if we can show that the mappings $j \circ
\phi_{r+1}$ and
\[
\psi\!: C^{r+1}(M,E) \times P \to C^r(M^{[1]},F),\;\;\;
\psi(\gamma,p):=\phi_{r+1}(\gamma,p)^{[1]}\]
are continuous.
We already know that $\phi_0$ is continuous,
whence
$j\circ \phi_{r+1}=\phi_0\circ
(i \times \id_P)$
is continuous.
Recall that $\phi_{r+1}(\gamma,p)(x)=\tilde{f}(\sigma(x),
\gamma(x),p)$ for $\gamma\in C^{r+1}(M,E)$,
$p\in P$ and $x\in M$.
The Chain Rule gives
\begin{eqnarray*}
\psi(\gamma,p)(x,y,t) & = & \phi_{r+1}(\gamma,p)^{[1]}(x,y,t)\\
& = &
\tilde{f}^{[1]}((\sigma(x),\gamma(x),p),\;(\sigma^{[1]}(x,y,t),
\gamma^{[1]}(x,y,t),0),\;t)
\end{eqnarray*}
for all $\gamma\in C^{r+1}(M,E)$, $p\in P$
and $(x,y,t)\in M^{[1]}$. Hence
\begin{equation}\label{jstarvar}
\psi(\gamma,p)
=g(\sbull,p)_*(\gamma\circ \pr_1,\gamma^{[1]})\,,
\end{equation}
where $\pr_1\!: M^{[1]}\to M$, $(x,y,t)\mto x$,
and $g:=\tilde{g}\circ \Big(
(\wh{T}\sigma)\times \id_{E^2}\times \id_P\Big)\!:
M^{[1]}\times E^2\times P\to F$ with
\[
\tilde{g}\!: \wt{M}^{[1]}\times E^2\times P\to F,\quad
\tilde{g}((x,y,t),\, (u,v),\, p):=
\tilde{f}^{[1]}((x,u,p),\, (y,v,0),\, t)
\]
of class $C^r_\K$
and $\wh{T}\sigma\!: M^{[1]}\to \wt{M}^{[1]}$,
$(\wh{T}\sigma)(x,y,t):=(\sigma(x),\sigma^{[1]}(x,y,t),t)$
of class $C^r_\F$.
By the induction hypothesis, the map
\[
C^r(M^{[1]},E^2)\times P\to C^r(M^{[1]},F),\quad
(\kappa,p)\mto g(\sbull,p)_*(\kappa)
\]
is continuous. As
$C^{r+1}(M,E)\to C^r(M^{[1]},E^2)\isom
C^r(M^{[1]},E)^2$,
$\gamma\mto (\gamma\circ \pr_1,\gamma^{[1]})$
is continuous as well (cf.\ Remark~\ref{simplobs}, Lemma~\ref{pullback},
Lemma~\ref{corprod}),
we deduce from (\ref{jstarvar})
that $\psi$ is continuous, and hence so is $\phi_{r+1}$.
\begin{center}
{\bf Induction step on~{\boldmath $k$}.}
\end{center}
Suppose the assertion of the lemma
is correct
for some $k\in \N_0$ and all $r\in \N_0$.
Let $\sigma$ and $\tilde{f}$ be given
which satisfy the hypotheses of the lemma
when $k$ is replaced with $k+1$.
Then $\phi\!: C^r(M,E)\times P\to C^r(M,F)$
is of class $C^k_\K$ (and thus continuous),
by induction.
Given $\gamma,\eta\in C^r(M,E)$,
$p,q\in H$ and $t\in \F$, we clearly have
$(\gamma,p,\eta,q,t)\in
(C^r(M,E)\times P)^{[1]}$
if and only if $(p,q,t)\in P^{[1]}$.
In this case, provided $t\not=0$ we calculate
\begin{eqnarray*}
\lefteqn{{\textstyle \frac{1}{t}(\phi(\gamma+t\eta,p+tq)-\phi(\gamma,p))\,(x)}}
\\
& = &
{\textstyle
\frac{1}{t}(\tilde{f}(\sigma(x),
\gamma(x)+t\eta(x),p+tq)-\tilde{f}(\sigma(x),\gamma(x),p))}\\
& = &
\tilde{f}^{[1]}((\sigma(x),\gamma(x),p),\;
(0,\eta(x),q),\; t)
\end{eqnarray*}
for all $x\in M$.
Hence
\begin{equation}\label{brf2}
\frac{1}{t}(\phi(\gamma+t\eta,p+tq)-\phi(\gamma,p))
=h(\sbull,(p,q,t))_*(\gamma,\eta)
\end{equation}
for all $\gamma,\eta\in C^r(M,E)$ and $(p,q,t)\in P^{[1]}$
such that $t\not=0$, where
$h:=\tilde{h}\circ (\sigma\times \id_{E^2}\times \id_{P^{[1]}})\!:
M\times E^2\times P^{[1]}\to F$
arises from the $C^{r+k}_\K$-map
\[
\tilde{h}\!: \wt{M}\times E^2\times P^{[1]}\to F,\quad
\tilde{h}(z,\, (u,v),\, (p,q,t)):=
\tilde{f}^{[1]}((z,u,p),\, (0,v,q),\, t)\,.
\] 
By the induction hypothesis,
the map
\[
\psi\!: C^r(M,E^2)\times P^{[1]}\to C^r(M,F),\quad
(\kappa,\, (p,q,t))\mto h(\sbull,(p,q,t))_*(\kappa)
\]
is of class $C^k_\K$ (and hence continuous).
In view of (\ref{brf2}), we see that $\phi$
is of class $C^1_\K$, with
$\phi^{[1]}$ given by $\phi^{[1]}((\gamma,p),\, (\eta,q),\, t)
=\psi((\gamma,\eta),\, (p,q,t))$
and thus of class $C^k_\K$.
Hence $\phi$
is of class $C^{k+1}_\K$, as required.
\end{proof}
It would not make sense
to omit the set of parameters~$P$
in the formulation of
Proposition~\ref{globcruc} (hoping to make the proof easier
this way).
In fact, even if~$P$ is a singleton,
a non-singleton set~$P_1$
will occur in the induction step on~$k$
of the preceding proof.
\begin{center}
{\bf The spaces {\boldmath $C^r_K(M,E)$} and mappings
between them}
\end{center}
If $\F$ is locally compact and $K\sub M$ is a compact subset,
we give\label{defcrsk}
\[
C^r_K(M,E):=\{\gamma\in C^r(M,E)\!:
\Supp(\gamma)\sub K\}
\]
the topology induced by $C^r(M,E)$.
The point evaluations $C^r(M,E)\to E$, $\gamma\mto\gamma(x)$
at the elements $x\in M$
being continuous linear maps, $C^r_K(M,E)$ is a closed
vector subspace of $C^r(M,E)$.
In the next proposition, we compile some useful properties
of function spaces. The simple proofs
are given in Appendix~\ref{appfun}.
Only part\,(c) is needed
for the Lie group constructions.
Part\,(d) serves to put our studies in perspective.
Before we can state the proposition,
let us recall various concepts.
\begin{numba}\label{defnkspa}
First, recall that a Hausdorff topological space~$X$
is called a {\em $k$-space\/}
if a subset $U\sub X$ is open precisely if
$U\cap K$ is open in $K$
for every compact subset $K\sub X$.
For example, every metrizable
topological space is a $k$-space.
\end{numba}
\begin{defn}
Let $E$ be a topological $\K$-vector space.
\begin{itemize}
\item[(a)]\label{defnseqc}
$E$ is called {\em sequentially complete\/}
if every Cauchy sequence in $E$ is convergent.
\item[(b)]\label{defnMC}
$E$ is called {\em Mackey complete\/}
if every Mackey-Cauchy sequence in~$E$
is convergent.
Here, a sequence $(x_n)_{n\in \N}$ in~$E$
is called a {\em Mackey-Cauchy sequence\/}\footnote{The
two concepts mainly are
of interest if $\K$ is a valued field.}
if there exists a bounded subset $B\sub E$
and elements $\mu_{n,m}\in \K$
such that $x_n-x_m\in \mu_{n,m}B$ for all $n,m\in \N$
and $\mu_{n,m}\to 0$ in~$\K$ as both
$n,m\to\infty$.
\end{itemize}
\end{defn}
Note that every Mackey-Cauchy sequence also is
a Cauchy sequence; hence every sequentially complete
topological $\K$-vector space is
Mackey complete.
\begin{prop}\label{propprop}
Let $M$ be a $C^r_\F$-manifold, modeled on a topological
$\F$-vector space~$Z$. Let $E$ be a topological $\K$-vector space.
Then the following holds:
\begin{itemize}
\item[\rm (a)]
Assume that
$Z^{[j]}$ is a $k$-space for
all $j\in \N_0$ such that $j\leq r$;
for example, this holds if both $\F$ and~$Z$ are metrizable.
Then $C^r(M,E)$ is complete $($resp.,
sequentially complete, resp., Mackey complete$)$
if $E$ is complete $($resp., sequentially
complete, resp., Mackey complete$)$.
\item[\rm (b)]
If $\K$ is $\R$, $\C$
or an ultrametric field and $E$ is locally convex,
then also $C^r(M,E)$ is locally convex.
\item[\rm (c)]
If $\F$ is locally compact, $E$ is metrizable
and $M$ is a $\sigma$-compact, finite-dimensional
$C^r_\F$-manifold, then $C^r(M,E)$ is metrizable.
\item[\rm (d)]
If $\F\in \{\R,\C\}$, $\K\in \{\F,\C\}$
and $E$ is locally convex, then the topology
on $C^r(M,E)$ is initial with respect
to the family $(D^j)_{r\geq j\in \N_0}$
of maps $D^j\!: C^r(M,E)\to C(T^jM,E)_{c.o.}$,
$\gamma\mto D^j\gamma$,
and hence it is the topology
traditionally considered on $C^r(M,E)$
$($see Appendix~{\rm \ref{appfun}}
for the notations$)$.
\end{itemize}
If $\F$ is locally compact, then
analogous conclusions hold
for the closed vector subspace
$C^r_K(M,E)$ of $\,C^r(M,E)$,
for every compact subset $K\sub M$.\Punkt
\end{prop}
The remainder of this section is
devoted to the following result
(and variants),
which will be needed, for example,
for the discussion of groups of $C^r$-maps.
Until the end of the section,
we now assume that the topological field $\F$ is locally compact.
\begin{prop}\label{pushforw2}
Let $E$, $F$, and $\wt{Z}$ be topological
$\K$-vector spaces, $U\sub E$
an open subset, $r,k\in \N_0\cup\{\infty\}$,
$\wt{M}$ be a $\K$-manifold of class $C^{r+k}_\K$
modeled on $\wt{Z}$,
and
$\tilde{f}\!: \wt{M}\times U \to F$
be a mapping of class $C^{r+k}_\K$.
Let $M$ be a finite-dimensional
$\F$-manifold
of class $C^r_\F$,
and $K\sub M$ be a compact
subset.
Given a mapping $\sigma\!:
M\to \wt{M}$ of class~$C^r_\F$,
we define
$f:=\tilde{f}\circ (\sigma\times \id_U)\!:
M\times U\to F$.
If $K\not=M$, we assume that
$0\in U$ and $f(x,0)=0$
for all $x\in M\,\take\, K$.
Then $C_K^r(M,U):=\{\gamma\in C^r_K(M,E)\!: \gamma(M)\sub U\}$
is an open subset of $C^r_K(M,E)$, and
\[
f_*\!:
C_K^r(M,U) \to C_K^r(M,F),\;\;\;\;
f_*(\gamma)(x):=f(x,\gamma(x))
\]
is a mapping of class $C^k_\K$.
\end{prop}
\begin{cor}\label{Cf2}
Let $E$ and $F$ be topological $\K$-vector
spaces
and $f\!:U\to F$ be a mapping of class~$C^{r+k}_\K$,
defined on an open subset~$U$ of~$E$.
Let $M$ be a finite-dimensional $\F$-manifold of class $C^r_\F$,
and $K\sub M$ a compact subset.
If $K\not=M$, we suppose $0\in U$ and $f(0)=0$.~Then
\[
C_K^r(M,f)\!: C_K^r(M,U)\to C_K^r(M,F),\;\;\;\;
\gamma\mto f\circ \gamma
\]
is a mapping of class~$C^k_\K$.
\end{cor}
\begin{proof}
Let $\wt{M}:=\{0\}$ be a singleton smooth $\K$-manifold,
and $\sigma\!: M\to \wt{M}$, $x\mto 0$,
which apparently is a smooth mapping.
Then $\tilde{g}\!: \wt{M}\times U\to F$,
$\tilde{g}(0,u):=f(u)$
is a mapping of class $C^{r+k}_\K$,
and
$C_K(M,f)=g_*$
for $g:=\tilde{g}\circ (\sigma\times \id_U)$.
By Proposition~\ref{pushforw2},
$g_*$ is $C^k_\K$.
\end{proof}
Instead of proving Proposition~\ref{pushforw2}
directly, we deduce it from a more flexible
technical result (Proposition~\ref{crucial} below),
which shall be re-used repeatedly
afterwards.
For its formulation,
we require certain sets $\lfloor K,U\rfloor_r\,$:
\begin{la}\label{flooropen}
For every compact subset $K\sub M$ and open subset $U\sub
E$, the set
\[
\lfloor K,U\rfloor_r:=\{\gamma\in C^r(M,E)\!:
\gamma(K)\sub U\}
\]
is open in $C^r(M,E)$.
\end{la}
\begin{proof}
There are compact subsets $A_1,\ldots, A_n\sub
K$ which cover $K$, and such that $A_i\sub U_i$
for some chart $\kappa_i\!: U_i\to V_i$
of~$M$ (cf.\ proof Lemma~\ref{cover}).
Set $K_i:=\kappa_i(A_i)$.
Then $\lfloor K_i,U\rfloor\sub C(V_i,E)$
is open by definition of the compact-open topology,
and $h_i\!: C^r(M,E)\to C(V_i,E)$,
$h_i(\gamma):=\gamma\circ \kappa_i^{-1}$
is continuous, for each $i\in \{1,\ldots, n\}$.
Thus $\lfloor K,U\rfloor_r=\bigcap_{i=1}^n
\lfloor A_i,U\rfloor_r
=\bigcap_{i=1}^n h_i^{-1}(\lfloor K_i,U\rfloor)$
is open in $C^r(M,E)$.
\end{proof}
We now formulate the main technical result of this
section. For the moment, only part\,(a)
of the proposition is needed,
but the more general part\,(b)
will become essential
in the proof of Proposition~\ref{nownow} below
(to tackle
also the case of infinite-dimensional~$M$ there).
\begin{prop}\label{crucial}
Let $E$, $F$, $H$ and $Z$ be topological
$\K$-vector spaces, $U\sub E$ and $P\sub H$
be open subsets, $r,k\in \N_0\cup\{\infty\}$,
and $N$ be a $\K$-manifold
of class $C^{r+k}_\K$ modeled on~$Z$.
Let $M$ be an $\F$-manifold
of class $C^r_\F$, modeled on a finite-dimensional topological
$\F$-vector space~$X$, $K\sub M$ be a compact
subset, $Y\sub K$ be an open, non-empty subset
of~$M$, and
$\sigma\!: Y\to N$
be a mapping of class $C^r_\F$.
Define $\lfloor K,U\rfloor_r\sub C^r(M,E)$ as above.
\begin{itemize}
\item[\n (a)]
If
$\,\tilde{g}\!: N \!\times \!U\!\times \!P\to F$
is a $C^{r+k}_\K$-map and
$g:=\tilde{g}\circ (\sigma\!\times\! \id_U\!\times
\!\id_P)\!: Y\!\times\! U\!\times\! P\to F$, then
\[
\lfloor K,U\rfloor_r
\times P \to C^r(Y,F),\;\;\;\;
(\gamma,p)\mto g(\sbull,p)_*(\gamma)
\]
is a mapping of class $C^k_\K$, where $g(\sbull,p)_*(\gamma)(x):=
g(x,\gamma(x),p)$ for $x\in Y$.
\item[\n (b)]
More generally,
let
$\wb{E}$ be a topological
$\K$-vector space,
$\wb{M}$ an $\F$-manifold of class $C^r_\F$,
modeled on a topological $\F$-vector space~$\wb{X}$,
and $\tilde{f}\!: N \times U\times \wb{E}\times P\to F$
be a $C^{r+k}_\K$-map.
Define
$f:=\tilde{f}\circ (\sigma\times \id_U\times \id_{\wb{E}}\times
\id_P)\!: Y\times U\times \wb{E}\times P\to F$.
Then the map
\[
\phi\!:
\lfloor K,U\rfloor_r \times C^r(\wb{M},\wb{E})
\times P \to C^r(Y\times \wb{M},F),\;\;\;\;
\phi(\gamma,\bar{\gamma},p):=f(\sbull,p)_*(\gamma\times \bar{\gamma})
\]
is of class $C^k_\K$,
where $f(\sbull,p)_*(\gamma\times \bar{\gamma})(x,\bar{x})
:=f(x,\gamma(x),\bar{\gamma}(\bar{x}),p)$
for $x\in Y$, $\bar{x}\in \wb{M}$.
\end{itemize}
\end{prop}
\begin{proof}
The proof
is similar to the one of Proposition~\ref{globcruc},
but longer and painfully technical in detail.
We defer it to Appendix~\ref{appcruc}.
\end{proof}
The following lemma helps
to deduce Proposition~\ref{pushforw2}
from Proposition~\ref{crucial}:
\begin{la}\label{restrK}
If $K\sub M$ is a compact subset and
$Y\sub M$ an open subset containing~$K$,
then the restriction map
\[
C^r_K(M,E)\to C^r_K(Y,E),\;\;\;\;
\gamma\mto\gamma|_Y
\]
is an isomorphism of topological $\K$-vector
spaces.
\end{la}
\begin{proof}
The restriction map clearly is an isomorphism
of $\K$-vector spaces.
Let $\cA_1$ be an atlas
for $Y$, and $\cA_2$ an atlas for
$M\,\take\, K$. Then $\cA:=\cA_1\cup\cA_2$
is an atlas for~$M$.
For each $\kappa\in \cA_2$,
we have $\theta_\kappa(\gamma)=0$
for each $\gamma\in C^r_K(M,E)$,
entailing that the initial topology
on $C^r_K(M,E)$ with respect to the mappings
$\theta_\kappa|_{C^r_K(M,E)}$, where $\kappa
\in \cA$, coincides with the
initial topology with respect to
the subset of mappings parametrized by $\kappa\in \cA_1$.
The assertion now readily follows with Lemma~\ref{atlas}.
\end{proof}
\noindent
{\bf Proof of Proposition~\ref{pushforw2}.\/}
Let $Z$ be the modeling space of~$M$.
Since $\F$ is locally compact, the canonical
Hausdorff vector topology on the finite-dimensional $\F$-vector
space~$Z$ is locally compact.
Hence $M$ is a locally compact topological space.
We therefore find a relatively compact
open neighbourhood $Y$ of $K$ in $M$.
The inclusion mappings
$i\!:$\linebreak
$C^r_K(M,E)\to C^r(M,E)$
and $j\!: C^r_K(Y,F)\to C^r(Y,F)$
are $\K$-linear and topological embeddings,
with closed image.
The restriction map
$\rho\!: C^r_K(M,F)\to C^r_K(Y,F)$
is an isomorphism of topological $\K$-vector
spaces by Lemma~\ref{restrK}.
Let $P:=H:=\{0\}$ (zero-dimensional $\K$-vector space),
and define $g\!: Y \times U\times P\to F$, $g(x,y,p):=f(x,y)$. 
Then, by Proposition~\ref{crucial}\,(a), the map
\[
\psi\!: \lfloor \wb{Y}, U\rfloor_r\times P
\to C^r(Y,F),
\;\;\;\;
\psi(\gamma,p):=g(\sbull,p)_*(\gamma)
\]
is of class~$C^k_\K$,
where $\lfloor \wb{Y},U\rfloor_r\sub C^r(M,E)$.
Note that $i(C^r_K(M,U))=\lfloor \wb{Y},U\rfloor_r\cap
C^r_K(M,E)$.
Thus $C^r_K(M,U)$ is open in $C^r_K(M,E)$,
and $i(C^r_K(M,U))\sub \lfloor \wb{Y}, U\rfloor_r$.
Since $j\circ \rho\circ f_*=\psi(\sbull,0)
\circ i|_{C^r_K(M,U)}^{\lfloor \wb{Y},
U\rfloor_r}$ apparently,
we see that $j\circ \rho\circ f_*$
is of class $C^k_\K$,
whence so is~$f_*$, by Lemma~\ref{corestr}.\Punkt
\section{Mapping groups and mapping algebras}\label{secmapgps}
In this section, we discuss mapping groups
and mapping algebras, based
on our studies in Sections~\ref{secveccts}
and \ref{secvecCr}.\\[3mm]
Throughout this section,
$r\in \N_0\cup\{\infty\}$.
If $r=0$,
we let $M$ be any
topological space,
and $\K$ any
topological field.
If $r>0$, we let $\F$ be a locally
compact
topological field, $M$ be a finite-dimensional
$\F$-manifold
of class~$C^r_\F$,
and~$\K$ be a topological field
possessing~$\F$ as a topological subfield.
In either case, we let $K\sub M$ be a compact subset.\pagebreak

\begin{center}
{\bf Mapping Groups}
\end{center}
Given a $\K$-Lie group~$G$,
we consider the set\label{defnmgp}
\[
C_K^r(M,G):=\{\gamma\in C^r(M,G)\!: \gamma|_{M\take K}=1\}
\]
of $G$-valued mappings of class $C^r_\F$ on~$M$
which are identically~$1$ off~$K$.\footnote{To harmonize
notation, we write $C^0(M,G):=C(M,G)$
now also in the case where~$M$ merely is a topological space,
and call continuous mappings $C^0$-maps.} 
It is clear that $C_K^r(M,G)$ is a group
under pointwise multiplication and
inversion. Then
$C_K^r(M,G)$ is a $\K$-Lie group
in a natural way:
\begin{prop}\label{propCXG}
On the group $C_K^r(M,G)$,
there is a uniquely determined
smooth $\K$-manifold structure
with the following
properties:
\begin{itemize}
\item[\n (a)]
it makes $C_K^r(M,G)$
a $\K$-Lie group; and:
\item[\n (b)]
There exists
a chart $\kappa\!: P\to Q$
from an open identity neighbourhood
$P\sub G$ onto an open zero-neighbourhood
$Q\sub L(G)$ such that $\kappa(1)=0$,
$T_1(\kappa)=\id_{L(G)}$,
and such that $C_K^r(M,P):=C_K^r(M,G)\cap P^M$
is open in $C_K^r(M,G)$ and
\[
C_K^r(M,\kappa)\!: C_K^r(M,P)\to C_K^r(M,Q)\sub C_K^r(M,L(G)),\;\;\;\;
\gamma\mto \kappa\circ \gamma
\]
is a diffeomorphism of smooth $\K$-manifolds.
\end{itemize}
Identifying
$L(C_K^r(M,G))$ with $T_0(C_K^r(M,L(G)))=C_K^r(M,L(G))$
via $T_1(C_K^r(M,\kappa))$,
the Lie bracket on $L(C_K^r(M,G))$
corresponds to the mapping
$C_K^r(M,\beta)\!: C_K^r(M,L(G)^2)\isom C_K^r(M,L(G))^2\to
C_K^r(M,L(G))$,
where $\beta\!: L(G)^2\to L(G)$ is
the Lie bracket of~$L(G)$
$($in other words,
$[\gamma,\eta](x)=[\gamma(x),\eta(x)])$.
\end{prop}
{\bf Proof.}
The following proof closely follows
the lines of \cite{GCX}, Section~3,
where only real and complex Lie groups
modeled on locally convex spaces are considered.
We proceed in steps.
\begin{numba}\label{newconstr1}
Let $\phi\!: U_1\to U$ be a chart of~$G$, defined on an
open identity neighbourhood~$U_1$ in~$G$,
with values in an open zero-neighbourhood~$U$ in~$L(G)$,
such that $\phi(1)=0$. Let~$V_1$ be an open, symmetric
identity neighbourhood in~$G$ such that $V_1V_1\sub U_1$,
and set $V:=\phi(V_1)$. Then the mappings
\[
\mu\!: V\times V\to U, \;\;\;\;\mu(x,y):=\phi(\phi^{-1}(x)\cdot
\phi^{-1}(y))
\]
and \hspace{3.7cm}$\iota\!: V\to V$,
$\quad\quad \iota(x):=\phi(\phi^{-1}(x)^{-1})$\\[4mm]
are smooth.
We equip $C_K^r(M,U_1):=\{\gamma\in C_K^r(M,G)\!: \gamma(M)\sub
U_1\}$ with the smooth $\K$-manifold structure
making the bijection
\[
C^r_K(M,\phi)\!: C_K^r(M,U_1)\to C_K^r(M,U),\;\;\;\;
\gamma\mto \phi\circ \gamma
\]
a diffeomorphism of smooth $\K$-manifolds
onto the open subset $C^r_K(M,U)\sub C^r_K(M,L(G))$.
\end{numba}
\begin{numba}\label{nconst2}
Since $C_K^r(M, V)\times C_K^r(M,V)\isom C_K^r(M,V\times V)$
as $C^\infty_\K$-manifolds
(cf.\ Lemma~\ref{corprod}),
and $C_K^r(M,\mu)\!: C_K^r(M,V\times V)\to C_K^r(M,U)$
is $C^\infty_\K$ by Corollary~\ref{Cf}
(resp., Corollary~\ref{Cf2}),
we deduce that the group multiplication of $C_K^r(M,G)$
induces a $C^\infty_\K$-mapping\linebreak
$C_K^r(M,V_1)\times C_K^r(M,V_1)\to C_K^r(M,U_1)$.
Similarly, inversion is $C^\infty_\K$
on $C_K^r(M,V_1)$.
\end{numba}
\begin{numba}\label{nconst3}
Let $\gamma\in C_K^r(M,G)$ now.
As $\gamma(M)\sub \gamma(K)\cup\{1\}$ is compact,
there is an open identity neighbourhood~$W_1\sub V_1$
in~$G$ and an open neighbourhood~$P$ of
$\gamma(M)$ in~$G$ such that $pW_1p^{-1}\sub U_1$
for all $p\in P$. Set $W:=\phi(W_1)$.
As $C_K^r(M, W)$ is open in $C_K^r(M,V)$,
we deduce that $C_K^r(M,W_1)$ is open in
$C_K^r(M,V_1)$. The mapping
$h\!: P\times W_1\to U_1$, $h(p,w):=pwp^{-1}$
being $C^\infty_\K$,
also
$\tilde{f}:=\phi\circ h\circ (\id_P\times \phi^{-1}|_W^{W_1})
\!: P\times W\to U$ is $C^\infty_\K$.
Then clearly the mapping
$f:= \tilde{f}\circ (\gamma \times \id_W)\!: M\times W\to U$,
$f(x,y)=\phi(\gamma(x)\phi^{-1}(y)\gamma(x)^{-1})$
satisfies the hypotheses
of Corollary~\ref{pushforw}
(resp., Proposition~\ref{pushforw2}),
with $k:=\infty$.
We deduce from
Corollary~\ref{pushforw}
(resp., Proposition~\ref{pushforw2})
that the mapping
$f_*\!: C_K^r(M, W)\to C_K^r(M,U)$
is $C^\infty_\K$.
Note that
\[
C^r_K(M,\phi)^{-1}\circ f_*\circ C^r_K(M,\phi)
\big|_{C^r_K(M,W_1)}^{C^r_K(M,W)}\;= \;
I_\gamma\big|_{C^r_K(M,W_1)}\,,
\]
where $I_\gamma\!: C^r_K(M,G)\to C^r_K(M,G)$,
$I_\gamma(\eta):=\gamma\eta\gamma^{-1}$.
Thus $I_\gamma(C^r_K(M,W_1))\sub C^r_K(M,U_1)$
and $I_\gamma\big|_{C^r_K(M,W_1)}^{C^r_K(M,U_1)}$
is $C^\infty_\K$
on the
open identity neighbourhood
$C_K^r(M,W_1)\sub C_K^r(M,V_1)$.
Now Proposition~\ref{locchar}
provides a unique smooth
$\K$-manifold structure on $C_K^r(M,G)$
such that $C_K^r(M,G)$ becomes
a $\K$-Lie group which possesses
$C_K^r(M,V_1)$ as an open submanifold.
\end{numba}
\begin{numba}\label{nconst4}
The Lie group $C_K^r(M,G)$ being modeled on $C_K^r(M,L(G))$,
its Lie algebra can be identified with $C_K^r(M,L(G))$ as
a topological vector space, by means of
$T_1(C^r_K(M,\phi|_{V_1}^V))$.
Let us show that the Lie bracket is the mapping
$C_K^r(M,[.,.])$ on $C_K^r(M,L(G)^2)\isom
C_K^r(M,L(G))^2$ (which is continuous
by Corollary~\ref{Cf}, resp., Corollary~\ref{Cf2}).
To this end, note first that the point evaluation
\mbox{$\pi_x\!:C_K^r(M,G)\to G$,} $\pi_x(\gamma):=\gamma(x)$
is a smooth homomorphism
for each $x\in M$, since $\pi_x\circ C_K^r(M,\phi^{-1}|_V)
=\phi^{-1}|_V\circ \Pi_x|^V_{C_K^r(M,V)}$ is smooth,
using that the point evaluation $\Pi_x\!: C_K^r(M,L(G))\to L(G)$
is a continuous linear map.
As we identify $T_1C_K^r(M,G)$ with $C_K^r(M,L(G))$
by means of $T_1C_K^r(M,\phi|_{V_1}^V)$,
and~$T_1\phi=\id_{L(G)}$ by hypothesis,
we clearly have $L(\pi_x)=T_1(\pi_x)=\Pi_x$.
As $L(\pi_x)$ is a Lie algebra homomorphism,
we deduce that $[\gamma,\eta](x)=[\gamma(x),\eta(x)]$
for all $\gamma,\eta\in C_K^r(M,L(G))$.
The assertion follows.
\end{numba}
\begin{numba}
The asserted uniqueness
of the Lie group structure on
$C_K^r(M,G)$ with the required properties
follows by standard arguments,
using that $C_K^r(M,\kappa_1\circ \kappa^{-1}_2)$
is a diffeomorphism (by Corollary~\ref{Cf},
resp., Corollary~\ref{Cf2})
if
both $\kappa_1$ and $\kappa_2$ are charts
of~$G$ with the described properties
(whose domains coincide, without loss
of generality).
This completes the proof of Proposition~\ref{propCXG}.\Punkt
\end{numba}
\begin{center}
{\bf Mapping Algebras}
\end{center}
Given an associative topological
$\K$-algebra~$A$ (possibly without an identity element),
we let $A_e$ be the associated
unital $\K$-algebra. Thus
$A_e=A\oplus \K e$ as a $\K$-vector space.
We give $A_e$ the product topology,
which makes it a unital, associative
topological $\K$-algebra.
\begin{prop}\label{mapcia}
If $A$ is a continuous inverse algebra
over~$\K$, then also $C^r_K(M,A)_e$
is a continuous inverse $\K$-algebra.
In the special case where~$M=K$
is compact,
also
$C^r(K,A)$ is a continuous inverse $\K$-algebra.
\end{prop}
\begin{proof}
Since we have Corollaries~\ref{Cf}
and \ref{Cf2} at our disposal,
the arguments used
in~\cite{ALG} to prove the analogous
result for locally convex real or
complex continuous inverse
algebras carry over to the present situation.
\end{proof}
\section{Mappings between direct sums}\label{secdirsum}
Throughout this section,
$(\K,|.|)$ denotes a valued field.
We study differentiability properties
of certain mappings between open subsets
of direct sums of topological $\K$-vector
spaces.\\[3mm]
Given a real number $\ve>0$,
we abbreviate $B_\ve(0):=B_\ve^\K(0)=
\{x\in \K\!: |x|<\ve\}$.
We recall that a subset $U\sub E$
of a $\K$-vector space~$E$
is called {\em balanced\/}
if $tU\sub U$ for all $t\in \K$ such that
$|t|\leq 1$. It is called
{\em absorbing\/} if, for $x\in E$,
there exists $\ve>0$
such that $B_\ve(0)\cdot x\sub U$
(see \cite{BTV}, Ch.\,I, \S1, no.\,5).
\begin{numba}\label{basicdsum}
Let $(E_i)_{i\in I}$
be a family of
topological $\K$-vector spaces,
and $E:=\bigoplus_{i\in I}E_i$
be its vector space direct sum.
Let $\cF$ be the set of all
sets~$U$ of the form
\[
{\textstyle U:=\bigoplus_{i\in I}U_i:=E\cap \prod_{i\in I}U_i}
\]
where $U_i$ is an open, balanced
zero-neighbourhood
in~$E_i$.
Then apparently every $U\in\cF$
is a balanced and absorbing subset
of~$E$, and $tU\in \cF$
for each $t\in \K^\times$.
It is also easy to find $V\in\cF$ such
that $V+V\sub U$.
As a consequence,
there is a unique topology on~$E$
turning~$E$ into a topological $\K$-vector
space, and such that $\cF$ is a basis
for the filter of zero-neighbourhoods
of~$E$
(see \cite{BTV}, Ch.\,I, \S1, no.\,5,
Prop.\,4).
Since $\bigcap\cF=\{0\}$,
this vector topology is Hausdorff.
\end{numba}
\begin{numba}
Let $x=(x_i)_{i\in I}\in E$,
and suppose that $U_i$ is an open neighbourhood
of~$x_i$ in~$E_i$, for all $i\in I$.
Then
$U=\bigoplus_{i\in I}U_i:=E\cap \prod_{i\in I}U_i$
is an open neighbourhood
of~$x$ in~$E$.
In fact, let $y=(y_i)_{i\in I}\in U$.
Then $U_i$ being a neighbourhood
of~$y_i$ in~$E_i$,
there exists a balanced, open zero-neighbourhood
$V_i$ in~$E_i$ such that $y_i+V_i\sub U_i$.
Then $V:=\bigoplus_{i\in I}V_i\in \cF$,
and thus $y+V\sub U$ shows that~$U$ is a neighbourhood
of~$y$. We have shown that~$U$ is open.\\[3mm]
In the preceding situation, we call~$U$\label{defnbxxx}
a {\em box neighbourhood\/} of~$x$.
Accordingly, the topology on~$E$
just defined will be called the {\em box topology}.
In this article, direct sums shall always
be equipped with the box topology.
\end{numba}
\begin{numba}\label{obvibox}
It is obvious from the definition that the box
topology on $E=\bigoplus_{i\in I}E_i$
is finer than the topology induced by the product
topology on $\prod_{i\in I}E_i$.
It is also obvious that the direct sum~$E$
induces the product topology on
$\prod_{i\in F}E_i=\bigoplus_{i\in F}E_i\sub E$,
for each finite subset $F\sub I$.
\end{numba}
\begin{numba}\label{lcx2}
Note that if $\K\in \{\R,\C\}$
and each $E_i$ is locally convex, then
also $\bigoplus_{i\in I}E_i$ is locally convex,
because $\bigoplus_{i\in I}U_i$ is convex for
any family $(U_i)_{i\in I}$
of convex open $0$-neighbourhoods $U_i\sub E_i$.
Likewise, if $\K$ is an ultrametric field
with valuation ring~$\bO$
and each $E_i$ is locally convex (see {\bf\ref{deflcx}}),
then
$\bigoplus_{i\in I}E_i$ is locally convex,
because $\bigoplus_{i\in I}U_i$ is an open $\bO$-submodule
of $\bigoplus_{i\in I}E_i$ for
any family $(U_i)_{i\in I}$
of open $\bO$-submodules $U_i\sub E_i$.
\end{numba}
\begin{numba}\label{boxisuniv}
We claim that $E$, equipped with the box topology,
is the direct sum
of the family $(E_i)_{i\in I}$
in the category of topological $\K$-vector spaces,
{\em provided that~$I$ is countable.}
Indeed, this assertion is trivial if $I$ is finite.
Otherwise, we may assume that $I=\N$.
In this case, suppose that $F$ is a topological
$\K$-vector space
and $\lambda_n\!: E_n\to F$ a continuous
linear mapping for each $n\in \N$.
As $E=\bigoplus_{n\in \N}E_n$ as
a $\K$-vector space, there is a uniquely
determined $\K$-linear map
$\lambda\!: E\to F$ such that $\lambda|_{E_n}=\lambda_n$
for each $n\in \N$.
Let $V_0$ be a zero-neighbourhood in~$F$.
Inductively, we find a sequence $(V_n)_{n\in \N}$
of open zero-neighbourhoods
$V_n\sub F$ such that $V_n+V_n\sub V_{n-1}$
for all $n\in \N$.
Then $U:=\bigoplus_{n\in \N}
\lambda_n^{-1}(V_n)$
is an open zero-neighbourhood in~$E$
such that $\lambda(U)\sub \sum_{n\in \N}
V_n\sub V_0$. We deduce that
$\lambda$ is continuous.
\end{numba}
We are primarily interested in the case of
countable direct sums, but our arguments will work
more generally.
\begin{rem}
If~$I$ is uncountable,
then the box topology on~$E$
need not make~$E$
the direct sum of the family $(E_i)_{i\in I}$
in the category of topological $\K$-vector spaces.
For example, if $\K=\R$ (or $\C$)
and $(E_i)_{i\in I}$ is an uncountable
family of non-zero locally convex topological
$\K$-vector
spaces, then the locally convex direct
sum topology is easily seen
to be properly finer than the box topology
(since this is so for
$\R^{(I)}$).\,\footnote{Note that the addition map
$\R^{(I)}\to\R$, $(r_i)_{i\in I}\mto \sum_{i\in I}r_i$
is discontinuous with respect to the box topology,
if $I$ is uncountable.}
\end{rem}
\begin{rem}\label{disclcx}
If $(E_i)_{i\in I}$
is any family of locally convex topological
vector spaces over an ultrametric field~$\K$,
then $E:=\bigoplus_{i\in I} E_i$, equipped with the box topology,
is locally convex (see {\bf\ref{lcx2}}),
and it is
the direct sum of the family
$(E_i)_{i\in I}$
in the category of locally
convex topological $\K$-vector spaces.
To see this,
let $M\sub E$ be an $\bO$-submodule
such that $M_i:=M \cap E_i$ is open in
$E_i$ for each~$i$.
Then $\bigoplus_{i\in I}M_i$ is a box neighbourhood
of~$0$ which is contained in~$M$
as~$M$ is an $\bO$-submodule
(and thus an additive subgroup) of~$E$.
Consequently, $M$ is open in~$E$.
Therefore the box topology is
the finest locally convex vector topology
on the direct sum $E$ which makes all
of the inclusion maps $E_i\to E$
continuous. Hence $E$, with the box topology, has
the universal property of the locally
convex direct sum: A linear map $f\!: E\to F$
in a locally convex space~$F$ is continuous
if and only if $f|_{E_i}\!: E_i\to F$ is
continuous for each $i\in I$.
\end{rem}
It is our goal now to explore
differentiability properties
of mappings between direct sums.
Our discussions will hinge on
symmetry properties
of the maps $f^{[k]}$.
In order to formulate
these symmetry properties conveniently,
we re-order the arguments of $f^{[k]}\!: U^{[k]}\to F$,
by grouping the variables in $E$ together
on the one hand, on the other hand those
in~$\K$.\\[3mm]
Given topological $\K$-vector spaces
$E$ and $F$ and a $C^r$-map $f\!: U\to F$
defined on an open subset of~$E$,
we let $U^{\{0\}}:=U$,
$f^{\{0\}}:=f$,
$U^{\{1\}}:=U^{[1]}$,
$f^{\{1\}}:=f^{[1]}$
and define mappings $f^{\{k+1\}}\!: U^{\{k+1\}}\to F$
for $k\in \N$, $k<r$ on the sets\label{verystrr}
\[
U^{\{k+1\}}:=\{(x,y,u,v,t)\in E^{2^k}\times E^{2^k}\times \K^{2^k-1}\times
\K^{2^k-1}\times \K\!:
(x,u,y,v,t)\in (U^{\{k\}})^{[1]}\}
\]
inductively via
\[
f^{\{k+1\}}(x,y,u,v,t):=
(f^{\{k\}})^{[1]}(x,u,y,v,t)\,.
\]
\begin{la}\label{firstsym}
Given $k\in \N$,
there exist $\ell\in \N$,
$i_\nu\in \N_0$ for $\nu=1,\ldots, 2^k$
and $j_\mu\in \N_0$ for $\mu=1,\ldots, 2^k-1$
with the following properties:
\begin{itemize}
\item[\n (a)]
Given an open subset $U$ of a topological
$\K$-vector space~$E$,
$x=(x_1,\ldots, x_{2^k})\in E^{2^k}$,
$p=(p_1,\ldots, p_{2^k-1})\in \K^{2^k-1}$
and $t\in \K^\times$, we have
$(x,tp)\in U^{\{k\}}$ if and only if
\[
(t^{i_1}x_1,\ldots, t^{i_{2^k}}x_{2^k},
t^{-j_1}p_1,\ldots, t^{-j_{2^k-1}}p_{2^k-1})
\in U^{\{k\}}\,.
\]
\item[\n (b)]
For any topological
$\K$-vector spaces~$E$, $F$,
any $C^k$-map $f\!: U\to F$
defined on an open subset of~$E$,
and
each $(x,p,t)\in E^{2^k}\times \K^{2^k-1}\times \K^\times$
such that $(x,tp)\in U^{\{k\}}$,
we have:
\begin{equation}\label{neuland}
f^{\{k\}}(x,tp)=
t^{-\ell}\cdot
f^{\{k\}}(t^{i_1}x_1,\ldots, t^{i_{2^k}}x_{2^k},
t^{-j_1}p_1,\ldots, t^{-j_{2^k-1}}p_{2^k-1})\,.
\end{equation}
\end{itemize}
\end{la}
\begin{proof}
The proof is by induction on $k\in \N$.
If $k=1$,
let $x_1,x_2\in E$, $p\in \K$,
$t\in \K^\times$.
Then
$(x_1,x_2,tp)\in U^{[1]}$
if and only if
$x_1\in U$ and $x_1+(tp)x_2=x_1+p(tx_2)\in U$,
which holds precisely if
$(x_1,tx_2,p)\in U^{[1]}$. Assume
that $(x_1,x_2,tp)\in U^{[1]}$.
If $p\not=0$,
we have
\[
{\textstyle
f^{[1]}(x_1,x_2,tp)=\frac{1}{tp}
(f(x_1+tpx_2)-f(x_1))=\frac{1}{t}f^{[1]}(x_1,tx_2,p)}.
\]
By continuity,
$f^{[1]}(x_1,x_2,tp)=\frac{1}{t}f^{[1]}(x_1,tx_2,p)$
then also holds if $p=0$.\\[3mm]
{\em Induction step}.
Suppose the lemma is correct for
a certain $k\in \N_0$;
let $\ell$ and $i_\nu$, $j_\mu$
be as described in the lemma.
Suppose further that $f\!: U\to F$ is of class~$C^{k+1}$.
Let $x,y\in E^{2^k}$, $u,v\in \K^{2^k-1}$, $s\in \K$
and $t\in \K^\times$ such that
$(x,y,tu,tv,ts)\in U^{\{k+1\}}$.
If $s\not=0$, we calculate
\begin{eqnarray}
\lefteqn{f^{\{k+1\}}(x,y,tu,tv,ts)}\label{hugearr}\\
& = & (f^{\{k\}})^{[1]}(x,tu;y,tv;ts)\nonumber\\
& = & {\textstyle \frac{1}{ts}
\big(f^{\{k\}}((x,tu)+ts(y,tv))-f^{\{k\}}(x,tu)\big)}\nonumber\\
& = & {\textstyle \frac{1}{ts}
\big(f^{\{k\}}(x+tsy,t^2(\frac{1}{t}u+sv))-f^{\{k\}}
(x,t^2(\frac{1}{t}u))\big)}\nonumber\\
& = & {\textstyle\frac{1}{ts\cdot t^{2\ell}}
[f^{\{k\}} (t^{2i_1}x_1+t^{2i_1+1}sy_1,\ldots,t^{2i_{2^k}}x_{2^k}
+t^{2i_{2^k}+1}sy_{2^k},
\;t^{-2j_1-1}u_1+t^{-2j_1}sv_1,}\nonumber\\
& & \hspace*{8cm}
\ldots, t^{-2j_{2^k-1}-1}u_{2^k-1}+t^{-2j_{2^k-1}}sv_{2^k-1})\nonumber\\
& & -
f^{\{k\}}(t^{2i_1}x_1,\ldots,t^{2i_{2^k}}x_{2^k},
t^{-2j_1-1}u_1,\ldots,t^{-2j_{2^k-1}-1}u_{2^k-1})]\nonumber\\
& = &
{\textstyle
\frac{1}{t^{2\ell+1}}
(f^{\{k\}})^{[1]}(t^{2i_1}x_1,\ldots,t^{2i_{2^k}}x_{2^k},
t^{-2j_1-1}u_1,\ldots, t^{-2j_{2^k-1}-1}u_{2^k-1},}\nonumber\\
& & \hspace*{4cm}
t^{2i_1+1}y_1,\ldots, t^{2i_{2^k}+1}y_{2^k},
t^{-2j_1}v_1,\ldots,t^{-2j_{2^k-1}}v_{2^k-1},\,s)\nonumber\\
& = &
{\textstyle
\frac{1}{t^{2\ell+1}}
f^{\{k+1\}}(t^{2i_1}x_1,\ldots,t^{2i_{2^k}}x_{2^k},
t^{2i_1+1}y_1,\ldots, t^{2i_{2^k}+1}y_{2^k},}\nonumber\\
& & \hspace*{4cm}t^{-2j_1-1}u_1,\ldots, t^{-2j_{2^k-1}-1}u_{2^k-1},
t^{-2j_1}v_1,\ldots,t^{-2j_{2^k-1}}v_{2^k-1},\,s)\,,\nonumber
\end{eqnarray}
where the calculation shows that the argument
of the function in the last line is
in $U^{\{k+1\}}$.
Here, the induction hypothesis was used to obtain
the fourth equality in~(\ref{hugearr}).
If $s=0$, there exists a zero-neighbourhood
$S$ in~$\K$ such that
$(x,y,tu,tv,ts')\in U^{\{k+1\}}$
for all $s'\in S$.
There exists $s'\in S\,\take\,\{0\}$.
By the above, we then have
\begin{eqnarray*}
& & (t^{2i_1}x_1,\ldots,t^{2i_{2^k}}x_{2^k},
t^{-2j_1-1}u_1,\ldots, t^{-2j_{2^k-1}-1}u_{2^k-1},
t^{2i_1+1}y_1,\ldots, t^{2i_{2^k}+1}y_{2^k},\\
& & \hspace*{4cm}
t^{-2j_1}v_1,\ldots,t^{-2j_{2^k-1}}v_{2^k-1},s')
\in (U^{\{k\}})^{[1]}
\end{eqnarray*}
and thus
$(t^{2i_1}x_1,\ldots,t^{2i_{2^k}}x_{2^k},
t^{-2j_1-1}u_1,\ldots, t^{-2j_{2^k-1}-1}u_{2^k-1})
\in U^{\{k\}}$,
entailing that
\begin{eqnarray*}
& & (t^{2i_1}x_1,\ldots,t^{2i_{2^k}}x_{2^k},
t^{-2j_1-1}u_1,\ldots, t^{-2j_{2^k-1}-1}u_{2^k-1},
t^{2i_1+1}y_1,\ldots, t^{2i_{2^k}+1}y_{2^k},\\
& & \hspace*{4cm}t^{-2j_1}v_1,\ldots,t^{-2j_{2^k-1}}v_{2^k-1},0)
\in (U^{\{k\}})^{[1]}
\end{eqnarray*}
and hence
\begin{eqnarray}
& & (t^{2i_1}x_1,\ldots,t^{2i_{2^k}}x_{2^k},
t^{2i_1+1}y_1,\ldots, t^{2i_{2^k}+1}y_{2^k},
t^{-2j_1-1}u_1,\ldots, t^{-2j_{2^k-1}-1}u_{2^k-1},\nonumber\\
& & \hspace*{4cm}t^{-2j_1}v_1,\ldots,t^{-2j_{2^k-1}}v_{2^k-1},\,s)
\in U^{\{k+1\}}\label{slyre}
\end{eqnarray}
with $s=0$.
By continuity, the first and final term in
display\,(\ref{hugearr})
also coincide when $s=0$.
To complete the proof of\,(a),
assume, conversely, that
$x,y\in E^{2^k}$, $u,v\in \K^{2^k-1}$, $s\in \K$
and $t\in \K^\times$ are given such that
(\ref{slyre}) holds.
If $s\not=0$, exploiting the induction
hypothesis we can go backwards
from bottom to top in the display (\ref{hugearr}),
and deduce that
$(x,y,tu,tv,ts)\in U^{\{k+1\}}$.
Arguing as above,
we see that this conclusion remains valid when $s=0$.
Thus (a) and (b) are established also for $k$
replaced with $k+1$.
\end{proof}
The proof shows that we can achieve $\ell=2^k-1$ here.\\[3mm]
We are now ready for the main result of this section.
\begin{prop}\label{mapsdirsums}
Suppose that $(E_i)_{i\in I}$
and $(F_i)_{i\in I}$
are families of topological $\K$-vector spaces
indexed by a set~$I$. Let
$k\in \N_0\cup\{\infty\}$,
and suppose that
$f_i\!: U_i \to F_i$
is a mapping of class $C^k$
for $i\in I$,
defined on an open non-empty subset
$U_i$ of~$E_i$.
Suppose that there is a finite subset $J\sub
I$ such that $0\in U_i$ and
$f_i(0)=0$, for all $i\in I\,\take\, J$.
Then $U:=\bigoplus_{i\in I}U_i$ is an open
subset of $E:=\bigoplus_{i\in I}E_i$,
and
\[
f:=\bigoplus_{i\in I}f_i\!: U \to F,\;\;\;\;
f((x_i)_{i\in I}):=(f_i(x_i))_{i\in I}
\]
is a mapping of class $C^k$ into
$F:=\bigoplus_{i\in I} F_i$.
For each $j\in \N$ such that $j\leq k$,
identifying $E^{2^j}$ with $\bigoplus_{i\in I}E_i^{2^j}$
in the natural way,
we have
\begin{equation}\label{obviousat2nd}
U^{\{j\}}=\{((x_i)_{i\in I},p)\in E^{2^j}\times \K^{2^j-1}\!:
\;\;(\forall i\in I)\;
(x_i,p)\in U_i^{\{j\}}\},\;\;\;\;\mbox{and}
\end{equation}
\[
f^{\{j\}}((x_i)_{i\in I},p)=(f_i^{\{j\}}(x_i,p))_{i\in I}\,.
\]
\end{prop}
\begin{proof}
We may assume that $k\in \N_0$;
the proof is by induction on~$k$.\\[2mm]
{\em The case $k=0$.}
Let $x=(x_i)_{i\in I}\in U$
and $V$ be a neighbourhood of $f(x)$.
Then $V$ contains a box-neighbourhood
$B=\bigoplus_{i\in I}V_i$ of $f(x)$, where $V_i$ is
an open neighbourhood of~$f_i(x_i)$.
As $f^{-1}(V)$
contains the box-neighbourhood $f^{-1}(B)=\bigoplus_{i\in I}
f_i^{-1}(V_i)$
of~$x$, the set $f^{-1}(V)$ is a neighbourhood of~$x$.
We have shown that $f$ is continuous at~$x$.\\[3mm]
{\em Induction step.}
Suppose the assertion holds for
a given $k\in \N_0$,
and suppose that each $f_i$ is a mapping of class
$C^{k+1}$. Then $f$ is of class~$C^k$,
and $f^{\{k\}}((x_i)_{i\in I},p)=(f_i^{\{k\}}(x_i,p))_{i\in I}$.
Equation\,(\ref{obviousat2nd}) holds for~$j\leq k$
by induction and thus also for $j=k+1$,
as an immediate consequence of the definitions.
We claim that $f^{\{k\}}$ is of class~$C^1$.
Let $x=(x_i)_{i\in I}$, $y=(y_i)_{i\in I} \in E^{2^k}\isom \bigoplus_{i\in I}
E_i^{2^k}$,
$u, v \in \K^{2^k-1}$
and $t\in \K$ such that
$(x,u,y,v,t)\in (U^{\{k\}})^{[1]}$.
If $t\not=0$, we have, by induction,
\begin{eqnarray*}
{\textstyle
\frac{1}{t} (f^{\{k\}}(x+ty,u+tv)-f^{\{k\}}(x,u))} & = &
{\textstyle
(\frac{1}{t}f_i^{\{k\}}(x_i+ty_i,u+tv)-f^{\{k\}}_i(x_i,u)))_{i\in I}}\\
& = & (f_i^{\{k+1\}} (x_i,y_i,u,v,t))_{i\in I}\,.
\end{eqnarray*}
Thus $f^{\{k\}}$ will be of class $C^1$
if we can show that
the mapping
\[
(U^{\{k\}})^{[1]}\to F,\quad
(x,u,y,v,t)\mto
(f_i^{\{k+1\}}(x_i,y_i,u,v,t))_{i\in I}
\]
is continuous,
or, equivalently, that
\[
g\!:
U^{\{k+1\}}\to F,\;\;\;
(x,y,u,v,t)\mto (f_i^{\{k+1\}}(x_i,y_i,u,v,t))_{i\in I}
\]
is continuous---this is our goal now.
We have
$\{0\}\times \K^{2^j-1}\sub U_i^{\{j\}}$
for all $i\in I\,\take\, J$
and
all $j\in \N$ such that $j\leq k+1$,
and
\begin{equation}\label{inotI}
f^{\{j\}}_i(0,p)=0 \quad \mbox{for all $\,p\in \K^{2^j-1}$,}
\end{equation}
by a simple induction.
Let $\bar{x}=(\bar{x}_i)_{i\in I}
\in E^{2^{k+1}}\isom \bigoplus_{i\in I}
E_i^{2^{k+1}}$,
$\bar{p}=(\bar{p}_\nu)_{\nu=1}^{2^{k+1}-1}\in \K^{2^{k+1}-1}$
such that
$(\bar{x}, \bar{p})\in U^{\{k+1\}}$.
Pick a real number
$r> \|\wb{p}\|_\infty$.
There is a finite subset $J_0\sub I$ such that
$J\sub J_0$ and such that $\bar{x}_i=0$
for all $i\in I\,\take\, J_0$.
Let $W$ be an open neighbourhood
of $g(\bar{x}, \bar{p})$
in~$F$; we may assume that $W=\bigoplus_{i\in I}W_i$,
where $W_i$ is an open neighbourhood of
$f^{\{k+1\}}_i(\bar{x}_i,\bar{p})$
in~$F_i$. For $i\in I\,\take\, J_0$,
we may assume that the zero-neighbourhood~$W_i$
is balanced.

Let $\ell\in \N$, $i_\mu\in \N_0$
for $\mu=1,\ldots, 2^{k+1}$
and $j_\nu\in \N_0$ for $\nu=1,\ldots, 2^{k+1}-1$
be as in the $C^{k+1}$-case of
Lemma~\ref{firstsym}.

For each $i\in I\,\take\, J_0$,
there exists $\ve_i>0$ and an
open balanced zero-neighbourhood $V_i\sub E_i$
such that
\[
f_i^{\{k+1\}}(V_i^{2^{k+1}}\times B_{\ve_i}(0)^{2^{k+1}-1})
\sub W_i\,.\]
There exists $\tau_i\in \K^\times$
such that $|\tau_i|>\max\{1,\frac{r}{\ve_i}\}$;
set $A_i:=\prod_{\mu=1}^{2^{k+1}}
\tau_i^{-i_\mu} V_i$.
Holding $i\in I\,\take\, J_0$ fixed for the moment,
let us write $\tau:=\tau_i$, for convenience.
For all $x=(x_\mu)_{\mu=1}^{2^{k+1}}\in A_i$ and
$p=(p_\nu)_{\nu=1}^{2^{k+1}-1}\in B_r(0)^{2^{k+1}-1}$,
we have
$\tau^{i_\mu}x_\mu\in V_i$ for $\mu=1,\ldots, 2^{k+1}$
and $|\tau^{-j_\nu-1}p_\nu|=|\tau|^{-j_\nu-1}\cdot |p_\nu|
<\frac{r}{|\tau|}<\ve_i$
for $\nu=1,\ldots, 2^{k+1}-1$,
{\em i.e.},
$\tau^{-j_\nu-1}p_\nu\in B_{\ve_i}(0)$.
Thus Lemma~\ref{firstsym}\,(a)
shows that $(x,p)\in U_i^{\{k+1\}}$
and
\begin{eqnarray*}
\lefteqn{f_i^{\{k+1\}}(x,p)-f_i^{\{k+1\}}(\wb{x}_i,\wb{p})}\\
& = & f_i^{\{k+1\}}(x,p)
= f_i^{\{k+1\}}(x,\tau(\tau^{-1}p))\\
& = &
{\textstyle \frac{1}{\tau^\ell}
f_i^{\{k+1\}}(\tau^{i_1}x_1,
\ldots,\tau^{i_{2^{k+1}}}x_{2^{k+1}},
\tau^{-j_1-1}p_1,\ldots,\tau^{-j_{2^{k+1}-1}-1}p_{2^{k+1}-1})
\in \frac{1}{\tau^\ell}W_i\sub W_i\,,}
\end{eqnarray*}
using (\ref{inotI}) to pass to the second line
and Lemma~\ref{firstsym}\,(b) to pass
to the third.
For each $i\in J_0$, on the other hand,
by continuity of $f_i^{\{k+1\}}$
there exists an open neighbourhood
$A_i\sub E_i$ of $\bar{x}_i$
and an open neighbourhood $Z_i$
of $\wb{p}$ in $\K^{2^{k+1}-1}$
such that $A_i\times Z_i\sub U_i^{\{k+1\}}$
and $f_i^{\{k+1\}}(A_i\times Z_i)\sub W_i$.
Then $Z:=B_r(0)^{2^{k+1}-1}\cap\bigcap_{i\in J_0}Z_i$
is an open neighbourhood of~$\bar{p}$
in $\K^{2^{k+1}-1}$. Let $A:=\bigoplus_{i\in I}A_i$.
Then
$A\times Z$
is an open neighbourhood
of $(\bar{x},\bar{p})$
in $U^{\{k+1\}}$ such that
\[
g(x,p) \in W\;\;\;
\mbox{for all $\, (x,p)\in A\times Z$.}
\]
We have shown that $g$ is continuous at~$(\bar{x},\bar{p})$.
Thus
$f^{\{k\}}$ is of class $C^1$
and hence
also $f^{[k]}$ is of class~$C^1$
(by the Chain Rule).
Hence $f$ is of class~$C^{k+1}$.
Furthermore, $f^{\{k+1\}}=g$ is of the asserted form.
\end{proof}
Results analogous to Proposition~\ref{mapsdirsums}
for mappings between
locally convex direct sums
of real or complex locally convex spaces
have first been established in~\cite{MEA};
the proofs are considerably easier in that case.

\begin{center}
{\bf Analogues for functions involving parameters}
\end{center}
When the ground field $\K$
is locally compact, Proposition~\ref{mapsdirsums}
can be generalized to functions
involving parameters (and its proof simplifies
substantially).
\begin{prop}\label{sumspara}
Let $(\K,|.|)$
be a valued field,
$P\not=\emptyset$ be a locally compact topological space,
$(E_i)_{i\in I}$
and $(F_i)_{i\in I}$
be families of topological $\K$-vector spaces
indexed by a set~$I$, and
$(f_i)_{i\in I}$ be a family of continuous mappings
$f_i\!: U_i\times P \to F_i$,
where
$U_i$ is a non-empty open subset of~$E_i$.
Suppose that there is a finite subset $J\sub
I$ such that $0\in U_i$ and
$f_i(0,p)=0$, for all $i\in I\,\take\, J$
and $p\in P$.
Then $U:=\bigoplus_{i\in I}U_i$ is an open
subset of $E:=\bigoplus_{i\in I}E_i$,
and
\[
f:=U \times P \to F,\;\;\;\;
f((x_i)_{i\in I},p):=(f_i(x_i,p))_{i\in I}
\]
is a continuous map into $F:=\bigoplus_{i\in I} F_i$.
If $\,\K$ is locally compact here,
$P$ an open subset of a finite-dimensional $\K$-vector
space~$Z$, and if there exists
$k\in \N_0\cup\{\infty\}$
such that $f_i$ is of class $C^k$ for all $i\in I$,
then also $f$ is of class~$C^k$.
\end{prop}
\begin{proof}
We may assume that $k\in \N_0$;
the proof is by induction on~$k$.\\[2mm]
{\em The case $k=0$.}
Let $x=(x_i)_{i\in I}\in U$,
$p\in P$,
and $V$ be a neighbourhood of $f(x,p)$ in~$F$.
Then $V$ contains a box-neighbourhood
$B=\bigoplus_{i\in I}V_i$ of $f(x,p)$, where $V_i$ is
an open neighbourhood of~$f_i(x_i,p)$ in~$F_i$.
There is a finite subset $J_0\sub I$
such that $J\sub J_0$
and such that $x_i=0$
for all $i\in I\setminus J_0$.
For each $i\in J_0$,
we find a compact neighbourhood
$K_i$ of $p$ in $P$ and
an open neighbourhood
$W_i\sub U_i$ of $x_i$
such that
$f_i(W_i \times K_i)\sub V_i$.
Then $K:=\bigcap_{i\in J_0} K_i$
is a compact neighbourhood of $p$ in~$P$.
For each $i\in I\setminus J_0$,
we have
\[
f_i(\{0\}\times K)=\{0\}\sub V_i\,.
\]
Using the compactness of~$K$,
we therefore
find an open zero-neighbourhood $W_i\sub E_i$
such that $f_i(W_i\times K)\sub V_i$.
Then $W:=\bigoplus_{i\in I}W_i\sub U$ is an open neighbourhood
of~$x$, and $f(W\times K)\sub B$
since $f_i(W_i\times K)\sub V_i$ for all~$i$.
Thus $f$ is continuous.\vspace{2mm}

{\em Induction step.}
Let $\K$ be locally compact now,
$k\in \N$, and suppose that the assertion
of the proposition holds when $k$ is replaced
with $k-1$. Let $P\sub Z$ and $C^k$-maps
$f_i\!: U_i\times P\to F_i$ be given.
Then $f\!: U\times P\to F$ is a $C^{k-1}$-map
(and thus continuous), by induction.
As
\[
(U\times P)^{[1]}=\left\{(x,p,y,q,t)\in (E\times Z)^2\times\K\!:
\,\mbox{$(x_i,p,y_i,q,t)\in (U_i\times P)^{[1]}$
for all $i\in I$}\,\right\}
\]
clearly (where $x=(x_i)_{i\in I}$,
$y=(y_i)_{i\in I}$),
we can define a mapping
\[
g\!: (U\times P)^{[1]}\to F,\quad
g(x,p,y,q,t)
:=(f_i^{[1]}(x_i,p, y_i,q,t))_{i\in I}\,.
\]
Let us show that $f$ is $C^1$, with $f^{[1]}=g$
of class $C^{k-1}$.
Since
\[
{\textstyle \frac{1}{t}}(f(x+ty,p+tq)-f(x,p))=\left(
{\textstyle \frac{1}{t}}
(f_i(x_i+ty_i,p+tq)-f_i(x_i,p))\right)_{i\in I}=g(x,p,y,q,t)
\]
for all $(x,p,y,q,t)\in (U\times P)^{[1]}$ such that
$t\not=0$,
it suffices to show that~$g$ is of class $C^{k-1}$.
Now $f$ being of class $C^{k-1}$, the map
$g$ is $C^{k-1}$
on the set
$\{(x,p,y,q,t)\in (U\times P)^{[1]}\!: t\not=0\}$.
It therefore only remains to show that $g$
is $C^k$ on some open neighbourhood of $(\bar{x},\bar{p},\bar{y},\bar{q},0)$,
for all $\bar{x}=(\bar{x}_i)\in U$, $\bar{p}\in P$, $\bar{y}=(\bar{y}_i)\in E$,
and $\bar{q}\in Z$.
For each $i\in I$, we find an open,
balanced zero-neighbourhood
$W_i\sub E_i$ such that $\bar{x}_i+W_i+W_i\sub U_i$.
Then $A_i:=\bar{x}_i+W_i\sub U_i$.
Since $\bar{y}_i=0$ for all but finitely many~$i$,
we find $r\in \; ]0,1]$ such that
$t \bar{y}_i\in W_i$ for all $i\in I$
and $t\in \K$ such that $|t|\leq r$.
Pick $\rho\in \K^\times$ such that $|\rho|\leq r$;
then $B_i:=\rho^{-1} W_i$ is an open neighbourhood
of~$\bar{y}_i$, for all $i\in I$.
There are
$s\in \,]0,r|\rho|]$
and open neighbourhoods $R\sub P$
of~$\bar{p}$ and $S\sub Z$ of~$\bar{q}$, 
such that $R\times S\times B_s(0)\sub P^{[1]}$,
where $B_s(0)\sub \K$.
Then also $A_i\times B_i\times B_s(0)\sub U_i^{[1]}$
for each $i\in I$ and hence
\[
A \times R \times B \times S \times B_s(0)\sub (U\times P)^{[1]}\,,
\]
where $A:=\bigoplus_{i\in I}A_i\sub E$
and $B:=\bigoplus_{i\in I} B_i\sub E$.
Let $Q:=R\times S\times B_s(0)$; then
\[
h_i\!: (A_i\times B_i)\times Q\to F,\quad h_i(x_i,y_i,p,q,t):=
f_i^{[1]}(x_i,p,y_i,q,t)
\]
is a $C^{k-1}$-map, for each $i\in I$.
Furthermore, $h_i|_{\{0\}\times Q}=0$ for all
$i\in I\setminus J$.
Define
\[
h\!: (A\times B) \times Q\to F,\quad
h(x,y,p,q,t):=(h_i(x_i,y_i,p,q,t))_{i\in I}=g(x,p,
y,q,t)\,.
\]
Here $A\times B\sub E\times E\isom \bigoplus_{i\in I}(E_i\times E_i)$,
and $Q$ is an open subset of the
finite-dimensional $\K$-vector space
$Z\times Z\times \K$, which is locally compact since
so is~$\K$.
By the induction hypothesis, $h$ is of class $C^{k-1}$.
Hence $g$ is $C^{k-1}$ on the open neighbourhood
$A\times R\times B\times S\times B_s(0)$
of $(\bar{x},\bar{p},\bar{y},\bar{q},0)$, which completes the proof.
\end{proof}
\begin{rem}
Let $\K$ be $\R$ or a local field,
$n\in \N$, and
$G:=\Diff^\infty_c(\K^n)\dsemi \Aff(\K^n)$
be the group of all smooth diffeomorphisms
of $\K^n$ which coincide with an affine isomorphism
of~$\K^n$ outside some compact set.
Using Proposition~\ref{sumspara},
it is possible to
make $G$ a $\K$-Lie group
modeled on the topological $\K$-vector space
$C^\infty_c(\K^n,\K^n)\times \aff(\K^n)$.
We omit the proof. A more profound application of Proposition~\ref{sumspara}
will be given in Theorem~\ref{OmegaP} below.
\end{rem}
\section{Weak direct products of Lie groups}\label{secweakprod}
The considerations in Section~\ref{secdirsum}
make it possible to construct Lie group
structures on weak direct products
of Lie groups.
\begin{prop}\label{propweakdp}
Let $(G_i)_{i\in I}$
be a family of Lie groups
over a valued field $(\K,|.|)$.
Then there exists a unique $\K$-Lie group
structure on
\[
{\textstyle \prod^*_{i\in I}G_i\!:=
\{
(g_i)_{i\in I}\in \prod_{i\in I}G_i\!:
\;
\mbox{$g_i=1$ for all but finitely many~$i$}\,\}\, ,}
\]
modeled on $\bigoplus_{i\in I}L(G_i)$,
equipped with the box topology,
such that, for certain charts
$\kappa_i\!: U_i\to V_i\sub L(G_i)$
of $G_i$ defined on an identity neighbourhood
$U_i\sub G_i$ and
taking $1$ to~$0$, the mapping
\[
{\textstyle
\bigoplus_{i\in I} V_i\to \prod_{i\in I}^*G_i,\;\;\;\;
(x_i)_{i\in I}\mto
(\kappa_i^{-1}(x_i))_{i\in I}}
\]
is a diffeomorphism of smooth $\K$-manifolds
onto an open subset of $\prod^*_{i\in I} G_i$.
\end{prop}
\begin{proof}
Using
Proposition~\ref{mapsdirsums}
instead of \cite[Prop.\,7.1]{MEA},
the proof of \cite[Prop.\,7.3]{MEA}
(devoted to weak direct products
of real or complex Lie groups
modeled on locally convex spaces)
carries over to the present situation
(see \cite[Thm.\,18.1]{INF} for further details).
\end{proof}
The following observations are immediate
from the construction of the Lie group
structure on weak direct products
and obvious properties of direct sums
of topological vector spaces:
\begin{la}\label{rearrange}
Let $\K$ be a valued field.
\begin{itemize}
\item[\rm (a)]
If $(G_i)_{i\in I}$ is a family of $\K$-Lie
groups and $(H_i)_{i\in I}$
a family of open subgroups $H_i\sub G_i$,
then $\prod_{i\in I}^*H_i$ is an open subset
of $\prod_{i\in I}^*G_i$.
The smooth manifold structure
making $\prod_{i\in I}^*H_i$ an open submanifold of
$\prod_{i\in I}^*G_i$ and the manifold structure
on the weak direct product of Lie groups
$\prod_{i\in I}^*H_i$ coincide.
\item[\rm (b)]
Assume that $I$ is a set, $J_i$ a finite set for each
$i\in I$, and $K:=\{(i,j)\!: i\in I, j\in J_i\}$.
Let $(G_{ij})_{(i,j)\in K}$ be a family
of $\K$-Lie groups.
Then the mapping
\[
{\textstyle \prod_{(i,k)\in K}^*}G_{ij}
\to
{\textstyle \prod_{i\in I}^*}\Big({\textstyle \prod_{j\in J_i} \,G_{ij}}
\Big), \,\quad (g_{ij})_{(i,j)\in K}\mto
\big( (g_{ij})_{j\in J_i} \big)_{i\in I}
\]
is an isomorphism of $\K$-Lie groups.
\item[\rm (c)]
If $(G_i)_{i\in I}$
and $(H_j)_{j\in J}$
are families of $\K$-Lie
groups, $\pi\!: J\to I$ is a bijection
and $\beta_j\!: G_{\pi(j)}\to H_j$
an isomorphism of $\K$-Lie groups for
each $j\in J$, then also
the map
\[
{\textstyle \prod_{i\in I}^*}\, G_i\to
{\textstyle \prod_{j\in J}^*}\, H_j\,, \quad
(g_i)_{i\in I}\mto \bigl( \beta_j(g_{\pi(j)}) \bigr)_{j\in J}
\]
is an isomorphism of $\K$-Lie groups.\Punkt
\end{itemize}
\end{la}
\section{\!\!Spaces of test functions and mappings
between them}\label{secmaptf}
In this section (and in Section~\ref{secalmloc}),
we study differentiability
properties of mappings between spaces of vector-valued
test functions on paracompact
finite-dimensional manifolds
over locally compact ground fields.
First, we collect some
properties of such
manifolds.
\begin{center}
{\bf Paracompact finite-dimensional manifolds over locally compact fields}
\end{center}
Throughout this subsection,
$\F$ is a (non-discrete),
locally compact topological field,
and $r\in \N_0\cup \{\infty\}$.\\[3mm]
Paracompact manifolds over locally compact
fields are amenable to investigation due
to the following well-known fact
(see \cite[Thm.\,5.1.27]{Eng}):
For every paracompact, locally compact topological space~$X$,
there exists a cover $(X_i)_{i\in I}$
of $X$ by mutually disjoint, $\sigma$-compact,
open (and closed) subsets $X_i\sub X$
(and thus $X=\coprod_{i\in I}X_i$).
As a special case, we obtain:
\begin{la}\label{paracamen}
Every paracompact, finite-dimensional
$C^r_\F$-manifold $M$
is a disjoint union $M=\coprod_{i\in I}M_i$
of a family $(M_i)_{i\in I}$
of $\sigma$-compact, open $($and closed\/$)$ submanifolds~$M_i\sub M$.\Punkt
\end{la}
\begin{numba}\label{defnuing}
If $\F$ is a local field, we fix the following
notation:
$|.|$ is an ultrametric absolute value on~$\F$
defining its topology,
$\bO$ the maximal compact subring of~$\F$,
and $\pi\in \F^\times$ a uniformizing element
(thus $|\pi|<1$ and $|\F^\times|=\langle|\pi|\rangle$).
Given $d\in \N$, we
let $\|\sbull\|_\infty$ be the maximum norm on~$\F^d$.
Given $a\in \K^d$ and $\ve>0$,
$B_\ve(a):=\{y\in \F^d\!: \|y-a\|_\infty<\ve\}$
denotes the ball with respect to the maximum norm.
Then $B:=\bO^d$ is an open compact $\bO$-submodule
of~$\F^d$, and it is easy to see that
each ball $B_\ve(a)$ is of the form $a+\pi^kB$
for some $k\in \Z$ and thus $C^\infty_\F$-diffeomorphic
to~$B$.
If~$M$ is a $d$-dimensional
$\F$-manifold of class $C^r_\F$,
we call an open subset of~$M$ a {\em ball\/}
if it is $C^r_\F$-diffeomorphic to~$B$.
It is clear that every point $x\in M$ is contained in some
ball. To avoid misunderstandings, the balls $B_\ve(a)\sub \F^d$ will
occasionally be called {\em metric balls\/} now.
\end{numba}
The following lemma assembles various useful
facts concerning paracompact manifolds
over local fields (cf.\ also \cite{Lu4}).
\begin{la}\label{onlyopen}
Let $\F$ be a local field,
$r\in \N_0\cup\{\infty\}$,
and $M$ be an $C^r_\F$-manifold over~$\F$,
of positive, finite dimension $d\in \N$.
Then the following holds:
\begin{itemize}
\item[\rm (a)]
If $M$ is $\sigma$-compact, then
$M$ is $C^r_\F$-diffeomorphic
to an open subset $U\sub \F^d$.
\item[\rm (b)]
If $M$ is paracompact, then $M$
is a disjoint union $M=\coprod_{i\in I}B_i$
of a family $(B_i)_{i\in I}$ of compact and
open balls $B_i\sub M$.
\end{itemize}
\end{la}
\begin{proof}
(a) Since $M$ is $\sigma$-compact,
there exists a sequence $(B_k)_{k\in \N}$
of balls covering~$M$.
We set $J_1:=\{B_1\}$.
Suppose that we have found
an open cover $J_k$
of $\bigcup_{j=1}^k B_j$ by disjoint balls
for $k=1,\ldots, n$,
such that $J_1\sub J_2\sub\cdots\sub J_n$.
Let $\psi\!: B_{n+1}\to B$ be a $C^r_\F$-diffeomorphism
onto $B:=\bO^d$.
Then $R:=B_{n+1}\,\take\, (\bigcup J_n)=B_{n+1}\,\take\, (\bigcup_{k=1}^n B_k)$
is an open,
compact subset of $B_{n+1}$
and thus $\psi(R)$ is an open, compact subset of $B$.
As $\psi(R)$ is open and compact,
there exists $\ve\in \; ]0,1]$ such that $\psi(R)+B_\ve(0)\sub \psi(R)$.
Since $B_\ve(0)$ is an open subgroup
of the compact additive group $B$,
we deduce that
$\psi(R)$ is the disjoint union of a finite
number of balls $B_\ve(a_1),\cdots,B_\ve(a_m)$
(i.e., cosets of $B_\ve(0)$),
for some $m\in \N_0$
and elements $a_1,\ldots, a_m\in B$.
Then $J_{n+1}:=J_n\cup\{\psi^{-1}(B_\ve(a_i))\!:i=1,\ldots,m\}$
is an open cover of~$B_{n+1}$ by mutually disjoint
balls, and $J_n\sub J_{n+1}$ by definition.
Proceeding in this way, we obtain an ascending
sequence $J_1\sub J_2\sub\cdots$,
where each $J_n$ is an open cover 
of $\bigcup_{k=1}^n B_k$ by mutually disjoint balls.
Thus $J:=\bigcup_{k\in \N}J_k$ is a countable cover of~$M$
by mutually disjoint balls.
Choose an injection $\kappa\!: J\to \N$.
For each ball $C_j:=j\in J$, there exists a $C^r_\F$-diffeomorphism
$\phi_j\!: C_j\to \pi^{-\kappa(j)}+B\sub \F^d$.
Then
$U:=\bigcup_{j\in J}(\pi^{-\kappa(j)}+B)$
is an open subset of~$\F^d$.
The union defining~$U$ is disjoint,
because $|\pi^{-\kappa(j)}+x|=\max\{|\pi^{-\kappa(j)}|,|x|\}=|\pi^{-\kappa(j)}|
=|\pi|^{-\kappa(j)}$ for each $j\in J$
and $x\in B$.
Hence $\phi:=\coprod_{j\in J}\phi_j\!: M\to U$
(the map determined by $\phi|_{C_j}=\phi_j$)
is a $C^r_\F$-diffeomorphism.\vspace{2mm}

(b) By Lemma~\ref{paracamen}, $M$ is a disjoint union
$M=\coprod_{i\in I}M_i$ of $\sigma$-compact,
open and closed submanifolds~$M_i$.
The proof of (a) shows that each $M_i$
is a disjoint union $M_i=\coprod_{j\in J_i}C_{i,j}$
of a countable family $(C_{i,j})_{j\in J_i}$
of balls $C_{i,j}\sub M_i$. Set
$K:=\{(i,j)\!: i\in I, j\in J_i\}$.
Then $M=\coprod_{(i,j)\in K}C_{i,j}$ is a disjoint union of balls.
\end{proof}
If $U$ is an open subset of $\F^d$, we can even
find partitions into metric balls
subordinate to any given open cover:
\begin{la}\label{disjballs}
Suppose that $\F$ is a local field,
$d\in \N$ and $U\sub \F^d$ a non-empty,
open subset. Let $(U_i)_{i\in I}$ be an open cover of~$U$.
Then there exist
families $(r_j)_{j\in J}$ and $(a_j)_{j\in J}$
of positive real numbers $r_j>0$,
resp., elements $a_j\in U$,
indexed by a countable set~$J$,
such that $(B_{r_j}(a_j))_{j\in J}$
is an open cover of~$U$ by mutually disjoint
sets and furthermore the open cover
$(B_{r_j}(a_j))_{j\in J}$ is subordinate to $(U_i)_{i\in I}$,
viz.\ for every $j\in J$, there exists $i(j)\in I$
such that $B_{r_j}(a_j)\sub U_{i(j)}$.
\end{la}
\begin{proof}
Since $U$ is $\sigma$-compact, we find a sequence
$(B_k)_{k\in \N}$ of metric balls covering~$U$
and which is subordinate to $(U_i)_{i\in I}$:
For each $k\in \N$, there exists $i_k\in I$ such that
$B_k \sub U_{i_k}$.
Adapting the proof of
Lemma~\ref{onlyopen}\,(a) in the obvious
way,\footnote{Thus, we choose each $\phi$ of the
form $\phi(z)=az+b$ with suitable $a\in \F^\times$, $b\in \F^d$.}
we arrive at a countable cover $J=\bigcup_{k\in \N} J_k$
of $U$ by mutually disjoint metric balls,
such that all balls $C \in J_1$ are subsets of $B_1\sub U_{i(1)}$
and all balls $C\in J_{k+1}\setminus J_k$
are subsets of $B_{k+1}\sub U_{i_{k+1}}$.
\end{proof}
Locally finite, relatively compact, open
covers can always be thickened.
\begin{la}\label{lathicken}
Let $\F$ be a locally compact field,
$M$ be a paracompact, finite-dimensional
$C^r_\F$-manifold, and $(U_i)_{i\in I}$ be a locally
finite cover of~$M$ by relatively compact, open subsets
$U_i\sub M$. Then there exists a locally finite
cover $(\wt{U}_i)_{i\in I}$ of $M$ by relatively compact, open
subsets $\wt{U}_i\sub M$, such that
for each $i\in I$ the closure $\wb{U}_i$ of~$U_i$
in~$M$ is contained in $\wt{U}_i$.
\end{la}
\begin{proof}
To reduce the assertion to the $\sigma$-compact case,
we first observe
that $U_i$ is $\sigma$-compact, for each $i\in I$
(using that $\wb{U}_i$ can be covered
by finitely many balls).
We now write $i\sim j$ for $i,j\in I$
if and only if there exists $n\in \N$
and $k_1,\ldots, k_n\in I$ such that
$k_1=i$, $k_n=j$ and
$U_{k_\nu}\cap U_{k_{\nu+1}}\not=\emptyset$
for $\nu=1,\ldots, n-1$.
Then $\sim$ is an equivalence relation.
Since $(U_i)_{i\in I}$ is a locally finite
cover and each $U_i$ is relatively compact,
we deduce that each equivalence class $C\in I/\! \sim\;\, =:J$
is countable; hence
\[
M_C:=\bigcup_{i\in C}U_i
\]
is a $\sigma$-compact open subset
of~$M$. By construction, $M=\coprod_{C\in J}M_C$
is a disjoint union of the
open (and hence also closed) sets~$M_C$. Since
$(U_i)_{i\in C}$ is a countable, locally finite
cover of~$M_C$ by relatively compact, open sets,
it suffices to prove our assertion for countable
covers of the $M_C$'s.
We may hence assume that $M$ is $\sigma$-compact
and that $I$ is countable. If $I$ is finite, then $M$ is compact and
the assertion is trivial. Thus $I=\N$ without loss of generality.\\[3mm]
To construct a suitable open cover $(\wt{U}_n)_{n\in \N}$,
choose a sequence $(K_n)_{n\in \N}$
of compact subsets $K_n\sub M$ such that $\bigcup_{n\in \N}K_n=M$
and such that $K_n$ is contained in the interior
$(K_{n+1})^0$, for each $n\in \N$.
Define $K_{-1}:=K_0:=\emptyset$ for convenience
of notation.
Then $I_m:=\{n\in \N\!: (K_m\,\take\, (K_{m-1})^0)
\cap \wb{U_n}\not=\emptyset\}$
is a finite set, for each $m\in \N$,
because also the sequence $(\wb{U_n})_{n\in \N}$
of the closures is locally finite.\footnote{Every
$x\in M$ has an open neighbourhood $U$
such that $\{n\in \N\!: U_n\cap U\not=\emptyset\}$
is finite. The set $U$ being open, we have
$\{n\in \N\!: U_n\cap U\not=\emptyset\}=
\{n\in \N\!: \wb{U_n}\cap U\not=\emptyset\}$.}
Also $J_n:=\{m\in\N\!: n \in I_m\}$
is a finite set for each $n\in \N$:
indeed, there is $m_0\in \N$ such that
$\wb{U_n}\sub K_{m_0}$;
then $\wb{U_n}\sub (K_m)^0$ for all $m\geq m_0+1$
and thus $\wb{U_n}\cap (K_m\,\take\, (K_{m-1})^0)=\emptyset$
for all $m\geq m_0+2$, entailing that $m\not\in J_n$
for all $m\geq m_0+2$. For each $m\in \N$
and $n\in I_m$, the set $V_m:=(K_{m+1})^0\,\take\, K_{m-2}$
is an open neighbourhood
of $\wb{U_n}\cap
(K_m\,\take\, (K_{m-1})^0)$,
which is contained in $K_{m+1}$ and therefore relatively
compact. We set $\wt{U}_n:=\bigcup_{m\in J_n} V_m$;
this is a relatively compact, open neighbourhood
of $\wb{U_n}$. Given $n,m\in \N$,
we have $K_m\cap \wt{U}_n=\bigcup_{m'\in J_n}
(K_m\cap V_{m'})$,
where $K_m\cap V_{m'}=\emptyset$ unless
$m'\leq m+2$. Let $m'\leq m+2$.
If $m'\in J_n$, then $n\in I_{m'}$.
Thus $K_m\cap \wt{U}_n=\emptyset$
unless $n\in \bigcup_{m'=1}^{m+2} I_{m'}$,
which is a finite set.
It now readily follows that the open cover $(\wt{U}_n)_{n\in \N}$
of~$M$ is locally finite.
In fact, given any $x\in M$ we find $m\in \N$ such that
$K_m^0$ is an open neighbourhood
of~$x$. By the preceding, $K_m$ (and hence $K_m^0$)
only meets $\wt{U}_n$ for finitely many~$n$.
\end{proof}
Cut-offs and
partitions
of unity on finite-dimensional
real manifolds are standard tools.
To enable unified proofs,
we now discuss
analogous concepts also over local fields.
\begin{defn}\label{defnparun}
Let $\F$ be a local field.
A {\em $C^r_\F$-partition of unity\/}
of a finite-dimensional
$C^r_\F$-manifold $M$
is a family $(h_i)_{i\in I}$
of continuous
mappings $h_i\!: M\to \{0,1\}\sub \F$,
such that the open and closed sets
$h_i^{-1}(\{1\})$ are mutually disjoint and
cover~$M$.
\end{defn}
Note that,
being locally constant, each $h_i$
is actually
$C^r_\F$.
\begin{la}\label{parloc}
Let $M$ be a $\sigma$-compact $C^r_\F$-manifold
over a local field~$\F$,
and $(U_i)_{i\in I}$
be an open cover of~$M$.
Then
there exists a partition of unity
$(h_i)_{i\in I}$
such that $\Supp(h_i)\sub U_i$.
\end{la}
\begin{proof}
The assertion is trivial if $\dim(M)=0$.
If $d:=\dim(M)>0$,
by Lemma~\ref{onlyopen}
we may assume that $M$ is an open subset of~$\F^d$.
Let $(B_{r_j}(a_j))_{j\in J}$
and $i(j)$ for $j\in J$
be as in Lemma~\ref{disjballs}.
The family of balls being locally finite,
the open sets $V_i:=\bigcup_{j\in J: i(j)=i}B_{r_j}(a_j)\sub U_i$
are also closed.
For each $i\in I$, define
$h_i\!: M\to \F$ via
$h_i(x):=1\in \F$ if $x\in V_i$,
$h_i(x):=0$ otherwise.
Then $(h_i)_{i\in I}$
is a partition of unity with the desired properties. 
\end{proof}
\begin{la}\label{cutthroat}
Let $\F$ be a local field,
$M$ be a finite-dimensional
$C^r_\F$-manifold, $K\sub M$ be compact, and
$U\sub M$ be an open subset containing~$K$.
Then there exists a $C^r_\F$-function
$h\!: M\to \{0,1\}\sub \F$ such that
$h|_K=1$ and $h|_{M\setminus U}=0$.
\end{la}
\begin{proof}
As each element $x\in K$ is contained
in some open and compact ball $B_x\sub U$,
exploiting the compactness of~$K$
we find finitely many open and compact
balls $C_1,\ldots, C_n\sub U$ such that
$K\sub \bigcup_{k=1}^n C_k=:W$.
Then $W$ is an open and closed neighbourhood
of $K$ such that $W\sub U$, and hence
$h\!: M\to \F$, $h(x):=1$ if $x\in W$, else $h(x):=0$
is a function with the desired properties.\vspace{-1mm}
\end{proof}
\begin{center}
{\bf Topologies on spaces of vector-valued test functions}
\end{center}
For the remainder of this section,
$\F$ denotes a locally compact, non-discrete topological
field, and $\K$ a topological extension field
of~$\F$, whose topology arises
from an absolute value $|.|\!: \K\to [0,\infty[$.
We let $r\in \N_0\cup\{\infty\}$.
\begin{defn}\label{defcpsupp}
Given a
paracompact $C^r_\F$-manifold~$M$,
modeled on a finite-dimensional
$\F$-vector space~$Z$,
and a (Hausdorff, not necessarily locally convex)
topological $\K$-vector space~$E$,
we let
\[
C^r_c(M,E):=\{\gamma\in C^r(M,E)\!: \;\mbox{$\Supp(\gamma)$
is compact}\,\}
\]
be the set of compactly supported
$E$-valued $C^r_\F$-functions
on~$M$.
Then $C^r_c(M,E)$ is a $\K$-vector subspace
of $C^r(M,E)$, and
$C^r_c(M,E)=\bigcup_{K\in \cK(M)}\, C^{\,r}_K(M,E)$,
where $\cK(M)$ denotes the set of all
compact subsets of~$M$.
In the following,
we consider three vector topologies on
$C^r_c(M,E)$:
\begin{itemize}
\item[(a)]\label{deftvstp}
We write
$C^r_c(M,E)_\tvs$
for $C^r_c(M,E)$, equipped
with the finest (a priori not necessarily
Hausdorff) vector topology
making the inclusion maps $\lambda_K\!:
C^r_K(M,E)\to C^r_c(M,E)$
continuous for each compact subset $K\sub M$.
Thus $C^r_c(M,E)_\tvs=\dl\, C^r_K(M,E)$\vspace{-.8mm}
in the category of
not necessarily Hausdorff topological $\K$-vector spaces
and continuous $\K$-linear maps.
\item[(b)]\label{deflcxtop}
If $E$ is locally convex, we
write
$C^r_c(M,E)_\lcx$
for $C^r_c(M,E)$, equipped
with the finest (a priori not necessarily
Hausdorff) locally convex vector topology
making the inclusion maps $\lambda_K\!:
C^r_K(M,E)\to C^r_c(M,E)$
continuous for each compact subset $K\sub M$.
Thus $C^r_c(M,E)_\lcx=\dl\, C^r_K(M,E)$\vspace{-.8mm}
in the category of
not necessarily Hausdorff, locally
convex topological $\K$-vector spaces
and continuous $\K$-linear maps.
\item[(c)]\label{defbxtop}
Given a locally finite cover $\cU=(U_i)_{i\in I}$
of $M$ by relatively compact, open subsets\footnote{Such a cover
always exists because $M$ is locally compact and paracompact.}
$U_i\sub M$, we let $\rho_i\!: C^r_c(M,E)\to C^r(U_i,E)$
be the restriction map for $i\in I$
and define
\[
\rho_\cU\!: C^r_c(M,E)\to \bigoplus_{i\in I}\, C^r(U_i,E),\quad
\rho_\cU(\gamma):=(\rho_i(\gamma))_{i\in I}=(\gamma|_{U_i})_{i\in I}\,.
\]
We write $C^r_c(M,E)_\bx$ for $C^r_c(M,E)$,
equipped with the topology $\cO_\cU$ induced by $\rho_\cU$,
where the direct sum is endowed with
the box topology.
\end{itemize}
\end{defn}
\begin{la}\label{boxindep}
In the situation of Definition~{\rm \ref{defcpsupp}\,(c)},
assume that both $\cU=(U_i)_{i\in I}$
and $\cV=(V_j)_{j\in J}$
are locally finite covers of~$M$
by relatively compact open subsets.
Then $\cO_\cU=\cO_\cV$.
In other words,
the box topology on $C^r_c(M,E)$
is independent of the choice of $\cU$.
\end{la}
\begin{proof}
The topologies $\cO_\cU$ and $\cO_\cV$ are induced by
$\rho_\cU\!:
C^r_c(M,E)\to \bigoplus_{i\in I}\,
C^r(U_i,E)$, $\rho_\cU(\gamma):=(\gamma|_{U_i})_{i\in I}$
and
$\rho_\cV\!:
C^r_c(M,E)\to \bigoplus_{j\in J}\,
C^r(V_j,E)$, $\rho_\cV(\gamma):=(\gamma|_{V_j})_{j\in J}$,
respectively.
\begin{numba}
Given $i\in I$,
the set $J_i:=\{j\in J\!: U_i\cap V_j\not=\emptyset\}$
is finite, as $U_i$ is relatively compact
and $\cV$ is a locally finite cover.
By Lemma~\ref{amend},
the topology on $C^r(U_i,E)$
is initial with respect to the family $(\mu_{i,j})_{j\in J_i}$
of restriction maps
\[
\mu_{i,j}\!: C^r(U_i,E)\to C^r(U_i\cap V_j,E),\quad
\mu_{i,j}(\gamma):=\gamma|_{U_i\cap V_j}\,.
\]
Likewise, the set $I_j:=\{i\in I\!: U_i\cap V_j\not=\emptyset\}$
is finite for each $j\in J$,
and the topology on $C^r(V_j,E)$
is initial with respect to the family $(\nu_{j,i})_{i\in I_j}$
of restriction mappings\linebreak
$\nu_{j,i}\!: C^r(V_j,E)\to
C^r(U_i\cap V_j,E)$.
\end{numba}
\begin{numba}\label{basisPB}
Let
$P_{i,j}$ be an open $0$-neighbourhood
of $C^r(U_i\cap V_j,E)$, for any $i\in I$, $j\in J_i$.
Then
$P_i:=\bigcap_{j\in J_i}\mu_{i,j}^{-1}(P_{i,j})$
is an open $0$-neighbourhood in $C^r(U_i,E)$
for each $i\in I$ and thus
\[
P:=\bigoplus_{i\in I} P_i
\]
is an open $0$-neighbourhood in
$\bigoplus_{i\in I}\, C^r(U_i,E)$.
It is clear
from the preceding that the set $\cB$
of such open
$0$-neighbourhoods $P$ is a basis
for the filter of $0$-neighbourhoods
of $\bigoplus_{i\in I}\, C^r(U_i,E)$,
and hence $\{\rho_\cU^{-1}(P)\!: P\in \cB\}$
is a basis for the filter of
$0$-neighbourhoods of $(C^r_c(M,E),\cO_\cU)$.
\end{numba}
To see that $\cO_\cU=\cO_\cV$,
it suffices
to show that $\cO_\cU\sub \cO_\cV$
(as we can interchange $\cU$ and $\cV$).
Since both $\cO_\cU$ and $\cO_\cV$
are vector topologies,
we only need to show that
$W \in \cO_\cV$ for $W$
ranging through a suitable basis of open
$0$-neighbourhoods of $(C^r_c(M,E),\cO_\cU)$.
It therefore suffices to consider
$W:=\rho_\cU^{-1}(P)$ for $P\in \cB$
as in {\bf \ref{basisPB}}.
Set $Q_j:=\bigcap_{i\in I_j}\nu_{j,i}^{-1}(P_{i,j})$
for $j\in J$;
then $Q:=\bigoplus_{j\in J}Q_j$
is an open $0$-neighbourhood in $\bigoplus_{j\in J}\,
C^r(V_j,E)$. For $\gamma\in C^r_c(M,E)$, we have
\begin{eqnarray*}
\gamma\in W=\rho_\cU^{-1}(P) \;\; &\aeq& \;\;\rho_\cU(\gamma)\in P
\quad \aeq \quad
(\forall i\in I)\;\gamma|_{U_i}\in P_i\\
&\aeq&\;\;
(\forall i\in I)\,(\forall j\in J_i) \quad
\gamma|_{U_i\cap V_j}\in P_{i,j}\\
&\aeq&\;\;
(\forall j\in J)\,(\forall i\in I_j)\quad
\gamma|_{U_i\cap V_j}\in P_{i,j}\\
&\aeq &\;\; 
(\forall j\in J)\; \gamma|_{V_j}\in Q_j
\quad\aeq\quad
\rho_\cV(\gamma)\in Q
\quad\aeq\quad
\gamma\in \rho_\cV^{-1}(Q)\,.
\end{eqnarray*}
Thus
$W=\rho_\cV^{-1}(Q)\in \cO_\cV$, which completes the proof.
\end{proof}
\begin{prop}\label{comparetop}
Let $M$ be a
paracompact $C^r_\F$-manifold,
modeled on a finite-dimensional
$\F$-vector space~$Z$,
and $E$ be a topological $\K$-vector space.
Then the following holds:
\begin{itemize}
\item[\rm (a)]
The box topology on $C^r_c(M,E)_\bx$
is Hausdorff.
For every locally finite cover $\cU=(U_i)_{i\in I}$
of~$M$ by relatively compact, open subsets $U_i\sub M$,
the map
\[
\rho_\cU\!: C^r_c(M,E)_\bx\to \bigoplus_{i\in I}\, C^r(U_i,E),\quad
\rho_\cU(\gamma):= (\gamma|_{U_i})_{i\in I}
\]
has closed image, and {\/\em $\rho_\cU|^{\text{im}\, \rho_\cU}$\/}
is an isomorphism of topological vector spaces.
The inclusion map $C^r_c(M,E)_\bx \to C^r(M,E)$
is continuous.
If $E$ is locally convex,
then $C^r_c(M,E)_\bx$ is locally convex.
\item[\rm (b)]
The inclusion map
$\lambda_K\!: C^r_K(M,E)\to C^r_c(M,E)_\bx$
is continuous and induces the given topology
on $C^r_K(M,E)$, for each compact subset
$K\sub M$.
\item[\rm (c)]\label{boxdirl1}
The map $\Phi\!: C^r_c(M,E)_\tvs\to C^r_c(M,E)_\bx$,
$\Phi(\gamma):=\gamma$
is continuous.
Thus $C^r_c(M,E)_\tvs$ is Hausdorff
and induces the given topology on each $C^r_K(M,E)$.
If $\F\not=\C$ and $M$ is $\sigma$-compact,
then $\Phi$ is an isomorphism of topological
$\K$-vector spaces.
\item[\rm (d)]
If $E$ is locally convex, then
$\Psi\!: C^r_c(M,E)_\lcx\to C^r_c(M,E)_\bx$,
$\Psi(\gamma):=\gamma$
is continuous.
Hence $C^r_c(M,E)_\lcx$ is Hausdorff
and induces the given topology on each $C^r_K(M,E)$.
If $\F\not=\C$ and $M$ is $\sigma$-compact,
then $\Psi$ is an isomorphism of topological
$\K$-vector spaces.
\item[\rm (e)]
If $\F$ is a local field
and $\cU=(U_i)_{i\in I}$
is a cover of~$M$ by mutually disjoint,
compact open sets $($cf.\ Lemma~{\rm \ref{onlyopen}\,(b)}$)$,
then
\[
\rho_\cU\!: C^r_c(M,E)_\bx \to \bigoplus_{i\in I}\, C^r(U_i,E),\quad
\rho_\cU(\gamma):=(\gamma|_{U_i})_{i\in I}
\]
is an isomorphism of topological vector spaces
onto the direct sum, equipped with the box topology.
\item[\rm (f)]
If $\F$ is a local field
and $E$ is locally convex, then $\Psi$
is an isomorphism of topological vector spaces,
i.e., $C^r_c(M,E)_\lcx=C^r_c(M,E)_\bx$.
\end{itemize}
In particular,
$C^r_c(M,E)_\bx=
C^r_c(M,E)_\tvs= C^r_c(M,E)_\lcx$
if $\F\not=\C$ and $M$
is $\sigma$-compact.
\end{prop}
\begin{proof}
(a) Let $\cU=(U_i)_{i\in I}$
be a locally finite cover of~$M$
by relatively compact, open sets.
The box topology on $\bigoplus_{i\in I} C^r(U_i,E)$
being Hausdorff and $\rho_\cU\!: C^r_c(M,E)\to
\bigoplus_{i\in I} C^r(U_i,E)$
being injective, the topology $\cO_\bx$ induced by
$\rho_\cU$ on $C^r_c(M,E)$ is Hausdorff
and $\rho_\cU|^{\text{im}\, \rho_\cU}$
is an isomorphism of topological vector spaces.
The box topology on $S:=\bigoplus_{i\in I}\, C^r(U_i,E)$
is properly finer than the topology
induced by the product
$P:=\prod_{i\in I}\, C^r(U_i,E)$.
The map
$\tau \!: C^r(M,E)\to P$, $\tau(\gamma):=(\gamma|_{U_i})_{i\in I}$
is a topological embedding
with closed image, by Lemma~\ref{amend}.
This entails, firstly,
that the inclusion map
$C^r_c(M,E)_\bx\to C^r(M,E)$ is continuous.
Secondly, it entails
that $\im(\tau)\cap S$ is closed
in~$S$. Note that $\tau(\gamma)\in S$ implies
that $\Supp(\gamma)$ is compact, i.e.,
$\gamma\in C^r_c(M,E)$.
Thus $\im(\rho_\cU)=\im(\tau)\cap S$ is closed
in~$S$.
If $E$ is locally convex, then
each of the spaces $C^r(U_i,E)$
is locally convex (Proposition~\ref{propprop}\,(b)),
whence so is the
direct sum $\bigoplus_{i\in I}\, C^r(U_i,E)$ (see
{\bf \ref{lcx2}})
and hence so is $\cO_\bx$.\vspace{2mm}

(b) Let $K\sub M$
be compact.
For $\cU=(U_i)_{i\in I}$
as before, there exists a finite subset
$J\sub I$ such that
$U_i\cap K=\emptyset$ for all $i\in I\setminus J$.
Thus
$K\sub \bigcup_{i\in J} U_i =:W$.
Then the composition $f\!: C^r_K(M,E)\to
\bigoplus_{i\in I} \, C^r(U_i,E)$ of the maps 
\[
C^r_K(M,E)\,\stackrel{\isom}{\to}\,
C^r_K(W,E)\,\emb\, C^r(W,E)\,\emb\, \prod_{i\in J}C^r(U_i,E)\,\emb\,
\bigoplus_{i\in I} \, C^r(U_i,E)
\]
is a topological embedding, where
the first map
and the coordinate functions
of the second map are the respective restriction
maps (see Lemma~\ref{restrK} and Lemma~\ref{amend}),
and the last map is inclusion (see {\bf \ref{obvibox}}).
Since $f=\rho_\cU\circ \lambda_K$,
where $\rho_\cU$ is a topological embedding,
we deduce that also $\lambda_K$
is a topological embedding.\vspace{2mm}

(c) Since $\lambda_K\!: C^r_K(M,E)\to C^r_c(M,E)_\bx$
is continuous for each $K\in \cK(M)$,
our definition of $C^r_c(M,E)_\tvs$
shows that the topology on
$C^r_c(M,E)_\tvs$ is finer than the one on
$C^r_c(M,E)_\bx$, and thus $\Phi$ is continuous.
Since $\lambda_K\!: C^r_K(M,E)\to C^r_c(M,E)_\tvs$
is continuous as a map into $C^r_c(M,E)_\tvs$ and $\Phi\circ \lambda_K
\!: C^r_K(M,E)\to C^r_c(M,E)_\bx$ is a topological embedding
by (b), also $\lambda_K\!: C^r_K(M,E)\to C^r_c(M,E)_\tvs$
is a topological embedding.\\[3mm]
We now assume that $M$ is $\sigma$-compact,
and we assume that $\F$ is not isomorphic to
$\C$ as a topological field; then $\F$ is a local
field or $\F\isom \R$ (see \cite{Wei}).
We have to show that
$C^r_c(M,E)_\bx=\dl\, C^r_K(M,E)$\vspace{-.8mm}
in the category of topological vector spaces,
with limit maps
$\lambda_K\!: C^r_K(M,E)\to C^r_c(M,E)_\bx$.
We already know from (b) that each $\lambda_K$
is continuous; thus $(C^r_c(M,E)_\bx,(\lambda_K)_{K\in \cK(M)})$
is a cone in the category of topological
$\K$-vector spaces and continuous $\K$-linear maps.
To see that it is a direct limit cone,
suppose that $(f_K)_{K\in \cK(M)}$
is a family of continuous linear maps
$f_K\!: C^r_K(M,E)\to F$
into a topological $\K$-vector space~$F$
such that $f_L|_{C^r_K(M,E)}=f_K$
whenever $K\sub L$.
Then $f\!: C^r_c(M,E)_\bx\to F$, $f(\gamma):=f_K(\gamma)$
if $\Supp(\gamma)\sub K$ is well-defined
and is the unique linear map $C^r_c(M,E)_\bx\to
F$ such that $f\circ \lambda_K=f_K$
for each~$K$. To establish
the desired direct limit property,
it only remains to show that $f$ is continuous.\\[3mm]
Let $\cU=(U_i)_{i\in I}$
be as before; $M$ being $\sigma$-compact,
we may assume that~$I$ is countable.
We pick a $C^r_\F$-partition of unity $(h_i)_{i\in I}$
of~$M$ such that
$\Supp(h_i)\sub U_i$
for each~$i\in I$ (see Lemma~\ref{parloc}
when $\F$ is a local field;
the real case is standard).
Since
$U_i$ is relatively compact, $K_i:=\Supp(h_i)\sub U_i$ is compact.
Let $e_i\!: C^{\,r}_{K_i}(U_i,E)\to C^{\,r}_{K_i}(M,E)$
be the isomorphism of topological
vector spaces which extends functions by~$0$
(cf.\ Lemma~\ref{restrK}).
Then $g_i\!: C^r(U_i,E)\to F$,
$g_i:= f_{K_i}\circ e_i \circ \mu_{h_i}$
is a continuous linear mapping,
where $\mu_{h_i}\!:$\linebreak
$C^r(U_i,E)\to C^{\,r}_{K_i}(U_i, E)$
is the multiplication operator
defined via
$\mu_{h_i}(\gamma):=h_i|_{U_i}\cdot \gamma$,
which is continuous linear
(in view of Lemma~\ref{corestr}
and Lemma~\ref{amend},
applied with a cover of coordinate neighbourhoods,
this assertion can be reduced to Lemma~\ref{multop}).
By the universal property
of the countable direct sum $S:=\bigoplus_{i\in I}\, C^r(U_i,E)$
(see {\bf \ref{boxisuniv}}),
the linear map
\[
g\!: S\to F,\;\;\; (\gamma_i)_{i\in I}\mto \sum_{i\in I}g_i(\gamma_i)
\]
(where $\gamma_i\in C^r(U_i,E)$)
is continuous,
because so is each~$g_i$.
Given $\gamma\in C^r_c(M,E)$,
we calculate
$g_i(\gamma|_{U_i})=f(e_i(h_i|_{U_i}\cdot \gamma|_{U_i}))
=f(e_i((h_i\cdot \gamma)|_{U_i}))=f(h_i\cdot \gamma)$,
whence $g(\rho_\cU(\gamma))=\sum_{i\in I}
g_i(\gamma|_{U_i})=\sum_{i\in I}f(h_i\cdot \gamma)=f(\sum_{i\in
I}h_i\cdot \gamma)=f(\gamma)$.
Thus $g\circ \rho_\cU=f$,
and so $f$ is continuous
on $C^r_c(M,E)_\bx$,
as required; the direct limit property
is fully established.\vspace{2mm}

(d) Since $C^r_c(M,E)_\bx$
is locally convex for locally convex~$E$ (see (a)),
we can repeat the proof of (c), except that
topological vector spaces have to be replaced
with locally convex spaces.\vspace{2mm}

(e) By the definition of the box topology,
$\rho_\cU$ is a topological embedding.
Each of the sets $U_i$ being compact and open,
given $(\gamma_i)_{i\in I}\in \bigoplus_{i\in I}\, C^r(U_i,E)$
the map $\gamma\!: M\to E$ defined via
$\gamma(x):=\gamma_i(x)$ for $x\in U_i$
is $C^r_\F$ and compactly supported.
Thus $\gamma\in C^r_c(M,E)$,
and $\rho_\cU(\gamma)=(\gamma_i)_{i\in I}$
by definition of~$\gamma$. Thus $\rho_\cU$
is also surjective, and thus $\rho_\cU$
is an isomorphism of topological vector spaces.\vspace{2mm}

(f) Assume that $f\!: C^r_c(M,E)_\bx\to F$
is a linear map into a locally convex topological
$\K$-vector space~$F$ such that $f\circ\lambda_K$
is continuous for each $K\in \cK(M)$.
Let $\cU=(U_i)_{i\in I}$
be as in~(e). Then $f|_{C^r_{U_i}(M,E)}$
is continuous in particular for each $i\in I$,
and hence so is $g_i:=f\circ e_i\!: C^r(U_i,E)\to F$,
where $e_i\!: C^r(U_i,E)\to C^r_{U_i}(M,E)$
is the isomorphism of topological vector spaces
obtained as the inverse of the restriction map\linebreak
$C^r_{U_i}(M,E)\to C^r(U_i,E)$
(see Lemma~\ref{restrK}).
Since the box topology makes
$\bigoplus_{i\in I}\, C^r(U_i,E)$
the category-theoretical
locally convex direct sum
in the present situation (see Remark~\ref{disclcx}),
the map
$g\!: \bigoplus_{i\in I}\, C^r(U_i,E)\to F$,
$g((\gamma_i)_{i\in I}):=\sum_{i\in I}g_i(\gamma_i)$
is continuous linear.
Hence also $f=g\circ \rho_\cU$ is continuous.
\end{proof}
\begin{cnv}\label{anoconv}
Throughout the following,
spaces of vector-valued
test functions
will always be equipped
with the box topology,
and we abbreviate $C^r_c(M,E):=C^r_c(M,E)_\bx$.
\end{cnv}
\begin{rem}\label{traditio}
If $\F=\K=\R$,
$M$ is $\sigma$-compact and $E$ is locally convex,
then the box topology on $C^r_c(M,E)$
coincides with the locally convex topology traditionally
considered on this space of test functions,
by Proposition~\ref{propprop}\,(d)
and Proposition~\ref{comparetop}\,(d).
\end{rem}
\begin{rem}\label{justffy}
As we shall mainly need spaces of test functions
with values in locally convex spaces
over local fields in the following
(for example, in our discussion of diffeomorphism
groups), we have chosen to work with the box topology,
which is the appropriate topology on $C^r_c(M,E)$
in this case, even for non-$\sigma$-compact~$M$
(see Proposition~\ref{comparetop}\,(f)).
If $M$ is a non-$\sigma$-compact, paracompact
finite-dimensional manifold over~$\R$
and $E$ a real locally convex space,
it is certainly more natural to
work with the (finer) locally convex direct limit topology
on $C^r_c(M,E)$, and study differentiability properties
of mappings between spaces of test functions
topologized in this way.
We did not find it advantageous
(nor necessary) to discuss
this situation in parallel here;
the interested reader can
find a separate discussion in \cite{SEC} and \cite{DIF}
(cf.\ also \cite{Mic}).
\end{rem}
\begin{rem}\label{remnosecs}
Spaces of compactly supported sections
in vector bundles can be treated
much in the same way as spaces of test functions.
However, the only vector bundles we shall really
have to work with
in our Lie group constructions (of diffeomorphism groups)
are the tangent bundles of paracompact finite-dimensional
smooth manifolds $M$ over local fields~$\K$.
Since any vector bundle over
such a manifold~$M$
is trivial (as a consequence
of Lemma~\ref{onlyopen}\,(b)
and Lemma~\ref{disjballs}),
it is not necessary for our purposes
to introduce the additional machinery required
to discuss spaces of sections and vector bundles,
and so we decided to defer their discussion to
an appendix (Appendix~\ref{appsections}).
The only facts we shall really use are the following:
1.\,For each paracompact, finite-dimensional
$C^r_\K$-manifold~$M$ and disjoint cover $(B_i)_{i\in I}$
by open and compact balls, the map
$C^r_c(M,TM)\to \bigoplus_{i\in I}C^r(B_i,TB_i)$,
$\sigma\mto (\sigma|_{B_i})_{i\in I}$
is an isomorphism of topological vector spaces
(Proposition~\ref{comparetop2}\,(e)).
2.\,If $\kappa\!: M\to B$ is a $C^r_\K$-diffeomorphism
from a $d$-dimensional
$C^r_\K$-manifold~$M$ onto a metric ball $B\sub \K^d$,
then $C^r(M,TM)\to C^r(B,\K^d)$,
$\sigma\mto (x\mto (d \kappa_i)\circ \sigma\circ \kappa^{-1})$
is an isomorphism of topological $\K$-vector spaces
(cf.\ Lemma~\ref{anyatlas} and Lemma~\ref{atlas}).\vspace{2mm}
\end{rem}
\begin{center}
{\bf Patched topological vector spaces and patched mappings}
\end{center}
To formalize the situation encountered
in Proposition~\ref{comparetop}\,(a),
we now introduce the notion of a
``patched'' topological vector space.
Roughly speaking, this is a topological
vector space, together with an embedding into a direct sum.
We then discuss differentiability properties
of mappings between patched topological vector spaces.
The general results obtained here shall allow us
to transfer our discussion of pushforwards
to the case of test functions (Proposition~\ref{pfwcsupp}).
We shall also derive certain very convenient
criteria ensuring differentiability properties
for mappings between spaces of test functions
(Section~\ref{secalmloc}).\\[3mm]
For analogous discussions
(and applications)
of ``patched locally convex spaces''
based on locally convex direct sums,
we refer to \cite{SEC} and \cite{DIF}.
\begin{defn}\label{defnpatched}
A {\em patched topological vector space\/}
over a valued field $(\K,|.|)$
is a pair\linebreak
 $(E,(p_i)_{i\in I})$, where
$E$ is a topological $\K$-vector space
and $(p_i)_{i\in I}$ a family of continuous linear
maps $p_i\!: E\to E_i$ to certain topological vector spaces~$E_i$,
such that
\begin{itemize}
\item[\n (a)]
For each $x\in E$, the set $\{i\in I\!: p_i(x)\not = 0\}$
is finite;
\item[\n (b)]
The linear map
\[
p\!: E\, \to \, \bigoplus_{i\in I} \, E_i\,,\;\;\;
x\mto (p_i(x))_{i\in I}=\sum_{i\in I}p_i(x)
\]
from $E$ to the direct sum $\bigoplus_{i\in I}E_i$
(equipped with the box topology)
is a topological embedding;
\item[\n (c)]
The image $p(E)$ is sequentially closed
in $\bigoplus_{i\in I}E_i$.
\end{itemize}
The mappings $p_i\!: E\to E_i$\label{pagerpa}
are called {\em patches}, and the family
$(p_i)_{i\in I}$ is called a {\em patchwork}.
\end{defn}
We retain the notation introduced earlier in this section.
\begin{example}\label{patchd}
Let $M$ be a paracompact,
finite-dimensional $C^r_\F$-manifold over~$\F$
(where $r\in \N_0\cup\{\infty\}$), and $E$ be a
topological $\K$-vector space.
Let $(U_i)_{i\in I}$
be a locally finite open cover
of~$M$ by relatively compact, open subsets
$U_i\sub M$
and $\rho_i\!: C^r_c(M,E)\to C^r(U_i,E)$,
$\rho_i(\gamma):=\gamma|_{U_i}$ be the restriction map
for $i\in I$. Then
\[
(C^r_c(M,E),(\rho_i)_{i\in I})
\]
is a patched topological vector space, by
Proposition~\ref{comparetop}\,(a).
\end{example}
We now discuss mappings between open subsets
of patched topological vector spaces.
\begin{defn}\label{defnpama}
Let $(E,(p_i)_{i\in I})$ and $(F,(q_i)_{i\in I})$
be patched topological $\K$-vector spaces over
the same index set~$I$.
Let $p\!: E\to \bigoplus_{i\in I} E_i$ and
$q\!: F\to \bigoplus_{i\in I}F_i$ be the canonical embeddings.
\begin{itemize}
\item[(a)]
A map $f\!: U\to F$,
defined on an open subset $U$ of~$E$,
is called a {\em patched mapping\/}
if there exists a family $(f_i)_{i\in I}$ of
mappings
$f_i\!: U_i\to F_i$ on certain
open neighbourhoods $U_i$ of $p_i(U)$ in~$E_i$,
which is {\em compatible with~$f$\/}
in the following sense:
we have $0\in U_i$ and $f_i(0)=0$ for all but finitely many~$i$,
and $q_i(f(x))=f_i(p_i(x))$ for all $i \in I$,
{\em i.e.}, $q\circ f= (\oplus f_i)\circ p|_U^{\oplus U_i}$.
\item[(b)]\label{defnpm2}
Given $k\in \N_0\cup\{\infty\}$,
we say that a patched mapping $f\!: U\to F$ as before
is {\em of class $C^k_\K$ on the patches\/}
if all of the mappings $f_i$ in (a)
can be chosen of class~$C^k_\K$.
\end{itemize}
\end{defn}
\begin{prop}\label{diffpatch}
Let $(E,(p_i)_{i\in I})$
and $(F,(q_i)_{i\in I})$
be patched topological $\K$-vector spaces over
the same index set~$I$.
Assume that $f\!: U \to F$
is a patched mapping from an open subset~$U\sub E$ to~$F$.
If $f$ is of class $C^k_\K$ on the patches,
then $f$ is of class $C^k_\K$.
\end{prop}
\begin{proof}
Let $p_i\!: E\to E_i$ and $q_i\!: F\to F_i$
be the patches of $E$, resp., $F$.
If $f$ is of class $C^k_\K$ on the patches,
then there exists a family $(f_i)_{i\in I}$
of $C^k_\K$-maps $f_i\!: U_i\to F_i$,
which is compatible with~$f$.
By Proposition~\ref{mapsdirsums},
the map $g:=\oplus_{i\in I}f_i\!:
\bigoplus_{i\in I}U_i \to \bigoplus_{i\in I}F_i$,
$g(\sum_{i\in I} u_i):=\sum_{i\in I}f_i(u_i)$
is of class~$C^k_\K$.
The linear map $p\!: E\to\bigoplus_{i\in I}
E_i$, $p(x):=\sum_{i\in I}p_i(x)$ being continuous,
the composition $g\circ p|_U^{\oplus\, U_i}$ is $C^k_\K$.
But $g\circ p|_U^{\oplus\, U_i}=q\circ f$,
where
$q\!: F\to\bigoplus_{i\in I}F_i$,\linebreak
$q(y):=\sum_{i\in I} q_i(y)$,
since
$g(p(x))=\sum_{i\in I}f_i(p_i(x))
=\sum_{i\in I}q_i(f(x))=q(f(x))$ for all $x\in U$.
By the preceding, the map $q\circ f$
is of class $C^k_\K$.
Its image being contained in
the sequentially closed vector subspace
$Q:=\im\, q$ of $\bigoplus_{i\in I}F_i$,
Lemma~\ref{corestr}
shows that also
the co-restriction $(q\circ f)|^Q$
is of class~$C^k_\K$.
As $q|^Q$ is an isomorphism of topological vector
spaces (by the axioms of a patched topological vector space),
we see that also $f=(q|^Q)^{-1}\circ (q\circ f)|^Q$
is~$C^k_\K$.
\end{proof}
\begin{center}
{\bf Example: Pushforwards of compactly supported functions}
\end{center}
We now
establish an analogue
of Proposition~\ref{pushforw2}
for pushforwards between
spaces of test functions.
The idea is to
use the technique of patched topological vector spaces
to reduce the assertion to Proposition~\ref{crucial}\,(a).
A further generalization
to mappings between
spaces of compactly supported sections
in vector bundles
is provided in Appendix~\ref{appsections}
(the ``$\Omega$-Lemma with parameters'').
In the real locally convex case,
stronger and much more refined
results are available: see \cite{SEC}.
\begin{prop}\label{pfwcsupp}
Let $E$, $F$ and $\wt{Z}$ be topological
$\K$-vector spaces, $U\sub E$
an open zero-neighbourhood, $r,k\in \N_0\cup\{\infty\}$,
$\wt{M}$ be a $\K$-manifold of class $C^{r+k}_\K$
modeled on $\wt{Z}$,
and
$\tilde{f}\!: \wt{M}\times U \to F$
be a mapping of class $C^{r+k}_\K$.
Let $M$ be a paracompact, finite-dimensional
$\F$-manifold of class $C^r_\F$.
Given a mapping $\sigma\!:
M\to \wt{M}$ of class~$C^r_\F$,
we define
$f:=\tilde{f}\circ (\sigma\times \id_U)\!:
M\times U\to F$.
We assume that
there exists a compact subset $K\sub M$
such that $f(x,0)=0$
for all $x\in M\,\take\, K$.
Then $C_c^r(M,U):=\{\gamma\in C^r_c(M,E)\!: \gamma(M)\sub U\}$
is an open subset of $C^r_c(M,E)$,
equipped with the box topology, and
\[
f_*\!:
C_c^r(M,U) \to C_c^r(M,F),\;\;\;\;
f_*(\gamma)(x):=f(x,\gamma(x))
\]
is a mapping of class $C^k_\K$.
\end{prop}
\begin{proof}
There exist locally finite
covers $\cV:=(V_i)_{i\in I}$
and $\cU:=(U_i)_{i\in I}$
of~$M$ by relatively compact, open sets $V_i$
(resp., $U_i$), such that $K_i:=\wb{V_i}\sub U_i$
for all $i\in I$ (cf.\ Lemma~\ref{lathicken}).\\[3mm]
To see that $C^r_c(M,U)$ is open in $C^r_c(M,E)$,
note that
$\lfloor K_i, U \rfloor_r\sub C^r(U_i,E)$
is an open $0$-neighbourhood, and thus
$Q:=\bigoplus_{i\in I}\, \lfloor K_i, U\rfloor_r$
is an open $0$-neighbourhood in
$\bigoplus_{i\in I} C^r(U_i,E)$.
By definition of the box topology on
$C^r(M,E)$, the map
$\rho_\cU\! : C^r_c(M,E) \to \bigoplus_{i\in I} C^r(U_i,E)$,
$\rho_\cU(\eta):=(\eta|_{U_i})_{i\in I}$
is continuous. Hence $\rho_\cU^{-1}(Q)$
is an open $0$-neighbourhood in $C^r_c(M,E)$.
Since
\[
\rho_\cU^{-1}(Q)\;=\; \{ \, \gamma\in C^r_c(M,E)\!:
\, (\forall i\in I)\;\, \gamma(K_i)\sub U\, \}\,,
\]
where $M=\bigcup_{i\in I}K_i$,
we have $\rho_\cU^{-1}(Q)=C^r_c(M,U)$.
Hence $C^r_c(M,U)$ is an open $0$-neighbourhood,
as required.\\[3mm]
To see that $f_*$ is $C^k_\K$, we
shall exploit that $(C^r_c(M,E),(\rho_i)_{i\in I})$
and $(C^r_c(M,F),(\tau_i)_{i\in I})$ are patched topological vector
spaces, with the patches
\[
\rho_i\!: C^r_c(M,E)\to C^r_c(U_i,E),\quad
\rho_i(\gamma):=\gamma|_{U_i} \quad \mbox{and}
\]
\[
\tau_i\!: C^r_c(M,F)\to C^r_c(V_i,F),\quad
\tau_i(\gamma):=\gamma|_{V_i}\,,
\]
respectively (see Example~\ref{patchd}).
The map $f_i:=f|_{U_i\times U}\!: U_i\times U\to F$
being of the form $\wt{f}\circ (\sigma|_{U_i} \times \id_U)$,
the pushforward
\[
(f_i)_*\!: \lfloor K_i, U\rfloor_r\to C^r(V_i,F)
\]
is $C^r_\K$ on the open subset
$\lfloor K_i, U\rfloor_r\sub C^r(U_i,E)$,
by Proposition~\ref{crucial}\,(a)
(applied with a singleton parameter set~$P$).
Apparently $\rho_i(C^r_c(M,U))\sub \lfloor K_i, U\rfloor_r$
for each $i\in I$, and the family of mappings
$((f_i)_*)_{i\in I}$ is compatible with $f_*$,
i.e., $\tau_i\circ f_*= (f_i)_*\circ \rho_i$
for all $i\in I$.
Thus $f_*$ is a patched mapping.
Every $(f_i)_*$ being $C^k_\K$,
the map $f_*$ is $C^k_\K$ on the patches
and hence $C^k_\K$, by Proposition~\ref{diffpatch}.
\end{proof}
\begin{cor}\label{functcsupp}
Let $E$ and $F$ be topological $\K$-vector
spaces
and $f\!:U\to F$ be a mapping of class~$C^{r+k}_\K$,
defined on an open zero-neighbourhood $U\sub E$,
such that $f(0)=0$.
Let $M$ be a paracompact,
finite-dimensional $\F$-manifold of class $C^r_\F$.
Then
\[
C_c^r(M,f)\!: C_c^r(M,U)\to C_c^r(M,F),\;\;\;\;
\gamma\mto f\circ \gamma
\]
is a mapping of class~$C^k_\K$.
\end{cor}
\begin{proof}
Let $\wt{M}:=\{0\}$ be a singleton smooth $\K$-manifold,
and $\sigma\!: M\to \wt{M}$, $x\mto 0$,
which apparently is a $C^r_\F$-map.
Then $\tilde{g}\!: \wt{M}\times U\to F$,
$\tilde{g}(0,u):=f(u)$
is a mapping of class $C^{r+k}_\K$,
and
$C_K^r(M,f)=g_*$
for $g:=\tilde{g}\circ (\sigma\times \id_U)$.
By Proposition~\ref{pfwcsupp},
$g_*$ is $C^k_\K$.
\end{proof}
\section{\!\!Test
function groups and algebras of test functions}\label{sectfgps}
As in Section~\ref{secmaptf},
let $\F$ be a locally compact topological field
and $\K$ be a valued field
which is a topological extension
field of~$\F$. Let $r\in \N_0\cup\{\infty\}$.
In view of
Proposition~\ref{pfwcsupp} and Corollary~\ref{functcsupp},
we can re-use the arguments from
Section~\ref{secmapgps}
to obtain the following:
\begin{prop}\label{disctefgp}
Let $M$ be a paracompact, finite-dimensional $\F$-manifold
of class $C^r_\F$.
\begin{itemize}
\item[\n (a)]
If $A$ is a topological $\K$-algebra, then also
$C^r_c(M,A)$ is a topological
$\K$-algebra $($using pointwise multiplication$)$.
\item[\n (b)]
If $A$ is an associative
topological $\K$-algebra and $E$ a topological $A$-module,
then $C^r_c(M,E)$ is a topological
$C^r_c(M,A)$-module.
\item[\n (c)]
If $G$ is a $\K$-Lie group modeled on a
topological $\K$-vector space~$E$, then
\[
C^r_c(M,G):=\{\gamma\in C^r(M,G)\!: \,
\mbox{$\wb{\gamma^{-1}(G\,\take\,\{1\})}$ is compact}\,\}
\]
can be given a $C^\infty_\K$-manifold structure modeled on
the topological $\K$-vector space $C^r_c(M,E)$
in one and only
one way, such that $C^r_c(M,G)$
becomes a $\K$-Lie group and
such that $C^r_c(M,U_\phi):=C^r_c(M,G)\cap
(U_\phi)^M$ is open in $C^r_c(M,G)$
and
\[
C^r_c(M,\phi)\!: C^r_c(M,U_\phi)\to C^r_c(M,V_\phi),\quad
\gamma\mto \phi\circ \gamma
\]
is a $C^\infty_\K$-diffeomorphism onto
the open subset $C^r_c(M,V_\phi)\sub C^r_c(M,E)$,
for some chart $\phi\!: U_\phi\to V_\phi\sub E$
of~$G$ around~$1$ such that $\phi(1)=0$.\Punkt
\end{itemize}
\end{prop}
The Lie groups $C^r_c(M,G)$ described in~(c)
are also called {\em test function groups\/}.
\section{Differentiability of almost local mappings}\label{secalmloc}
We describe a criterion
ensuring differentiability properties
for mappings between open subsets of spaces
of vector-valued test functions (equipped with the
box topology).
Cf.\
\cite{SEC}, \cite{DIF}, \cite{INF}
and their precursor \cite{DRn}
for analogous results
in the real locally
convex case, based on the
locally convex direct limit topology.
\begin{numba}\label{settall}
Our general setting is the following:
$\F$ is the field of real numbers
or a local field, and $\K$ a valued field
which is a topological extension field of~$\F$.\,\footnote{The
case $\F=\C$ has to be excluded now, since we have to use
compactly supported cut-off functions.}
For $r,s,k \in \N_0\cup\{\infty\}$,
we are given a paracompact, finite-dimensional
$\F$-manifold~$M$ of class $C^r_\F$;
a paracompact,
finite-dimensional
$\F$-manifold $N$ of class $C^s_\F$;
and topological $\K$-vector spaces
$E$ and $F$.
We consider a mapping $f\!: P \to C^s_c(N,F)$,
defined on an open subset $P\sub C^r_c(M,E)$.
\end{numba}
Our investigations are stimulated by the following question:
\begin{numba}\label{questi} {\bf Question.\/}
If $f|_{P\cap C^r_K(M,E)}$
is of class $C^k_\K$ for all compact
subsets $K\sub M$, does it follow
that $f$ is $C^k_\K$\,?
\end{numba}
The answer is 
{\em negative}. For
example, the self-map
\[
f\!: C^\infty_c(\R,\R)\to C^\infty_c(\R,\R),\quad
\gamma\; \mto \;\, \gamma\circ \gamma\, - \,\gamma(0)
\]
of the space of real-valued test functions on the line
is discontinuous
at $\gamma=0$, although $f|_{C^\infty_K(\R,\R)}$
is smooth, for all compact subsets
$K\sub \R$ (see \cite{DIS}).\\[3mm]
The goal of this section is to describe
a simple {\em additional condition\/}
which prevents
the type of pathology just described.
As we shall see, Question~\ref{questi}
has an affirmative answer if we require
in addition that $f$ be ``almost local.''
Being almost local is a rather mild
condition, which is satisfied by most
of the mappings of relevance, for example
by all mappings encountered
in the construction of the Lie group
structure on groups of compactly supported
diffeomorphisms of finite-dimensional smooth
manifolds over the reals (see \cite{DIF}).
\begin{defn}\label{defnal}
(a) A map
$f\!: P \to C^s_c(N,F)$ (as in {\bf \ref{settall}})
is called {\em almost local\/}
if there exist locally finite covers
$(U_i)_{i\in I}$
of $M$ and $(V_i)_{i\in I}$
of~$N$
by relatively
compact, open sets
such that, for all $i\in I$
and $\gamma,\eta\in P$
with $\gamma|_{U_i}=\eta|_{U_i}$,
we have $f(\gamma)|_{V_i}=f(\eta)|_{V_i}$.\,\footnote{In other words,
$f(\gamma)|_{V_i}$ only depends
on $\gamma|_{U_i}$.}\\[2mm]
(b) A map $f\!: P\to C^s_c(N,F)$
is called {\em locally almost
local\/} if every $\gamma\in P$ has an open neighbourhood
$Q\sub P$ such that $f|_Q$ is almost local.\\[2mm]
(c) In the special case where $M=N$,
we call $f\!: P\to C^s_c(M,F)$
a {\em local\/} mapping if, for all $x\in M$
and $\gamma\in P$,
the element $f(\gamma)(x)$ only depends on the
germ of~$\gamma$ at~$x$.\footnote{More precisely,
we require
$f(\gamma)(x)=f(\eta)(x)$
for all $x\in M$ and $\gamma,\eta\in P$ with
the same germ at~$x$.}
It is easy to see that every local mapping
is almost local.
\end{defn}
Cf.\ already \cite[Defn.\,14.13]{KMS}
for the related notion of a ``local operator.''
\begin{thm}[Smoothness Theorem]\label{smoothy}
Let $f\!:
C^r_c(M,E)\supseteq P\to C^s_c(N,F)$
be a map as described in {\bf \ref{settall}}.
If
$f_K:=f|_{P\cap C^r_K(M,E)}$
is of class $C^k_\K$
for every compact subset $K\sub M$
and~$f$ is locally almost local,
then~$f$ is of class $C^k_\K$.
\end{thm}
{\bf Proof.} We proceed in steps.
\begin{numba}
Given $\gamma\in P$,
there exists an open neighbourhood~$Q$
of $\gamma$ in $P$ such that $f|_Q$
is almost local.
As $\gamma$ was arbitrary,
the assertion
will follow if we can show that $f|_W$
is of class $C^k_\K$
for some open neighbourhood $W$ of~$\gamma$
in~$Q$.
To this end, it suffices to show
that the mapping
$g\!: Q-\gamma\to C^s_c(N, F)$,
$g(\eta):=f(\gamma+\eta)-f(\gamma)$
is of class $C^k_\K$
on some open
zero-neighbourhood. As $f|_Q$ is
almost local, we find
locally finite covers $(U_i)_{i\in I}$
of $M$
and $(V_i)_{i\in I}$ of $N$,
with each $U_i$ and $V_i$ relatively compact
and open,
such that $f(\eta)|_{V_i}$ only depends
on $\eta|_{U_i}$, for all $\eta\in Q$.
Then apparently also
$g(\eta)|_{V_i}=g(\xi)|_{V_i}$
for all $\eta,\xi\in Q-\gamma$ such that
$\eta|_{U_i}=\xi|_{U_i}$,
showing that also~$g$ is almost local.
Furthermore, given a compact subset
$K\sub M$, the map
$g|_{(Q-\gamma)\cap C^r_K(M,E)}$
is of class $C^k_\K$,
since so is the restriction of~$f$
to $Q\cap C^r_{K\cup\sSup(\gamma)}(M,E)$.
We abbreviate $R:= Q-\gamma$.
\end{numba}
\begin{numba}
We pick a locally finite open cover $(\wt{U}_i)_{i\in I}$
of~$M$
such that
$\wb{U_i}\sub\wt{U}_i$
holds for the compact
closures, for all $i\in I$;
such a ``thickening''
exists by Lemma~\ref{lathicken}. 
For each $i\in I$,
we pick
a mapping $h_i\in C^r(\wt{U}_i,\F)$,
with compact support $K_i:=\Supp(h_i)$,
which is constantly~$1$ on $U_i$
(see Lemma~\ref{cutthroat} if $\F$ is a local field;
the real case is standard).
\end{numba}
\begin{numba}
By Example~\ref{patchd}, the family
$(\rho_i)_{i\in I}$ of restriction maps
$\rho_i\!:
C^r_c(M,E)\to C^r(\wt{U}_i, E)$
is a patchwork for $C^r_c(M,E)$.
We let $\rho\!: C^r_c(M, E)\to \bigoplus_{i\in I}
C^r(\wt{U}_i,E)=:S$
be the corresponding embedding
taking $\eta$ to $\sum_{i\in I} \rho_i(\eta)$.
Similarly, the family $(\sigma_i)_{i\in I}$
of restriction maps
$\sigma_i\!: C^s_c(N,F)\to C^s(V_i,F)$
is a patchwork for $C^s_c(N,F)$.
\end{numba}
\begin{numba}
The mapping $\rho$ being a topological embedding,
we find an open $0$-neighbourhood $H\sub S$
such that $\rho^{-1}(H)\sub R$.
The direct sum being
equipped with the box topology,
after shrinking~$H$ we may assume that $H=\bigoplus_{i\in I}
A_i$ for a family
$(A_i)_{i\in I}$ of open $0$-neighbourhoods
$A_i\sub C^r(\wt{U}_i, E)$.
The multiplication operator $\mu_{h_i}\!:
C^r(\wt{U}_i, E)
\to C^r_{K_i}(\wt{U}_i, E)$,
$\eta \mto h_i\cdot \eta$
is continuous linear
(in view of Lemma~\ref{corestr}
and Lemma~\ref{amend},
applied with a cover of coordinate neighbourhoods,
this assertion can be reduced to Lemma~\ref{multop}).
Hence, we find an open zero-neighbourhood
$W_i\sub A_i$ such that $h_i\cdot W_i\sub R$,
where we identify
$C^{\, r}_{K_i}(\wt{U}_i, E)$
with $C^{\, r}_{K_i}(M, E)\sub C^r_c(M,E)$
as a topological $\K$-vector space
in the natural way, extending functions by~$0$
(cf.\ Lemma~\ref{restrK}).
Then $W:=\rho^{-1}(\bigoplus_{i\in I} W_i)\sub R$
is an open zero-neighbourhood in
$C^r_c(M,E)$ such that $\rho_i(W)\sub W_i$
for each $i\in I$.
We define
\[
g_i\!: W_i \to C^s(V_i,F),
\;\;\; g_i:= \sigma_i \circ
g|_{R\cap C^r_{K_i}(M,E)}
\circ \mu_{h_i}|_{W_i}^R\, .
\]
Then $g_i$ is of class $C^k_\K$,
being a composition of $C^k_\K$-maps.
Note that $\sigma_i(g(\eta))
=g(\eta)|_{V_i}=g(h_i \cdot \eta)|_{V_i}
=g_i(\eta|_{\wt{U}_i})$ for each $\eta\in W$ and $i\in I$.
Thus $(g_i)_{i\in I}$ is compatible with $g|_W$
in the sense of Definition~\ref{defnpama}.
We have shown that $g|_W$
is a patched mapping which is of class $C^k_\K$ on the patches.
By Proposition~\ref{diffpatch},
$g|_W$ is of class $C^k_\K$.\Punkt
\end{numba}
\section{Smoothness of evaluation and composition}\label{seccompo}
We discuss
differentiability properties of
evaluation and composition of maps.
\begin{prop}\label{evalCk}
Let $\K$ be a locally compact
topological field, $k\in \N_0\cup\{\infty\}$,
$M$ a finite-dimensional
$C^k_\K$-manifold,
and $E$ a topological 
$\K$-vector space. Then the ``evaluation map''
\[
\ve\!: C^k(M,E)\times M\to E,\;\;\;\;
\ve(\gamma,x):=\gamma(x)
\]
is of class~$C^k_\K$.
\end{prop}
\begin{proof}
Given $x\in M$, let $\kappa\!: U\to V$
be a chart of~$M$ around~$x$,
where $V$ is an open subset
of the modeling space~$Z$ of~$M$.
Then $\ve(\gamma,\kappa^{-1}(y))=(\gamma\circ \kappa^{-1})(y)
=\tilde{\ve}(C^k(\kappa^{-1},E)(\gamma),y)$
for all $y\in V$,
in terms of the evaluation map
$\tilde{\ve}\!:
C^k(V,E)\times V\to E$
and the pullback
$C^k(\kappa^{-1},E)\!: C^k(M,E)\to C^k(V,E)$
which is continuous linear and thus smooth
(Lemma~\ref{pb2}).
It therefore suffices to consider
the case where $M=V$ is an open subset
of a finite-dimensional $\K$-vector space~$Z$.
The inclusion map $C^\infty(V,E)\to C^k(V,E)$
being continuous linear for all
$k\in \N_0$ (Remark~\ref{simplobs}\,(a)),
it also suffices to consider
finite~$k$. We proceed by induction.\\[3mm]
{\em The case $k=0$} is well known
(see, {\em e.g.}, \cite{Eng}, Thm.\,3.4.3 and Prop.\,2.6.11).\\[3mm]
{\em Induction step.} Given $k\in \N$,
suppose that the assertion of the lemma
holds if $k$ is replaced with $k-1$.
Given $(\gamma,x,\eta,y,t)\in (C^k(V,E)\times V)^{[1]}$
such that $t\not=0$, we calculate
\begin{eqnarray}
{\textstyle \frac{1}{t}}(\ve(\gamma+t\eta,x+ty)-\ve(\gamma,x))
& = &
{\textstyle\frac{1}{t}}(\gamma(x+ty)-\gamma(x))+\eta(x+ty)\nonumber\\
& = & \gamma^{[1]}(x,y,t)\,+\, \eta(x+ty)\,.\label{os}
\end{eqnarray}
Let $\ve_1\!: C^{k-1}(V^{[1]},E)\times V^{[1]}\to E$
denote the evaluation map, which is of class
$C^{k-1}_\K$ by induction.
Then, using Remark~\ref{simplobs}\,(b),
\[
\theta\!: (C^k(V,E)\times V)^{[1]}\to E,\;\;\;\;
\theta(\gamma,x,\eta,y,t):=\ve(\eta,x+ty)+\ve_1(\gamma^{[1]},(x,y,t))
\]
is a mapping of class $C^{k-1}_\K$.
In view of (\ref{os}),
we deduce that $\ve$ is of class~$C^1_\K$,
with $\ve^{[1]}=\theta$ of class~$C^{k-1}_\K$,
and thus $\ve$ is of class~$C^k_\K$,
which completes the inductive proof.
\end{proof}
Let us turn to the composition map
now. We shall show:
\begin{prop}\label{compcomp}
Let $\K$ be a
locally compact topological field,
$E$ a topological $\K$-vector space,
$r,k\in \N_0\cup\{\infty\}$,
$M$ be a finite-dimensional $C^r_\K$-manifold,
$F$ a finite-dimensional $\K$-vector space,
$U\sub F$ be open,
and $K\sub M$ compact.
Then the composition map
\[
\Gamma\!:
C^{r+k}(U,E)\times C^r_K(M,U)\to
C^r(M,E),\;\;\;
\Gamma(\gamma,\eta):=\gamma\circ \eta
\]
is of class~$C^k_\K$.
If $k\geq 1$,
then
\begin{equation}\label{formdGK}
d\Gamma(\gamma,\eta,\gamma_1,\eta_1)=d\gamma\circ (\eta,\eta_1)
+\gamma_1\circ \eta
\end{equation}
for all $\gamma,\gamma_1\in C^{r+k}(U,E)$,
$\eta\in C^r_K(M,U)$, and
$\eta_1\in C^r_K(M,F)$.
\end{prop}
If $U=F$, then
we need not assume that~$M$ be finite-dimensional.
In this case, we have:
\begin{prop}\label{comparbit}
Let $\K$ be a
locally compact topological field,
$E$ a topological $\K$-vector space,
$r,k\in \N_0\cup\{\infty\}$,
$M$ a $C^r_\K$-manifold,
modeled on an arbitrary topological $\K$-vector space,
and $F$ be a finite-dimensional $\K$-vector space.
Then the composition map
\[
\Gamma\!:
C^{r+k}(F,E)\times C^r(M,F)\to
C^r(M,E),\;\;\;
\Gamma(\gamma,\eta):=\gamma\circ \eta
\]
is of class~$C^k_\K$.
If $k\geq 1$,
then
\begin{equation}\label{formdGarbit}
d\Gamma(\gamma,\eta,\gamma_1,\eta_1)=d\gamma\circ (\eta,\eta_1)
+\gamma_1\circ \eta
\end{equation}
for all $\gamma,\gamma_1\in C^{r+k}(F,E)$ and
$\eta,\eta_1 \in C^r(M,F)$.
\end{prop}
For finite-dimensional~$M$, both propositions are
immediate
consequences of the following
technical result, which we prove now.
A direct proof for Proposition~\ref{comparbit}
(including the case of infinite-dimensional~$M$)
is given in Appendix~\ref{appcomparbit}.
\begin{la}\label{compsmooth}
Let $\K$ be a
locally compact topological field,
$E$ a topological $\K$-vector space,
$r,k\in \N_0\cup\{\infty\}$,
$M$ a finite-dimensional $C^r_\K$-manifold,
$F$ a finite-dimensional $\K$-vector space,
$U$ an open subset of $F$,
$K$ a compact subset of~$M$,
and $Y\sub K^0$ be a non-empty, open subset.
Let $H$ be a finite-dimensional
$\K$-vector space, and $P\sub H$
be open.
Then
\[
\Theta\!:
C^{r+k}(U\times P,E)\times \lfloor K,U\rfloor_r\times P\to
C^r(Y,E),\;\;\;
\Theta(\gamma,\eta,p):=\gamma(\sbull,p)\circ \eta|_Y\,,
\]
where $\lfloor K,U\rfloor_r\sub C^r(M,F)$,
is a mapping of class~$C^k_\K$.
If $k\geq 1$,
then
\begin{eqnarray}
\lefteqn{\!\!\!\!\!\!\!\!\!\!\!\!\!\!\!\!\Theta^{[1]}((\gamma,\eta,p),
\,(\gamma_1,\eta_1,p_1),\,t)\quad\quad\quad\quad}\nonumber \\
&=&
\gamma^{[1]}((\sbull,p),\, (\sbull,p_1),\, t)
\circ (\eta,\eta_1)|_Y
+\gamma_1(\sbull,p+tp_1)\circ (\eta+t\eta_1)|_Y\label{formcomp3}
\end{eqnarray}
for all $((\gamma,\eta,p),\, (\gamma_1,\eta_1,p_1),\,t)\in
(C^{r+k}(U\times P,E)\times \lfloor K,U\rfloor_r\times P)^{[1]}$.\\[3mm]
Hence, as a special case, the map
\[
\Gamma\!:
C^{r+k}(U,E)\times \lfloor K,U\rfloor_r\to
C^r(Y,E),\;\;\;
\Gamma(\gamma,\eta):=\gamma\circ \eta|_Y
\]
is of class~$C^k_\K$.
If $k\geq 1$,
then
\begin{equation}\label{formcomp1}
\Gamma^{[1]}((\gamma,\eta),\,(\gamma_1,\eta_1),\,t)=
\gamma^{[1]}(\sbull,t)\circ (\eta,\eta_1)|_Y+\gamma_1\circ (\eta+t\eta_1)|_Y
\end{equation}
for all $((\gamma,\eta),\,(\gamma_1,\eta_1),\,t)\in
(C^{r+k}(U,E)\times \lfloor K,U\rfloor_r)^{[1]}$.
In particular,
\begin{equation}\label{formdG}
d\Gamma((\gamma,\eta),\,(\gamma_1,\eta_1))=d\gamma\circ (\eta,\eta_1)|_Y
+\gamma_1\circ \eta|_Y
\end{equation}
for all $\gamma,\gamma_1\in C^{r+k}(U,E)$,
$\eta\in \lfloor K,U\rfloor_r$, and
$\eta_1\in C^r(M,F)$.
\end{la}
\begin{rem}
In the real or complex locally convex case,
the desired properties of $\Gamma$ can be established
directly,
without recourse to
parameters.
In the general case envisaged here,
a direct induction without parameter sets
(based on Proposition~\ref{evalCk} and Lemma~\ref{halfcartesian})
would only show that $\Gamma$ is $C^k_\K$ when
$C^{r+k}(U,E)$ is replaced with
$C^{r+\frac{1}{2}k(k+1)}(U,E)$,
due to the loss in the order of differentiability
in Lemma~\ref{halfcartesian}.
\end{rem}
{\bf Proof of Lemma~\ref{compsmooth}.}
Clearly, we only need to prove the assertions
concerning~$\Theta$: then also $\Gamma$ will have the
asserted properties.
It suffices to consider finite
$k\in \N_0$ (cf.\ Remark~\ref{simplobs}\,(a)).
We may also assume that $r\in \N_0$ (cf.\
proof of Proposition~\ref{globcruc}).
Thus, we assume that both $r$ and~$k$ are
finite, and prove the assertion
by induction on~$k$.
\begin{center}
{\bf The case {\boldmath $k=0$}}
\end{center}
We proceed by induction on~$r$.
Let us suppose that $r=0$ first.
We recall that the topology we have defined on spaces
of $C^0$-maps coincides with the compact-open topology
(Remark~\ref{compareco}).
For $(\gamma,\eta,p)\in
C(U\times P,E)\times \lfloor K,U\rfloor \times P$,
we have
\begin{equation}\label{gtcts}
\Theta(\gamma,\eta,p)=
\wt{\Gamma}\left(\gamma,\;\eta|_Y\!\times \!\id_P\right)(\sbull,p)\, ,
\quad \mbox{where}
\end{equation}
\[
\wt{\Gamma}\!:
C(U\!\times \!P,\, E)_{c.o.}\times C(Y\!\times\! P,\,U\!\times \!P)_{c.o.}\to
C(Y\!\times\! P,\, E)_{c.o.},\quad \wt{\Gamma}(\sigma,\tau):=\sigma\circ \tau
\]
is the composition map, which
is continuous since $U\times P$
is locally compact \cite[Thm.\,3.4.2]{Eng}.
It easily follows from the definition of
the compact-open topology that the mapping\linebreak
$\lfloor K,U\rfloor\to
C(Y\times P,U\times P)_{c.o.}$, $\eta\mto \eta|_Y\times \id_P$
is continuous.
The map
\[
f^\vee : \, P\to C(Y,E)\, , \qquad f^\vee(p):=f(\sbull,p)
\]
is continuous for $f\in C(Y\times P,E)$, and also the map
$C(Y\times P,E)\to C(P,C(Y,E))$, $f\mto (f^\vee\!: p\mto f(\sbull,p))$
is continuous
\cite[Thm.\,3.4.7]{Eng}.
Furthermore, $P$ being locally compact, the evaluation map
$\ve\!: C(P,C(Y,E))\times P\to
C(Y,E)$ is continuous
(cf.\ \cite{Eng}, Thm.\,3.4.3
and Prop.\,1.6.11).
Reading
(\ref{gtcts}) as
$\Theta(\gamma,\eta,p)=\ve(\wt{\Gamma}(\gamma,\,
\eta|_Y\!\times \!\id_P)^\vee,\,
p)$,
we see that
$\Theta$ is continuous.\\[3mm]
{\em Induction step on $r$.}
Let $r\in \N$, and suppose that the proposition
holds for $k=0$, when $r$ is replaced with
$r-1$.
It then suffices to show continuity of
$\Theta$ in the case
where $M$ is an open subset of its
modeling space~$Z$.
In fact, suppose that
$M$ is a $C^r_\K$-manifold.
For each $y\in Y$, there exists
a chart $\kappa_y\!: W_y\to V_y\sub Z$ of~$Y$
around~$y$. Let $L_y\sub W_y$ be a compact
neighbourhood of~$y$, and $K_y:=\kappa_y(L_y)$;
let $Y_y:=K_y^0$ be the interior of~$K_y$.
Since $(L_y^0)_{y\in Y}$
is an open cover of~$Y$,
we deduce with Lemma~\ref{amend}
that $\Theta$ will be continuous
if we can show that
\[
h_y\!:
C^r(U\times P,E)\times \lfloor K,U\rfloor_r\times P\to
C^r(Y_y,E),\quad
h_y(\gamma,\eta,p):=\Theta(\gamma,\eta,p)\circ \kappa_y^{-1}|_{Y_y}
\]
is continuous, for all $y\in Y$.
But
\begin{equation}\label{noend}
h_y(\gamma,\eta,p) =\gamma(\sbull,p)\circ (\eta\circ \kappa^{-1}_y)|_{Y_y}
=\Theta_y(\gamma,\eta\circ \kappa^{-1}_y,p)
\end{equation}
with $\Theta_y\!:
C^r(U\times P,E)\times \lfloor K_y,U\rfloor_r\times P\to
C^r(Y_y,E)$, $\Theta_y(\gamma,\sigma,p):=\gamma(\sbull,p)\circ
\sigma|_{Y_y}$,
where $\lfloor K_y,U\rfloor_r\sub C^r(V_y,F)$.
Note that
the pullback $C^r(M,F)\to C^r(V_y,F)$, $\eta\mto\eta\circ\kappa_y^{-1}$
is continuous linear (Lemma~\ref{pb2})
and takes the open set $\lfloor K,U\rfloor_r$ into
$\lfloor K_y,U\rfloor_r$.
Thus (\ref{noend}) shows that $h_y$ will be continuous
if each $\Theta_y$ is continuous. Since $V_y$ is open in~$Z$,
this completes the reduction step to the case
where $M$ is open in~$Z$.\\[3mm]
To complete the induction step on~$r$ in the case $k=0$,
by the preceding we may assume now that $M$ is an
open subset of~$Z$.
The map
$\Theta\!:
C^r(U\times P,E)\times \lfloor K,U\rfloor_r\times P\to
C^r(Y,E)$ is continuous as a map into $C(Y,E)$,
by the case $r=0$ already settled and Remark~\ref{simplobs}\,(a).
Hence, in view of Remark~\ref{simplobs}\,(b),
$\Theta$ will be continuous if we can show
that the map
\[
C^r(U\times P,E)\times \lfloor K,U\rfloor_r\times P\to
C^{r-1}(Y^{[1]},E),\quad (\gamma,\eta,p)\mto 
\Theta(\gamma,\eta,p)^{[1]}
\]
is continuous at each given element
$(\gamma_0,\eta_0,p_0)$ in its domain, where
\begin{equation}\label{formkl1}
\Theta(\gamma,\eta,p)^{[1]}(x,y,t)=\gamma^{[1]}((\eta(x),p),\,
(\eta^{[1]}(x,y,t),0),\, t)
\end{equation}
for all $(x,y,t)\in Y^{[1]}$, by the Chain Rule.
Let $(x_0,y_0,t_0)\in Y^{[1]}$ be given.
There exist open neighbourhoods
$U_1\sub U$ of~$\eta_0(x_0)$,
$U_2\sub F$ of $\eta_0^{[1]}(x_0,y_0,t_0)$
and $U_3\sub \K$ of~$t_0$,
such that $U_1\times P\times U_2 \times
\{0\}\times U_3\sub (U\times P)^{[1]}$.
Then
\[
\rho\!: C^{r-1}((U\times P)^{[1]},E)
\to C^{r-1}(U_1\times U_2\times U_3\times P,E),\quad
\rho(\xi)(x,y,t,p):=\xi(x,p,y,0,t)
\]
is a continuous linear map (Lemma~\ref{pb2}).
There exist open neighbourhoods
$V_1\sub Y$ of $x_0$, $V_2\sub Z$ of~$y_0$, and
$V_3\sub U_3$ of~$t_0$ such that
$\eta_0(V_1)\sub U_1$, $V_1\times V_2\times V_3\sub Y^{[1]}$, and
$\eta_0^{[1]}(V_1\times V_2\times V_3)\sub U_2$.
There exist compact neighbourhoods
$K_1\sub V_1$ of~$x_0$, $K_2\sub V_2$ of~$y_0$,
and $K_3\sub V_3$ of~$t_0$.
Set $Y_i:=K_i^0$ for $i=1,2,3$.
By induction, the map
\[
\wt{\Theta}\!:\!
C^{r-1}(U_1\!\times\!U_2\!\times U_3\!\times \!P,E)
\!\times \!\lfloor K_1\!\times \!K_2\times\!K_3,
U_1\!\times\!U_2\!\times\!U_3 \rfloor_{r-1}
\!\!\times \!P
\!\to C^{r-1}(Y_1\!\times \!Y_2 \!\times\!Y_3,E)
\]
taking $(\sigma,\tau,p)$ to
$\sigma(\sbull,p)\circ \tau|_{Y_1\times Y_2\times Y_3}$
is continuous; here
$\lfloor K_1\times K_2\times K_3,U_1\times U_2 \times U_3 \rfloor_{r-1}
\sub C^{r-1}(V_1\times V_2\times V_3,F\times F\times \K)$.
Note that
\[
\Omega:=\left\{
\eta\in \lfloor K,U\rfloor_r\cap \lfloor K_1,U_1\rfloor_r
\!: \eta^{[1]}|_{V_1\times V_2\times
V_3}\in \lfloor K_1\times K_2\times K_3, U_2\rfloor_{r-1}
\right\}
\]
is an open neighbourhood of~$\eta_0$ in
$\lfloor K,U\rfloor_r$.
Furthermore, by Lemma~\ref{pb2}, the map
\[
h\!:\Omega\to
\lfloor K_1\times K_2\times K_3,U_1\times U_2\times U_3\rfloor_{r-1},\quad
h(\eta):=\left( \eta\circ \pi_1,\, \eta^{[1]}|_{V_1\times V_2\times V_3},\,
\pi_3\right)
\]
is continuous, where $\pi_1\!: V_1\times V_2\times V_3\to V_1$
and $\pi_3\!: V_1\times V_2\times V_3\to V_3\sub U_3$
are the coordinate projections.
Since,
by (\ref{formkl1}) and the definition of~$\rho$ and
$h$, we have
\[
\Theta(\gamma,\eta,p)^{[1]}|_{Y_1\times Y_2\times Y_3}
=\wt{\Theta}(\rho(\gamma^{[1]}),h(\eta))
\]
for all $p\in P$,
$\gamma\in C^r(U\times P,E)$, and $\eta\in \Omega$,
we see that $(\gamma,\eta,p)\mto
\Theta(\gamma,\eta,p)^{[1]}|_{Y_1\times Y_2\times Y_3}\in C^{r-1}(Y_1\times Y_2\times Y_3,E)$
is continuous at $(\gamma_0,\eta_0,p_0)$.
Since $Y^{[1]}$ can be covered by sets
of the form $Y_1\times Y_2\times Y_3$ as before,
using Lemma~\ref{amend} we now deduce that
the mapping $(\gamma,\eta,p)\mto \Theta(\gamma,\eta,p)^{[1]}\in
C^{r-1}(Y^{[1]},E)$ is continuous at $(\gamma_0,\eta_0,p_0)$,
as desired.
\begin{center}
{\bf Induction step on~{\boldmath $k$}}
\end{center}
Let $k\in \N$, and
suppose that the assertion of
the lemma holds when $k$ is replaced with
$k-1$, for all $r\in \N_0$. Let $r\in \N_0$.
Given an element
$((\gamma,\eta,p),\,(\gamma_1,\eta_1,p_1),\, t)\in
\Omega:=(C^{r+k}(U\times P,E)\times \lfloor K,U\rfloor_r\times P)^{[1]}$
such that $t\not=0$,
we calculate for $x\in Y$:
\begin{eqnarray}
\lefteqn{\!\!\!\!\!\!\!\!\!{\textstyle \frac{1}{t}}
\left(\Theta(\gamma+t\gamma_1,\eta+t\eta_1,p+tp_1)-\Theta(\gamma,\eta,p)\right)
(x)\quad\quad\quad}\nonumber\\
&=&
{\textstyle \frac{1}{t}}\Big(
\gamma(\eta(x)+t\eta_1(x), p+tp_1)-\gamma(\eta(x),p)\Big)
\,+\,\gamma_1(\eta(x)+t\eta_1(x),p+tp_1)\nonumber\\
&=&
\gamma^{[1]}((\eta(x),p),\, (\eta_1(x),p_1),\, t)
\,+\, \gamma_1(\eta(x)+t\eta_1(x),p+tp_1)\,,\label{hbreak}
\end{eqnarray}
in accordance with (\ref{formcomp3}).
Since $\Theta$ is $C^{k-1}_\K$ and hence continuous
as a consequence of the induction hypothesis,
in order that $\Theta$ be $C^k_\K$,
it therefore only remains to show that the mapping
$\Omega\to C^r(Y,E)$
described in (\ref{formcomp3}),
let us call it~$g$,
is of class $C^{k-1}_\K$
(then $g=\Theta^{[1]}$).
Since $\Theta$ is $C^{k-1}_\K$,
the map $g$ is $C^{k-1}_\K$
on an open neighbourhood of
each given element $((\bar{\gamma},\bar{\eta},\bar{p}),\,
(\bar{\gamma}_1,\bar{\eta}_1,\bar{p}_1),\,
\bar{t})\in \Omega$,
provided $\bar{t}\not=0$. It remains to consider the case
where $\bar{t}=0$.
There is a balanced, open zero-neighbourhood $W\sub F$ such that
$\bar{\eta}(K)+W+W+W\sub U$.
Next, there are open neighbourhoods
$P_0\sub P$ of~$\bar{p}$, $P_1\sub H$
of $\bar{p}_1$, and $r\in \;]0,1]$ such that
\[
P_0+P_2P_1\sub P\quad \mbox{and thus}\quad
P_0\times P_1\times P_2\sub P^{[1]}\,,
\]
where $P_2:=\{t\in \K\!: |t|\leq r\}$.
After shrinking~$r$, we may assume that
furthermore $P_2\cdot \bar{\eta}_1(K) \sub W$.
We let $U_0:=\bar{\eta}(K)+W\sub U$ and
$U_1:=\bar{\eta}_1(K)+W$.
Then
\[
U_0+P_2U_1\sub \bar{\eta}(K)+W+P_2\bar{\eta}_1(K)+P_2W
\sub \bar{\eta}(K)+W+W+W\sub U
\]
and hence $U_0\times U_1\times P_2\sub U^{[1]}$.
Furthermore, we have
\[
(\bar{\eta},\bar{\eta}_1)\in
\lfloor K,U_0\times U_1\rfloor_r\sub C^r(M,F\times F)\,.
\]
Then
$U_0\times P_0\times U_1\times P_1\times P_2\sub (U\times P)^{[1]}$,
and the map
\begin{eqnarray*}
\rho \!: C^{r+k}(U,E) \quad& \to &
C^{r+(k-1)}((U_0\times U_1)\times (P_0\times P_1\times P_2),E),\\
\rho(\gamma)((u_0,u_1),\, (p,p_1,t))&:=& \gamma^{[1]}((u_0,p),\,
(u_1,p_1),\, t)
\end{eqnarray*}
is continuous linear by Remark~\ref{simplobs}\,(b) and
Lemma~\ref{pb2}.
Hence $\rho$ is $C^{k-1}_\K$.
By the induction hypothesis, the map
\[
\wt{\Theta}\!:
C^{r+(k-1)}((U_0\!\times\!U_1)\times (P_0\!\times\!P_1\times\!P_2),E)
\times \lfloor K,U_0\!\times \!U_1\rfloor_r
\times (P_0\!\times \!P_1\!\times \!P_2)
\to C^r(Y,E)
\]
taking $(\xi,\zeta,(p,p_1,t))$ to $\xi(\sbull,(p,p_1,t))\circ \zeta$
is of class $C^{k-1}_\K$.
The set
\[
C^{r+k}(U\times P,E)\times \lfloor K,U_0\rfloor_r
\times P_0\times C^{r+k}(U\times P,E)\times \lfloor K,U_1\rfloor_r
\times P_1\times P_2
\]
is an open neighbourhood of
$(\bar{\gamma},\bar{\eta},\bar{p},\bar{\gamma}_1,\bar{\eta}_1,\bar{p}_1,0)$
in the domain $\Omega$ of~$g$.
For all elements\linebreak
$(\gamma,\eta,p,\gamma_1,\eta_1,p_1,t)$
in this open neighbourhood, we
have
\[
g(\gamma,\eta,p,\gamma_1,\eta_1,p_1,t)=\wt{\Theta}(\rho(\gamma),\,
(\eta,\eta_1),\, (p,p_1,t))\;+\;
\Theta(\gamma_1,\eta+t\eta_1,p+tp_1)\,,
\]
showing that $g$ is $C^{k-1}_\K$
on this open neighbourhood.
This completes the proof.\vspace{2mm}\Punkt

\noindent
Proposition~\ref{compcomp} and Proposition~\ref{comparbit}
(for finite-dimensional~$M$) now readily follow:\\[3mm]
{\bf Proof of Proposition~\ref{comparbit}
for finite-dimensional~{\boldmath $M$.}}
Every $x\in M$ has a compact neighbourhood~$K_x$;
let $Y_x:=K_x^0$ be its interior.
Then $(Y_x)_{x\in M}$ is an open cover of~$M$.
Let $\rho_x\!: C^r(M,E)\to C^r(Y_x,E)$, $\gamma\mto \gamma|_{Y_x}$
be the restriction map.
Then, as a consequence
of Lemma~\ref{corestr} and Lemma~\ref{amend},
the composition map
$\Gamma\!: C^{r+k}(F,E)\times C^r(M,F)\to C^r(M,E)$
will be of class $C^k$ if we can show that
\[
\rho_x\circ \Gamma\!: C^{r+k}(F,E)\times C^r(M,F)\to C^r(Y_x,E)
\]
is of class $C^s$, for each $x\in M$. However, we have
$\rho_x\circ \Gamma=\Gamma_x$, where
\[
\Gamma_x\!: C^{r+k}(F,E)\times \lfloor K_x,F\rfloor_r
\to C^r(Y_x,E),\quad \Gamma_x(\gamma,\eta):=
\gamma\circ (\eta|_{Y_x})
\]
is of class $C^k$ by Lemma~\ref{compsmooth}.
Hence $\Gamma$ is of class~$C^k$. Suppose that $k\geq 1$ now.
The mapping $\rho_x$ being continuous linear,
we have
$d\Gamma_x=d(\rho_x\circ \Gamma)=\rho_x\circ d\Gamma$.
Hence (\ref{formdG}) implies that, for all $\gamma,\gamma_1\in C^{r+k}(F,E)$,
$\eta,\eta_1\in C^r(M,F)$:
\[
d\Gamma((\gamma,\eta),\, (\gamma_1,\eta_1))|_{Y_x}=
d\gamma\circ (\eta,\eta_1)|_{Y_x}+\gamma_1\circ \eta|_{Y_x}
=\left(d\gamma\circ (\eta,\eta_1)+\gamma_1\circ \eta\right)|_{Y_x}
\]
for all $x\in M$, entailing that
$d\Gamma((\gamma,\eta),\, (\gamma_1,\eta_1))=
d\gamma\circ (\eta,\eta_1)+\gamma_1\circ \eta$,
as asserted.\vspace{3mm}\Punkt

\noindent
{\bf Proof of Proposition~\ref{compcomp}.}
Let $K_x$, $Y_x$ and
$\rho_x\!: C^r(M,E)\to C^r(Y_x,E)$
be as in the preceding proof.
In order that the composition map
\[
\Gamma\!:
C^{r+k}(U,E)\times C^r_K(M,U) \to
C^r(M,E)
\]
be of class $C^k$, we only need to show
that $\rho_x\circ \Gamma$ is of class $C^k$ for all
$x\in M$. But $\rho_x\circ \Gamma=\Gamma_x\circ \lambda_x$,
where
\[
\Gamma_x\!: C^{r+k}(U,E)\times \lfloor K_x, U\rfloor_r
\to C^r(Y_x,E),\quad \Gamma_x(\gamma,\eta):=
\gamma\circ (\eta|_{Y_x}^U)
\]
is of class~$C^k$ by Lemma~\ref{compsmooth}, and
\[
\lambda_x\!: C^{r+k}(U,E)\times C^r_K(M,U)\to
C^{r+k}(U,E)\times \lfloor K_x,U\rfloor_r,\quad
(\gamma,\eta)\mto (\gamma,\eta)
\]
is obtained by restricting and co-restricting
a continuous linear map to open sets
and therefore smooth. Hence $\rho_x\circ \Gamma=\Gamma_x\circ \lambda_x$
is of class~$C^k$, being a composition of $C^k$-maps.
The desired formula for $d\Gamma$
(if $k\geq 1$) can now be deduced
as in the preceding proof.\Punkt
\section{Basic exponential law for smooth mappings}\label{secexplaw}
In this section, we establish an exponential
law for smooth mappings on products
of suitable manifolds, and related results.
\begin{la}\label{halfcartesian}
Let $\K$ be a topological field,
$r,k\in \N_0\cup\{\infty\}$,
$M$ and $N$ be $C^{r+k}_\K$-manifolds
modeled on topological $\K$-vector spaces,
and $E$ be a topological
$\K$-vector space.
Then the following holds:
\begin{itemize}
\item[\n (a)]
For each mapping
$f\!: M\times N\to E$ of class
$C^{r+k}_\K$, the associated mapping
\[
f^\vee \!: M\to C^r(N,E),\;\;\;
f^\vee(x):=f(x,\sbull)
\]
is of class~$C^k_\K$.
\item[\n (b)]
The linear map
$\Phi\!: C^{r+k}(M\times N,E)\to C^k(M,C^r(N,E))$,
$\Phi(f):=f^\vee$
is continuous.
\end{itemize}
\end{la}
\begin{proof}
The lemma will hold in general if we can prove
the case where $M$ and $N$ are open subsets
of topological $\K$-vector spaces $X$ and $Y$,
respectively. In fact, suppose that $M$ and $N$ are
$C^{r+k}_\K$-manifolds.
Let $f\!: M\times N\to E$ be a $C^{r+k}_\K$-map.
The mapping $f^\vee$ will be of class $C^k_\K$
if we can show that it is $C^k_\K$ on some open neighbourhood
of each given point $x_0\in M$.
Given $x_0$, we let $\phi\!: U_\phi\to V_\phi\sub X$
be a chart of~$M$ around~$x_0$. Let $\cA$ be an atlas
for $N$, of charts $\psi\!: U_\psi\to V_\psi\sub Y$ of~$N$.
As a consequence of Lemma~\ref{corestr},
Lemma~\ref{amend} and Lemma~\ref{pb2},
the map $f^\vee|_{U_\phi}$
is $C^k_\K$ if and only if
\begin{equation}\label{ifff}
C^r(\psi^{-1},E)\circ f^\vee|_{U_\phi}\!: U_\phi\to C^r(V_\psi,E),\quad
x\mto f^\vee(x)\circ \psi^{-1}
\end{equation}
is $C^k_\K$ for each $\psi\in \cA$.
This holds if and only if
\[
C^r(\psi^{-1},E)\circ f^\vee\circ \phi^{-1}\!:
V_\phi\to C^r(V_\psi,E)
\]
is $C^k_\K$, for each $\psi\in \cA$.
Now, for given $\psi$, the latter
map coincides with $g^\vee$,
where $g:=f\circ (\phi^{-1}\times \psi^{-1})\!: V_\phi\times V_\psi\to E$,
and here $V_\phi\sub X$ and $V_\psi\sub Y$ are open subsets
of topological $\K$-vector spaces. It therefore suffices
to show that each $g^\vee$ is of class $C^k_\K$.\\[3mm]
To see that also (b) can be reduced to the case
of open subsets
of topological vector spaces,
note that,
as a consequence of
Lemma~\ref{pb2},
Lemma~\ref{amend},
Lemma~\ref{linearcase}
and
Lemma~\ref{Ctopinit}, 
the topology on $C^k(M,C^r(N,E))$
is initial with respect to the family
of mappings
\[
h_{\phi,\psi}:=C^k(V_\phi,C^r(\psi^{-1},E))\circ
C^k(\phi^{-1},C^r(N,E))\!: C^k(M,C^r(N,E))\to C^k(V_\phi,C^r(V_\psi,E))
\]
taking $g\in C^k(M,C^r(N,E))$ to $C^r(\psi^{-1},E)\circ g\circ \phi^{-1}$,
where $\phi$ and $\psi$
range through the charts
of $M$ and $N$, respectively.
Since $h_{\phi,\psi}(f^\vee)=(f\circ (\phi^{-1}\times \psi^{-1}))^\vee
=\left(C^{r+k}(\phi^{-1}\times \psi^{-1},E)(f)\right)^\vee$,
where $C^{r+k}(\phi^{-1}\times \psi^{-1},E)\!:
C^{r+k}(M\times N,E)\to C^{r+k}(V_\phi\times V_\psi,E)$
is continuous and takes $f$
to a mapping defined on the product
$V_\phi\times V_\psi$ of open subsets of~$X$ and $Y$,
it suffices to prove (b) for mappings on such
products $V_\phi\times V_\psi$.\\[3mm]
By the preceding, we may assume for the rest of the proof
that
$U:=M \sub X$ and $V:=N \sub Y$ are open subsets of topological
vector spaces.
Recall from Remark~\ref{simplobs}\,(a)
that $C^\infty(V,E)={\pl}_{r\in \N_0}C^r(V,E)$.
Accordingly, $C^k(U,C^\infty(V,E))={\pl}_{r\in \N_0}
C^k(U,C^r(V,E))$\vspace{-.8mm}
(Lemma~\ref{inpl}, Lemma~\ref{Ctopinit}).
It therefore suffices to prove the assertions
when $r\in \N_0$. By a similar argument,
we may assume that~$k$ is finite.
The proof is by induction on~$k\in \N_0$.\\[3mm]
{\em The case $k=0$.}
If $r=0$, then (a) and (b) are special
cases of \cite{Eng}, Thm.\,3.4.1 and 3.4.7,
respectively. To proceed by induction on~$r$,
suppose that
$r\in \N$, and suppose the assertion of the lemma
holds when $r$ is replaced with~$r-1$.
The topology on $C^r(V,E)$ is initial
with respect to the maps $\alpha\!:
C^r(V,E)\to C^{r-1}(V,E)$, $\gamma\mto \gamma$
and $\beta\!: C^r(V,E)\to C^{r-1}(V^{[1]},E)$,
$\beta(\gamma):=\gamma^{[1]}$
(Remark~\ref{simplobs}).
Hence the topology on $C(U,C^r(V,E))$
is initial with respect to the mappings
$C(U,\alpha)$ and $C(U,\beta)$ (Lemma~\ref{Ctopinit}).\\[3mm]
(a) Let $f\!: U\times V\to E$ be a $C^r_\K$-map,
and $f^\vee\!: U\to C^r(V,E)$ be as above.
By the induction hypothesis,
$U\to C^{r-1}(V,E)$, $x\mto f(x,\sbull)=\alpha\circ f^\vee$
is a continuous mapping.
In view of the preceding, $f^\vee\!: U\to C^r(V,E)$
will be continuous if we can show that also
\[
\beta\circ f^\vee\!: U\to C^{r-1}(V^{[1]},E),\;\;\;
x\mto (f^\vee(x))^{[1]}
\]
is continuous.
However,
\begin{equation}\label{eqnpsif}
\psi_f\!: U\times V^{[1]}\to E,\;\;
\psi_f(x,(v,h,t)):=f^{[1]}((x,v),\,(0,h),\, t)
\;\mbox{for $x\in U$, $(v,h,t)\in V^{[1]}$}
\end{equation}
is of class~$C^{r-1}_\K$,
being a partial map of $f^{[1]}$.
Since $(\beta\circ f^\vee)(x)=\psi_f(x,\sbull)=(\psi_f)^\vee(x)$,
the map $\beta\circ f^\vee=(\psi_f)^\vee$
is continuous by the induction hypothesis.
Thus $f^\vee$ is continuous, and thus (a) holds in the
$C^r_\K$-case, when $k=0$.\\[3mm]
(b) As an immediate consequence of the induction hypothesis,
the mapping
$C(U,\alpha)\circ \Phi\!:
C^r(U\times V,E)\to C(U,C^{r-1}(V,E))$,
$f\mto (x\mto f(x,\sbull))$ is continuous.
The mapping
\[
\Psi\!: C^r(U\times V,E)\to C^{r-1}(U\times V^{[1]},E),\;\;\;
\Psi(f):=\psi_f
\]
(with $\psi_f$ as in (\ref{eqnpsif}))
is continuous by Remark~\ref{simplobs}\,(b)
and Lemma~\ref{pullback}. Furthermore,
\[
\Xi\!: C^{r-1}(U\times V^{[1]},E)\to C(U,C^{r-1}(V^{[1]},E)),
\;\;\;\;\Xi(g)(x):=g(x,\sbull)\;\;\mbox{for $x\in U$}
\]
is continuous, by induction.
Thus $C(U,\beta)\circ \Phi=\Xi\circ \Psi$
is continuous. The topology on $C(U,C^r(V,E))$
being initial with respect to
$C(U,\alpha)$ and $C(U,\beta)$,
we deduce from the preceding
that $\Phi\!: C^r(U\times V,E)\to C(U,C^r(V,E))$
is continuous. Thus also the $C^r_\K$-case of (b)
is established, when $k=0$.\\[3mm]
{\em Induction step on~$k$.}
Let $k\in \N$, and suppose that the assertions
of the lemma hold for all~$r$,
when $k$ is replaced with~$k-1$.
Let $f\!: U\times V\to E$ be a mapping of class~$C^{r+k}_\K$.
As a consequence of the induction hypothesis
and Remark~\ref{simplobs}\,(a),
the map
$f^\vee\!: U\to C^r(V,E)$
is of class $C^{k-1}_\K$, and $C^{r+k}(U\times V,E)\to C^{k-1}(U,
C^r(V,E))$, $f\mto f^\vee$ is a continuous linear map.
We now observe that
\begin{equation}\label{tnot0}
{\textstyle\frac{1}{t}}(f^\vee(x+ty)-f^\vee(x))(v)=
{\textstyle \frac{1}{t}}(f(x+ty,v)-f(x,v))=f^{[1]}(x,v,y,0,t)
\end{equation}
for all $v\in V$ and
$(x,y,t)\in U^{[1]}$ such that $t\not=0$.
The mapping $\delta(f)\!: U^{[1]}\times V\to E$,
$\delta(f)(x,y,t,v):=f^{[1]}(x,v,y,0,t)$
is $C^{r+k-1}_\K$,
and $\delta\!: C^{r+k}(U\times V,E)\to C^{r+k-1}(U^{[1]}\times V,E)$
is a continuous linear map (see Remark~\ref{simplobs}\,(b),
Lemma~\ref{pullback}). By the induction hypothesis,
for any $g\in C^{r+k-1}(U^{[1]}\times V,E)$,
the map
\[
\Psi(g):=g^\vee \! : U^{[1]}\to C^r(V,E),\quad
\Psi(g)(x,y,t)(v):=g((x,y,t),\,v)
\quad\mbox{for $(x,y,t)\in U^{[1]}$, $v\in V$}
\]
is of class $C^{k-1}_\K$,
and the map
$\Psi\!: C^{r+k-1}(U^{[1]}\times V,E)\to C^{k-1}(U^{[1]},C^r(V,E))$
so obtained is continuous and linear.
By the preceding, given
$f\in C^{r+k}(U\times V,E)$,
we have $\Psi(\delta(f))\in C^{k-1}(U^{[1]},C^r(V,E))$.
In particular, $\Psi(\delta(f))=\delta(f)^\vee$
is continuous.
Note that (\ref{tnot0}) can be read as
\[
(f^\vee)^{]1[}(x,y,t)\; =\; (\delta(f))^\vee(x,y,t)\qquad
\mbox{for all $\,(x,y,t)\in U^{]1[}$.}
\]
Thus $f^\vee$ is of class
$C^1_\K$ with $(f^\vee)^{[1]}=\delta(f)^\vee=\Psi(\delta(f))$.
Now $f^\vee$ being of class $C^1_\K$ with $(f^\vee)^{[1]}$
of class~$C^{k-1}_\K$, the mapping $f^\vee$
is of class~$C^k_\K$.
Since $\Phi$ is continuous when
considered as a mapping
$C^{r+k}(U\times V,E)\to C(U,C^r(V,E))$
as a consequence of the
case $k=0$,
and the map $C^{r+k}(U\times V,E)\to C^{k-1}(U^{[1]},C^r(V,E))$,
$f\mto (f^\vee)^{[1]}=(\Psi\circ\delta)(f)$
is continuous by the preceding,
we deduce with Remark~\ref{simplobs}\,(b)
that $\Phi\!:
C^{r+k}(U\times V,E)\to C^k(U,C^r(V,E))$,
$\Phi(f)=f^\vee$ is continuous.
This completes the proof.
\end{proof} 
\begin{prop}\label{nownow}
Let $\K$ be a locally compact
topological field,
$E$ be a topological $\K$-vector space,
$M$ be a $C^\infty_\K$-manifold modeled
on a topological $\K$-vector space,
and
$N$ be a finite-dimensional
$C^\infty_\K$-manifold.
Then the following holds:
\begin{itemize}
\item[\n (a)]
A mapping $g\!: M\to C^\infty(N,E)$
is of class $C^\infty_\K$ if and only if
\[
g^\wedge\!: M\times N\to E,\quad g^\wedge(x,y):=g(x)(y)
\]
is of class $C^\infty_\K$.
\item[\n (b)]
The mapping $\quad \Phi\!:
C^\infty(M\times N, E)\to C^\infty(M,C^\infty(N,E)),\;\;\;
\Phi(f):=f^\vee$\\[2mm]
is an isomorphism of topological $\K$-vector spaces,
with inverse given by
\[
\Phi^{-1}\!: C^\infty(M,C^\infty(N,E))\to C^\infty(M\times N,E),\quad
\Phi^{-1}(g)=g^\wedge\,.
\]
\end{itemize}
\end{prop}
\begin{proof}
(a) By Proposition~\ref{evalCk},
the evaluation map $\ve\!: C^\infty(N,E)\times N\to E$
is smooth.
The formula $g^\wedge=\ve\circ (g\times \id_N)$
for $g\in C^\infty(M,C^\infty(N,E))$
shows that $g^\wedge$ is smooth whenever
so is~$g$. If, on the other hand,
$g\!: M\to C^\infty(N,E)$
is a mapping such that $f:=g^\wedge$ is smooth,
then $g=f^\vee$ is smooth, by Lemma~\ref{halfcartesian}\,(a).

(b) As a consequence of Lemma~\ref{halfcartesian}
and Part\,(a) of the present proposition,
the mapping $\Phi$ is an isomorphism of vector spaces
and continuous,
with inverse given by $\Phi^{-1}(g)=g^\wedge$
for $g\in C^\infty(M,C^\infty(N,E))$.
In order that $\Phi^{-1}$ be continuous, in view of
Lemma~\ref{amend} and Lemma~\ref{pb2},
we only need to show that
\[
C^\infty(\id_M\times \psi^{-1}, E)\circ \Phi^{-1}\!:
C^\infty(M,C^\infty(N,E))\to C^\infty(M\times V_\psi,E),\]
\[
f\mto \Phi^{-1}(f)\circ (\id_M\times \psi^{-1})=
f^\wedge \circ (\id_M\times \psi^{-1})
\]
is continuous, for each chart $\psi\!: U_\psi\to V_\psi\sub Y$
of~$N$, where $Y$ is the modeling space of~$N$.
Note that
$f^\wedge\circ(\id_M\times \psi^{-1})
=(C^\infty(\psi^{-1},E)\circ f)^\wedge
=\left(C^\infty(M,C^\infty(\psi^{-1},E))(f)\right)^\wedge$
for $f\in C^\infty(M,C^\infty(N,E))$, and thus
\[
C^\infty(\id_M\times \psi^{-1}, E)\circ \Phi^{-1}
=\Psi\circ C^\infty(M,C^\infty(\psi^{-1},E))\,,
\]
where $\Psi\!: C^\infty(M, C^\infty(V_\psi,E))
\to C^\infty(M\times V_\psi,E)$, $g\mto g^\wedge$.
The map $C^\infty(M,C^\infty(\psi^{-1},E))$
being continuous (Lemma~\ref{pb2}, Lemma~\ref{linearcase}),
it only remains to prove that $\Psi$ is continuous.
We fix $\psi$ now, and write $V:=V_\psi$
for brevity.\\[3mm]
Let $(W_i)_{i\in I}$ be an open cover of~$V$,
where $W_i\sub V$
is relatively compact for each $i\in I$,
with compact closure $K_i:=\wb{W_i}\sub V$.
As a consequence of
Lemma~\ref{amend},
the map $\Psi$ will be continuous
if we can show that $\rho_i\circ \Psi$
is continuous for each $i\in I$, where\linebreak
$\rho_i\!: C^\infty(M\times V,E)\to C^\infty(M\times W_i,E)$
is the restriction map. Hold $i\in I$ fixed.
We have
\begin{equation}\label{tiresome}
\rho_i(\Psi(g))=\rho_i(g^\wedge)=(\ve\circ
(g\times \id_V))|_{M\times W_i}\,,
\end{equation}
where $\ve\!: C^\infty(V,E)\times V\to E$ is evaluation
(which is $C^\infty_\K$ by
Proposition~\ref{evalCk}).
We want to re-write (\ref{tiresome})
further in order to be able to apply Proposition~\ref{crucial}\,(b).
To this end, we let $\sigma\!: W_i\to V$ be inclusion.
We define
\[
\tilde{h}\!:
V\times V \times C^\infty(V,E) \to E,\quad
\tilde{h}(v,y,\gamma ):=\ve(\gamma,y)=\gamma(y)
\]
and $h:=\tilde{h}\circ (\sigma\times \id_V\times
\id_{C^\infty(V,E)})\!: W_i\times V\times
C^\infty(V,E)\to E$,
$h(v,y,\gamma)=\gamma(y)$.
Then (\ref{tiresome})
can be re-written as
$\rho_i(\Psi(g))=C^\infty(\tau,E)\left(\phi(\id_V, g)\right)$,
where $\tau\!: M\times W_i\to W_i\times M$ is the coordinate
flip and where
\[
\phi\!: \lfloor K_i, V\rfloor_\infty\times C^\infty(M,C^\infty(V,E))
\to C^\infty(W_i\times
M, E),\quad
\phi(f,g):=h_*(f,g)
\]
is smooth by Proposition~\ref{crucial}\,(b);
here
$h_*(f,g)(y,x)=h(y,f(y),g(x))
=g(x)(f(y))$ for $y\in W_i$, $x\in M$,
and $\lfloor K_i,V\rfloor_\infty\sub C^\infty(V,Y)$.
Hence $\rho_i\circ \Psi$ is smooth and thus
continuous, which completes the proof of Part\,(b).
\end{proof}
The remainder of this section
is devoted to a variant of Proposition~\ref{nownow}
for manifolds modeled on metrizable topological
vector spaces.
In order to prove the result efficiently,
we introduce as an auxiliary concept
the notion of {\em conveniently $\K$-smooth\/} maps,
inspired
by the convenient differential calculus
of Fr\"{o}licher, Kriegl and Michor
(devoted
to the real or complex locally convex case).
\begin{defn}\label{bscsconv}
Given a topological field~$\K$
and topological $\K$-vector space~$E$,
a subset $U\sub E$ will be called {\em $c^\infty$-open\/}
if $U$ is open in the final topology on~$E$
with respect to the set of all $C^\infty_\K$-maps
$\gamma\!: I\to E$,
where $I$ is an open subset of $\K^n$ for some
$n\in \N$.\footnote{Thus, we require that $\gamma^{-1}(U)$ be open in
$I$ for any $\gamma$.}
A mapping $f\!: U\to F$
from a $c^\infty$-open subset $U\sub E$ to a topological
$\K$-vector space $F$ is called
{\em conveniently $\K$-smooth\/}
(or also a {\em $c^\infty_\K$-map\/})
if $f\circ \gamma\!: I\to F$ is $C^\infty_\K$,
for every $n\in \N$, open subset $I\sub \K^n$, and
$C^\infty_\K$-map $\gamma\!: I\to E$ such that $\gamma(I)\sub U$.
\end{defn}
Apparently every open subset $U\sub E$
is $c^\infty$-open,
and every $C^\infty_\K$-map is also $c^\infty_\K$.
Furthermore, it is obvious that
compositions of composable $c^\infty_\K$-maps
are $c^\infty_\K$-maps.
It does not pose any problems
to develop a theory of $c^\infty_\K$-manifolds,
along the lines of convenient differential calculus,
but we refrain from doing so here, as we wish
to focus on the $C^r_\K$-theory.
For the present purposes,
the following limited definition is sufficient:
We call a map $f\!: M\to E$
from a $C^\infty_\K$-manifold
to a topological $\K$-vector space
{\em conveniently $\K$-smooth\/}
(or a $c^\infty_\K$-map)
if $f\circ \gamma$ is $C^\infty_\K$
for $C^\infty_\K$-maps $\gamma\!: \K^n\supseteq I\to M$,
or, equivalently, if $f\circ \kappa^{-1}$
is conveniently $\K$-smooth
for every chart $\kappa$ in an $C^\infty_\K$-atlas
for~$M$.
\begin{rem}\label{specialconv}
Throughout this remark,
suppose that $\K=\R$, or that $\K$ is an ultrametric field.
Then, as a consequence of \cite[Thm.\,11.3\,(a)]{Ber}
(applied to subsets of $\K^n$),
a subset $U\sub E$ is $c^\infty$-open
if and only if $\gamma^{-1}(U)$ is open in $\K$
for every $C^\infty_\K$-curve $\gamma\!: \K\to E$.
Furthermore, a map $f\!: U\to F$ is
conveniently $\K$-smooth if and only
if $f\circ \gamma$ is $C^\infty_\K$
for every $n\in \N$ and every
$C^\infty_\K$-map $\gamma\!: \K^n\to E$
with image in~$U$ (defined on all of~$\K^n$).
If $E$ is a metrizable topological $\K$-vector space,
then a subset $U\sub E$ is open if and only if
it is $c^\infty$-open \cite[Thm.\,11.3\,(a)]{Ber};
in this case, a mapping $f\!: U\to F$
into a topological $\K$-vector space~$F$
is $C^\infty_\K$ if and only if it
is conveniently $\K$-smooth (cf.\
\cite[Thm.\,12.4]{Ber}).
\end{rem}
Before we can formulate the exponential
law, we need to have a second look at the evaluation map.
\begin{la}\label{seclook}
Let $\K$ be a topological field
which is metrizable $($or, more generally,
a topological field such that
$\K^n$ is a $k$-space for all $n\in \N)$.
Let $M$ be a $C^\infty_\K$-manifold modeled on a
topological $\K$-vector space,
and $E$ a topological $\K$-vector space.
Then the evaluation~map
\[
\ve\!: C^\infty(M,E)\times M\to E,\quad
\ve(\gamma,x):=\gamma(x)
\]
is conveniently $\K$-smooth.
\end{la}
\begin{proof}
Arguing similarly as in the proof of Proposition~\ref{evalCk},
we reduce to the case where $M$ is an open subset
of its modeling space~$Z$, which we assume now.
To establish the lemma,
we show by induction on $k\in \N_0$
that $\ve\circ c$ is of class $C^k_\K$,
for every $C^\infty_\K$-map $c=(c_1,c_2)\!: I\to C^\infty(M,E)\times M$
defined on an open subset $I\sub \K^n$
for some $n\in \N$.\\[3mm]
{\em The case $k=0$.}
Let $c$ be as before.
Since $\K^n$ is a $k$-space, so is its open
subset $I$. Hence $\ve\circ c$
will be continuous
if we can show that $\ve\circ c|_K\!: K\to E$
is continuous, for every compact subset $K\sub I$.
As $c_2\!: I\to M$ is continuous,
the set $L:=c_2(K)\sub M$ is compact.
Since
\[
\ve(c(x))=c_1(x)(c_2(x))=\tilde{\ve}(c_1(x)|_L,c_2(x))
\]
for all $x\in K$, where the restriction map
$C^\infty(M,E)\to C(L,E)_{c.o.}$, $\eta\mto \eta|_L$
is continuous (cf.\ Remark~\ref{compareco}
\& \cite[p.\,157, Eqn.\,(2)]{Eng})
and the
evaluation map
$\tilde{\ve}\!: C(L,E)_{c.o.}\times L\to E$
is continuous, we see that $\ve\circ c|_K$
is continuous, as desired.\\[3mm]
{\em Induction step.}
Suppose that $k\in \N_0$
and suppose that $\ve\circ c$ is of class $C^k_\K$,
for all $c$ as before.
For all $(x,y,t)\in I^{]1[}$,
we calculate
\begin{eqnarray}
(\ve\circ c)^{]1[}(x,y,t)\!&=& \!{\textstyle \frac{1}{t}}((\ve\circ c)(x+ty)
-(\ve\circ c)(x))\nonumber\\
&=& \!{\textstyle \frac{1}{t}(c_1(x+ty)-c_1(x))(c_2(x+ty))
+ \frac{1}{t} (c_1(x)(c_2(x+ty))-c_1(x)(c_2(x)))}\nonumber\\
&=& \!c_1^{[1]}(x,y,t)(c_2(x+ty))\,+\, c_1(x)^{[1]}(c_2(x),
c_2^{[1]}(x,y,t),t)\nonumber\\
&=& \!\ve(c_1^{[1]}(x,y,t),c_2(x+ty))\,+\,
(\tilde{\ve}\circ \tilde{c})(x,y,t)\,,\label{lastfm}
\end{eqnarray}
where $\tilde{c}\!: I^{[1]}\to C^\infty(M^{[1]},E)\times M^{[1]}$,
$\tilde{c}(x,y,t):=(c_1(x)^{[1]},\, (c_2(x), c_2^{[1]}(x,y,t),t))$
is smooth (cf.\ Remark~\ref{simplobs}),
and where evaluation
$\tilde{\ve}\!: C^\infty(M^{[1]},E)\times M^{[1]}\to E$
takes $C^\infty_\K$-maps on open subsets of $\K^m$
(for any $m\in \N$) to $C^k_\K$-maps, by induction.
Since, trivially, also
$I^{[1]}\to C^\infty(M,E)\times M$, $(x,y,t)\mto (c_1^{[1]}(x,y,t),
c_2(x+ty))$
is $C^\infty_\K$ (cf.\ {\bf \ref{usefulsimp}}),
we deduce from the induction hypothesis
that the map
\[
g\!: I^{[1]}\to E,\quad g(x,y,t):=
\ve(c_1^{[1]}(x,y,t),c_2(x+ty))\,+\,
(\tilde{\ve}\circ \tilde{c})(x,y,t)
\]
is of class $C^k_\K$ and thus continuous.
Since $g|_{I^{]1[}}=(\ve\circ c)^{]1[}$ by (\ref{lastfm}),
we deduce that $\ve\circ c$ is $C^1_\K$,
with $(\ve\circ c)^{[1]}=g$, and thus $\ve\circ c$
is of class $C^{k+1}_\K$, which completes the induction.
\end{proof}
\begin{prop}\label{expmetriz}
Let $\K$ be a metrizable topological field,
$E$ be a topological $\K$-vector space,
and $M$, $N$ be $C^\infty_\K$-manifolds modeled
on arbitrary topological $\K$-vector spaces.
Let $g\!: M\to C^\infty(N,E)$ be a map.
Then the following holds:
\begin{itemize}
\item[\n (a)]
If $g$ is
a $c^\infty_\K$-map, then also
$g^\wedge\!: M\times N\to E$, $g^\wedge(x,y):=g(x)(y)$
is a $c^\infty_\K$-map.
\item[\n (b)]
Assume that $\,\K=\R$ or $\K$ is an ultrametric field.
If both $M$ and $N$ are modeled on metrizable
topological $\K$-vector spaces,
then $g$ is $C^\infty_\K$
if and only if $g^\wedge$ is $C^\infty_\K$.
\item[\n (c)]
$\Phi\!: C^\infty(M\times N,E)\to C^\infty(M,C^\infty(N,E))$,
$f\mto f^\vee$
is a continuous isomorphism of vector spaces
in the situation of {\rm (b)},
whose inverse $g\mto g^\wedge$ is a $c^\infty_\K$-map.
\end{itemize}
\end{prop}
\begin{proof}
(a) Suppose that $g$ is a $c^\infty_\K$-map.
Then $g^\wedge=\ve\circ (g\times \id_N)$,
where the evaluation map $\ve\!: C^\infty(N,E)\times N\to E$
is $c^\infty_\K$ by Lemma~\ref{seclook}.
If $\gamma=(\gamma_1,\gamma_2)\!: \K^n\supseteq I\to M\times N$
is a $C^\infty_\K$-map, then
$g^\wedge\circ \gamma=\ve\circ (g\circ \gamma_1, \gamma_2)$
is of class $C^\infty_\K$
since $(g\circ \gamma_1, \gamma_2)\!: I\to C^\infty(N,E)\times N$
is a $C^\infty_\K$-map and $\ve$ is $c^\infty_\K$.
Thus $g^\wedge$ is $c^\infty_\K$.\\[3mm]
(b) Since $M$ and $M\times N$ are $C^\infty_\K$-manifolds
modeled on metrizable topological
vector spaces, where $\K$ is $\R$ or an ultrametric field,
mappings on these manifolds are $C^\infty_\K$
if and only if they are $c^\infty_\K$ (cf.\
Remark~\ref{specialconv}).
Hence (b) readily follows from (a) and
Proposition~\ref{halfcartesian}\,(a).\\[3mm]
(c) It is immediate from (b) and Proposition~\ref{halfcartesian}\,(b)
that $\Phi$ is a continuous linear bijection, with
inverse $\Psi\!: C^\infty(M,C^\infty(N,E))\to C^\infty(M\times N,E)$,
$\Psi(g)=g^\wedge$.
To see that $\Psi$ is a $c^\infty_\K$-map,
let $\gamma\!: \K^n\supseteq I\to C^\infty(M,C^\infty(N,E))$
be a $C^\infty_\K$-map. We have to show that
$\Psi\circ \gamma\!: I\to C^\infty(M\times N,E)$
is $C^\infty_\K$.
By Prop.\,\ref{halfcartesian}\,(a),
this will hold if we can show that
\[
f:=(\Psi\circ \gamma)^\wedge\!: I\times M\times N\to E
\]
is of class $C^\infty_\K$
(since then $f^\vee=\Psi\circ \gamma$ apparently).
This in turn holds if and only if $f$ is a $c^\infty_\K$-map,
the manifold $I\times M\times N$ being modeled
on a metrizable topological vector space.
However, given
$\eta=(\eta_1,\eta_2,\eta_3)\!: \K^m\supseteq J\to
I\times M\times N$ of class $C^\infty_\K$,
we have
\[
(f\circ \eta)(z)=\ve_2\Big(\ve_1\big(
\gamma(\eta_1(z)),\eta_2(z)\big),\, \eta_3(z)\Big)\,,
\]
where
$\ve_1\!: C^\infty(M,C^\infty(N,E))\times M\to C^\infty(N,E)$
and $\ve_2\!: C^\infty(N,E)\times N\to E$
are the respective evaluation maps,
which are $c^\infty_\K$-maps by Lemma~\ref{seclook}.
Consequently, $f\circ \eta$ is of class $C^\infty_\K$
and thus $f$ a $c^\infty_\K$-map, which completes
the proof.
\end{proof}
Note that, in the real or complex case,
none of the topological vector spaces involved
in Lemma~\ref{halfcartesian},
Proposition~\ref{nownow} and Proposition~\ref{expmetriz}
need to be locally convex.
Cf.\ \cite{Bil}
for a careful discussion
of the exponential law for maps $M\times N\to E$,
when $E$ is a real locally convex space
and both $M$ and $N$ are open subsets of
real locally convex spaces.
In the real and complex locally convex case,
the exponential law for conveniently smooth
maps plays a central role in convenient differential
calculus (see \cite{FaK}, \cite{KaM}).
The reader should be aware that
the locally convex topology on spaces
of conveniently smooth functions
primarily used in convenient differential calculus
(initial with respect to pullbacks
along smooth curves)
is in general properly coarser than
the topology we use here,
already for $C^\infty(\R^2,\R)$
(cf.\ \cite{Bil}).
%
%
%
%
%
%
%
%
%
%
%
%
%
%
%
\section{Diffeomorphism groups of finite-dimensional,\\
paracompact smooth manifolds over local fields}\label{secdiffeos}
Let $\K$ be a local field (of arbitrary characteristic),
and $M$ be a paracompact, finite-dimensional
smooth manifold over~$\K$.
In this section, we turn the group
$\Diff^\infty(M)$ of smooth diffeomorphisms
of~$M$ into a $\K$-Lie group,
modeled on the space $C^\infty_c(M,TM)$
of compactly supported smooth vector fields,
equipped with the box topology.\footnote{By
Prop.\,\ref{comparetop}\,(f),
this is the locally convex
direct limit topology on
$C^\infty_c(M,TM)=\dl_K \, C^\infty_K(M,TM)$.}
Since $M=\coprod_{i\in I} B_i$
is a disjoint union of balls,
we first turn the diffeomorphism group
$\Diff^\infty(B)$ of a ball
$B$ into a Lie group, which is quite easy.
Then also the weak direct product
$\prod_{i\in I}^*\Diff^\infty(B_i)$
is a Lie group by our general construction
from Section~\ref{secweakprod}.
As this weak direct product
can be identified with a subgroup
of $\Diff^\infty(M)$ in an obvious way,
it only remains to show in a final step
that $\Diff^\infty(M)$ can be given
a Lie group structure
making
$\prod_{i\in I}^*\Diff^\infty(B_i)$
an open subgroup; this will complete our construction.\\[3mm]
In the next section,
which can be read independently of the present one,
we describe an alternative (slightly more complicated)
construction, which is restricted to
$\sigma$-compact manifolds but provides
information also on the
groups $\Diff^r(M)$ of $C^r$-diffeomorphisms
for finite~$r$.
\begin{numba}\label{comebck}
Throughout this section and the next,
$\K$ denotes a local field.
We let $|.|$ be an ultrametric absolute value on~$\K$
defining its topology, and
$\bO$ be the maximal compact subring of~$\K$.
Given $d\in \N$, we
let $\|\sbull\|_\infty$ be the maximum norm on~$\K^d$.
Given $a\in \K^d$ and $\ve>0$,
$B_\ve(a):=\{y\in \K^d\!: \|y-a\|_\infty<\ve\}$
denotes the open ball with respect to the maximum norm.
\end{numba}
\begin{center}
{\bf The diffeomorphism group of a ball}
\end{center}
Given $d\in \N$, consider the
ball $B:=\bO^d\sub \K^d$.
We show:
\begin{prop}\label{diffball}
The set $\Diff^\infty(B)$ of all
$C^\infty_\K$-diffeomorphisms of $B$
is an open subset of $C^\infty(B,\K^d)$.
Consider $\Diff^\infty(B)$
as an open smooth submanifold
of $C^\infty(B,\K^d)$.
Then $\Diff^\infty(B)$, with composition
of mappings as the group operation, is a $\K$-Lie group.
\end{prop}
{\bf Proof.} We prove the proposition in various
steps.
\begin{numba}\label{gammult}
Define $\End^\infty(B)=C^\infty_B(B,B)=
\{\gamma\in C^\infty(B,\K^d)\!:
\gamma(B)\sub B\}$.
Since $B$ is both open and compact,
Proposition~\ref{pushforw2}
shows that $\End^\infty(B)$
is an open subset of $C^\infty(B,\K^d)=C^\infty_B(B,\K^d)$.
By Proposition~\ref{compcomp},
the composition map
\[
\Gamma\!: \End^\infty(B)\times \End^\infty(B)\to
\End^\infty(B)\, ,\quad \Gamma(\gamma,\eta):=\gamma\circ \eta
\]
is of class $C^\infty_\K$.
In particular, $\Gamma$ is continuous and thus
$\End^\infty(B)$ is a topological monoid,
with identity element~$\id_B$.
Hence, by standard arguments,
the unit group
$\Diff^\infty(B)=\End^\infty(B)^\times$
of the topological monoid $\End^\infty(B)$
will be open in
$\End^\infty(B)$ (and hence in
$C^\infty(B,\K^d)$) if we can show
that it contains an
identity neighbourhood.
\end{numba}
\begin{numba}
Given $\gamma\in C^\infty(B,\K^d)$ and $x\in B$,
we abbreviate $\gamma'(x):=d\gamma(x,\sbull)$.
We let
$\Omega\sub \End^\infty(B)$ be the
set of all
$\gamma\in \End^\infty(B)$ such that
$(\gamma-\id_B)^{[1]}(B\times B \times \bO)
\sub B_{\frac{1}{2}}(0)$.
Then $\Omega$
is an open identity neighbourhood
in $\End^\infty(B)$
and $\|\gamma'(x)-\id\|<\frac{1}{2}$
for all $\gamma\in \Omega$ and $x\in B$
(using the operator norm with
respect to the maximum norm $\|.\|_\infty$).
We now show
that $\Omega\sub \Diff^\infty(B)$.
To this end, let $\gamma\in \Omega$
and abbreviate $\sigma:=\gamma-\id_B$.
Then $\gamma'(x)\in \GL_d(\bO)=
\Iso(\K^d,\|\sbull\|_\infty)$ is a linear isometry
for all $x\in B$,
because $\|\gamma'(x)-\id\|<\frac{1}{2}$
(cf.\ \cite{Ser}, Chapter~IV, Appendix~1).
Furthermore, $\|\gamma'(x)\|=\|\gamma'(x)^{-1}\|=1$.
We conclude that
\begin{eqnarray}\label{prsntly0}
\|\gamma(z)-\gamma(y)-\gamma'(x).(z-y)\|_\infty
&=& \|\sigma(z)-\sigma(y)-\sigma'(x).(z-y)\|_\infty\nonumber\\
&\leq & \min \left\{\, \|\sigma(z)-\sigma(y)\|_\infty ,\,
\|\sigma'(x).(z-y)\|_\infty \right\} \nonumber \\
&<& \frac{1}{2}\, \|z-y\|_\infty \;=\;
\frac{1}{2\, \|\gamma'(x)^{-1}\|}\, \|z-y\|_\infty\label{horrib0}
\end{eqnarray}
for all $x,y,z\in B$ such that $y\not=z$.
Indeed, because we are using the maximum norm here,
given $x,y,z$ as before there exists $0\not= \xi\in \K$ such that
$|\xi|=\|z-y\|_\infty \leq 1$. Then
$\|\sigma'(x).(z-y)\|_\infty
\leq \|\sigma'(x)\|\, \|z-y\|_\infty
< \frac{1}{2} \|z-y\|_\infty$
and
\[
\sigma(z)-\sigma(y)
=\xi\, {\textstyle \frac{1}{\xi}}
\big(\sigma\big(y+\xi \textstyle{\frac{z-y}{\xi}}\big)-\sigma(y)\big)
=\xi\, \sigma^{[1]}\big(y,\textstyle{\frac{z-y}{\xi}},\xi\big)
\]
with $\frac{z-y}{\xi}\in \bO^d$ and $\xi\in \bO$,
entailing that
$\|\sigma(z)-\sigma(y)\|_\infty
\leq |\xi|\cdot \|\sigma^{[1]}(y,\frac{z-y}{\xi},\xi)\|_\infty
<\frac{1}{2}|\xi|=\frac{1}{2}\|z-y\|_\infty$.
Thus (\ref{horrib0}) holds.
Using (\ref{horrib0}) with $x=0$,
\cite[Lemma\,6.1\,(b)]{IMP} shows that
$\gamma$
is an isometry from $B$
onto $\gamma(0)+\gamma'(0).B
=\gamma(0)+B  = B$.
Since $\gamma'(x)\in \GL_d(\bO)$
is invertible for all~$x$, we deduce from the Inverse Function
Theorem~\cite[Thm.\,7.3]{IMP} that $\gamma$ is a $C^\infty_\K$-diffeomorphism
and thus $\gamma\in \Diff^\infty(B)$.
We have shown that $\Omega\sub \Diff^\infty(B)$.
Hence $\Diff^\infty(B)$ is open.
\end{numba}
\begin{numba}
The group multiplication on $\Diff^\infty(B)$ being smooth
by {\bf \ref{gammult}}, it only remains to show
smoothness of the inversion map $\iota\!: \Diff^\infty(B)\to
\Diff^\infty(B)$, $\iota(\gamma):=\gamma^{-1}$.
We only need to prove that
\[
\iota^\wedge\!: \Diff^\infty(B)\times B \to \K^d\,,\quad
\iota^\wedge(\gamma,x):=
\iota(\gamma)(x)=\gamma^{-1}(x)
\]
is smooth; then $\iota=(\iota^\wedge)^\vee\!: \Diff^\infty(B) \to
C^\infty(B,\K^d)$ will be
smooth,
by Lemma~\ref{halfcartesian}\,(a).
By Lemma~\ref{evalCk},
the evaluation map
\[
\ve\!: \Diff^\infty(B)\times B\to \K^d\,,\quad
\ve(\gamma,x):=\gamma(x)
\]
is smooth.
Note that $\ve(\gamma,\sbull)=\gamma$
is a diffeomorphism of $B$
for each $\gamma\in \Diff^\infty(B)$,
where $\Diff^\infty(B)$ is an open subset
of the metrizable topological vector space
$C^\infty(B, \K^d)$
(see Proposition~\ref{propprop}\,(c)).
Therefore
the Inverse Function
Theorem with Parameters \cite[Thm.\,8.1\,(c)$'$]{IMP}
can be applied to the map~$\ve$,
using the diffeomorphism $\gamma\in \Diff^\infty(B)$
itself as the parameter.
The theorem
shows that $\Diff^\infty(B)\times B\to \K^d$, $(\gamma,x)\mto
(\ve(\gamma,\sbull))^{-1}(x)=\gamma^{-1}(x)=\iota^\wedge(\gamma,x)$
is smooth. As just observed, this entails smoothness
of~$\iota$.\vspace{2mm}\Punkt
\end{numba}
\begin{numba}\label{sltly}
Slightly more generally,
let us consider
a $C^\infty_\K$-manifold~$M$ now
which is isomorphic to $B=\bO^d$
as a $C^\infty_\K$-manifold.
Given a $C^\infty_\K$-diffeomorphism
$\psi\!: M\to B$,
we simply give $\Diff^\infty(M)$
the uniquely determined $\K$-Lie group
structure modeled on $C^\infty(M,TM)\isom C^\infty(B,\K^d)$
which makes the isomorphism of groups
\[
\Theta_\psi\!: \Diff^\infty(M)\to\Diff^\infty(B)\,,
\quad   \Theta_\psi(\gamma):=\psi\circ \gamma\circ \psi^{-1}
\]
an isomorphism of Lie groups.
\end{numba}
\begin{la}\label{strindep}
For $M\isom B$ as before,
the Lie group structure on $\Diff^\infty(M)$
just defined is independent of the choice
of $C^\infty_\K$-diffeomorphism $\psi\!: M\to B$.
\end{la}
\begin{proof}
If both $\phi$ and $\psi$ are $C^\infty_\K$-diffeomorphism
$M\to B$, then the composition
$\Theta_\phi\circ (\Theta_\psi)^{-1}\!:
\Diff^\infty(B)\to \Diff^\infty(B)$,
$\gamma\mto (\phi\circ \psi^{-1}) \circ \gamma\circ (\phi\circ \psi^{-1})^{-1}$
is an inner automorphism of the Lie group $\Diff^\infty(B)$
and hence a $C^\infty_\K$-diffeomorphism.
The assertion follows.
\end{proof}
A very similar argument shows:
\begin{la}\label{conju}
Suppose that $M$ and $N$ are finite-dimensional
$C^\infty_\K$-manifolds such that $M\isom N\isom B$.
Let $\phi\!: M\to N$ be a $C^\infty_\K$-diffeomorphism.
Then
\[
\Lambda\!: \Diff^\infty(M)\to \Diff^\infty(N)\,,\qquad \gamma\mto
\phi\circ \gamma\circ \phi^{-1}
\]
is an isomorphism of Lie groups.
\end{la}
\begin{proof}
Let $\psi\!: N\to B$ be a $C^\infty_\K$-diffeomorphism.
Then also $\psi\circ \phi\!: M\to B$ is a $C^\infty_\K$-diffeomorphism
and hence
$\Theta_\psi\!: \Diff^\infty(N)\to\Diff^\infty(B)$,
$\gamma\mto \psi\circ \gamma\circ \psi^{-1}$
and
$\Theta_{\psi\circ\phi}\!: \Diff^\infty(M)\to\Diff^\infty(B)$,
$\gamma\mto (\psi\circ \phi) \circ \gamma\circ (\psi\circ \phi)^{-1}$
are isomorphisms of Lie groups.
Hence also
$\Lambda=\big(\Theta_\psi \big)^{-1}\circ
\Theta_{\psi\circ \phi}$ is an isomorphism of Lie groups.
\end{proof}
Another technical lemma is useful:
\begin{la}\label{nonmetric}
If $M$ is an open submanifold of $\K^d$
such that $M\isom B$, then $\Diff^\infty(M)$
is an open subset of $C^\infty(M,\K^d)$.
The manifold structure making
$\Diff^\infty(M)$ an open submanifold
of $C^\infty(M,\K^d)$ coincides with the
manifold structure underlying the Lie group
$\Diff^\infty(M)$, as defined in
{\bf \ref{sltly}}.
\end{la}
\begin{proof}
Let $\psi\!: M\to B$ be a $C^\infty_\K$-diffeomorphism.
As $M\sub \K^d$ is open and compact,
$C^\infty(M,M)$ is open in $C^\infty(M,\K^d)$.
The pullback $C^\infty(\psi^{-1},\K^d)\!:
C^\infty(M,\K^d)\to C^\infty(B,\K^d)$
is a linear isomorphism which takes
$C^\infty(M,M)$ $C^\infty_\K$-diffeomorphically
onto $C^\infty(B,M)$ (cf.\ Lemma~\ref{pb2}).
The map $C^\infty(B,\psi)\!:
C^\infty(B,M)\to C^\infty(B,B)$
is a $C^\infty_\K$-diffeomorphism,
since so is~$\psi$ (cf.\ Corollary~\ref{Cf2}).
Hence also
\[
\Phi\!: C^\infty(M,M)\to C^\infty(B,B),\quad \gamma\mto \psi
\circ \gamma\circ\psi^{-1}
\]
is a $C^\infty_\K$-diffeomorphism.
Since $\Diff^\infty(B)$ is an open submanifold
of $C^\infty(B,B)$, the set
$\Phi^{-1}(\Diff^\infty(B))=\Diff^\infty(M)$
is open in $C^\infty(M,M)$,
and $\Phi$ induces a $C^\infty_\K$-diffeomorphism
$\Theta_\psi$ from the open submanifold $\Diff^\infty(M)\sub C^\infty(M,M)$
onto $\Diff^\infty(B)$.
But the same map $\Theta_\psi$
also is a $C^\infty_\K$-diffeomorphism
from the Lie group $\Diff^\infty(M)$
onto $\Diff^\infty(B)$, by definition of the Lie group
structure in {\bf \ref{sltly}}.
\end{proof}
$\;$\pagebreak

\begin{center}
{\bf Passage to paracompact manifolds}
\end{center}
Let $M$ be a paracompact, finite-dimensional
smooth $\K$-manifold now,
of dimension~$d$, say.
Then $M$ is a disjoint union
$M=\coprod_{i\in I} B_i$ of a family
of open and compact balls $B_i\sub M$
(see Lemma~\ref{onlyopen}\,(b)).
For each $i\in I$, we equip $\Diff^\infty(B_i)$
with the $\K$-Lie group structure
modeled on $C^\infty(B_i,TB_i)$
described in {\bf \ref{sltly}} and Lemma~\ref{strindep}.
We then endow the weak direct product
\[
{\textstyle \prod^*_{i\in I}\Diff^\infty(B_i)}
\]
with a Lie group structure, as described in
Proposition~\ref{propweakdp}.
Consider the mapping
\begin{equation}\label{intrpsi}
\Psi : \, {\textstyle \prod^*_{i\in I}\Diff^\infty(B_i)}
\,\to \, \Diff^\infty(M)\,\quad
(\gamma_i)_{i\in I}\mto \coprod_{i\in I}\gamma_i
\end{equation}
taking
$(\gamma_i)_{i\in I}$ to the map
$\gamma\!: M\to M$ determined by $\gamma|_{B_i}=\gamma_i$
for each $i\in I$.
Then indeed $\Psi$ takes its
values in $\Diff^\infty(M)$, and apparently
$\Psi$ is injective and a homomorphism of groups.
Throughout the following, using $\Psi$
we shall
identify $\prod^*_{i\in I}\Diff^\infty(B_i)$
with the corresponding subgroup $\im\,\Psi\sub \Diff^\infty(M)$.
We shall also identify
the modeling space $\bigoplus_{i\in I}C^\infty(B_i,TB_i)$
with $C^\infty_c(M,TM)$ in the obvious way
(Proposition~\ref{comparetop2}\,(e));
cf.\ Remark~\ref{remnosecs}).
\begin{la}\label{prerefine}
Assume that $M\isom B$ and assume
that $M=\coprod_{j\in J}C_j$ for a finite family
$(C_j)_{j\in J}$ of balls.
Then $\prod_{j\in J}\Diff^\infty(C_j)$
is open in $\Diff^\infty(M)$ and $\Diff^\infty(M)$
induces the given manifold structure
on the product
$\prod_{j\in J}\Diff^\infty(C_j)$
of the Lie groups $\Diff^\infty(C_j)$.
\end{la}
\begin{proof}
{\em Reduction to the case where $M$ is a metric ball\/}:
Let $\psi\!: M\to B$ be a $C^\infty_\K$-diffeomorphism.
Then $\Theta_\psi\!: \Diff^\infty(M)\to \Diff^\infty(B)$,
$\gamma\mto \psi\circ \gamma\circ \psi^{-1}$ is an
isomorphism of Lie groups, by {\bf \ref{sltly}}.
Set $B_j:=\psi(C_j)\sub B$. Then $\Lambda_j\!:
\Diff^\infty(C_j)\to \Diff^\infty(B_j)$,
$\Lambda_j(\gamma):=\psi|_{C_j}^{B_j}\circ \gamma\circ
\big(\psi|_{C_j}^{B_j}\big)^{-1}$ is an isomorphism of
Lie groups for each $j\in J$, by Lemma~\ref{conju}.
Since the restriction of $\Theta_\psi$ to
$\prod_{j\in J}\Diff^\infty(C_j)$
coincides with the map $\prod_{j\in J}\Lambda_j$
onto $\prod_{j\in J}\Diff^\infty(B_j)$,
we clearly only need to prove the assertion for
$B=\coprod_{j\in J}B_j$ (then the assertion concerning
$M=\coprod_{j\in J }C_j$ follows).\\[3mm]
By the preceding, we may assume now that $M=B$.
By Lemma~\ref{amend},
the map $\rho\!: C^\infty(M,\K^d)\to\prod_{j\in J}C^\infty(C_j,\K^d)$,
$\rho(\gamma):=(\gamma|_{C_j})_{j\in J}$
is an isomorphism of topological vector spaces.
By Lemma~\ref{nonmetric},
the Lie group $\Diff^\infty(C_j)$ is an open submanifold
of $C^\infty(C_j,\K^d)$. Hence $\rho^{-1}$
induces an isomorphism $\Psi$ from the direct
product of Lie groups
$P:=\prod_{j\in J}\Diff^\infty(C_j)$
onto the open subset $\Psi(P)=\rho^{-1}(P)\sub \Diff^\infty(M)$
of $C^\infty(M,\K^d)$. Here $\Psi$ is the map
from (\ref{intrpsi}).
\end{proof}
\begin{la}\label{refineballs}
Assume that $M=\coprod_{i\in I}B_i$
is a disjoint union of a family
$(B_i)_{i\in I}$ of balls
$B_i\sub M$, and assume that, for each $i\in I$,
the ball $B_i=\coprod_{j\in J_i}C_{ij}$
is a disjoint union of balls
$C_{ij}\sub B_i$ for some finite set $J_i$.
Abbreviate $K:=\{(i,j)\!: i\in I, j\in J_i\}$.
Then
$\prod^*_{(i,j)\in K}\Diff^\infty(C_{ij})$
is an open subgroup of
$\prod^*_{i\in I}\Diff^\infty(B_i)$, and
$\prod^*_{i\in I}\Diff^\infty(B_i)$
induces the given manifold structure
on the weak direct product
$\prod^*_{(i,j)\in K}\Diff^\infty(C_{ij})$.
\end{la}
\begin{proof}
By Lemma~\ref{rearrange}\,(b),
there is a natural isomorphism
$\prod^*_{(i,j)\in K}\Diff^\infty(C_{ij})\isom
\prod_{i\in I}^* H_i$,
with $H_i:=\prod_{j\in J_i}\Diff^\infty(C_{ij})$.
Here $H_i$ is an open subgroup (and submanifold)
of $\Diff^\infty(B_i)$,
by Lemma~\ref{prerefine}.
Hence, by Lemma~\ref{rearrange}\,(a),
also the weak direct product
$\prod_{i\in I}^* H_i$
is an open subgroup and submanifold of
$\prod_{i\in I}\Diff^\infty(B_i)$,
as asserted.
\end{proof}
\begin{thm}[The Lie group structure on
{\boldmath $\mbox{\bf Diff}^\infty(M)$}]\label{thmdiffpara}
$\;$Let $M$ be a paracompact,\linebreak
finite-dimensional smooth
manifold over a local field~$\K$.
Then
there exists a uniquely\linebreak
determined $C^\infty_\K$-manifold
structure on $\Diff^\infty(M)$,
modeled on the space
$C^\infty_c(M,TM)$ of compactly supported smooth
vector fields, such that
$\prod^*_{i\in I}\Diff^\infty(B_i)$
is an open subgroup of $\Diff^\infty(M)$
and $\Diff^\infty(M)$ induces the given manifold
structure on
$\prod^*_{i\in I}\Diff^\infty(B_i)$,
for every cover $(B_i)_{i\in I}$
of~$M$ by mutually disjoint balls~$B_i$.
\end{thm}
\begin{proof}
For the moment, we fix a cover
$(B_i)_{i\in I}$ of $M$ by mutually disjoint
balls. Let $U:=\prod^*_{i\in I}\Diff^\infty(B_i)\sub
\Diff^\infty(M)$,
equipped with its natural Lie group structure
introduced above. Suppose we can show
the following claim:\\[3mm]
{\bf Claim.}
{\em For every $\gamma\in \Diff^\infty(M)$,
there exists an open identity neighbourhood
$V\sub U$ such that $I_\gamma(V)\sub U$,
and such that $I_\gamma|_V^U\!: V\to U$ is smooth,
where $I_\gamma\!: \Diff^\infty(M)\to \Diff^\infty(M)$,
$I_\gamma(\eta):=\gamma\circ \eta\circ \gamma^{-1}$.}\\[3mm]
Then there exists a uniquely determined
Lie group structure on $\Diff^\infty(M)$ with $U$
as an open submanifold, by Proposition~\ref{locchar}.
\begin{numba}
To prove the claim, let
$\gamma\in \Diff^\infty(M)$.
As a consequence
of Lemma~\ref{disjballs}, for each $i\in I$
there exists a finite cover
$(C_{ij})_{j\in J_i}$ of $B_i$
by mutually disjoint balls
that is subordinate
to the open cover
$\{B_i\cap \gamma^{-1}(B_k)\!: k\in I\}$ of~$B_i$.
Let $J:=\{(i,j)\!: i\in I, j\in J_i\}$.
Given $k\in I$, define $L_k:=\{(i,j)\in J\!: \gamma(C_{ij})\sub B_k\}$.
Let $L:=\{(k,\ell)\!: k\in I, \ell\in L_k\}$
and $D_k:=\gamma(C_\ell)$ for $(k,\ell)\in L$.
Then $(D_{k\ell})_{\ell\in L_k}$
is a finite
partition of~$B_k$ into balls,
for each $k\in I$.
The map $\pi\!: L\to J$, $\pi(k,\ell):=\ell$
is a bijection, and we have
\begin{equation}\label{infogamma}
\gamma(C_{\pi(k,\ell)})=\gamma(C_\ell)=
D_{k\ell}\quad\mbox{for all $(k,\ell)\in L$.}
\end{equation}
Define
$V:=\prod^*_{(i,j)\in J}\Diff^\infty(C_{ij})$
and $W:=\prod^*_{(k,\ell)\in L}\Diff^\infty(D_{k\ell})$;
then $V$ and $W$ are open subgroups (and submanifolds)
of $U$, by Lemma~\ref{refineballs}.
Let $\eta\in V$.
Given $x\in D_{k\ell}$, we have
$\gamma^{-1}(x)\in C_\ell$ by (\ref{infogamma})
and hence $\eta(\gamma^{-1}(x))\in C_\ell$.
Thus
$I_\gamma(\eta)(x)=
\gamma(\eta(\gamma^{-1}(x)))\in D_{k\ell}$,
using (\ref{infogamma}) again.
Therefore $I_\gamma(V)\sub W\sub U$,
as desired. Interpreting $\eta$ as
the corresponding element $(\eta_{ij})_{(i,j)\in J}
\in \prod^*_{(i,j)\in J}\Diff^\infty(C_{ij})$
here with $\eta_{ij}:=\eta|_{C_{ij}}$,
by the preceding we have
\[
I_\gamma(\eta)_{k\ell}:=I_\gamma(\eta)|_{D_{k\ell}}
\, =\, \Lambda_{k\ell}(\eta_{\pi(k,\ell)})
\]
for all $(k,\ell)\in L$,
where $\Lambda_{k\ell}\!: \Diff^\infty(C_{\pi(k,\ell)})\to
\Diff^\infty(D_{k\ell})$, $\Lambda_{k\ell}(\sigma):=
\gamma|_{C_\ell}\circ \sigma\circ \gamma^{-1}|_{D_{k\ell}}^{C_\ell}$
is an isomorphism of $\K$-Lie groups by Lemma~\ref{conju}.
Thus $I_\gamma|_V^W\!: V\to W$
is a mapping of the type discussed
in Lemma~\ref{rearrange}\,(c),
and thus $I_\gamma|_V^W$ is an isomorphism
of $\K$-Lie groups (and hence a $C^\infty_\K$-diffeomorphism).
Thus, the above claim is established. 
\end{numba}
By the preceding, $\Diff^\infty(M)$ admits a unique
$\K$-Lie group structure making
the Lie group $\prod^*_{i\in I}\Diff^\infty(B_i)$
a subgroup and open submanifold.
To complete the proof of Theorem~\ref{thmdiffpara},
it only remains to show that the Lie group structure
on $\Diff^\infty(M)$ so obtained is independent
of the choice of the partition $M=\coprod_{i\in I} B_i$
of $M$ into balls. To this end, suppose that
$M=\coprod_{j\in J} C_j$ is a second partition of~$M$
into balls.
For each $(i,j)\in I\times J$
such that $B_i\cap C_j\not=\emptyset$,
the open, compact submanifold $C_i\cap C_j$
can be partitioned into finitely many balls
(cf.\ Lemma~\ref{onlyopen}\,(b)).
As a consequence, we find a partition
$(D_k)_{k\in K}$ of~$M$ into balls
that is subordinate to the disjoint open cover
$\{B_i\cap C_j\!: i\in I, j\in J\}$ of~$M$.
Given $i\in I$, the set $K_i:=\{k\in K\!: D_k\sub B_i\}$
is finite, and $B_i=\coprod_{k\in K_i}D_k$.
Likewise, for any $j\in J$ the set $L_j:=\{k\in K\!: D_k\sub
C_j\}$ is finite, and $C_j=\coprod_{k\in L_j}D_k$.
Hence, by Lemma~\ref{refineballs},
the weak direct product $\prod_{k\in K}^*\Diff^\infty(D_k)$
is a subgroup and open submanifold of
both $U:=\prod_{i\in I}^*\Diff^\infty(B_i)$
and $V:=\prod_{j\in J}^*\Diff^\infty(C_j)$.
Therefore $\prod_{k\in K}^*\Diff^\infty(D_k)$
is a subgroup and open submanifold of both
$\Diff^\infty(M)$,
equipped with the unique Lie group structure
making $U$ an open submanifold,
and of
$\Diff^\infty(M)$,
equipped with the unique Lie group structure
making $V$ an open submanifold.
As a consequence, the two
Lie group structures on $\Diff^\infty(M)$ coincide.
\end{proof}
\begin{rem}
It is not hard to see that the natural action
$\Diff^\infty(M)\times M\to M$ (the evaluation map)
is smooth, entailing
that every smooth homomorphism $\pi\!: G\to \Diff^\infty(M)$
from a $\K$-Lie group~$G$ to $\Diff^\infty(M)$
gives rise to a smooth action $\pi^\wedge\!: G\times M\to M$.
If~$G$ is finite-dimensional here
or modeled on a metrizable topological vector space,
then $G$ has an open subgroup fixing any element
outside a compact subset of~$M$.
Conversely, given a smooth action $\sigma\!: G\times M\to M$,
the associated homomorphism
$\sigma^\vee\!: G\to \Diff^\infty(M)$ is smooth, provided
there exists an open subgroup $U\sub G$
and a compact subset $K\sub M$ such that $\sigma(x,y)=y$ for all
$x\in U$ and $y\in M\setminus K$.
This condition is, of course, automatically satisfied
if $M$ is compact; in this special case,
smooth actions of Lie groups on~$M$
are in one-to-one correspondence
with smooth homomorphisms into $\Diff^\infty(M)$.
Compare \cite{DIF}
for full proofs in the real case;
they are easily adapted to the present situation.
\end{rem}
In the real case,
it is well known that every continuous action
of a finite-dimensional Lie group on a manifold
by $C^\infty$-diffeomorphisms is automatically
smooth \cite[Thm., p.\,212]{MaZ}.
It is also known that every locally compact group
acting effectively on a connected
finite-dimensional smooth manifold by diffeomorphisms
is a Lie group (see \cite[Ch.\,I, Thm.\,4.6]{KaN}
and \cite[\S5.2]{MaZ}),
whence every locally compact subgroup of $\Diff^\infty(M)$
is a Lie group in particular.
It is natural to ask for analogues in the
$p$-adic case.
The following problems are open and deserve to be
investigated:
\begin{problem}
Is every compact subgroup~$G$
of the diffeomorphism group $\Diff^\infty(M)$
of a paracompact,
finite-dimensional smooth $p$-adic manifold
$M$ a $p$-adic Lie group\,?
Does this hold at least if~$G$ is topologically
finitely generated\,?
\end{problem}
One strategy might be to try to verify
the hypotheses of Lazard's characterization
of finite-dimensional $p$-adic Lie groups
(see
\cite[A1, Thm.\,(1.9)]{Laz}, \cite[Thm.\,8.32]{DSM}; cf.\
\cite[p.\,157]{Ser}).
However, the author suspects that
counterexamples can be found.
Compare also \cite{Lu10} for some related
studies of subgroups of diffeomorphism groups.
\begin{problem}
Are continuous actions of finite-dimensional
$p$-adic Lie groups by smooth
diffeomorphisms on paracompact,
finite-dimensional smooth $p$-adic
manifolds always\linebreak
smooth\,?
Arguing as in the real case (see \cite{MiP} or \cite{DIF}),
this would entail that every continuous homomorphism
from a finite-dimensional $p$-adic Lie
group to $\Diff^\infty(M)$ is smooth.
\end{problem}
\section{\!\!\!The diffeomorphism groups {\boldmath
$\text{Diff}^r\!(M)$\/} and {\boldmath $\text{Diff}^\infty\!(M)\wt{\,}$}}\label{secdiff2}
Let $\K$ be a totally disconnected local field
and $M$ be a $\sigma$-compact smooth
manifold over~$\K$,
of positive finite dimension~$d$.
Given $r\in \N \cup\{\infty\}$,
let $\Diff^r(M)$ be the group of all $C^r_\K$-diffeomorphisms
of~$M$.
In this section,
we equip $\Diff^r(M)$ with a $C^\infty_\K$-manifold structure
modeled on the space $C^r_c(M,TM)$
of compactly supported
vector fields of class $C^r_\K$ on~$M$
which makes $\Diff^r(M)$ a topological group,
with smooth right translation maps.
For $r=\infty$,
the preceding smooth manifold structure
makes $\Diff^\infty(M)$
a Lie group,
modeled on the LF-space
$C^\infty_c(M,TM)$;
it coincides with the Lie group
constructed in the preceding section.
However,
we shall also define
a {\em second\/}
smooth manifold structure
on $\Diff^\infty(M)$
making it a $\K$-Lie group, denoted
$\Diff^\infty(M)\wt{\,}$;
it is modeled
on the projective limit\label{defspatil}
\[
C^\infty_c(M,TM)\wt{\,}
\;:=\; \pl_{k\in \N_0}C^k_c(M,TM)\;=\;\bigcap_{k\in \N_0} C^k_c(M,TM)
\]
of LB-spaces.
Note that $C^\infty_c(M,TM)\wt{\,}$
coincides with $C^\infty_c(M,TM)$
as a vector space, but its topology is
coarser (and can be properly coarser).
Since $M$ is diffeomorphic
to an open subset $U$ of $\K^d$ (see Lemma~\ref{onlyopen}\,(a)),
we first discuss
$\Diff^r(U)$ and only pass to general~$M$
at the very end. Throughout this section,
we retain the conventions from {\bf \ref{comebck}}.
\begin{center}
{\bf The monoids {\boldmath ${\rm End}_c^r(U)$} and
{\boldmath ${\rm End}_c^\infty(U)\wt{\,}$}}
\end{center}
\begin{numba}\label{dfri}
Let $d\in \N$ and $U\sub \K^d$ be a non-empty, open
subset. By Lemma~\ref{disjballs},
there exist a countable set~$I$,
positive real numbers $r_i$ for $i\in I$
and elements $a_i\in U$ such that
$U=\bigcup_{i\in I}B_{r_i}(a_i)$
as a disjoint union.
Abbreviate $U_i:=B_{r_i}(a_i)$.
Then, for every $r\in \N_0\cup\{\infty\}$, by
Proposition~\ref{comparetop}\,(e) the map
\[
C^r_c(U,\K^d)\to \bigoplus_{i\in I}C^r(U_i,\K^d),\;\;\;
\gamma\mto (\gamma|_{U_i})_{i\in I}
\]
is an isomorphism of topological $\K$-vector
spaces (when the box topology is used
on the direct sum);
we use it to identify $C^r_c(U,\K^d)$ and $\bigoplus_{i\in I}
C^r(U_i,\K^d)$.
We define\label{defendom}
\[
\End_c^r(U):=\{\gamma\in C^r(U,\K^d)\!: \;
\mbox{$\gamma(U)\sub U$ and $\gamma-\id_U\in C^r_c(U,\K^d)$}\}
\;\;\;\mbox{and}
\]
\[
\cE_c^r(U):=\{\gamma\in C^r_c(U,\K^d)\!: \id_U+\gamma\in \End_c^r(U)\}\,.
\]
Then\label{defbetr}
\[
\beta_r\!: \cE_c^r(U)\to \End_c^r(U),\;\;\; \beta_r(\gamma):=\id_U+\gamma
\]
is a bijection.
\end{numba}
\begin{numba}
Every $\gamma\in \End_c^r(U)$ is a proper
map, since $\gamma^{-1}(K)\sub K\cup \Supp(\gamma-\id_U)$
for every compact subset~$K$ of~$U$.
Given $\gamma,\eta\in \cE_c^r(U)$,
we have
\[
(\id_U+\gamma)\circ (\id_U+\eta)\; =\; \id_U+\eta+\gamma\circ (\id_U+\eta)\,,
\]
where $\gamma\circ (\id_U+\eta)\in C^r_c(U,\K^d)$
since $\id_U+\eta$ is proper,
and thus $\eta+\gamma\circ (\id_U+\eta)\in C^r_c(U,\K^d)$.
Therefore $(\id_U+\gamma)\circ (\id_U+\eta)\in \End_c^r(U)$.
We have shown that $\End_c^r(U)$ is closed under composition of
maps, and thus $\End_c^r(U)$ is a monoid under composition,
with $\id_U$ as the identity element.
We give $\cE_c^r(U)$ the monoid structure which makes $\beta_r$
an isomorphism of monoids. Thus $0$ is the identity
in $\cE_c^r(U)$, and the monoid multiplication $\mu_r\!:
\cE_c^r(U)\times \cE_c^r(U)\to \cE_c^r(U)$ is given
by
\begin{equation}\label{defnmu}
\mu_r(\gamma,\eta)=\eta+\gamma\circ (\id_U+\eta)\,.
\end{equation}
\end{numba}
\begin{numba}\label{firsttime}
We claim that $\cE_c^r(U)$ is open
in~$C^r_c(U,\K^d)$.
In fact, let $\gamma\in \cE_c^r(U)$
be given. For every $x\in U$,
there exists
$i(x)\in I$ such that $x+\gamma(x)\in U_{i(x)}=B_{r_{i(x)}}(a_{i(x)})
=x+\gamma(x)+B_{r_{i(x)}}(0)\sub U$.
There exists $s(x)\in \,]0,r_{i(x)}]$ such that
$B_{s(x)}(x)\sub U$ and
$\gamma(B_{s(x)}(x))\sub
\gamma(x)+B_{r_{i(x)}}(0)$. By Lemma~\ref{disjballs},
there exists a countable cover
$(V_j)_{j\in J}$ of~$U$
by mutually disjoint compact open sets~$V_j$,
which is subordinate to $(B_{s(x)}(x))_{x\in U}$.
Given $j\in J$, choose $x_j\in U$ such that
$V_j\sub B_{s(x_j)}(x_j)$, and abbreviate $i(j):=i(x_j)$.
If $\eta\in C^r_c(U,\K^d)$ such that
$\eta(V_j) \sub B_{r_{i(j)}}(\gamma(x_j))$,
then
\begin{equation}\label{doublestar}
y+\eta(y)\in x_j+B_{s(x_j)}(0)
+ \gamma(x_j)+B_{r_{i(j)}}(0)=x_j+\gamma(x_j)+B_{r_{i(j)}}(0)=U_{i(j)}\sub U
\end{equation}
for all $y\in V_j$.
We have shown that the open neighbourhood
$\bigoplus_{j\in J} C^r(V_j,B_{r_{i(j)}}(\gamma(x_j)))$
of~$\gamma$ in $C^r_c(U,\K^d)$
is contained in $\cE_c^r(U)$. Thus $\cE_c^r(U)$ is indeed
open in $C^r_c(U,\K^d)$.
\end{numba}
\begin{numba}\label{transstruc}
We consider $\cE_c^r(U)$ as an open $C^\infty_\K$-submanifold
of $C^r_c(U,\K^d)$.
We equip $\End_c^r(U)$ with the smooth $\K$-manifold structure
making $\beta_r\!: \cE_c^r(U)\to \End_c^r(U)$
a $C^\infty_\K$-diffeomorphism.
\end{numba}
\begin{numba}\label{wildetilde}
Define $C_c^\infty(U,\K^d)\wt{\,}
:=\bigcap_{k\in \N_0} C^k_c(U,\K^d)=
\pl_{k\in \N_0}C_c^k(U,\K^d)$\vspace{-.8mm}.
Then the vector space underlying $C_c^\infty(U,\K^d)\wt{\,}$ is
$C^\infty_c(U,\K^d)$, but the projective limit topology
on $C_c^\infty(U,\K^d)\wt{\,}$
can be properly
coarser than the topology on $C^\infty_c(U,\K^d)$
if~$U$ is non-compact.
Since
\begin{equation}\label{usfl}
\cE_c^r(U)=\cE_c^0(U)\cap C^r(U,\K^d)
\end{equation}
for each $r\in \N_0\cup\{\infty\}$,
we have $\cE^\infty_c(U)=\cE^0_c(U)\cap C^\infty_c(U,\K^d)$
in particular. As a consequence,
$\cE^\infty_c(U)$ is an open subset of $C^\infty_c(U,\K^d)\wt{\,}$.
When equipped with the topology induced by $C^\infty_c(U,\K^d)\wt{\,}$,
we write $\cE_c^\infty(U)\wt{\,}$ for $\cE_c^\infty(U)$.
We write $\End_c^\infty(U)\wt{\,}$
for the monoid $\End_c^\infty(U)$, equipped with the $C^\infty_\K$-manifold
structure
making $\tilde{\beta}\!: \cE_c^\infty(U)\wt{\,}\to
\End_c^\infty(U)\wt{\,}$, $\gamma\mto \id_U+\gamma$
a $C^\infty_\K$-diffeomorphism and an isomorphism of monoids.
\end{numba}
In various steps, we now show:\pagebreak

\begin{prop}\label{cpendsmooth}
In the preceding situation, we have:
\begin{itemize}
\item[\n (a)]
For every $r,k\in \N_0\cup\{\infty\}$, the mapping
\[
m_{r,k}\!: \End_c^{r+k}(U)\times
\End_c^r(U)\to \End_c^r(U),\;\;\;
m_{r,k}(\gamma,\eta):=\gamma\circ \eta
\]
is of class $C^k_\K$.
In particular, for each $r\in \N_0$
the composition
map
$m_r := m_{r,0}\!: \End_c^r(U)\times \End_c^r(U)\to \End_c^r(U)$
is continuous,
and the composition maps
$m_\infty:=m_{\infty,\infty}\!: \End_c^\infty(U)\times
\End_c^\infty(U)\to \End_c^\infty(U)$
and
\[
\wt{m}\!: \End_c^\infty(U)\wt{\,}\times
\End_c^\infty(U)\wt{\,}\to \End_c^\infty(U)\wt{\,}\,,\quad
\wt{m}(\gamma,\eta):=\gamma\circ \eta
\]
are smooth.
\item[\n (b)]
For every $r\in \N_0\cup\{\infty\}$
and $\eta\in \cE^r_c(U)$, the right translation map
$\rho_{r,\eta}:=m_{r,0}(\sbull,\eta)\!: \End^r_c(U)\to \End^r_c(U)$
is of class~$C^\infty_\K$.
\item[\n (c)]
For every $r\in \N\cup \{\infty\}$,
the group of invertible elements
$\Diff_c^r(U):=\End_c^r(U)^\times$ is open in $\End_c^r(U)$,
and $\Diff_c^r(U)=\Diff_c^1(U)\cap \End_c^r(U)$.
Also $\Diff^\infty_c(U)\wt{\,}:=(\End_c^\infty(U)\wt{\,}\; )^\times$
is open in $\End_c^\infty(U)\wt{\,}$.
\item[\n (d)]
Given $r\in \N\cup\{\infty\}$, let
$\gamma$ be a $C^r_\K$-diffeomorphism of~$U$.
Then $\gamma\in \Diff^r_c(U)$ if and only if there
exists a compact subset $K\sub U$ such that
$\gamma(x)=x$ for all $x\in X\,\take\, K$.
\item[\n (e)]
For each $r\in \N\cup\{\infty\}$ and $k\in \N_0\cup\{\infty\}$,
the map
$\iota_{r,k}\!:
\Diff_c^{r+k}(U)\to \Diff_c^r(U)$, $\gamma\mto\gamma^{-1}$
is~$C^k_\K$. In particular,
inversion $\iota_r:=\iota_{r,0}\!: \Diff_c^r(U)\to \Diff_c^r(U)$ is
continuous for
each $r\in \N$,
and the inversion maps $\iota_\infty
:=\iota_{\infty,\infty}\!:\Diff^\infty_c(U)\to\Diff^\infty_c(U)$
and $\tilde{\iota}\!:
\Diff^\infty_c(U)\wt{\,}
\to \Diff^\infty_c(U)\wt{\,}$ are smooth.
\end{itemize}
Thus
$\Diff_c^\infty(U)$ and $\Diff^\infty_c(U)\wt{\,}$
are $\K$-Lie groups when considered as open
$C^\infty_\K$-
submanifolds of~$\End_c^\infty(U)$, resp., $\End_c^\infty(U)\wt{\,}$.
Furthermore,
$\Diff^r_c(U)$ is a topological group
with respect to the topology underlying
the smooth manifold $\Diff^r_c(U)$,
for each $r\in \N$.
\end{prop}
{\bf Proof.} We begin with the proof of (a) and (b).
\begin{numba}
The maps $\beta_r$ and $\beta_{r+k}$
being $C^\infty_\K$-diffeomorphisms and
isomorphism of monoids, in view of (\ref{defnmu})
apparently $m_{r,k}$ will be of class $C^k_\K$ 
if we can show
that the mapping $\mu_{r,k}\!:
\cE_c^{r+k}(U)\times \cE_c^r(U)\to C^r_c(U,\K^d)$,
$(\gamma,\eta)\mto \eta+\gamma\circ (\id_U+\eta)$
is of class~$C^k_\K$.
The map $(\gamma,\eta)\mto \eta$ involved
being continuous linear and thus smooth,
it suffices to show that
\[
f\!: C^{r+k}_c(U,\K^d)\times \cE_c^r(U)\to C^r_c(U,\K^d),
\;\;\; f(\gamma,\eta):=\gamma\circ (\id_U+\eta)
\]
is of class~$C^k_\K$.
To see this, let $\gamma\in C^{r+k}_c(U,\K^d)$,
$\eta\in \cE_c^r(U)$
be given. As in {\bf \ref{firsttime}},
we find a countable open cover
$(V_j)_{j\in J}$ of~$U$ by mutually disjoint
compact open sets~$V_j$, elements $x_j\in U$, and a mapping
$i\!: J\to I$ such that $\eta\in
\bigoplus_{j\in J} C^r(V_j,B_{r_{i(j)}}(\eta(x_j)))\sub
\cE_c^r(U)$, and such that
$(\id_U+\tau)(V_j)\sub U_{i(j)}$
for all $j\in J$ and $\tau\in
\bigoplus_{j\in J} C^r\big(V_j,B_{r_{i(j)}}(\eta(x_j))\big)\,$
(cf.\ (\ref{doublestar})).
Abbreviate $B_j:=B_{r_{i(j)}}(\eta(x_j))$ and
$Q:=\bigoplus_{j\in J}C^r(V_j,B_j)$.
Then $C^{r+k}_c(U,\K^d)\times Q$
is an open neighbourhood of $(\gamma,\eta)$.
By the preceding,
$f(\sigma,\tau)|_{V_j}=\sigma \circ (\id_U+\tau)|_{V_j}
=\sigma|_{U_{i(j)}}\circ
(\id_{V_j}+\tau|_{V_j})|^{U_{i(j)}}$
for all $\sigma\in C^{r+k}_c(U,\K^d)$
and $\tau\in Q$.
Thus $f|_{C^{r+k}_c(U,\K^d)\times Q}$
can be written as the composition
\begin{eqnarray*}
{\textstyle C^{r+k}_c(U,\K^d)\times Q} & \stackrel{\isom}{\longrightarrow}&
{\textstyle (\bigoplus_{i\in I} C^{r+k}(U_i,\K^d))\times
(\bigoplus_{j\in J}C^r(V_j,B_j))}\\
&\stackrel{\isom}{\longrightarrow} & {\textstyle \bigoplus_{i\in I, j\in J}
(C^{r+k}(U_i,\K^d)\times C^r(V_j,B_j))}\\
&\stackrel{p}{\longrightarrow}&
{\textstyle
\bigoplus_{j\in J}(C^{r+k}(U_{i(j)},\K^d)\times C^r(V_j,B_j))}\\
&\stackrel{\oplus_{j\in J}f_j}{\longrightarrow}&
{\textstyle \bigoplus_{j\in I}C^r(V_j,\K^d)
\stackrel{\isom}{\longrightarrow}C^r_c(U,\K^d)}
\end{eqnarray*}
where
``$\isom$'' are the obvious isomorphisms of topological
vector spaces (or their restrictions to
$C^\infty_\K$-diffeomorphisms of open sets),
$p$ is the map $(\sigma_i,\tau_j)_{i\in I, j\in J}\mto
(\sigma_{i(j)},\tau_j)_{j\in J}$ which~is~$C^\infty_\K$ as the
restriction of a continuous linear map,
and $f_j\!: C^{r+k}(U_{i(j)},\K^d)\!\times \!C^r(V_j,B_j)\to
C^r(V_j,\K^d)$, $f_j(\sigma,\tau):=
\sigma\circ(\id_{V_j}+\tau)|^{U_{i(j)}}$.
Then $f_j=\Gamma_j\circ (\id_{C^{r+k}(U_{i(j)},\K^d)}\times g_j)$, 
where the composition map
\[
\Gamma_j\!: C^{r+k}(U_{i(j)},\K^d)\times C^r(V_j,U_{i(j)})
\to C^r(V_j,\K^d)
\]
is~$C^k_\K$ by Proposition~\ref{compsmooth},\,\footnote{We apply
the proposition with $K:=Y:=V_j$; note that
$C^r(V_j,U_{i(j)})=\lfloor V_j,U_{i(j)}\rfloor_r$
here.}
and $g_j\!: C^r(V_j,B_j)\to C^r(V_j, U_{i(j)})$,
$g_j(\tau):=\id_{V_j}+\tau$ is smooth, being a restriction
of an affine-linear map. Thus each $f_j$ is of class~$C^k_\K$
and hence so is $\oplus_{j\in J}f_j$, by Proposition~\ref{mapsdirsums}.
Being a composition of $C^k_\K$-maps,
$f|_{C^{r+k}_c(U,\K^d)\times Q}$ is~$C^k_\K$.\\[3mm]
If $k=0$ here, then $f(\sbull,\eta)\!: C^r_c(U,\K^d)
\to C^r_c(U,\K^d)$ is a continuous map by the preceding.
Since it is also linear, we deduce that the map $f(\sbull,\eta)$ is
smooth. As a consequence,
$\rho_{r,\eta}=m_{r,0}(\sbull,\eta)$ is smooth,
establishing~(b).
\end{numba}
\begin{numba}\label{onceenough}
Clearly $\wt{m}$ will be smooth if we can show
that $\wt{\mu}\!: \cE^\infty_c(U)\wt{\,}\times \cE^\infty_c(U)\wt{\,}\to
C^\infty_c(U,\K^d)\wt{\,}$,
$\wt{\mu}(\gamma,\eta):=\eta+\gamma\circ (\id_U+\eta)$
is smooth.  
By Lemma~\ref{inpl}, $\wt{\mu}$ will be smooth if
$\lambda_r\circ \wt{\mu}$ is smooth
for every $r\in \N_0$,
where $\lambda_r\!: C^\infty_c(U,\K^d)\wt{\,}\to C^r_c(U,\K^d)$
is the inclusion map.
But this is the case.
In fact, given any $k\in \N_0$,
we have
$(\lambda_r\circ \wt{\mu})|^{\cE^r_c(U)}=
\mu_{r,k}\circ (\lambda_{r+k}\times
\lambda_r)|_{(\cE^\infty_c(U)\wt{\,}\,)^2}^{\cE^{r+k}_c(U)
\times \cE^r_c(U)}$,
where $\mu_{r,k}$ is of class~$C^k_\K$ and $\lambda_{r+k}$
and $\lambda_r$
are continuous linear and thus smooth.
Thus $\lambda_r\circ \wt{\mu}$ is of class~$C^k_\K$.
This completes the proof of Part\,(a) of Proposition~\ref{cpendsmooth}.
\end{numba}
\begin{numba} To prove (d), let $r\in \N
\cup\{\infty\}$. If $\gamma\in \Diff_c^r(U)$,
then $\gamma\in \End_c^r(U)$, whence there exists a compact subset
$K\sub U$ such that $(\gamma-\id_U)|_{U\take K}=0$
and thus $\gamma|_{U\take K}=\id_U|_{U\take K}$.
Conversely, if $\gamma$ is a $C^r_\K$-diffeomorphism
of~$X$ such that $\gamma$ coincides with $\id_U$
off some compact set~$K$,
then apparently also the inverse map
$\gamma^{-1}$ satisfies $\gamma^{-1}|_{X\take K}=\id_U|_{X\take K}$,
whence $\gamma^{-1}\in \End_c^r(U)$.
Thus $\gamma$ is invertible in the monoid
$\End_c^r(U)$; (c) is established.
\end{numba}
\begin{numba}
To prove (c), note first that
\begin{equation}\label{invertC1Cr}
\End^r_c(U)\cap \End_c^1(U)^\times =\End_c^r(U)^\times\;\;\;
\mbox{for all $r\in \N\cup\{\infty\}$.}
\end{equation}
In fact, clearly $\End_c^r(U)^\times\sub
\End^r_c(U)\cap \End_c^1(U)^\times$.
If, conversely, $\gamma\in \End^r_c(U)\cap \End_c^1(U)^\times$,
then $\gamma$ is a $C^1_\K$-diffeomorphism
and thus $d\gamma(x,\sbull)$ is invertible
for all $x\in U$. Since, furthermore, $\gamma$ is
of class~$C^r_\K$, the Ultrametric Inverse Function
Theorem~\cite[Thm.\,7.3]{IMP}
entails that $\gamma^{-1}$ is of class~$C^r_\K$
(cf.\ also \cite[Rem.\,5.4]{IMP}).
Thus $\gamma$ is a $C^r_\K$-diffeomorphism,
and, in view of (d), apparently $\gamma\in \Diff^r_c(U)=\End^r_c(U)^\times$.
\end{numba}
\begin{numba}\label{Omi}
Given $r\in \N\cup\{\infty\}$,
by the preceding
$\End_c^r(U)$ is a topological monoid.
Therefore its unit group
$\End_c^r(U)^\times$ will be open if we can show that it is
an identity neighbourhood.
The inclusion map $\End_c^r(U)\to \End_c^1(U)$ being continuous,
in view of (\ref{invertC1Cr}), we only need to show
that $\Diff_c^1(U)$ is open in $\End_c^1(U)$,
or equivalently, that
$\cE_c^1(U)^\times$ is a $0$-neighbourhood in $\cE_c^1(U)$.
Let $B_i:=B_{r_i}(0)\sub \K^d$
and $D_i:=\{t\in \K\!: |t|<r_i\}$
for $i\in I$,
with $r_i$ as in {\bf \ref{dfri}}.
Then $W:=\bigoplus_{i\in I}C^1(U_i,B_i)$ is an open zero-neighbourhood
in $C^1_c(U,\K^d)$, and we have $W\sub \cE_c(U)$
since
$x+\gamma(x)\in U_i+B_i=a_i+B_{r_i}(0)+B_{r_i}(0)=
a_i+B_{r_i}(0)=U_i\sub U$
for all $\gamma\in W$,
$i\in I$, and $x\in U_i$.
Define
\[
\Omega_i:=\{\sigma\in C^1(U_i,B_i)\!:
\mbox{$\sigma^{[1]}(U_i\times \bO^d \times D_i)
\sub B_{\frac{1}{2}}(0)$
and $d\sigma(U_i\times \bO^d)\sub B_{\frac{1}{2}}(0)$}\}\,.
\]
Then $\Omega:=\bigoplus_{i\in I}\Omega_i\sub W\sub \cE_c^1(U)$
is an open zero-neighbourhood.
We claim that $\Omega\sub \cE_c^1(U)^\times$,
or equivalently,
$\beta_1(\Omega)\sub \Diff_c^1(U)$.
To see this, let $\sigma=(\sigma_i)_{i\in I}\in \Omega$,
where $\sigma_i=\sigma|_{U_i}\in \Omega_i$ for $i\in I$.
Define $\gamma:=\beta_1(\sigma)=\id_U+\sigma$
and $\gamma_i=\id_{U_i}+\sigma_i=\gamma|_{U_i}$.
Then $\gamma_i'(x):=d\gamma_i(x,\sbull)=
\one+d\sigma_i(x,\sbull)\in \GL_d(\bO)=
\Iso(\K^d,\|\sbull\|_\infty)$ for all $x\in U_i$
(cf.\ \cite{Ser}, Chapter~IV, Appendix~1)
and $\|\gamma_i'(x)\|=\|\gamma_i'(x)^{-1}\|=1$,
because $\|\gamma_i'(x)-\one\|=\|\sigma_i'(x)\|<\frac{1}{2}$.
We conclude that
\begin{eqnarray}\label{prsntly}
\|\gamma_i(z)-\gamma_i(y)-\gamma_i'(x).(z-y)\|_\infty
&= &
\|\sigma_i(z)-\sigma_i(y)-\sigma_i'(x).(z-y)\|_\infty \nonumber \\
&\leq & \min \left\{\, \|\sigma_i(z)-\sigma_i(y)\|_\infty ,\,
\|\sigma_i'(x).(z-y)\|_\infty \right\} \nonumber \\
&<& \frac{1}{2}\, \|z-y\|_\infty \;=\;
\frac{1}{2\, \|\gamma'_i(x)^{-1}\|}\, \|z-y\|_\infty\label{horrib}
\end{eqnarray}
for all $x,y,z\in U_i$ such that $y\not=z$.
Indeed, because we are using the supremum norm here,
given $x,y,z$ as before there exists $0\not= \xi\in \K$ such that
$|\xi|=\|z-y\|_\infty <r_i$. Then
$\|\sigma_i'(x).(z-y)\|_\infty\leq \|\sigma_i'(x)\|\,
\|z-y\|_\infty< \frac{1}{2}\|z-y\|_\infty$
and
\[
\sigma_i(z)-\sigma_i(y)
=\xi\, {\textstyle \frac{1}{\xi}}
\big(\sigma_i(y+\xi {\textstyle \frac{z-y}{\xi}})-\sigma_i(y)\big)
=\xi\, \sigma_i^{[1]}(y,{\textstyle \frac{z-y}{\xi}},\xi)
\]
with $\frac{z-y}{\xi}\in \bO^d$ and $\xi\in D_i$,
entailing that
$\|\sigma_i(z)-\sigma_i(y)\|_\infty
\leq |\xi|\cdot \|\sigma_i^{[1]}(y,\frac{z-y}{\xi},\xi)\|_\infty
<\frac{1}{2}|\xi|=\frac{1}{2}\|z-y\|_\infty$.
Thus (\ref{horrib}) holds.
Using (\ref{horrib}) with $x=a_i$,
\cite[Lemma\,6.1\,(b)]{IMP} shows that
$\gamma_i$
is an isometry from $U_i=B_{r_i}(a_i)$
onto $\gamma_i(a_i)+\gamma_i'(a_i).B_{r_i}(0)
=\gamma_i(a_i)+B_{r_i}(0) =a_i+\sigma_i(a_i)+B_{r_i}(0)=
a_i+B_{r_i}(0)=B_{r_i}(a_i)=U_i$.
As a consequence, $\gamma$ is an
isometry from $U$ onto~$U$. Since $\gamma'(x)=\one+\sigma'(x)\in \GL_d(\bO)$
is invertible for all~$x$, we deduce from the Inverse Function
Theorem~\cite[Thm.\,7.3]{IMP} that $\gamma$ is a $C^1_\K$-diffeomorphism
and thus $\gamma\in \Diff_c^1(U)$, using~(d).
\end{numba}
\begin{numba}
In view of (\ref{invertC1Cr}),
the openness of $\End^1_c(U)^\times$ in $\End_c^1(U)$ just
established entails that $\Diff_c^\infty(U)\wt{\,}=
\End^\infty_c(U)\wt{\,}\cap \End^1_c(U)^\times$
is open in $\End^\infty_c(U)\wt{\,}$. This completes the proof
of Part~(c) of Proposition~\ref{cpendsmooth}.
\end{numba}
\begin{numba}
To prove (e), we first observe that, for given $r\in \N\cup\{\infty\}$
and $k\in \N_0\cup\{\infty\}$, the map $\iota_{r,k}$
will be of class~$C^k_\K$ if its restriction to some
open identity neighbourhood $Y\sub \Diff^{r+k}_c(U)$ is of class~$C^k_\K$.
In fact, if $\gamma\in \Diff_c^{r+k}(U)$ is given,
then $Y\circ \gamma$ is an open neighbourhood
of $\gamma$ in $\Diff_c^{r+k}(U)$ since right translation by
$\gamma$ is a $C^\infty_\K$-diffeomorphism of
$\End^{r+k}_c(U)$,
as a consequence of~(b).
In view of (a) and (b), the formula
$\iota_{r,k}|_{Y\circ \gamma}(\eta)=
\eta^{-1}=\gamma^{-1}\circ (\eta\circ \gamma^{-1})^{-1}
=m_{r,k}\Big(\gamma^{-1},\iota_{r,k}|_Y(\rho_{r+k,\gamma^{-1}}(\eta))
\Big)$
for $\eta\in Y\circ \gamma$
shows that $\iota_{r,k}$ will be $C^k_\K$ on $Y\circ \gamma$
if it is $C^k_\K$ on~$Y$.
\end{numba}
\begin{numba}
By {\bf \ref{transstruc}} and the preceding, $\iota_{r,k}$ will
be of class $C^k_\K$ if we can show that
the map
\[
j_{r,k}\!: \cE^{r+k}_c(U)^\times\to C^r_c(U,\K^d),\;\;\;
j_{r,k}(\gamma):=\gamma^*:= (\id_U+\gamma)^{-1}-\id_U
\]
is $C^k_\K$ on some open $0$-neighbourhood
in $\cE^{r+k}_c(U)$.
Note that $S:=\bigoplus_{i\in I}
\cE^{r+k}_c(U_i)^\times$\linebreak
is an open subset of $\cE^{r+k}_c(U)^\times \sub C^{r+k}_c(U,\K^d)
=\bigoplus_{i\in I} C^{r+k}(U_i,\K^d)$,
and $j_{r,k}|_S=\oplus_{i\in I} \,j_{i,r,k}\!:
S\to \bigoplus_{i\in I} C^r(U_i,\K^d)=C^r_c(U,\K^d)$,
where $j_{i,r,k}\!:
\cE^{r+k}_c(U_i)^\times\to C^r_c(U_i,\K^d)$,
$j_{i,r,k}(\gamma):=(\id_{U_i}+\gamma)^{-1}-\id_{U_i}$.
By Proposition~\ref{mapsdirsums},
$\oplus_{i\in I}\, j_{i,r,k}$ will be of class~$C^k_\K$
if each $j_{i,r,k}$ is of class~$C^k_\K$.
Summing up, in order that $\iota_{r,k}$ be
of class~$C^k_\K$, for all $r$ and~$k$,
we only need to establish the following claim:\\[3mm]
{\em Claim.}
{\em For each $i\in I$, $r\in \N\cup\{\infty\}$
and $k\in \N_0\cup\{\infty\}$,
the map 
$j_{i,r,k}$ is of class~$C^k_\K$
on some open zero-neighbourhood in $\cE^{r+k}(U_i)^\times$,
where $\cE^{r+k}(U_i):=\cE^{r+k}_c(U_i)$.}
\end{numba}
\begin{numba}\label{proofincts}
Fix $i$, $r$ and $k$.
Since $U_i$ is compact, $\Diff^r_c(U_i)=\Diff^r(U_i)$
is the set of all $C^r_\K$-diffeomorphisms
of~$U_i$, and the map $C^r(U_i,\K) \to C^r(U_i,\K^d)$,
$\gamma\mto \gamma+\id_{U_i}$ is an affine-linear
homeomorphism and hence a $C^\infty_\K$-diffeomorphism,
which takes $\cE^r(U_i)^\times$
diffeomorphically onto the open subset $\Diff^r(U_i)
\sub C^r(U_i,\K^d)$.
Thus $\Diff^r(U_i)$ simply is an open
$C^\infty_\K$-submanifold of $C^r(U_i,\K^d)$.
Likewise for $\Diff^{r+k}(U_i)$.
In order that $j_{i,r,k}$ be $C^k_\K$
on some open zero-neighbourhood,
it therefore suffices to show that
inversion $h \!: P \to C^r(U_i,\K^d)$,
$h (\gamma):=\gamma^{-1}$
is $C^k_\K$ on the open identity neighbourhood
\[
P :=\{\id_{U_i}+ \sigma\!:\sigma \in \Omega_i\cap C^{r+k}(U_i,\K^d)\}
\]
of $\Diff^{r+k}(U_i)$,
where $\Omega_i\sub C^1(U_i,\K^d)$ is as in
{\bf \ref{Omi}}.\footnote{The discussion in {\bf \ref{Omi}}
shows that indeed $P \sub \Diff^{r+k}(U_i)$.}
We only need to prove that\linebreak
$h^\wedge\!: P \times U_i\to \K^d$, $h^\wedge(\gamma,x):=
h(\gamma)(x)=\gamma^{-1}(x)$
is $C^{r+k}_\K$; then $h=(h^\wedge)^\vee\!: P \to
C^r(U_i,\K^d)$ will be $C^k_\K$,
by Lemma~\ref{halfcartesian}\,(a).
By Lemma~\ref{evalCk},
the evaluation map
\[
\ve\!: P \times U_i\to \K^d\,,\quad
\ve(\gamma,x):=\gamma(x)
\]
is $C^{r+k}_\K$.
Since $\ve(\gamma,\sbull)=\gamma$
is a diffeomorphism of $U_i$
for each $\gamma\in P$,
and $C^{r+k}(U_i,\K^d)$
is metrizable (Proposition~\ref{propprop}\,(c)),
the Inverse Function
Theorem with Parameters \cite[Thm.\,8.1\,(c)$'$]{IMP}
can be applied to $\ve$,
with the diffeomorphism $\gamma\in P$ as the parameter.
The theorem
shows that $P \times U_i\to \K^d$, $(\gamma,x)\mto
(\ve(\gamma,\sbull))^{-1}(x)=\gamma^{-1}(x)=h^\wedge(\gamma,x)$
is of class $C^{r+k}_\K$. The claim is established.
\end{numba}
\begin{numba}
Arguing as in {\bf \ref{onceenough}},
we deduce from the fact that $\iota_{r,k}$ is of class~$C^k_\K$
for all $r,k\in \N$ that
$\wt{\iota}$ is smooth.
This completes the proof of Proposition~\ref{cpendsmooth}.\Punkt
\end{numba}
\begin{la}\label{pbproper}
Let $n,m\in \N$, $r\in \N_0\cup\{\infty\}$,
$U\sub \K^n$, $V\sub\K^m$ be open subsets,
$\phi\!: U\to V$ be a proper mapping of class $C^r_\K$,
and $E$ be a topological $\K$-vector space.
Then
\[
C^r_c(\phi,E)\!: C^r_c(V,E)\to C^r_c(U,E),\;\;\;
\gamma \mto \gamma\circ \phi
\]
is a continuous $\K$-linear map.
\end{la}
\begin{proof}
Indeed, the mapping $C^r_c(\phi,E)$ being linear,
by Proposition~\ref{comparetop}\,(c) we only
need to show that its restriction
to $C^r_K(V,E)$ is continuous, for
every compact open subset $K$ of~$V$.
Since~$\phi$ is assumed to be proper,
$L:=\phi^{-1}(K)\sub U$ is a compact, open subset.
It is clear that $C^r_c(\phi,E)$ takes
$C^r_K(V,E)$ into $C^r_L(U,E)$.
Using the obvious identifications
$C^r_K(V,E)\isom C^r(K,E)$ and $C^r_L(U,E)\isom
C^r(L,E)$, the map
$C^r_c(\phi,E)|_{C^r_K(V,E)}^{C^r_L(U,E)}$
corresponds to the pullback
$C^r(\phi|_L^K,E)\!: C^r(K,E)\to C^r(L,E)$,
which is a continuous linear map by Lemma~\ref{pullback}.
\end{proof}
\begin{la}\label{whynomatt}
Let $U\sub \K^d$, $V\sub \K^d$ be open subsets
and $\phi\!:
U\to V$ be a bijection. Then the following holds:
\begin{itemize}
\item[\n (a)]
Let $r,k\in
\N_0\cup\{\infty\}$,
and suppose that $\phi$ is a $C^{r+k}_\K$-diffeomorphism.
Then
\[
\Phi\!: \End_c^r(U)\to \End_c^r(V),\;\;\;
\Phi(\gamma):=\phi\circ \gamma\circ \phi^{-1}
\]
is a $C^k_\K$-diffeomorphism and an isomorphism of monoids.
\item[\n (b)]
If $\phi$ is a $C^\infty_\K$-diffeomorphism,
then $\End_c^\infty(U)\wt{\,}\to \End_c^\infty(V)\wt{\,}$,
$\gamma\mto \phi\circ \gamma\circ \phi^{-1}$
is a $C^\infty_\K$-diffeomorphism and an isomorphism
of monoids.
\end{itemize}
\end{la}
\begin{proof}
(a) It is obvious that~$\Phi$ is a bijection,
whose inverse
$\gamma\mto \phi^{-1}\circ \gamma\circ\phi$
also is a map as described in the lemma.
Furthermore, clearly $\Phi$ is a homomorphism
of monoids.
In view of the preceding it only remains to
show that $\Phi$ is a $C^k_\K$-map,
or equivalently, that
\[
\Psi_{r,k}\!: \cE_c^r(U)\to \cE_c^r(V),\;\;\;
\Psi_{r,k}(\gamma):=\phi\circ (\id_U+\gamma)\circ \phi^{-1}-\id_V
=\phi\circ (\phi^{-1}+\gamma\circ \phi^{-1})-\id_V
\]
is of class~$C^k_\K$.
To see this, we note first that $\Psi_{r,k}$
is almost local.
Indeed: Choose any locally finite
cover $(U_\ell)_{\ell\in L}$ of $U$ by relatively compact,
open sets $U_\ell$.
Then $V_\ell:=\phi(U_\ell)$
defines a locally finite cover
$(V_\ell)_{\ell\in L}$
of~$V$ by relatively compact, open subsets
$V_\ell\sub V$.
Given any $\ell\in L$, for every
$x\in V_\ell$ and $\gamma\in \cE_c^r(U)$
we have
\[
\Psi_{r,k}(\gamma)(x)
=\phi\big(\phi^{-1}(x)+\gamma(\phi^{-1}(x))\big)
\;-\, x\,,
\]
which only depends on the value of $\gamma$ at
$\phi^{-1}(x)\in U_\ell$.
Thus $\Psi_{r,k}$ is almost local.\\[3mm]
By the Smoothness Theorem (Theorem~\ref{smoothy}),
$\Psi_{r,k}$ will be $C^k_\K$ if we can show that
the restriction $f$ of $\Psi_{r,k}$
to the open subset $\cE^r_c(U)\cap C^r_K(U,\K^d)$
of $C^r_K(U,\K^d)$
is $C^k_\K$,
for every compact subset $K\sub U$.
It suffices to show this for $K$ open and compact,
which we assume now.
The image of $f$ is contained in $C^r_{\phi(K)}(V,\K^d)$.
The inclusion mapping\linebreak
$C^r_{\phi(K)}(V,\K^d)\to C^r_c(V,\K^d)$
being continuous linear and hence smooth,
it therefore suffices to prove that
$f$ is $C^k_\K$ as a map into $C^r_{\phi(K)}(V,\K^d)$.
Since $\phi(K)$ is compact and open, the restriction map
\[
C^r_{\phi(K)}(V,\K^d)\to C^r_{\phi(K)}(\phi(K),\K^d)\;=\;
C^r(\phi(K),\K^d)
\]
is an isomorphism of topological vector spaces
(Lemma~\ref{restrK}).
In order that $f$ be $C^k_\K$, we therefore
only need to show that
\[
g\!: \cE^r_c(U)\cap C^r_K(U,\K^d)\to C^r(\phi(K),\K^d),\quad
g(\gamma):=f(\gamma)|_{\phi(K)}
\]
is $C^k_\K$. Note that
\begin{eqnarray}
g(\gamma) &= &
\phi\circ (\id_K+\gamma|_K)\circ \phi^{-1}|_{\phi(K)}
\;-\,\id_{\phi(K)}\nonumber\\
&=& \big(C^r(\phi^{-1}|_{\phi(K)},\K^d)\circ
C^r(K,\phi)\big)\,(\,\underbrace{\id_K+\gamma|_K}_{\in C^r(K,U)}\,)
\;-\,\id_{\phi(K)}\,. \label{seesm}
\end{eqnarray}
The pullback $C^r(\phi^{-1}|_{\phi(K)},\K^d)\!:
C^r(K,\K^d)\to C^r(\phi(K),\K^d)$
is continuous linear and hence smooth,
by Lemma~\ref{pb2};
the map $C^r_c(U,\K^d)\to C^r(U,\K^d)\to C^r(K,\K^d)$,
$\gamma\mto \gamma|_K$
composed of inclusion and restriction is
continuous linear and hence smooth (cf.\ Proposition~\ref{comparetop}\,(a)
and Lemma~\ref{pb2});
and $C^r(K,\phi)\!: C^r(K,U)\to C^r(K,V)$
is $C^k_\K$
as $\phi$ is $C^{r+k}_\K$, by Corollary~\ref{Cf2}.
Hence (\ref{seesm}) shows that $g$ is $C^k_\K$,
as required.\vspace{1mm}

(b) Apparently (b) will hold if we can show that
the mapping $\cE_c^\infty(U)\wt{\,}\to \cE_c^\infty(V)\wt{\,}$,
$\gamma\mto \phi\circ (\id_U+\gamma)\circ \phi^{-1}-\id_V$
is~$C^k_\K$ for all $k\in \N$.
But this readily follows from the fact that $\Psi_{k,k}$
is of class~$C^k_\K$ for all $k\in \N$
(cf.\ {\bf\ref{onceenough}}).
\end{proof}
\begin{prop}\label{DiffU}
Let $d\in \N$
and $U\sub \K^d$
be a non-empty open subset. Then the following holds:
\begin{itemize}
\item[\n (a)]
For each $r\in \N\cup\{\infty\}$,
there is a uniquely determined
$C^\infty_\K$-manifold structure on the group
$\Diff^r(U)$ of all $C^r_\K$-diffeomorphisms of~$U$
such that $\Diff^r(U)$ becomes a topological group,
the right translation maps
$R_\gamma\!: \Diff^r(U)\to \Diff^r(U)$,
$R_\gamma(\eta):=\eta\circ \gamma$ are $C^\infty_\K$
for each $\gamma\in \Diff^r(U)$,
and such that $\Diff^r_c(U)$ is an open $C^\infty_\K$-submanifold
of $\Diff^r(U)$.
\item[\rm (b)]
For any $r\in \N\cup\{\infty\}$
and $k\in \N_0\cup\{\infty\}$,
the composition map
\begin{equation}\label{anothcomp}
\Diff^{r+k}(U)\times \Diff^r(U)\to \Diff^r(U),\quad
(\gamma,\eta)\mto \gamma\circ \eta
\end{equation}
and the inversion map
\begin{equation}\label{anothinv}
\Diff^{r+k}(U)\to \Diff^r(U),\quad \gamma\mto \gamma^{-1}
\end{equation}
are of class $C^k_\K$,
with respect to the smooth manifold
structures from {\rm (a)}.
\item[\n (c)]
There is a uniquely determined smooth
manifold structure on $\Diff^\infty(U)$
turning it into a $\K$-Lie group
and such that $\Diff_c^\infty(U)$ is an open
smooth submanifold of $\Diff^\infty(U)$.
\item[\n (d)]
There is a uniquely determined smooth
manifold structure on $\Diff^\infty(U)$
turning it into a $\K$-Lie group
$($which we denote by $\Diff^\infty(U)\wt{\,}\,)$,
such that $\Diff_c^\infty(U)\wt{\,}$ is an open
smooth submanifold of $\Diff^\infty(U)\wt{\,}$.
\end{itemize}
\end{prop}
\begin{proof}
It readily follows from
Proposition~\ref{cpendsmooth}\,(d)
that
$\Diff_c^r(U)$ is
a normal subgroup of $\Diff^r(U)$,
for each $r\in \N\cup\{\infty\}$.
Given $\gamma\in \Diff^r(U)$, let $I_\gamma\!:
\Diff_c^r(U)\to \Diff_c^r(U)$ denote the automorphism
of groups
$\eta\mto \gamma\circ \eta\circ \gamma^{-1}$.\vspace{1mm}

(a) Given $\gamma\in \Diff^r(U)$,
consider the map
\[
\kappa_\gamma\!: \Diff^r_c(U)\circ \gamma\to \cE^r_c(U)^\times\,,
\quad \kappa_\gamma(\eta):= \beta_r^{-1}(\eta\circ \gamma^{-1})\,.
\]
We claim that $\cA:=\{\kappa_\gamma\!: \gamma\in \Diff^r(U)\}$
is an atlas defining a $C^\infty_\K$-manifold
structure on $\Diff^r(U)$ (equipped with the final
topology with respect to the family of the maps $\kappa_\gamma^{-1}\!:
\cE^r_c(U)^\times \to \Diff^r(U)$).
The domains of the maps $\kappa_\gamma$
cover $\Diff^r(U)$. Let us prove compatibility
of the charts.
If $\gamma,\wb{\gamma}\in \Diff^r(M)$
such that $\Diff^r_c(U)\circ \gamma$
and $\Diff^r_c(U)\circ \wb{\gamma}$ have non-empty
intersection, then $\wb{\gamma}\circ \gamma^{-1}\in \Diff^r_c(U)$
and the two cosets coincide.
For $\eta\in \cE^r_c(U)^\times$, we have
\[
(\kappa_\gamma \circ \kappa_{\wb{\gamma}}^{-1})(\eta)
=\beta_r^{-1}(\beta_r(\eta)\circ \wb{\gamma}\circ \gamma^{-1})
=\beta_r^{-1}\big(
\rho_{r,\wb{\gamma}\circ \gamma^{-1}}
(\beta_r(\eta))\big)
\]
using right translation
$\rho_{r,\wb{\gamma}\circ \gamma^{-1}}$ on $\Diff^r_c(U)$,
which is smooth.
Hence $\kappa_\gamma \circ \kappa_{\wb{\gamma}}^{-1}$ is smooth,
as required for compatibility. Now standard
arguments provide a smooth manifold
structure on $\Diff^r(U)$ with atlas $\cA$.
Since $\kappa_{\text{id}}=\beta_r^{-1} \!:
\Diff^r_c(U)\to \cE^r_c(U)^\times$,
we see that $\Diff^r_c(U)$ is an open submanifold
of $\Diff^r(U)$. Given $\gamma\in \Diff^r(U)$,
for each $\eta\in \Diff^r(U)$ we have\linebreak
$\big(\kappa_{\eta\circ \gamma}\big)^{-1}\circ R_\gamma
\circ \kappa_\eta^{-1}=\id$ on $\cE^r_c(U)^\times$,
entailing that $R_\gamma$ is smooth.
The topology underlying
$\Diff^r(U)$ makes it a topological group,
because it has the following properties
(cf.\ \cite[Thm.\,4.5]{HaR}):
the topological group
$\Diff^r_c(U)$ is an open subgroup
of $\Diff^r(U)$; all right translations are homeomorphisms
of $\Diff^r(U)$;
and $I_\gamma$ is continuous for each $\gamma\in \Diff^r(U)$,
by Lemma~\ref{whynomatt}.\vspace{1mm}

(b) Let $\gamma\in \Diff^{r+k}(U)$,
$\eta\in \Diff^r(U)$.
For all $\wb{\gamma}\in \Diff^{r+k}_c(U)$ and
$\wb{\eta}\in \Diff^r_c(U)$,
we have
\[
(\wb{\gamma}\circ \gamma)\circ (\wb{\eta}\circ \eta)
=\wb{\gamma}\circ (\gamma\circ \wb{\eta}\circ \gamma^{-1})\circ (\gamma\circ
\eta)\,.
\]
Right translation by $\gamma\circ \eta$
being smooth,
$I_\gamma\!: \Diff^r_c(U)\to\Diff^r_c(U)$
being $C^k_\K$ (Lemma~\ref{whynomatt})
and composition
$\Diff^{r+k}_c(U)\times\Diff^r_c(U)\to \Diff^r_c(U)$
being $C^k_\K$ (Proposition~\ref{cpendsmooth}),
the preceding formula
defines a $C^k_\K$-function $\Diff^{r+k}_c(U)\times \Diff^r_c(U)\to
\Diff^r(U)$ of $(\bar{\gamma},\bar{\eta})$.
Hence the composition map
(\ref{anothcomp})
is $C^k_\K$
on the open neighbourhood
$(\Diff^{r+k}_c(U)\circ \gamma)\, \times \,(\Diff^r_c(U)\circ \eta)$
of $(\gamma,\eta)$.
Similarly, the inversion map (\ref{anothinv})
is $C^k_\K$ on the open neighbourhood
$\Diff^{r+k}_c(U)\circ \gamma$ of $\gamma$ because
\[
(\wb{\gamma}\circ \gamma)^{-1}= \gamma^{-1}\circ  \wb{\gamma}^{\,-1}=
\gamma^{-1}\circ  \wb{\gamma}^{\,-1}\circ \gamma\circ \gamma^{-1}
=(I_{\gamma^{-1}}( \wb{\gamma}^{\,-1} ))\circ \gamma^{-1}\,,
\]
where inversion $\Diff^{r+k}_c(U)\to \Diff^r_c(U)$
is $C^k_\K$ by Proposition~\ref{cpendsmooth}\,(e)
and $I_{\gamma^{-1}}$ is $C^k_\K$ by Lemma~\ref{whynomatt}\,(a).\vspace{1mm}

(c) By (b), the $C^\infty_\K$-manifold
structure
from (a) makes $\Diff^\infty(U)$
a Lie group.\vspace{1mm}

(d)
By Lemma~\ref{whynomatt}\,(b),
the automorphism $I_\gamma$ of $\Diff^r_c(U)\wt{\,}$
is $C^\infty_\K$,
for any $\gamma\in \Diff^r(U)$.
Therefore Part\,(d)
directly follows from Proposition~\ref{locchar}.
\end{proof}
\begin{defn}\label{defnDM}
Let $\K$ be a local field,
$r\in \N\cup\{\infty\}$,
and $M$ be a $\sigma$-compact
$\K$-manifold of class~$C^\infty_\K$,
of finite, positive dimension~$d$.
By Lemma~\ref{onlyopen}\,(a), there exists a $C^\infty_\K$-diffeomorphism
$\psi\!: M\to U_\psi$ from~$M$ onto
an open subset $U_\psi\sub \K^d$.
Then
\[
\Theta_\psi\!: \Diff^r(M)\to \Diff^r(U_\psi),\;\;\;
\Theta_\psi(\xi):= \psi\circ \xi \circ \psi^{-1}
\]
is an isomorphism of groups.
The map
\[
\Xi\!: C^r_c(U_\psi,\K^d)\to C^r_c(M,TM),\quad
\Xi(f)(x):=T\psi^{-1}\big(\psi(x),f(\psi(x))\big)
\]
being an isomorphism of topological vector
spaces, there exists a uniquely
determined
$C^\infty_\K$-manifold structure
on $\Diff^r(M)$,
modeled on $C^r_c(M,TM)$,
which makes the bijection
$\Theta_\psi$ a $C^\infty_\K$-diffeomorphism.
The charts of $\Diff^r(M)$
are of the form
\[
\kappa_\psi\!: \Diff^r(M)\to C^r_c(M,TM),\quad
\kappa_\psi(\gamma):=\Xi(\kappa(\Theta_\psi(\gamma)))\,,
\]
for
$\kappa\!: P_\kappa\to Q_\kappa$
ranging through the charts of
$\Diff^r(U_\psi)$.
If $r=\infty$ here,
then apparently $\Diff^\infty(M)$
is a Lie group, and $\Theta_\psi$
is an isomorphism of Lie groups.
Analogously, we make
$\Diff^\infty(M)\wt{\,}$ a Lie group modeled
on $C^\infty_c(M,TM)\wt{\,}:=\pl_{k\in \N_0} C^k_c(M,TM)$.
\end{defn}
\begin{prop}
Let $\K$ be a local field
and $M$ be a $\sigma$-compact
$C^\infty_\K$-manifold
of finite, positive dimension~$d$.
\begin{itemize}
\item[\n (a)]
For each $r\in \N\cup\{\infty\}$,
the $C^\infty_\K$-manifold structure
on $\Diff^r(M)$
is independent of the choice
of~$\psi$ in Definition~{\n \ref{defnDM}}.
It makes $\Diff^r(M)$ a topological group
and the right translation maps $\Diff^r(M)\to \Diff^r(M)$,
$\eta\mto\eta\circ \gamma$
are smooth for each $\gamma\in \Diff^r(M)$.
Furthermore, for any $k\in \N_0\cup\{\infty\}$,
both the composition map
\[
\Diff^{r+k}(M)\times \Diff^r(M)\to \Diff^r(M)
\]
and the inversion map $\Diff^{r+k}(M)\to \Diff^r(M)$ are $C^k_\K$.
\item[\n (b)]
The $C^\infty_\K$-manifold structure
on $\Diff^\infty(M)$ $($resp., $\Diff^\infty(M)\wt{\,}\,)$
is independent of the choice
of~$\psi$ in Definition~{\n \ref{defnDM}};
it makes $\Diff^\infty(M)$ $($resp., $\Diff^\infty(M)\wt{\,}\,)$
a $\K$-Lie group.
\end{itemize}
\end{prop}
\begin{proof}
We only need to show that the manifold structures
are independent of the choice of~$\psi$;
all other assertions are immediate consequences
of Proposition~\ref{DiffU}.

(a) If both $\phi\!: M\to U_\phi$ and
$\psi\!: M\to U_\psi$ are $C^\infty_\K$-diffeomorphisms
onto open subsets of $\K^d$,
then $\Theta_\phi\circ (\Theta_\psi)^{-1}\!:
\Diff^r(U_\psi)\to \Diff^r(U_\phi)$,
$\xi\mto (\phi\circ \psi^{-1}) \circ \xi\circ (\phi\circ \psi^{-1})^{-1}$
is an isomorphism of groups
which takes $\Diff^r_c(U_\psi)$
$C^\infty_\K$-diffeomorphically
onto $\Diff^r_c(U_\phi)$,
by Lemma~\ref{whynomatt}\,(a).
Since right translations
in the groups $\Diff^r(U_\phi)$ and $\Diff^r(U_\psi)$
are smooth and the homomorphism
$f:=\Theta_\phi\circ (\Theta_\psi)^{-1}$
is smooth on an open identity neighbourhood,
the usual argument shows that the homomorphism
$f$ is smooth.
Interchanging the roles of $\phi$ and $\psi$,
we see that also $f^{-1}$ is smooth.\vspace{1mm}

(b) In view of Lemma~\ref{whynomatt}\,(b),
the same argument applies to $\Diff^\infty(M)\wt{\,}$.
\end{proof}
\begin{rem}
If $M$ is a $\sigma$-compact,
finite-dimensional $C^r_\K$-manifold
for some $r\in \N$
(but not smooth),
we can still use the same arguments to make
$\Diff^r(M)$ a topological group
and $C^0_\K$-manifold
modeled on $C^0_c(M,TM)$.\,\footnote{Although
we obtain an isomorphism of topological
groups onto the $C^\infty_\K$-manifold
$\Diff^r(U_\psi)$ for each
$C^r_\K$-diffeomorphism
$\psi\!: M\to U_\psi\sub \K^d$,
we cannot expect anymore
that the corresponding
$C^\infty_\K$-manifold structures
on $\Diff^r(M)$ are
independent of the choice of $\psi$.}
\end{rem}
\appendix
\section{Proof of Proposition~4.19}\label{appfun}
In this section, we prove
the properties of function spaces
asserted in Proposition~\ref{propprop}.
In view of Remark~\ref{simplobs}\,(a)
and Lemma~\ref{amend} (applied with a countable cover
of coordinate neighbourhoods in case of (c)),
it suffices to prove assertions (a), (b) and (c) of
Proposition~\ref{propprop} when $r\in \N_0$
and $M=U$ is an open subset of~$Z$,
which we assume now.\footnote{For (a),
note that cartesian products
and closed vector subspaces of Mackey complete
topological vector spaces are Mackey complete.}\vspace{1.5mm}

(a) The proof is by induction on $r\in \N_0$.
Let us assume first that $E$ is complete.
If $r=0$, suppose that $(\gamma_\alpha)$ is
a Cauchy net in $C(U,E)_{c.o.}$.
Then $(\gamma_\alpha(x))$ is a
Cauchy net in~$E$ for each fixed element $x\in U$
and hence convergent, to $\gamma(x)\in E$, say.
For each compact subset $K\sub U$,
the restrictions $\gamma_\alpha|_K$
converge uniformly to $\gamma|_K$,
whence $\gamma|_K$ is continuous.
Hence $\gamma$ is continuous, using that~$U$
(being open in the $k$-space~$Z$) is a $k$-space.
Furthermore, $\gamma_\alpha\to \gamma$ in $C(U,E)_{c.o.}$.\\[3mm]
Induction step: Assume the assertion is correct
for some~$r$. Then both $C(U,E)$ and $C^r(U^{[1]},E)$
are complete and hence so is the topological vector space
$C^{r+1}(U,E)$,
being isomorphic to a closed vector subspace of
$C(U,E)\times C^r(U^{[1]},E)$ by Lemma~\ref{closdinprod}.
This completes the induction.\\[3mm]
If $E$ is merely sequentially complete (resp., Mackey
complete)
the sequential completeness of $C^r(U,E)$ can be proved in the same
way, replacing Cauchy nets by Cauchy sequences
(resp., Mackey-Cauchy sequences).\footnote{Note
that continuous linear maps take Mackey-Cauchy sequences to
Mackey-Cauchy sequences because they take
bounded sets to bounded sets.}\vspace{2mm}

(b) Assume that $\K\in \{\R,\C\}$
(resp., that $\K$ is an ultrametric field,
with valuation ring~$\bO$)
and $E$ is locally convex.
If $K\sub U$ is compact and $V\sub E$
is a convex, open $0$-neighbourhood
(resp., an open $\bO$-submodule submodule),
then apparently also the open $0$-neighbourhood
$\lfloor K, V\rfloor \sub C(U,E)$
is convex (resp., an $\bO$-submodule).
Hence $C(U,E)$ is locally convex.
Likewise, $C(U^{[j]},E)$ is locally convex for
each~$j$ and hence so is $C^r(U,E)$,
its topology being initial with respect to
linear maps into the spaces $C(U^{[j]},E)$
for $j\leq r$ (by definition).\vspace{2mm}

(c) Case~$r=0$: There exists an ascending sequence
$(K_j)_{j\in \N}$
of compact subsets $K_j$ of~$U$
such that $K_j$ is contained in the interior
of $K_{j+1}$ for each $j\in \N$, and $U=\bigcup_{j\in \N}K_j$.
Then every compact subset of~$U$ is contained
in $K_j$ for some~$j$, entailing that
the topology of uniform convergence
on compact sets on $C(U,E)$ is the topology
making the map $C(U,E)\to \prod_{j\in \N}
C(K_j,E)$, $\gamma\mto (\gamma|_{K_j})_{j\in \N}$
a topological embedding,
where $C(K_j,E)$ is equipped with the topology of uniform
convergence.
If $(V_n)_{n\in \N}$ is a countable basis of open
$0$-neighbourhoods in~$E$, then
$(\lfloor K_j, V_n\rfloor)_{n\in \N}$ is a countable
basis of open $0$-neighbourhoods
for $C(K_j,E)$, entailing that this space is
metrizable. We readily deduce that also $C(U,E)$ is metrizable.\\[3mm]
Suppose that $r\in \N$ now, and suppose that
the assertion of the lemma is correct
for $r-1$. By Lemma~\ref{closdinprod},
$C^r(U,E)$ is isomorphic to a topological vector subspace
of $C(U,E)\times C^{r-1}(U^{[1]},E)$.
The factors of the product
being metrizable by induction, also the product is
metrizable and hence so is
$C^r(U,E)$. This completes the proof of~(c).\vspace{2mm}

(d) Given $r$, $M$ and $E$ as described in the proposition,
let us write $C^r(M,E)_D$ for $C^r(M,E)$,
equipped with the initial topology with respect to
the family $(D^j)_{\N\ni j\geq r}$
of the mappings $D^j\!: C^r(M,E)\to C(T^jM,E)_{c.o.}$,
$\gamma\mto D^j\gamma$,
defined as follows:\label{defnDjs}
we set $D^0\gamma:=\gamma$,
let $D^1\gamma:=D\gamma:=d\gamma \!: TM\to E$ be the second component of
the tangent map $T\gamma\!: TM\to TE=E\times E$,
$T_xM\ni v\mto (\gamma(x),d\gamma(v))$,
and define $D^j \gamma:=D(D^{j-1}\gamma)\!:
T^jM:=T(T^{j-1}M)\to E$ recursively.
Before we establish Proposition~\ref{propprop}\,(d),
let us first recall some useful
properties of the spaces $C^r(M,E)_D$\,:
\begin{la}\label{appalem} In the situation of Proposition
{\rm\ref{propprop}\,(d)}, we have:
\begin{itemize}
\item[\rm (a)]
If $r=\infty$, then $C^\infty(M,E)_D=\pl_{k\in \N_0}C^k(M,E)_D$
as a topological vector space,
with the inclusion maps $C^\infty(M,E)_D\to C^k(M,E)_D$
as the limit maps.
\item[\rm (b)]
If $(U_i)_{i\in I}$ is an open cover
of~$M$, then the topology on $C^r(M,E)_D$ is initial
with respect to the family $(\rho_i)_{i\in I}$
of restriction maps $\rho_i\!: C^r(M,E)\to C^r(U_i,E)_D$,
$\rho_i(\gamma):=\gamma|_{U_i}$.
\item[\rm (c)]
If $\phi\!: M\to U \sub Z$ is a $C^r_\R$-diffeomorphism,
then $C^r(\phi,E)\!: C^r(U,E)_D\to C^r(M,E)_D$,
$\gamma\mto \gamma\circ \phi$ is an isomorphism
of topological $\K$-vector spaces.
\item[\rm (d)]
If $r\in \N_0$, then the topology on
$C^{r+1}(M,E)_D$ is initial with respect to
the maps\linebreak
$\beta_1\!: C^{r+1}(M,E)\to C(M,E)_{c.o.}$,
$\beta_1(\gamma):=\gamma$ and
$\beta_2\!: C^{r+1}(M,E)\to
C^r(TM,E)_D$, $\beta_2(\gamma):=D \gamma$.
\end{itemize}
\end{la}
\begin{proof}
(a) and (d) are immediate from the
definition of the topologies.
For a proof of (b), see
\cite[Prop.\,24.10]{INF} 
(for example).
For (c), see \cite[Prop.\,24.8]{INF}.\footnote{The cited
results are formulated in \cite{INF}
only for real manifolds, but they carry over
to complex manifolds, with identical proofs.}
See also~\cite{SEC}.
\end{proof}
In view of
Lemma~\ref{amend}, Lemma~\ref{pb2}
and their analogues compiled
in Lemma~\ref{appalem}\,(b) and\,(c),
Proposition~\ref{propprop}\,(d)
will hold in general if we can prove it in the special case
when $U:=M\sub Z$ is an open subset of the modeling space,
which we assume now.
In view of Remark~\ref{simplobs}\,(a)
and Lemma~\ref{appalem}\,(a), we may also assume
that $r\in \N_0$.
The following four lemmas will enable
us to complete the proof:
\begin{la}\label{itervsusual}
For $Z$, $E$, $U\sub Z$ as before and each $r\in \N_0$,
the topology on $C^r(U,E)_D$ is initial with respect to
the family $(d^j)_{r\geq j\in \N_0}$,
where $d^j\!: C^r(U,E)\to C(U\times Z^j,E)_{c.o.}$,
$\gamma\mto d^j\gamma$ is as in {\rm\bf\ref{differentials}}.
\end{la}
\begin{proof}
Note first that $d^j\gamma$
is a partial map of $D^j\gamma$:
There is an injective map $\kappa\!: \{1,\ldots,j\}\to
\{1,\ldots, 2^j-1\}$ (independent of $Z$, $U$ and $E$),
such that
\[
d^j\gamma(x,y)\,=\, D^j\gamma(x,\phi(y))\quad
\mbox{for all $\gamma\in C^r(U,E)$, $x\in U$ and $y\in Z^j$,}
\]
where $\phi\!: Z^j\to Z^{2^j-1}$
is the (continuous linear) map with $k$th
component
\[
\pr_k(\phi(y_1,\ldots, y_j))\;=\;
\left\{
\begin{array}{cl}
y_i & \;\mbox{if \,$\kappa(i)=k$}\\
0 &\;\mbox{else,}
\end{array}
\right.
\]
for $k=1,\ldots, 2^j-1$
and $y_1,\ldots,y_j\in Z$
\cite[Claim~2, p.\,50]{RES}.
Accordingly, the map
$d^j=C(\id_U\times \phi,E)\circ D^j\!: C^r(U,E)_D\to
C(U\times Z^j,E)_{c.o.}$
is a composition of $D^j$ and
a pullback along
a continuous map,
and thus $d^j$ is continuous
on $C^r(U,E)_D$, for
each $j\leq r$.\\[3mm]
We now show that, conversely,
each $D^j$ is continuous on $C^r(U,E)_d$,
i.e., on $C^r(U,E)$,
equipped with the topology initial with respect to
the family $(d^j)_{j\leq r}$.
To this end, we recall that
\[
D^j\gamma(x,y_1,\ldots,y_{2^j-1})
\; =\; \sum_{\ell=1}^j\,
\sum_{1\leq i_1<i_2<\cdots<i_\ell\leq 2^j-1}
c_{i_1,\ldots,i_\ell}\, d^\ell\gamma(x,y_{i_1},\ldots,y_{i_\ell})
\]
for all $\gamma\in C^r(U,E)$, $x\in U$
and $y_1,\ldots,y_{2^j-1}\in Z$,
for suitable numbers $c_{i_1,\ldots, i_\ell}\in \N_0$
which are independent of $Z$, $U$, $E$, $\gamma$,
$x$, and $y_1,\ldots,y_{2j-1}$
(cf.\ \cite[Eqn.\,(3)]{RES}).
Hence $D^j$ is a sum of terms
of the form $C(\id_U\times \psi_{i_1,\ldots,i_\ell},E)\circ d^\ell$
with a suitable continuous (linear) map
$\psi_{i_1,\ldots,i_\ell}\!: E^{2^j-1}\to E^\ell$,
and hence $D^j$ is continuous on $C^r(U,E)_d$.
Thus $C^r(U,E)_D=C^r(U,E)_d$.
\end{proof}
\begin{la}\label{complspace}
Let $Z$, $U\sub Z$ and $E$ be
as before, and $F$ be a locally convex topological
$\K$-vector space such that $E$ is a topological vector
subspace of~$F$. The the inclusion maps
\[
C^r(U,E)\to C^r(U,F)\quad\mbox{and}\quad
C^r(U,E)_D\to C^r(U,F)_D
\]
are topological embeddings.
\end{la}
\begin{proof}
In view of the definition of the topologies,
the assertion readily follows from the well-known
(and apparent) fact that the inclusion map
$C(Y,E)_{c.o.}\to C(Y,F)_{c.o.}$
is a topological embedding, for any topological space~$Y$.
(Apply this with $Y:=U^{[j]}$, resp., $Y:=U\times Z^j$,
for all $j\in \N_0$ such that $j\leq r$).
\end{proof}
\begin{la}\label{intcts}
Let $Z$, $U\sub Z$ and $E$ be
as before.
Assume that $E$ is complete.
Define
\[
\lambda : \, C(U\times [0,1],E)\to C(U,E)
\]
via
$\lambda(\gamma)(x):=
\int_0^1 \gamma(x,t)\, dt$
for $\gamma\in
C(U\times [0,1],E)$ and $x\in U$.
Then $\lambda$
is a continuous $\K$-linear map.
\end{la}
\begin{proof}
The maps $\lambda(\gamma)\!: U\to E$
are in fact continuous,
being parameter-dependent integrals with continuous
integrands (see, for example, \cite[La.\,6.15]{INF}
or \cite{GaN}).
As clearly $\lambda$ is linear,
it only remains to show that $\lambda$
is continuous at~$0$.
To this end, let $V\sub C(U,E)$ be
a $0$-neighbourhood.
Then there exists a compact subset
$K\sub U$
and a closed,
convex $0$-neighbourhood $W\sub E$
such that
$\lfloor K,W\rfloor:=\{\gamma\in C(U,E)\!: \gamma(K)\sub W\}
\sub V$.
Set $I:=[0,1]$.
Then 
$\lfloor K\times I,W\rfloor:=\{
\gamma\in C(U\times I,E)\!: \gamma(K\times I)\sub W\}$
is a $0$-neighbourhood such that
$\lambda(\lfloor K\times I,W\rfloor)\sub \lfloor K, W\rfloor\sub
V$ (cf.\ \cite[La.\,1.7]{RES}). Hence $\lambda$ is continuous.
\end{proof}
\begin{la}\label{paradep}
Given $r\in \N_0$ and an open neighbourhood
$I$ of $[0,1]$ in~$\F$, consider the set
$\Omega\sub C^r(U\times I,E)_D$ of all
$\gamma\in C^r(U\times I,E)_D$ such that the weak integrals
\[
\iota(\gamma)(x)\;  :=\, \int_0^1 \gamma(x,t)\; dt
\]
exist in~$E$ for all $x\in U$,
as well as the weak integrals $\int_0^1 d_1^j\gamma(x,t,y)\, dt$,
for all $j\leq r$ and
$(x,y)\in U\times E^j$.
Here
$d_1^j\gamma(x,t,y):=d^j(\gamma(\sbull,t))(x,y)$
for all $x\in U$, $t\in I$, and $y=(y_1,\ldots,y_j)\in Z^j$,
or, more explicitly:
\begin{equation}\label{eqcrss}
d_1^j\gamma(x,t,y_1,\ldots,y_j)\,=\, d^j\gamma((x,t),\,
(y_1,0), \ldots, (y_j,0))\,.
\end{equation}
Then $\iota(\gamma)\in C^r(U,E)$ for all $\gamma\in \Omega$, and
\begin{equation}\label{propiota}
\iota : \; C^r(U\times I, E)_D\,\supseteq \,\Omega\;
\to \; C^r(U,E)_D
\end{equation}
is a continuous $\K$-linear map.
\end{la}
\begin{proof}
(e) By \cite[Prop.\,8.7]{INF} (or \cite[La.\,7.5]{Ber})
and \cite[Prop.\,7.4]{Ber},\footnote{Erratum:
In the complex case,
$I$ should be an open neighbourhood of~$[0,1]$ in~$\C$
in \cite[La.\,7.5]{Ber}. Similar apparent adaptations
are necessary in the proof of \cite[Prop.\,7.4]{Ber}.}
the map $\iota(\gamma)$ is $C^r_\F$,
with
\begin{equation}\label{shape}
d^j\big(\iota(\gamma)\big)(x,y)\;=\; \int_0^1 d_1^j\gamma(x,t,y)\; dt\quad
\mbox{for all $j\leq r$, $x\in U$ and $y\in E^j$.}
\end{equation}
Clearly $\iota$ is linear.
To prove the continuity of~$\iota$,
after passing to the completion of~$E$ we may assume
that $E$ is complete, and hence $\Omega=C^r(U\times I,E)$,
for convenience (cf.\ Lemma~\ref{complspace}).
The topology on $C^r(U,E)_D$
being initial with respect to the maps
$d^j\!: C^r(U,E)\to C(U\times Z^j,E)_{c.o.}$
by Lemma~\ref{itervsusual},
the mapping $\iota$ will be continuous if
$d^j\circ \iota\!:$\linebreak
$C^r(U\times I,E)_D \to C(U\times Z^j ,E)$, $\gamma\mto
d^j\iota(\gamma)$
is continuous for $j\leq r$.
By (\ref{shape}) and (\ref{eqcrss}),
$d^j\circ \iota$ is
a composition of the continuous
map $d^j\!:
C^r(U\times I,E)_D\to C((U\times I)\times (Z\times \F)^j,E)$,\linebreak
the continuous pullback $C(f,E)\!:
C((U\times I)\times (Z\times\F)^j,E)\to C(U\times I\times Z^j,E)$,\linebreak
$\eta\mto\eta\circ f$ with $f(x,t,y_1,\ldots,y_j):=
(x,t,\, (y_1,0),\ldots, (y_j,0))$,
and the integration map
$\lambda_j\!:$\linebreak
$C(U\times I\times Z^j, E)\to C(U\times Z^j,E)$,
$\lambda_j(\eta)(x,y):=\int_0^1 \eta(x,t,y)\; dt$,
which is continuous as a consequence
of Lemma~\ref{intcts}.
Hence $d^j\circ \iota$ is
continuous for each $j\leq r$,
and hence so is $\iota$,
as asserted.
\end{proof}
We are now in the position to complete the proof
of Proposition~\ref{propprop}\,(d).
The proof is by induction on~$r\in \N_0$.
The case $r=0$ is trivial:
By definition, both $C^0(U,E)$ and $C^0(U,E)_D$
are equipped with the compact-open topology,
and hence coincide.\\[3mm]
Induction step: Assume that Proposition~\ref{propprop}\,(d) holds
for some $r\in \N_0$.
Consider the map
$f\!: C^{r+1}(U,E)\to C^{r+1}(U,E)_D$, $f(\gamma):=\gamma$.
For $\beta_1$ and $\beta_2$ as in Lemma~\ref{appalem}\,(d),
the composition
\[
\beta_1\circ f \!: C^{r+1}(U,E)\to C(U,E)_{c.o.}\,,
\quad \gamma\mto \gamma
\]
is continuous by Remark~\ref{simplobs}\,(b),
and also
\[
\beta_2\circ f\!: C^{r+1}(U,E)\to C^r(TU,E)_D\,,
\quad \gamma\mto d\gamma=\gamma^{[1]}(\sbull,0)
\]
is continuous as it is a composition of the continuous map
$(\sbull)^{[1]}\!: C^{r+1}(U,E)\to C^r(U^{[1]},E)$
(Remark~\ref{simplobs}\,(b)),
the restriction map
$C^r(U^{[1]},E)\to C^r(TU,E)$ which is a pullback
and hence continuous (Lemma~\ref{pullback}),
and the identity map $C^r(TU,E)\to C^r(TU,E)_D$
which is an isomorphism of topological vector spaces
by induction.
The topology on $C^{r+1}(U,E)_D$ being initial with respect to
the maps $\beta_j$ ($j\in \{1,2\}$),
the continuity of the compositions $\beta_j\circ f$
entails that~$f$ is continuous.\\[3mm]
It remains to show that also $f^{-1}$ is continuous.
In view of Remark~\ref{simplobs}\,(b),
we only need to show that $\alpha_1\circ f^{-1}$
and $\alpha_2\circ f^{-1}$ are continuous,
where $\alpha_1\!: C^{r+1}(U,E)\to C(U,E)$
is the inclusion map and $\alpha_2\!: C^{r+1}(U,E)\to C^r(U^{[1]},E)$,
$\alpha_2(\gamma):=\gamma^{[1]}$.
Here $\alpha_1\circ f^{-1}=\beta_1$ is continuous.
To see that $\alpha_2\circ f^{-1}$ is continuous, note that
for any $x\in U$ and $y\in Z$, we find an open balanced $0$-neighbourhood $J_{(x,y)}\sub \F$
and open neighbourhoods $V_{(x,y)}\sub U$ of~$x$ and $W_{(x,y)}\sub Z$
of $y$ such that
\begin{equation}\label{VJW}
V_{(x,y)}+2\, J_{(x,y)}\, W_{(x,y)}\;\sub\; U
\end{equation}
(and thus $V_{(x,y)}\times W_{(x,y)}\times 2J_{(x,y)}\sub U^{]1[}$).
Note that $U^{[1]}$ is covered by $U^{]1[}$,
together with the sets $Y_{(x,y)}:=V_{(x,y)}\times
W_{(x,y)}\times J_{(x,y)}$ for $(x,y)\in U\times Z$.
Hence, by Lemma~\ref{amend}, the map $\alpha_2\circ
f^{-1}$ will be continuous if we can show that
\[
\phi : \, C^{r+1}(U,E)_D\, \to \, C^r(U^{]1[},E)\, ,\quad
\phi(\gamma)\, :=\, \gamma^{[1]}|_{U^{]1[}}\, =\,
\gamma^{]1[}
\]
is continuous, as well as the mappings
\[
\psi_{(x,y)} : \,
C^{r+1}(U,E)_D \, \to\, C^r(Y_{(x,y)},E)\,,
\quad
\psi_{(x,y)}(\gamma)\, :=\, \gamma^{[1]}|_{Y_{(x,y)}}\,,
\]
for all $(x,y)\in U\times Z$.
The formula $\phi(\gamma)(x,y,t)=\frac{1}{t}
(\gamma(x+ty)-\gamma(x))$ shows that $\phi$ is built up from
the following continuous maps and hence continuous:
1. The map $C^{r+1}(U,E)_D\to C^r(U^{]1[}, E)$,
$\gamma\mto [(x,y,t)\mto \gamma(x+ty)]$.
This map is the composition of the inclusion map
$C^{r+1}(U,E)_D\to C^r(U,E)_D=C^r(U,E)$
(which is continuous by definition of the $D$-topologies
and the induction hypothesis) and the pullback
$C^r(g,E)\!: C^r(U,E)\to C^r(U^{]1[},E)$
along $g\!: U^{]1[}\to U$, $g(x,y,t):=
x+ty$, which is continuous by Lemma~\ref{pullback}.
\,2. The map $C^{r+1}(U,E)_D\to C^r(U^{]1[},E)$,
$\gamma\mto [(x,y,t)\mto \gamma(x)]$,
which is continuous by the same argument.
\,3. The addition map of the topological vector
space $C^r(U^{]1[},E)$.
\,4. The multiplication operator
$m_\tau\!: C^r(U^{]1[},E)\to C^r(U^{]1[},E)$,
$\gamma\mto \tau\cdot \gamma$ with $\tau\!: U^{]1[}\to \F$,
$\tau(x,y,t):=\frac{1}{t}$, which is continuous by
Lemma~\ref{multop}. Hence $\phi$ is continuous.\\[3mm]
In view of Lemma~\ref{complspace},
it suffices to prove continuity of the maps
$\psi_{(x,y)}$ when $E$ is complete, which we assume
now (we can always replace $E$ by its completion).
For the proof, fix $(x,y)\in U\times Z$; we abbreviate
$V:=V_{(x,y)}$, $W:=W_{(x,y)}$, $J:=J_{(x,y)}$,
$Y:=Y_{(x,y)}=V\times W\times J$,
and $\psi:=\psi_{(x,y)}$. Let $B_2(0):=B_2^\F(0)\sub \F$.
Then
\[
\psi(\gamma)(u,v,t)\, =\,
\gamma^{[1]}(u,v,t)\, =\,
\int_0^1 d\gamma(u+stv,v)\, ds\quad
\mbox{for all $(u,v,t)\in Y$}
\]
(see \cite[Prop.\,7.4]{Ber} and its proof).
We can therefore interpret $\psi$ as a composition
\begin{eqnarray*}
C^{r+1}(U,E)_D &\to &
C^r(U\times Z,E)_D\, =\,
C^r(U\times Z,E)\, \to\,
C^r(Y\times B_2(0),E)\\
&=& C^r(Y\times B_2(0),E)_D\, \to\,
C^r(Y,E)_D\, =\, C^r(Y,E)
\end{eqnarray*}
of the continuous map $D\!: C^{r+1}(U,E)_D\to C^r(U\times Z,E)_D$
(Lemma~\ref{appalem}\,(d)), the continuous pullback
$C^r(h,E)\!: C^r(U\times Z,E)\to C^r(Y\times B_2(0),E)$
along $h\!: Y\times B_2(0)\to U\times Z$,
$h(u,v,t,s):=(u+stv,v)$ (Lemma~\ref{pullback}), and the integration map
\[
\iota :\,
C^r(Y\times B_2(0),E)_D\, \to\, C^r(Y,E)_D
\]
taking $\eta$ to the map $(u,v,t)\mto
\int_0^1 \eta(u,v,t,s)\, ds$.\,\footnote{Note that all
of the required weak integrals exist, by completeness of~$E$.}
Here $\iota$ is continuous by Lemma~\ref{paradep}.
Hence $\psi$ is continuous, and hence so is $f^{-1}$,
which completes the proof of Proposition~\ref{propprop}.\Punkt
\section{Proof of Proposition 4.23}\label{appcruc}
Part\,(a) follows from Part\,(b)
in an obvious way.
We therefore only need to prove\,(b).
Furthermore, as in the proof of Proposition~\ref{globcruc},
we may assume that $r,k\in \N_0$.\\[1mm]
\begin{center}
{\bf Reduction to open subsets of {\boldmath $X$, $\wb{X}$ and $Z$}}
\end{center}
As $\F$ is locally compact and $X$ finite-dimensional,
we deduce that $X$ is locally compact
(cf.\ \cite[I, \S2, No.\,3, Thm.\,2]{BTV}), and hence so
are~$M$ and~$Y$.
For each $i:=(x,\bar{x})\in Y\times \wb{M}=:I$, we find a chart
$\tau_i\!: W_i\to N_i$
of~$N$
around $\sigma(x)$,
onto an open subset $N_i\sub Z$.
Let $\bar{\kappa}_i\!: \wb{S}_i\to \wb{M}_i\sub \wb{X}$
be a chart of $\wb{M}$ around~$\bar{x}$, and
$\kappa_i\!: S_i\to
M_i\sub X$
be a chart of~$M$ around~$x$
such that $S_i\sub \sigma^{-1}(W_i)$.
As $Y$
is locally compact,
we find a compact neighbourhood
$C_i$ of $x$ in~$Y$
such that $C_i\sub S_i$.
We let $U_i:=C_i^0$ be the interior of~$C_i$.
Define
$K_i:=\kappa_i(C_i)$,
$Y_i:=\kappa_i(U_i)=K_i^0$,
and $\eta_i:=\kappa_i|_{U_i}^{Y_i}$.
Then $Y=\bigcup_{i\in I} U_i$,
and $\{\eta_i\times \bar{\kappa}_i\!: i\in I\}$
is an atlas for~$Y\times \wb{M}$.
Abbreviate $\Theta_i:=\theta_{\eta_i\times\bar{\kappa}_i}\!:
C^r(Y\times \wb{M},F)\to C^r(Y_i\times \wb{M}_i,F)$,
$\theta_i:=\theta_{\kappa_i}\!:
C^r(M,E)\to C^r(M_i,E)$, and
$\bar{\theta}_i:=\theta_{\bar{\kappa}_i}\!:
C^r(\wb{M},E)\to C^r(\wb{M}_i,E)$
(see (\ref{dag}) in {\bf\ref{top2}} for the notations).
Then
\[
\lambda\!: C^r(Y\times \wb{M},F)
\to \prod_{i\in I} C^r(Y_i\times \wb{M}_i,F),\;\;\;
\gamma\mto (\Theta_i(\gamma))_{i\in I}=
(\gamma\circ (\eta_i^{-1}\times \wb{\kappa}_i^{\,-1}))_{i\in I}
\]
is a topological embedding
(Lemma~\ref{atlas})
whose image is easily seen to be closed
(cf.\ Lemma~\ref{amend}).
Hence, by Lemma~\ref{corestr} and Lemma~\ref{inprod},
the map $\phi$ will be of class~$C^k_\K$
if we can show that $\Theta_i\circ \phi$
is of class $C^k_\K$ for each $i\in I$.
Using the $C^r_\F$-map
$\sigma_i:=\tau_i\circ \sigma|_{S_i}^{W_i}\circ \eta_i^{-1}\!:
Y_i\to N_i$
and the $C^{r+k}_\K$-map
$\tilde{f}_i:= \tilde{f}\circ (\tau_i^{-1}\times\id_U\times\id_{\wb{E}}
\times\id_P)\!:
N_i \times U\times \wb{E}\times P\to F$, we define
$f_i:=\tilde{f}_i\circ (\sigma_i\times \id_U\times \id_{\wb{E}}
\times \id_P)\!:
Y_i\times U\times \wb{E} \times P\to F$.
For $\gamma\in \lfloor K,U\rfloor_r$,
$\bar{\gamma}\in C^r(\wb{M},\wb{E})$,
$x\in Y_i$, $\bar{x}\in \wb{M}_i$, and $p\in P$, we have
\begin{eqnarray*}
\Theta_i(\phi(\gamma,\bar{\gamma},p))(x,\bar{x}) &=&
f(\eta_i^{-1}(x),\gamma(\eta_i^{-1}(x)),
\bar{\gamma}(\bar{\kappa}_i^{-1}(\bar{x})),p)
= f_i(x,\theta_i(\gamma)(x),\bar{\theta}_i(\bar{\gamma})(\bar{x}),p)\\
&=& \phi_i(\theta_i(\gamma),\bar{\theta}_i(\bar{\gamma}),p)(x,\bar{x})\, ,
\end{eqnarray*}
where
\[
\phi_i\!: \lfloor K_i,U\rfloor_r\times C^r(\wb{M}_i,\wb{E})\times
P\to C^r(Y_i\times \wb{M}_i,F),\;\;\;
\phi_i(\eta,\bar{\eta},p)(x,\bar{x}):=f_i(x,\eta(x),\bar{\eta}(\bar{x}),p)\]
with $\lfloor K_i,U\rfloor_r\sub C^r(M_i,E)$.
Thus
$\Theta_i\circ\phi=\phi_i\circ
(\theta_i|_{\lfloor K,U\rfloor_r}^{\lfloor K_i,U\rfloor_r}\times
\bar{\theta}_i\times \id_P)$,
showing that $\Theta_i\circ \phi$ will be of class~$C^k_\K$
if~$\phi_i$ is of class~$C^k_\K$.
Thus, each $M_i$ being an open subset
of~$X$, $\wb{M}_i$ open in $\wb{X}$,
and $N_i$ an open subset of~$Z$,
the proposition will hold if
we can prove (b) in the case where~$M$,
$\wb{M}$ and $N$ are
open subsets of~$X$, $\wb{X}$ and $Z$,
respectively, which we shall assume
for the rest of the proof. We proceed by induction on~$k$.\pagebreak

\begin{center}
{\bf The case {\boldmath $k=0$}.}
\end{center}
The proof is by induction on~$r$.
If $r=0$, then the topology on $C^0(M,E)$,
$C^0(\wb{M},\wb{E})$ and $C^0(Y\times \wb{M},F)$
is the topology of uniform convergence on compact
sets (see {\bf \ref{newold}}).
Let $\gamma \in \lfloor K,U \rfloor\sub C^r(M,E)$, $\bar{\gamma}\in
C^r(\wb{M},\wb{E})$,
$p\in P$, $L$ be a compact subset of~$Y\times \wb{M}$,
and $V\sub F$ be an open zero-neighbourhood.
Let $W\sub F$ be an open zero-neighbourhood
such that $W-W\sub V$.
For each $i:=(x,\bar{x})\in L$,
there are open neighbourhoods $A_i\sub Y$
and $\wb{A}_i\sub \wb{M}$ of $x$, resp., $\bar{x}$,
and open zero-neighbourhoods
$B_i\sub E$, $\wb{B}_i\sub \wb{E}$
and $C_i\sub H$
such that $\gamma(A_i)+B_i\sub \gamma(K)+B_i\sub U$,
$p+C_i\sub P$,
and
\[
f(y,u,\bar{u},q)-f(x,\gamma(x),\bar{\gamma}(x),p)\in W
\]
for all $y\in A_i$, $u\in \gamma(A_i)+B_i$,
$\bar{u}\in \bar{\gamma}(\wb{A}_i)+\wb{B}_i$,
and $q\in p+C_i$.
By compactness, $L\sub \bigcup_{i\in I}(A_i\times \wb{A}_i)$
for some finite subset $I\sub L$.
Then $B:=\bigcap_{i\in I}B_i\sub E$,
$\wb{B}:=\bigcap_{i\in I}\wb{B}_i\sub \wb{E}$
and $C:=\bigcap_{i\in I}C_i\sub H$
are open zero-neighbourhoods,
and $\wb{K}:=\pr_2(L)\sub \wb{M}$ is compact,
where $\pr_2\!:Y\times \wb{M}\to\wb{M}$
is the coordinate projection.
Let $\xi\in \gamma+\lfloor K,B\rfloor\sub C^0(M,E)$,
$\bar{\xi}\in \bar{\gamma}+\lfloor\wb{K},\wb{B}\rfloor\sub
C^0(\wb{M},\wb{E})$, 
and $q\in p+C\sub P$.
Given $(y,\bar{y})\in L$,
there is $i=(x,\bar{x})\in I$ such that
$(y,\bar{y})\in A_i\times \wb{A}_i$.
Then $\xi\in \lfloor K,U\rfloor$, and
\begin{eqnarray*}
\lefteqn{f(y,\xi(y),\bar{\xi}(\bar{y}), q)
-f(y,\gamma(y),\bar{\gamma}(\bar{y}),p)}\\
& = & \!f(y,\xi(y),\bar{\xi}(\bar{y}),q)-
f(x,\gamma(x),\bar{\gamma}(\bar{x}),p)
\,-\,(f(y,\gamma(y),\bar{\gamma}(\bar{y}),p)-
f(x,\gamma(x),\bar{\gamma}(\bar{x}),p))\\
&\in& \!W-W\, \sub \,V\,.
\end{eqnarray*}
We have shown that
$\phi(\xi,\bar{\xi},q)-\phi(\gamma,\bar{\gamma},p)
\in \lfloor L,V\rfloor\sub
C(Y\times \wb{M},F)$
for all $(\xi,\bar{\xi},q)$
in the open neighbourhood
$(\gamma+\lfloor K,B\rfloor)\times (\bar{\gamma}+\lfloor \wb{K},\wb{B}\rfloor)
\times (p+C)$
of $(\gamma,\bar{\gamma},p)$.
Thus $\phi$ is continuous.\\[3mm]
{\em Induction step on~$r$.}
We write $\phi_r$ for~$\phi$,
to emphasize its dependence on~$r$.
Suppose the assertion of the lemma is correct
for $k=0$ and some $r\in \N_0$.
Suppose that the hypotheses of the lemma
are satisfied by $\tilde{f}$ and $\sigma$,
with $r$ replaced by $r+1$.
The map $\phi_r$ being continuous,
we see as in the proof of Proposition~\ref{globcruc}
that $\phi_{r+1}$ will be continuous if we can show that the map
\[
\psi\!: \lfloor K,U\rfloor_{r+1} \times C^{r+1}(\wb{M},\wb{E})\times
P \to C^r((Y\times \wb{M})^{[1]},F),\;\;\;
\psi(\gamma,\bar{\gamma},p):=\phi_{r+1}(\gamma,\bar{\gamma},p)^{[1]}
\]
is continuous.
Here we have, by the Chain Rule,
\begin{eqnarray}
\lefteqn{\!\!\!\!\!\!\!\!\psi(\gamma,\bar{\gamma},p)
(x,\bar{x},y,\bar{y},t)\quad\quad}\nonumber\\
& = &
\phi_{r+1}(\gamma,\bar{\gamma},p)^{[1]}(x,\bar{x},y,\bar{y},t)\nonumber\\
& = &
\tilde{f}^{[1]}((\sigma(x),\gamma(x),\bar{\gamma}(\bar{x}),p),\;
(\sigma^{[1]}(x,y,t),
\gamma^{[1]}(x,y,t),\bar{\gamma}^{[1]}(\bar{x},\bar{y},t),0),\;t)\label{jstar}
\end{eqnarray}
for $(x,\bar{x},y,\bar{y},t)\in (Y\times \wb{M})^{[1]}$.
Let $(\gamma_0,\bar{\gamma}_0,p_0)\in \lfloor
K,U\rfloor_{r+1}\times C^{r+1}(\wb{M},\wb{E})\times P$ be given;
our goal is to show that~$\psi$
is continuous at~$(\gamma_0,\bar{\gamma}_0,p_0)$.
We set $X_1:=X\times X\times \F$,
$\wb{X}_1:=\wb{X}\times \wb{X}\times \F$,
$Z_1:=Z\times Z \times \K$,
$E_1:=E\times E$, and $\wb{E}_1:=\wb{E}\times \wb{E}$.
Given $(x_0,\bar{x}_0,y_0,\bar{y}_0,t_0)\in (Y\times \wb{M})^{[1]}$,
we have $(\sigma(x_0),\sigma^{[1]}(x_0,y_0,t_0),t_0)\in N^{[1]}$
and $(\gamma_0(x_0),\gamma_0^{[1]}(x_0,y_0,t_0),t_0)\in U^{[1]}$.
There are open neighbourhoods
$R_1 \sub N$ of $\sigma(x_0)$,
$R_2\sub Z$ of $\sigma^{[1]}(x_0,y_0,t_0)$,
$R_3\sub \K$ of~$t_0$,
$V_1\sub U$ of $\gamma_0(x_0)$,
and $V_2\sub E$ of $\gamma_0^{[1]}(x_0,y_0,t_0)$
such that
\[
R_1\times R_2\times R_3\sub N^{[1]}\;\;\;\;
\mbox{and}\;\;\;\;
V_1\times V_2\times R_3\sub U^{[1]}\,.
\]
Then $R_1\times V_1\times \wb{E}\times P\times R_2\times V_2\times\wb{E}
\times \{0\}\times
R_3\sub (N\times U\times \wb{E} \times P)^{[1]}$.
Abbreviate $N_1:=R_1\times R_2\times R_3$
and $U_1:=V_1\times V_2$.
Then
\[
\tilde{f}_1\!: N_1\times U_1\times \wb{E}_1 \times P\to F,\;\;\;\;
\tilde{f}_1(x,y,t,u,v,\bar{u},\bar{v},p):=
\tilde{f}^{[1]}((x,u,\bar{u},p),\, (y,v,\bar{v},0),\,t)
\]
for $(x,y,t)\in N_1$, $(u,v)\in U_1=V_1\times V_2$,
$(\bar{u},\bar{v})\in \wb{E}_1$,
$p\in P$ defines a mapping of class
$C^{r+1+k-1}_\K=C^{r+k}_\K$ on the open subset
$N_1\times U_1\times \wb{E}_1 \times P$
of $Z_1\times E_1\times \wb{E}_1\times H$.
There exists an open neighbourhood
$A_1\sub Y$ of~$x_0$ with compact closure $C_1\sub Y$,
and open neighbourhoods
$A_2\sub X$ of~$y_0$,
$\wb{A}_1\sub \wb{M}$ of $\bar{x}_0$,
$\wb{A}_2\sub \wb{X}$ of $\bar{y}_0$,
and
$A_3\sub \F\cap R_3$ of~$t_0$
such that $M_1:=A_1\times A_2\times A_3\sub Y^{[1]}$,
$\wb{M}_1:=\wb{A}_1\times \wb{A}_2\times A_3\sub \wb{M}^{[1]}$,
$\sigma(A_1)\sub R_1$,
$\sigma^{[1]}(A_1\times A_2\times A_3)\sub R_2$,
$\gamma_0(C_1)\sub V_1$,
and $\gamma_0^{[1]}(A_1\times A_2\times A_3)\sub V_2$.
Let $K_1\sub M_1$ be a compact neighbourhood
of $(x_0,y_0,t_0)$, with interior $Y_1:=K_1^0$.
Define $\sigma_1\!: Y_1\to N_1$,
$\sigma_1(x,y,t):=(\sigma(x),\sigma^{[1]}(x,y,t),t)$
for $x\in A_1$, $y\in A_2$, $t\in A_3$.
Then $\sigma_1$ is a $C^r_\F$-map.
We set
\[
f_1:=\tilde{f}_1\circ (\sigma_1\times \id_{U_1}\times
\id_{\wb{E}_1}\times\id_P)\!:
Y_1 \times U_1\times \wb{E}_1\times P\to F.
\]
Abbreviate $z:=(x_0,\bar{x}_0,y_0,\bar{y}_0,t_0)$.
Then
\[
\psi_1\!:
\lfloor K_1,U_1\rfloor_r\times C^r(\wb{M}_1,\wb{E}_1)
\times P\to C^r(Y_1\times \wb{M}_1,F),\;\;\;\;
\psi_1(\gamma,\bar{\gamma},p):=f_1(\sbull,p)_*(\gamma\times \bar{\gamma})
\]
is a continuous
mapping by induction,
where
$\lfloor K_1,U_1\rfloor_r\sub C^r(M_1,E_1)$.
Let $B_1\sub A_1$, $B_2\sub A_2$ and $B_3\sub A_3$
be open neighbourhoods of~$x_0$, $y_0$
and $t_0$, respectively, such that
$B_1\times B_2\times B_3\sub Y_1$.
Then
$Q_z:=B_1\times \wb{A}_1\times B_2\times \wb{A}_2\times B_3$
is an open neighbourhood of~$z$
in $(Y\times \wb{M})^{[1]}$.
We define
$\delta\!:
Q_z\to
Y_1\times \wb{M}_1$, $\delta(x,\bar{x},y,\bar{y},t):=
(x,y,t,\bar{x},\bar{y},t)$.
Since $\psi_1$ is continuous, also
\[
\psi_2:=C^r(\delta,F)\circ \psi_1\!:
\lfloor K_1,U_1\rfloor_r\times C^r(\wb{M}_1,\wb{E}_1)
\times P\to C^r(Q_z,F)
\]
is continuous (by Lemma~\ref{pullback}).
Note that
\[
\Omega:=\{\gamma\in \lfloor K,U\rfloor_{r+1}\cap
\lfloor C_1,V_1\rfloor_{r+1}\!:
\gamma^{[1]}|_{M_1}\in \lfloor K_1,V_2\rfloor_r \} 
\]
is an open neighbourhood
of $\gamma_0$ in $\lfloor K,U\rfloor_{r+1}$
(cf.\ Remark~\ref{simplobs}\,(b),
Lemma~\ref{pullback} and Lemma~\ref{flooropen}).
The linear maps
\[
\;\;\;\;\pi\!: C^{r+1}(M,E)\to C^r(M_1,E),\;\;\;\;
\pi(\gamma)(x,y,t):=\gamma(x)\quad\mbox{and}
\]
\[
\!\!\!\!\bar{\pi}\!: C^{r+1}(\wb{M},\wb{E})\to C^r(\wb{M}_1,\wb{E}),\;\;\;\;
\bar{\pi}(\bar{\gamma})(\bar{x},\bar{y},t):=\bar{\gamma}(\bar{x})
\]
are continuous (see Remark~\ref{simplobs}\,(a) and Lemma~\ref{pullback}),
and we have $\pi(\Omega)\sub \lfloor K_1,V_1\rfloor_r$.
Also
\[
\;\;\;\;C^{r+1}(M,E)\to C^r(M_1,E),\;\;\;\;
\gamma\mto \gamma^{[1]}|_{M_1}\quad\mbox{and}
\]
\[
\!\!\!\!C^{r+1}(\wb{M},\wb{E})\to C^r(\wb{M}_1,\wb{E}),\;\;\;\;
\bar{\gamma}\mto \bar{\gamma}^{[1]}|_{\wb{M}_1}
\]
are continuous linear mappings (see Remark~\ref{simplobs}
and Lemma~\ref{pullback}),
and the first of these maps $\Omega$ into $\lfloor K_1,V_2\rfloor_r$.
Let $\rho_{z}\!:
C^r((Y\times \wb{M})^{[1]},F)\to C^r(Q_z,F)$,
$\rho_z(\eta):=\eta|_{Q_z}$
be the restriction map.
In view of (\ref{jstar}) and the definition of $\psi_2$,
we have
\[
\psi(\gamma,\bar{\gamma},p)(x,\bar{x},y,\bar{y},t)=\psi_2
\left((\pi(\gamma),\gamma^{[1]}|_{M_1}),\,
(\bar{\pi}(\bar{\gamma}),\bar{\gamma}^{[1]}|_{\wb{M}_1}),\,
p\right)(x,\bar{x},y,\bar{y},t)
\]
for $(\gamma,\bar{\gamma},p)\in \Omega\times C^{r+1}(\wb{M},\wb{E})
\times P$,
$(x,\bar{x},y,\bar{y},t)\in Q_z$,
showing that
\[
\rho_z\circ \psi|_{\Omega\times C^{r+1}(\wb{M},\wb{E})\times P}
\]
is continuous
on $\Omega\times C^{r+1}(\wb{M},\wb{E})\times P$, which is a neighbourhood
of $(\gamma_0,\bar{\gamma}_0,p_0)$.
Thus, we have achieved the following:
given any $z=(x_0,\bar{x}_0,y_0,\bar{y}_0,t_0)\in (Y\times \wb{M})^{[1]}$,
we have found an open neighbourhood
$Q_z$ of $z$ in $(Y\times \wb{M})^{[1]}$
such that
$\rho_z\circ \psi$ is continuous
at $(\gamma_0,\bar{\gamma}_0,p_0)$.
In view of Lemma~\ref{cover},
this entails that $\psi$ is
continuous at $(\gamma_0,\bar{\gamma}_0,p_0)$,
as desired.
\begin{center}
{\bf Induction step on~{\boldmath $k$}.}
\end{center}
Suppose the assertion of the lemma
is correct
for some $k\in \N_0$ and all $r\in \N_0$.
Let $\sigma$ and $\tilde{f}$ be given
which satisfy the hypotheses of the lemma
when $k$ is replaced with $k+1$.
Then $\phi\!: \lfloor K,U\rfloor_r\times C^r(\wb{M},\wb{E})\times
P\to C^r(Y\times \wb{M},F)$
is of class $C^k_\K$ (and thus continuous),
by induction.
For all $(\gamma,\bar{\gamma},p,\eta,\bar{\eta},q,t)\in
(\lfloor K,U\rfloor_r\times C^r(\wb{M},\wb{E})\times P)^{]1[}\sub
\lfloor K,U\rfloor_r\times C^r(\wb{M},\wb{E})\times
P\times C^r(M,E)\times C^r(\wb{M},\wb{E})\times H\times\K$,
we calculate
\begin{eqnarray}
\lefteqn{{\textstyle \!\!\!\!\!\!
\frac{1}{t}(\phi(\gamma+t\eta,\bar{\gamma}+t\bar{\eta},
p+tq)-\phi(\gamma,\bar{\gamma},p))\,(x,\bar{x})}\quad\quad}
\nonumber\\
& = &
{\textstyle
\frac{1}{t}(\tilde{f}(\sigma(x),
\gamma(x)+t\eta(x),\bar{\gamma}(\bar{x})+t\bar{\eta}(\bar{x}),
p+tq)-\tilde{f}(\sigma(x),\gamma(x),\bar{\gamma}(\bar{x}),p))}\nonumber\\
& = &
\tilde{f}^{[1]}((\sigma(x),\gamma(x),\bar{\gamma}(\bar{x}),p),\;
(0,\eta(x),\bar{\eta}(\bar{x}),q),\; t)\label{badref}
\end{eqnarray}
for all $x\in Y$ and $\bar{x}\in \wb{M}$.
On the open subset
\[
W\sub N \times U\times E \times \wb{E}\times \wb{E}
\times P\times H\times \K
\]
consisting of those
$(x,y,z,\bar{y},\bar{z},p,q,t)$
such that
$(y,z,t)\in U^{[1]}$ and
$(p,q,t)\in P^{[1]}$,
we define a $C^{r+k+1-1}_\K=C^{r+k}_\K$-map
\[
\tilde{h}\!: W\to F,\;\;\;\;
\tilde{h}(x,y,z,\bar{y},\bar{z},p,q,t):=\tilde{f}^{[1]}((x,y,\bar{y},p),\;
(0,z,\bar{z},q),\; t)\,.
\]
Motivated by (\ref{badref}), we consider the map
\[
\begin{array}{c}
\psi\!: (\lfloor K,U\rfloor_r\times C^r(\wb{M},\wb{E})\times
P)^{[1]}\to C^r(Y\times \wb{M},F)\,,\\
\psi((\gamma,\bar{\gamma},p),\; (\eta,\bar{\eta},q),\; t)(x,\bar{x}):=
\tilde{h}(\sigma(x),\gamma(x),\eta(x),\bar{\gamma}(\bar{x}),
\bar{\eta}(\bar{x}),p,q,t)
\end{array}
\]
which extends $\phi^{]1[}$.
If we can show that~$\psi$
is $C^k_\K$,
then $\psi$ is continuous and thus
$\psi=\phi^{[1]}$.
Thus $\phi$ will be of class~$C^1_\K$,
with $\phi^{[1]}=\psi$ of class~$C^k_\K$,
entailing that $\phi$ is of class $C^{k+1}_\K$ (see
{\bf \ref{usefulsimp}}),
as required.\\[3mm]
Let $\cA$ be the set of all relatively
compact, non-empty, open subsets of~$Y$.
As a consequence of
Lemmas \ref{corestr}, \ref{inprod},
\ref{cover} and \ref{amend},
the mapping $\psi$ will be of class $C^k_\K$
if we can show that
$\rho_Q\circ \psi$
is of class $C^k_\K$
for all $Q\in \cA$,
where $\rho_Q\!: C^r(Y\times \wb{M},F)\to C^r(Q\times \wb{M},F)$,
$\rho_Q(\gamma):=\gamma|_{Q\times \wb{M}}$.
Here
$\rho_Q\circ \psi$ will be of class~$C^k_\K$
if we can show that
every
$(\gamma_0,\bar{\gamma}_0,p_0,\eta_0,\bar{\eta}_0,
q_0,t_0)\in (\lfloor K, U\rfloor_r\times C^r(\wb{M},\wb{E})\times P)^{[1]}$
has an open neighbourhood~$T$
such that
$(\rho_Q\circ \psi)|_T$ is of class $C^k_\K$
(see Lemma~\ref{Crlocal}).\\[3mm]
{\em Case $t_0\not=0$}:
Then we can take $T:=
(\lfloor K,U\rfloor_r\times C^r(\wb{M},\wb{E})\times P)^{]1[}$.
In fact, $\psi|_T=\phi^{]1[}$ (and thus $\rho_Q\circ \psi|_T$)
is of class $C^k_\K$
since so is~$\phi$.\\[3mm]
{\em Case $t_0=0$}:
In view of the compactness
of $\gamma_0(\wb{Q})\sub U$ and $\eta_0(\wb{Q})\sub E$,
where $\wb{Q}$ denotes the closure of~$Q$ in~$Y$,
there exist open neighbourhoods
$B_1$ of $\gamma_0(\wb{Q})$ in~$U$,
$B_2$ of $\eta_0(\wb{Q})$ in~$E$
and
an open zero-neighbourhood $B_3\sub \K$
such that $B_1+B_3\cdot B_2\sub U$
and thus $B_1\times B_2\times B_3\sub U^{[1]}$.
Shrinking $B_3$ if necessary,
we find
open neighbourhoods $B_4$ of $p_0$ in $P$
and $B_5$ of~$q_0$ in~$H$
such that $B_4+B_3\cdot B_5\sub P$
and thus $B_4\times B_5\times B_3\sub P^{[1]}$.
Define $U_1:=B_1\times B_2\sub U\times E$
and $P_1:=B_4\times B_5\times B_3\sub P\times H\times \K$.
Then $N\times U_1\times \wb{E}\times \wb{E} \times P_1\sub W$, and
$\tilde{g}:=\tilde{h}|_{N\times U_1\times \wb{E}^2\times P_1}\!:
N\times U_1\times \wb{E}^2\times P_1\to F$
is a mapping of class~$C^{r+k}_\K$. We define
\[
g:=\tilde{g}\circ (\sigma|_Q \times \id_{U_1}\times \id_{\wb{E}^2}\times
\id_{P_1})
\!: Q\times U_1\times \wb{E}^2\times P_1 \to F\, .
\]
By induction hypothesis,
the mapping
\[
\Phi\!: \lfloor \wb{Q},U_1\rfloor_r\times C^r(\wb{M},\wb{E}^2)\times
P_1\to C^r(Q\times \wb{M},F),\;\;\;
\Phi(\gamma,\bar{\gamma},p):=g(\sbull,p)_*(\gamma\times \bar{\gamma})
\]
is of class~$C^k_\K$,
where $\lfloor \wb{Q}, U_1\rfloor_r\sub C^r(M,E\times E)$.
Note that
\[
\begin{array}{c}
\theta\!:
(C^r(M,E)\times C^r(\wb{M},\wb{E})\times H)^2\times
\K\to C^r(M,E\times E)\times C^r(\wb{M},\wb{E}\times \wb{E})\times
H\times H\times \K,\\
\! (\gamma,\bar{\gamma},p,\eta,\bar{\eta},q,t)\mto ((\gamma,\eta),
(\bar{\gamma},\bar{\eta}),p,q,t)\quad\quad
\end{array}
\]
is an isomorphism of topological
$\K$-vector spaces
such that $(\gamma_0,\bar{\gamma}_0,p_0,\eta_0,\bar{\eta}_0,q_0,t_0)\in
T:=\theta^{-1}(\lfloor \wb{Q},U_1 \rfloor_r\times
C^r(\wb{M},\wb{E}^2)\times P_1)
\cap (\lfloor K , U \rfloor_r\times C^r(\wb{M},\wb{E})\times P)^{[1]}$.
To complete the proof,
it only remains to observe that
$\rho_Q\circ \psi|_T=
\Phi\circ \theta|_T^{\theta(T)}$
is a mapping of class~$C^k_\K$.
\section{Proof of Proposition~11.3}\label{appcomparbit}
In the situation of Proposition~\ref{comparbit},
let $H$ be a finite-dimensional
$\K$-vector space, and $P\sub H$
be open. If we can prove the following lemma,
then apparently Proposition~\ref{comparbit} will follow:
\begin{la}\label{lacompa}
$\;\;$The mapping
$\;\,\Theta\!:
C^{r+k}(F\times P,E)\times C^r(M,F)\times P\to
C^r(M,E)$,\\
$\Theta(\gamma,\eta,p):=\gamma(\sbull,p)\circ \eta\,$
is of class~$C^k_\K$.
If $k\geq 1$,
then
\begin{eqnarray}
\lefteqn{\!\!\!\!\!\!\!\!\!\!\!\!\!\!\!\!\Theta^{[1]}((\gamma,\eta,p),
\,(\gamma_1,\eta_1,p_1),\,t)\quad\quad\quad\quad}\nonumber \\
&=&
\gamma^{[1]}((\sbull,p),\, (\sbull,p_1),\, t)
\circ (\eta,\eta_1)
+\gamma_1(\sbull,p+tp_1)\circ (\eta+t\eta_1)\label{formcomp5}
\end{eqnarray}
for all $((\gamma,\eta,p),\, (\gamma_1,\eta_1,p_1),\,t)\in
(C^{r+k}(F\times P,E)\times C^r(M,F)\times P)^{[1]}$.
\end{la}
\begin{proof}
As in the proof of Lemma~\ref{compsmooth},
we may assume that
$k,r\in \N_0$.
The proof is by induction on~$k$.
\begin{center}
{\bf The case {\boldmath $k=0$}}
\end{center}
We proceed by induction on~$r$.
If $r=0$, then a trivial variation of the argument
used in the proof of Lemma~\ref{compsmooth} shows
that $\Theta$ is continuous.\\[3mm]
{\em Induction step on $r$.}
Let $r\in \N$, and suppose that the lemma
holds for $k=0$, when $r$ is replaced with
$r-1$.
It then suffices to show continuity of
$\Theta$ in the case
where $M$ is an open subset of its
modeling space~$Z$.
In fact, if
$M$ is a $C^r_\K$-manifold,
Lemma~\ref{amend} entails that $\Theta$ will be continuous
if, for any chart $\kappa\!: U\to V\sub Z$
of~$M$, the map
\[
h\!:
C^r(F\times P,E)\times C^r(M,F)\times P\to
C^r(V,E),\quad
h(\gamma,\eta,p):=\Theta(\gamma,\eta,p)\circ \kappa^{-1}
\]
is continuous.
But
\begin{equation}\label{noend2}
h(\gamma,\eta, p)\,=\,\gamma(\sbull,p)\circ (\eta\circ \kappa^{-1})
\,=\, \Xi(\gamma,\eta\circ \kappa^{-1},p)\,,
\end{equation}
where $\Xi\!:
C^r(F\times P,E)\times C^r(V,F)\times P\to
C^r(V,E)$, $\Xi(\gamma,\sigma,p):=\gamma(\sbull,p)\circ
\sigma$.
Recall that
the pullback $C^r(M,F)\to C^r(V,F)$, $\eta\mto\eta\circ\kappa^{-1}$
is continuous linear (Lemma~\ref{pb2}).
Thus (\ref{noend2}) shows that $h$ will be continuous
if $\Xi$ is continuous. Since $V$ is open in~$Z$,
this completes the reduction step to the case
where $M$ is open in~$Z$.\\[3mm]
To complete the induction step on~$r$ in the case $k=0$,
by the preceding we may assume now that $M$ is an
open subset of~$Z$.
The map
$\Theta\!:
C^r(F\times P,E)\times C^r(M,F) \times P\to
C^r(M,E)$ is continuous as a map into $C(M,E)$,
by the case $r=0$ already settled.
Hence, in view of Remark~\ref{simplobs}\,(b),
$\Theta$ will be continuous if we can show
that the map
\[
C^r(F\times P,E)\times C^r(M,F) \times P\to
C^{r-1}(M^{[1]},E),\quad (\gamma,\eta,p)\mto 
\Theta(\gamma,\eta,p)^{[1]}
\]
is continuous, where
\begin{equation}\label{formkl3}
\Theta(\gamma,\eta,p)^{[1]}(x,y,t)=\gamma^{[1]}((\eta(x),p),\,
(\eta^{[1]}(x,y,t),0),\, t)
\end{equation}
for all $(x,y,t)\in M^{[1]}$, by the Chain Rule.
By (\ref{formkl3}), we have
\begin{equation}\label{fm10}
\Theta(\gamma,\eta,p)^{[1]}=
\wt{\Theta}(\gamma^{[1]}\circ \rho,\,
(\eta\circ \pr_1,\eta^{[1]},\pr_3),\, p)
\end{equation}
for all $(\gamma,\eta,p)$ in the domain of~$\Theta$,
where $\pr_1\!: M^{[1]}\to M$, $(x,y,t)\mto x$
and $\pr_3\!: M^{[1]}\to\K$, $(x,y,t)\mto t$
are the coordinate projections,
\[
\rho\!: F^{[1]}\times P\to (F\times P)^{[1]},\quad
\rho((u,v,t),\, p):=((u,p),\, (v,0),\, t)
\]
is continuous linear, and
$\wt{\Theta}\!: C^{r-1}(F^{[1]}\!\times\! P,E)\times
C^{r-1}(M^{[1]},F^{[1]})\times P\to C^{r-1}(M^{[1]},E)$,
\[
\wt{\Theta}(\sigma,\tau,p):=\sigma(\sbull,p)\circ \tau
\]
is continuous by induction.
Because the mapping $C^r(M,F)\to C^{r-1}(M^{[1]},F)$, $\eta\mto \eta^{[1]}$
and both of the pullbacks
$C^r(M,F)\to C^{r-1}(M^{[1]},F)$, $\eta\mto \eta\circ \pr_1$
and $C^{r-1}(\rho,E)\!:$\linebreak
$C^{r-1}((F\times P)^{[1]},E)\to C^{r-1}(F^{[1]}\times P,E)$
are continuous,
we deduce from (\ref{fm10}) that\linebreak
$(\gamma,\eta,p)\mto
\Theta(\gamma,\eta,p)^{[1]}$ is continuous.
Hence so is $\Theta$.
\begin{center}
{\bf Induction step on~{\boldmath $k$}}
\end{center}
Let $k\in \N$, and
suppose that the assertion of
the lemma holds when $k$ is replaced with
$k-1$, for all integers $r\in \N_0$. Let $r\in \N_0$.
Given an element
$((\gamma,\eta,p),\,(\gamma_1,\eta_1,p_1),\, t)\in
(C^{r+k}(F\times P,E)\times C^r(M,F) \times P)^{]1[}$,
we calculate for $x\in M$:
\begin{eqnarray}
\lefteqn{\!\!\!\!\!\!\!\!\!{\textstyle \frac{1}{t}}
\left(\Theta(\gamma+t\gamma_1,\eta+t\eta_1,p+tp_1)-\Theta(\gamma,\eta,p)\right)
(x)\quad\quad\quad}\nonumber\\
&=&
{\textstyle \frac{1}{t}}\Big(
\gamma(\eta(x)+t\eta_1(x), p+tp_1)-\gamma(\eta(x),p)\Big)
\; +\,\gamma_1(\eta(x)+t\eta_1(x),p+tp_1)\nonumber\\
&=&
\gamma^{[1]}((\eta(x),p),\, (\eta_1(x),p_1),\, t)
\; +\, \gamma_1(\eta(x)+t\eta_1(x),p+tp_1)\, .\label{hbreak2}
\end{eqnarray}
Thus
$\Theta^{]1[}$ coincides with the restriction
to
$(C^{r+k}(F\times P,E)\times C^r(M,F) \times P)^{]1[}$
of the mapping
$g\!: (C^{r+k}(F\times P,E)\times C^r(M,F) \times P)^{[1]}\to C^r(M,F)$,
\begin{equation}\label{tnot}
g((\gamma,\eta,p),\, (\gamma_1,\eta_1,p_1),\,t)
:=
\wt{\Theta}(\gamma^{[1]}\circ \rho,\, (\eta,\eta_1),\, (p,p_1,t))
\,+\,
\Theta(\gamma_1,\eta+t\eta_1,p+tp_1)\,,
\end{equation}
where $\rho\!: F^2 \times P^{[1]}\to
(F\times P)^{[1]}$, $\rho(u,v,p,p_1,t):=((u,p),\, (v,p_1),\, t)$
is continuous linear and
$\wt{\Theta}\!: C^{r+k-1}(F^2\!\times \!P^{[1]}, E)\times
C^r(M,F^2)\times P^{[1]}\to C^r(M,E)$,
\[
\wt{\Theta}(\sigma,\tau,\,(p,p_1,t)):=
\sigma(\sbull,(p,p_1,t))\circ \tau
\]
is of class $C^{k-1}_\K$, by induction.
Since $\Theta$ is $C^{k-1}_\K$ and hence continuous
as a consequence of the induction hypothesis,
in order that $\Theta$ be $C^k_\K$,
it only remains to show that $g$ is of class
$C^{k-1}_\K$ (then $\Theta^{[1]}=g$,
which also establishes (\ref{formcomp5})).
Now, the second summand in (\ref{tnot})
clearly describes a $C^{k-1}_\K$-map.
Also the first summand is $C^{k-1}_\K$,
as pullbacks are continuous linear
and $\wt{\Theta}$ is $C^{k-1}_\K$.
\end{proof}
\section{Smoothness vs.\ weak smoothness over local fields}\label{grothend}
We vary a classical result of A. Grothendieck concerning mappings on
open sets in $\R^n$:
\begin{thm}\label{thmgroth}
Let $\K$ be a local field,
$E$ and $F$ be topological $\K$-vector spaces,
$f\!: U\to F$ a map on an open
subset of~$E$, and $k\in \N_0$.
If $E$ is metrizable and $F$ is Mackey complete and locally convex,
then we have:
\begin{itemize}
\item[\rm (a)]
If $f$ is $C^k$, then $f$ is weakly $C^k$,
viz.\ $\lambda\circ f$ is $C^k$ for each $\lambda\in F'$.
\item[\rm (b)]
If $f$ is weakly $C^{k+1}$, then $f$ is $C^k$.
\end{itemize}
In particular, $f$ is smooth if and only if $f$ is weakly smooth.
\end{thm}
\begin{proof}
(a) is a trivial consequence of the Chain Rule.\vspace{.7mm}

(b) If we can prove (b) in the special case where
$U=E$ is finite-dimensional,
then for general $f$ the composition
$f\circ \gamma$ will be weakly $C^{k+1}$ and thus $C^k$,
for every smooth map $\gamma\!: \K^{k+1}\to U$.
Hence $f$ will be $C^k$ by \cite[Thm.\,12.4]{Ber}.
We may therefore assume that $U=E=\K^\ell$ for some~$\ell$.
The proof is by induction on~$k$.\vspace{1mm}

If $k=0$ and $f\!: E=\K^\ell\to F$ is weakly $C^1$,
let $x\in U$ and
$W\sub F$ be a $0$-neighbourhood. Define
\[
B\; :=\; \{f^{]1[}(x,y,t)\!:y\in K, t\in \bO\setminus\{0\}\}\,,
\]
where $\bO\sub \K$ is the valuation ring,
$K$ a compact $0$-neighbourhood in~$E$, and $f^{]1[}(x,y,t):=
t^{-1}(f(x+ty)-f(x))$.
Then $\lambda(B)\sub (\lambda\circ f)^{[1]}(\{x\}\times K\times\bO)$
is compact and thus bounded, for each $\lambda\in F'$,
whence $B$ is bounded in~$F$ by \cite[Thm.\,4.21]{Thi}.
Thus there is $t\in \bO\setminus \{0\}$
such that $tB\sub W$. Then
$f(y)-f(x)\in tB\sub W$
for every $y\in x+tK$. Hence $f$ is continuous
at~$x$.\vspace{1mm}

{\em Induction step.\/} If $k\geq 1$ and $f\!: E=\K^\ell\to F$
is weakly $C^{k+1}$,
given $x,y\in E$
choose a sequence $(t_n)_{n\in \N}$ of pairwise distinct elements
in $\bO\setminus \{0\}$
such that $t_n\to 0$. Set
\[
B\; :=\; \{(t_n-t_m)^{-1}(f^{]1[}(x,y,t_m)
-f^{]1[}(x,y,t_n)): n,m\in \N\}\,.
\]
Then $\lambda(B)=\{
(\lambda\circ f)^{[2]}((x,y,t_n),\,(0,0,1),\,t_m-t_n)\!: n,m\in \N\}$
is contained in the compact set
$(\lambda\circ f)^{[2]}(\{x\}\times\{y\}\times \bO\times \{(0,0,1)\}
\times \bO)$ and thus
bounded, for each $\lambda\in F'$, whence $B$ is bounded
by \cite[Thm.\,4.21]{Thi}.
Since $f^{]1[}(x,y,t_m)
-f^{]1[}(x,y,t_n) \in (t_m-t_n) B$,
we deduce that $(f^{]1[}(x,y,t_n))_{n\in \N}$
is a Mackey-Cauchy sequence in~$F$ and thus convergent;
we let $g(x,y,0)$ be its limit.
Then $\lambda(g(x,y,0))=\lim_{n\to\infty}(\lambda\circ f)^{]1[}(x,y,t_n)=
(\lambda\circ f)^{[1]}(x,y,0)$
for each $\lambda$. Furthermore, trivially
$\lambda(g(x,y,t))=(\lambda\circ f)^{[1]}(x,y,t)$
for $g(x,y,t):=t^{-1}(f(x+ty)-f(x))$ for all
$(x,y,t)\in E\times E\times \K^\times$.
Thus $\lambda\circ g=(\lambda\circ f)^{[1]}$
is $C^k$ for each $\lambda$, whence
$g$ is $C^{k-1}$ (and thus $C^0$),
by induction. Hence $f$ is $C^1$
with $f^{[1]}=g$ of class $C^{k-1}$, whence
$f$ is $C^k$.\vspace{-6mm}
\end{proof}
\section{Towards a {\boldmath $p$}-adic analogue
of Boman's Theorem}\label{appboman}
Boman's Theorem asserts that a mapping $f\!: \R^n\to\R$
is smooth if and only if $f\circ \gamma\!: \R\to\R$
is smooth for each smooth curve $\gamma\!: \R\to
\R$ (see \cite[Thm.\,1]{Bom}).
As is well known,
this implies that a mapping $f\!: U\to F$
from an open subset $U\sub E$ of a metrizable
real locally convex space~$E$ to a real locally convex
space~$F$ is smooth if and only
if $f\circ \gamma$ is smooth for each smooth curve
$\gamma\!: \R\to U$ (cf.\ \cite[Thm.\,12.8]{KaM}).
It is natural to ask whether
versions of Boman's theorem remain valid
over arbitrary locally compact topological fields,
or at least in the $p$-adic case.
A positive answer to this question (which is still open)
would be quite useful.\footnote{For example,
combining the results of this section and \cite[Cor.\,4.3]{ANA}, this would
entail that the exponential map of a smooth
(not necessarily analytic) $p$-adic Banach-Lie
group is automatically a $C^\infty$-diffeomorphism.
This would be an important step towards proving the conjecture
that every smooth $p$-adic Banach-Lie group admits a
smoothly compatible $p$-adic analytic Lie group structure.}
In the present section,
we show that, as in the real
locally convex case, the validity of Boman's Theorem
for functions $\K^2\to\K$ 
entails its validity for functions
between open subsets of suitable
topological vector spaces.
Our main tools are
the exponential laws from Section~\ref{secexplaw},
the adaptations of Grothendieck's Theorem
from Appendix~\ref{grothend},
and the characterization of
smooth maps on metrizable topological vector spaces
from \cite[Thm.\,12.4]{Ber}.
\begin{la}\label{preBom}
Let $\K$ be the field of real numbers or an
ultrametric field.
Let $\cA$ be one of the following
classes of
topological $\K$-vector spaces:
\begin{itemize}
\item
The class of all topological $\K$-vector spaces;
\item
The class of all sequentially complete
topological $\K$-vector spaces;
\item
The class of all Mackey complete
topological $\K$-vector spaces;
\item
The class of all locally convex topological $\K$-vector spaces;
\item
The class of all sequentially complete,
locally convex topological $\K$-vector spaces; or:
\item
The class of all Mackey complete,
locally convex topological $\K$-vector spaces.
\end{itemize}
Suppose that
Boman's Theorem holds for mappings
from $\K^2$ to topological vector spaces in~$\cA$,
i.e., assume the validity of the following statement:\\[3mm]
\hspace*{6mm}If $f\!: \K^2\to F$ is a map into a topological
$\K$-vector space~$F\in \cA$,
and $f\circ \gamma\!: \K\to F$\linebreak
\hspace*{6mm}is $\,C^\infty_\K$
for each $C^\infty_\K$-curve $\gamma\!: \K\to \K^2$,
then $f$ is of class $C^\infty_\K$.\\[3mm]
Then Boman's Theorem holds for mappings $f\!: U\to F$
from open subsets $U\sub E$ of metrizable
topological $\K$-vector spaces $E$
to topological $\K$-vector spaces $F\in \cA$,
i.e., smoothness of $f\circ \gamma\!: \K\to F$
for each smooth curve $\gamma\!: \K\to U$
entails smoothness of~$f$.
\end{la}
\begin{proof}
Assuming that the hypothesis of the lemma is correct,
we first show by induction on $n\in \N$
that Boman's Theorem holds for mappings
$f\!: \K^n\to F$, for all $F\in \cA$.
The case $n=1$ is trivial, and the case $n=2$
holds by the hypothesis of the lemma.
Thus, assume that $n\geq 2$ now,
assume that
Boman's Theorem holds for functions on $\K^n$,
and assume that
$f\!: \K^{n+1}\to F$ is a function into some
$F\in \cA$ such that $f\circ \gamma$ is smooth for
each smooth curve $\gamma\!: \K\to \K^{n+1}$.
Then $f(x,\sbull)\!: \K^{n-1}\to F$
is smooth along smooth curves, for any $x\in \K^2$,
and thus $f(x,\sbull)$ is smooth, by induction.
Therefore\linebreak
$f^\vee\!: \K^2\to C^\infty(\K^{n-1},F)$,
$f^\vee(x):=f(x,\sbull)$
is correctly defined.
By Proposition~\ref{expmetriz}\,(b),
the map $f=(f^\vee)^\wedge$
will be smooth if $f^\vee$
is smooth. Since $C^\infty(\K^{n-1},F)\in \cA$
by Proposition~\ref{propprop}, $f^\vee$
will be smooth if $f^\vee\circ \gamma\!:
\K\to C^\infty(\K^{n-1},F)$
is smooth for each smooth curve $\gamma=(\gamma_1,\gamma_2)\!: \K\to \K^2$,
by the hypothesis of the lemma.
By Proposition~\ref{expmetriz}\,(b),
$f^\vee\circ \gamma$ will be smooth
if $g:=(f^\vee\circ \gamma)^\wedge\!: \K^n\to F$
is smooth. Note that $g(t,y)=f(\gamma_1(t),\gamma_2(t),y)$
for $t\in \K$, $y\in \K^{n-1}$,
and thus $(g\circ \eta)(t)=f(\gamma_1(\eta_1(t)),\gamma_2(\eta_1(t)),
\eta_2(t),\ldots, \eta_n(t))$
is smooth for each smooth curve $\eta=(\eta_1,\ldots,\eta_n)\!:
\K\to\K^n$, as $f$ is smooth along smooth curves.
Being smooth along smooth curves,
$g\!: \K^n\to F$ is smooth, by the induction hypothesis.
In view of our reduction steps, this means that also
$f$ is smooth. This finishes the induction.\\[3mm]
To complete the proof, let $U$ be an open subset of a metrizable
topological $\K$-vector space, $F\in \cA$
and $f\!: U\to F$ be a mapping which is smooth along smooth curves.
Then the composition $f\circ \eta\!: \K^n\to F$
is smooth along smooth curves and hence smooth by what
has already been shown,
for any $n\in \N$ and smooth map $\eta\!: \K^n\to U$.
Hence $f$ is smooth, by \cite[Thm.\,12.4]{Ber}.
\end{proof}
\begin{thm}\label{thmBom}
Let $\K$ be a local field.
Suppose that
Boman's Theorem holds for functions
from $\K^2$ to $\K$,
i.e., assume the validity of the following statement:\\[3mm]
\hspace*{6mm}If $f\!: \K^2\to \K$ is a function such
that $f\circ \gamma\!: \K\to \K$ is $C^\infty_\K$
for each $C^\infty_\K$-curve\linebreak
\hspace*{6mm}$\gamma : \, \K\to \K^2$,
then $f$ is of class $C^\infty_\K$.\\[3mm]
Then Boman's theorem holds for mappings $f\!: U\to F$
from open subsets $U\sub E$ of metrizable
topological $\K$-vector spaces $E$
to Mackey complete,
locally convex topological $\K$-vector spaces~$F$,
i.e., smoothness of $f\circ \gamma\!: \K\to F$
for each smooth curve $\gamma\!: \K\to U$
entails smoothness of~$f$.
\end{thm}
\begin{proof}
Let $\cA$ be the class of Mackey complete,
locally convex topological $\K$-vector spaces.
We only need to show that the present hypothesis
entails the hypothesis of Lemma~\ref{preBom}.
To this end, let $f\!: \K^2\to F$
be a map into a Mackey
complete, locally convex topological $\K$-vector space which
is smooth along each smooth curve.
If $\lambda\!: F\to\K$ is a continuous linear functional,
then $\lambda\circ f\!: \K^2\to\K$ is smooth along
each smooth curve and hence smooth,
by the hypothesis of the present theorem.
Thus $f$ is weakly smooth and hence smooth,
by our analogue of Grothendieck's Theorem
(Theorem~\ref{thmgroth}).
Thus, the hypothesis of Lemma~\ref{preBom} is verified.
\end{proof}
In the real case, we do have Boman's
Theorem available. Thus, we arrive at the conclusion:
\begin{prop}\label{Bomreal}
Let $E$ be a
metrizable real topological
vector space, $U\sub E$ be an open subset
and $F$ a locally convex real topological vector space.
Then a mapping $f\!: U\to F$ is smooth if and only if
$f\circ \gamma$ is smooth for each smooth curve
$\gamma\!: \R\to U$.
\end{prop}
\begin{proof}
We only need to verify the hypothesis
of Lemma~\ref{preBom} for $\K=\R$ and the class~$\cA$
of all locally convex, real
topological vector spaces.
Thus, suppose that $f\!: \R^2\to F$
is a map into a real locally convex space which is
smooth along smooth curves.
Then $\lambda\circ f\!: \R^2\to\R$
is smooth along smooth curves and hence
smooth by Boman's Theorem \cite[Thm.\,1]{Bom},
for each continuous linear functional
$\lambda$ on the completion $\wb{F}$
of~$F$. Hence $f\!: \R^2\to \wb{F}$ is weakly smooth
and hence smooth,
by Grothendieck's classical theorem.
For any $n\in \N$ and
$x,y\in \R^2$, we have $d^nf(x,y,\ldots,y)=\frac{d^n}{dt^n}\big|_{t=0}f(x+ty)
\in F$, as $f$ is smooth along the smooth curve
$\R\to\R^2$, $t\mto x+ty$.
Hence
$d^nf(x,\sbull)$ has image in~$F$,
by polarization (cf.\ \cite[Thm.\,A]{BS1}).
Now, $f\!: \R^2\to F$ being a smooth map
into $\wb{F}$ with \,$\im(d^nf)\sub F$
for each $n\in \N_0$, we deduce that $f$
is smooth as a map into~$F$,
as required.
\end{proof}
\begin{rem}
Note that $E$ need not be locally convex in
Proposition~\ref{Bomreal}, and that
no completeness properties whatever
are presumed for $E$ nor $F$.
Therefore, the result is slightly more general
than the results in the literature
(and those in the folklore).
\end{rem}
\section{Spaces of sections in vector bundles
and mappings between them}\label{appsections}
In this appendix,
we define
bundles of topological vector spaces,
topologize their spaces of sections,
and study differentiability properties
of mappings between open subsets of spaces of sections.\\[3mm]
In the following,
$\K$ denotes an arbitrary
topological field.
Additional properties of
$\K$ (local compactness)
will be spelt out explicitly where they are needed.
Unless specified otherwise,
$r\in \N_0\cup\{\infty\}$.
\begin{defn}\label{defnbdle}
Let $M$ be a $C^r_\K$-manifold,
modeled on a topological $\K$-vector space~$Z$,
and $F$ be a topological~$\K$-vector space.
A $C^r_\K$-{\em vector bundle\/} over $M$, with typical fibre~$F$,
is a $C^r_\K$-manifold~$E$,
together with a $C^r_\K$-surjection $\pi\!: E\to M$
and equipped with a
$\K$-vector space structure on each fibre
$E_x:=\pi^{-1}(\{x\})$,
such that for each $x_0\in M$, there exists
an open neighbourhood $M_\psi$ of~$x_0$ in~$M$
and a $C^r_\K$-diffeomorphism
\[
\psi\!: \pi^{-1}(M_\psi)\to M_\psi\times F\]
(called a ``local trivialization of $E$ about $x_0$'')
such that
$\psi(E_x)=\{x\}\times F$ for each $x\in M_\psi$
and $\pr_F\circ \psi|_{E_x}\!: E_x\to
F$ is $\K$-linear
(and thus an isomorphism of topological $\K$-vector spaces
with respect to the topology on $E_x$ induced by~$E$),
where $\pr_F\!: M_\psi\times F\to F$ is the projection on
the second coordinate.
\end{defn}
\begin{rem}\label{chnge}
In the situation of Definition~\ref{defnbdle},
suppose we are given two local trivializations
$\psi\!: \pi^{-1}(M_\psi)\to M_\psi\times F$
and $\phi\!: \pi^{-1}(M_\phi)\to M_\phi\times F$.
Then $\phi(\psi^{-1}(x,v))=(x, g_{\phi,\psi}(x).v)$
for some function $g_{\phi,\psi}\!: M_\phi\cap M_\psi\to
\GL(F)\sub L(F,F)$ (the space of continuous
linear self-maps),
and $G_{\phi,\psi}\!:(M_\phi\cap M_\psi)\times F\to F$,
$(x,v)\mto g_{\phi,\psi}(x).v$
is a $C^r_\K$-map (since $\phi\circ \psi^{-1}$
is so).
\end{rem}
\begin{defn}\label{bscssec}
A {\em $C^r_\K$-section\/} of a $C^r_\K$-vector bundle
$\pi\!: E\to M$ is a $C^r_\K$-map $\sigma\!:
M\to E$ such that $\pi\circ \sigma=\id_M$.
Its {\em support\/} $\Supp(\sigma)$
is the closure of $\{x\in M\!:
\sigma(x)\not=0_x \}$.
We let $C^r(M,E)$ be the set of all
$C^r_\K$-sections in~$E$.
If $\K$ is locally compact and $M$ is finite-dimensional,
we let
$C^r_K(M,E)$ be the set of all
$C^r_\K$-sections with support
contained
in a given compact subset $K\sub M$.
\end{defn}
Making use of scalar multiplication
and addition in the individual fibres,
we obtain natural vector space structures
on $C^r(M,E)$ and $C^r_K(M,E)$.
The zero-element is
the {\em zero-section\/} $0_\sbull\!:
M\to E$, $x\mto 0_x\in E_x$.
\begin{defn}\label{sigsi}
If $\pi\!: E\to M$ is a $C^r_\K$-vector bundle
with typical fibre~$F$,
$\sigma\!:M\to E$ a $C^r_\K$-section,
and $\psi\!:\pi^{-1}(M_\psi)\to M_\psi\times F$ a local
trivialization, we define
$\sigma_\psi:=\pr_F\circ \psi\circ \sigma|_{M_\psi}\!: M_\psi\to F$.
Thus $\psi(\sigma(x))=(x,\sigma_\psi(x))$
for all $x\in M_\psi$.
\end{defn}
Note that $\sigma_\psi$ is a mapping of class~$C^r_\K$ here.
The symbols $g_{\phi,\psi}$,
$G_{\phi,\psi}$, and $\sigma_\psi$
will always be used with the meanings just described,
without further explanation.
\begin{defn}
If $\pi\!: E\to M$ is a vector bundle
and $\cA$ a set of local trivializations $\psi$ of~$E$
whose domains cover~$E$,
then we call $\cA$ an {\em atlas\/} of local trivializations.
\end{defn}
\begin{la}\label{famgivessec}
If $\pi\!: E\to M$ is a $C^r_\K$-vector bundle
with typical fibre~$F$, and $\cA$ an atlas
of local trivializations for~$E$,
then
\[
\Gamma\!: C^r(M,E)\to \prod_{\psi\in\cA} C^r(M_\psi,F),\;\;\;\;
\sigma\mto (\sigma_\psi)_{\psi\in\cA}
\]
is an injection, whose image is the closed vector subspace
\[
H\, :=\,
\big\{(f_\psi)\in \prod_{\psi\in \cA}
C^r(M_\psi,F)\!: \;(\forall \phi,\psi\in\cA,\;
\forall x\in M_\phi\cap M_\psi)\;
f_\phi(x)=g_{\phi, \psi}(x).f_\psi(x)\,\big\}
\]
of $\, \prod_{\psi\in\cA} C^r(M_\psi,F)$.
\end{la}
\begin{proof}
It is obvious that $\, \im\,\Gamma\sub H$,
and clearly~$\Gamma$ is injective.
If now $(f_\psi)_{\psi\in\cA}\in H$,
we define $\sigma\!: M\to E$ via
$\sigma(x):=\psi^{-1}(x,f_\psi(x))$ if $x\in M_\psi$.
By definition of~$H$, $\sigma(x)$ is independent of the choice of~$\psi$.
As $\psi\circ \sigma|_{M_\psi}=(\id_{M_\psi},f_\psi)$,
the mapping $\sigma\!: M\to E$ is of class~$C^r_\K$.
Thus $\sigma$ is a $C^r_\K$-section,
and $\Gamma(\sigma)=(f_\psi)_{\psi\in\cA}$
by definition of~$\sigma$.
We deduce that $\im\,\Gamma=H$.
The closedness of~$H$ follows
from the continuity of the
point evaluations $C^r(M_\psi,F)\to F$,
$\gamma\mto \gamma(x)$ for $x\in M_\psi$.
\end{proof}
\begin{defn}\label{dfntops}
Let $\pi\!:E\to M$ be a $C^r_\K$-vector bundle, with typical fibre~$F$,
and $\cA$ be the atlas of all local trivializations
of~$E$.
We give $C^r(M,E)$ the vector topology making the
linear mapping
\[
\Gamma\!: C^r(M,E)\to\prod_{\psi\in\cA}C^r(M_\psi, F),\;\;\;
\sigma\mto (\sigma_\psi)_{\psi\in\cA}\]
a topological embedding.
\end{defn}
By the preceding definition,
the topology on $C^r(M,E)$ is initial
with respect to the family
$(\theta_\psi)_{\psi\in\cA}$,
where $\theta_\psi\!: C^r(M,E)\to C^r(M_\psi,F)$,
$\sigma\mto \sigma_\psi$.
\begin{rem}\label{seclcx}
In the situation of
Definition~\ref{dfntops},
assume that $\K$ is $\R$, $\C$ or
an ultrametric field, and assume that $F$ is locally convex.
Then $C^r(M_\psi,F)$ is locally convex for
each $\psi\in \cA$, by Proposition~\ref{propprop}\,(b).
Hence also
$C^r(M,E)$ is locally convex.
\end{rem}
It suffices to work with any atlas
of local trivializations.
\begin{la}\label{anyatlas}
The topology on $C^r(M,E)$ described in Definition~{\n \ref{dfntops}\/}
is initial with respect to $(\theta_\psi)_{\psi\in \cB}$,
for any atlas $\cB\sub \cA$ of local trivializations
for~$E$.
\end{la}
\begin{proof}
We let $\cO$ be the initial topology
on $C^r(M,E)$ with respect to $(\theta_\psi)_{\psi\in \cB}$,
which apparently is coarser than initial topology
with respect to $(\theta_\psi)_{\psi\in\cA}$.
Fix a local trivialization
$\phi\in \cA$.
Then $\{M_\phi\cap M_\psi\!:\psi\in \cB\}$
is an open cover for $M_\phi$.
In view of Lemma~\ref{amend}, the map
$\theta_\phi$
will be continuous on $(C^r(M,E),\cO)$ if the map
$(C^r(M,E),\cO)\to
C^r(M_\phi\cap M_\psi, F)$,
$\sigma\mto \theta_\phi(\sigma)|_{M_\phi\cap M_\psi}$
is continuous for all $\psi\in \cB$.
But, with $G_{\phi,\psi}$
as in Remark~\ref{chnge}, the latter mapping is the
composition of $(G_{\phi,\psi})_*\!: C^r(M_\phi\cap M_\psi,F)\to
C^r(M_\phi\cap M_\psi,F)$ (which is continuous
by Proposition~\ref{globcruc})
and $C^r(M,E)\to C^r(M_\phi\cap M_\psi,F)$,
$\sigma\mto \theta_\psi(\sigma)|_{M_\phi\cap M_\psi}$,
which is continuous by Lemma~\ref{pb2}.
Thus, $\theta_\phi$
being continuous on $(C^r(M,E),\cO)$
for each $\phi\in\cA$,
the topology~$\cO$ is finer than the initial
topology with respect to the family $(\theta_\phi)_{\phi\in \cA}$,
which completes the proof.
\end{proof}
\begin{defn}\label{defnparvect}
Suppose that $\pi_1\!:E_1\to M$ and $\pi_2\!: E_2\to M$
are $C^r_\K$-vector bundles over the same base,
with typical fibres $F_1$ and $F_2$, respectively.
A mapping $f\!: U \to E_2$,
defined on an open subset~$U$
of $E_1$,
is called a {\em bundle map\/}
if it preserves fibres, {\em i.e.},
$f(U\cap (E_1)_x)\sub (E_2)_x$ for all $x\in M$.
Then, given local trivializations
$\psi\!: \pi_1^{-1}(M_\psi)\to M_\psi\times F_1$
and $\phi\!:\pi_2^{-1}(M_\phi)\to M_\phi\times F_2$
of $E_1$, resp., $E_2$, we have
\[
\phi(f(\psi^{-1}(x,v)))=(x, f_{\phi,\psi}(x,v))
\]
for all $(x,v)\in U_{\phi,\psi}:=\psi(U\cap E_1|_{M_\psi\cap M_\phi})
\sub (M_\psi\cap M_\phi)\times F_1$,
for a uniquely determined mapping
\[
f_{\phi,\psi}\!: U_{\phi,\psi}\to F_2\, .
\]
\end{defn}
\begin{thm}\label{pushfbdl}
Let
$r,k\in \N_0\cup\{\infty\}$,
$M$ be a $($not necessarily finite-dimensional$)$
$C^{r+k}_\K$-manifold and
$\pi_1\!: E_1\to M$,
$\pi_2\!: E_2\to M$ be $C^{r+k}_\K$-vector bundles
over~$M$,
whose typical fibres $F_1$,
resp., $F_2$ are arbitrary topological
$\K$-vector spaces. Then the following holds:
\begin{itemize}
\item[\rm (a)]
If $H$ is a topological $\K$-vector space,
$P\sub H$ an open subset
and
$f\!: E_1\times P \to E_2$ a $C^{r+k}_\K$-map
such that $f_p:=f(\sbull,p)\!: E_1\to E_2$
is a bundle map, for each $p\in P$,
then
\[
\Theta : \, C^r(M,E_1)\times P\to C^r(M,E_2)\,,
\quad \Theta(\sigma,p)\, :=\, C^r(M,f_p)(\sigma)\, =\, f(\sbull,p)\circ \sigma
\]
is a $C^k_\K$-map.
\item[\rm (b)]
If $f \!: E_1\to E_2$ is a bundle map of class $C^{r+k}_\K$,
then
\[
C^r(M,f) : \, C^r(M,E_1)\to C^r(M,E_2),\quad
\sigma \mto f\circ \sigma
\]
is a mapping of class~$C^k_\K$.
\end{itemize}
\end{thm}
\begin{proof}
(b) directly follows from (a)
by taking $P:=H:=\{0\}$.
It thus suffices to prove~(a).
We let $(U_j)_{j\in J}$ be an open cover of
$M$ such that for every $j\in J$,
there are local trivializations
$\psi_j\!:\pi_1^{-1}(U_j)\to U_j \times F_1$
and $\phi_j\!: \pi_2^{-1}(U_j)\to U_j\times F_2$.
For each $j\in J$,
the mapping $h_j\!: U_j  \times F_1\times P\to F_2$,
$h_j(x,y,p):=(f_p)_{\phi_j,\psi_j}(x,y)$ is $C^{r+k}_\K$.
We have
\[
\beta_j(\Theta(\sigma))  \,= \,((f_p)_{\phi_j,\psi_j})_* (\alpha_j(\sigma))
\, =\, h_j(\sbull,p)_*(\alpha_j(\sigma))
\]
for all $p\in P$ and $\sigma\in C^r(M,E_1)$,
where $\alpha_j\!: C^r(M,E_1)\to C^r(U_j,F_1)$,
$\sigma\mto \sigma_{\psi_j}$
and $\beta_j\!: C^r(M,E_2)\to C^r(U_j,F_2)$,
$\sigma\mto \sigma_{\phi_j}$
are continuous linear maps.
By Proposition~\ref{globcruc}, the map
$C^r(U_j,F_1)\times P\to C^r(U_j,F_2)$,
$(\gamma,p)\mto
h(\sbull,p)_* (\gamma)$
is~$C^k_\K$.
In view of
Lemma~\ref{famgivessec} and Lemma~\ref{anyatlas},
Lemma~\ref{corestr}
shows that $\Theta$ is
of class~$C^k_\K$.
\end{proof}
The most interesting cases of Theorem~\ref{pushfbdl}
are (i) $r=k=\infty$;
\,(ii) $r\in \N_0\cup\{\infty\}$, $k=0$.\\[3mm]
To illustrate the results,
we show that spaces of sections
are topological modules
over the corresponding function algebras.
Recall the notion of (Whitney) sums
of vector bundles:
\begin{numba}\label{sumsvble}
If $\pi_j\!: E_j\to M$ are $C^r_\K$-vector bundles
with fibre $F_j$
for $j\in\{1,2\}$,
over the same base~$M$, then $E_1\oplus E_2:=
\bigcup_{x\in M} (E_1)_x\times (E_2)_x$
is a $C^r$-vector bundle over $M$ in a natural way,
with projection $\pi\!: E_1\oplus E_2\to M$, $v\mto x$
if $v\in (E_1)_x\times
(E_2)_x $. Let $\phi_j\!: \pi_j^{-1}(M_{\phi_j})\to M_{\phi_j}\times F_j$
be local trivializations of~$E_j$ for $j\in \{1,2\}$,
and
$M_{\phi_1\oplus \phi_2}:=M_{\phi_1}\cap M_{\phi_2}$.
Then
\[
\phi_1\oplus \phi_2\!: \pi^{-1}(M_{\phi_1\oplus\phi_2}) \to
M_{\phi_1\oplus\phi_2}\times (F_1\times F_2),
\]
$(\phi_1\oplus\phi_2)(v,w) :=(x, \pr_{F_1}(\phi_1(v)),
\pr_{F_2}(\phi_2(w)))$
for $(v,w)\in (E_1)_x\times (E_2)_x$,
is a local trivialization of~$E_1\oplus E_2$.
It is easy to see that the linear mapping
\[
C^r(M,E_1)\times C^r(M,E_2)\to C^r(M,E_1\oplus E_2),\;\;\;
(\sigma_1,\sigma_2)\mto (x\mto (\sigma_1(x),\sigma_2(x)))
\]
is an isomorphism of topological $\K$-vector spaces.
\end{numba}
As an immediate consequence of Theorem~\ref{pushfbdl},
we have:
\begin{cor}\label{aretopmod}
Let $\K$ be a topological field
and $\pi\!: E\to M$ be a $C^r_\K$-vector bundle,
whose fibre is a topological $\K$-vector
space~$F$.
Then
$C^r(M,E)$ is a topological $C^r(M,\K)$-module.
\end{cor}
\begin{proof}
The function space
$C^r(M,\K)$ can be identified with
the space $C^r(M,M \times \K)$
of $C^r_\K$-sections
of the trivial bundle
$\pr_1\!: M\times \K\to M$
with fibre~$\K$ (cf.\ Lemma~\ref{anyatlas}).
Thus $C^r(M,\K)\times C^r(M,E)\isom C^r(M, (M\times \K)\oplus E)$.
Using this identification,
the multiplication map $C^r(M,\K)\times C^r(M,E)\to
C^r(M,E)$ has the form
\[
C^r(M,\mu)\!: C^r(M, (M\times \K)\oplus E)\to C^r(M,E),
\]
where $\mu\!: (M\times \K)\oplus E \to E$
is the bundle map defined via $\mu((x,z),v):=zv\in E_x$ (scalar multiplication)
for all $x\in M$, $z\in \K$, and $v\in E_x$.
Given any local trivialization $\psi\!:\pi^{-1}(M_\psi)\to M_\psi\times
F$ of~$E$,
using the global trivialization $\phi:=\id\!:M\times \K\to
M\times \K$
we have $\mu_{\psi, \phi\oplus \psi}(x,z,v)=
zv\in F$, for all $(x,z,v)\in M_\psi\times\K\times F$,
showing that the map $\mu_{\psi, \phi \oplus \psi}\!:
M_\psi\times\K\times F\to F$ is of class~$C^r_\K$.
Thus $\mu$ is a $C^r_\K$-bundle map.
By Theorem~\ref{pushfbdl}\,(b) (applied with $k=0$),
$C^r(M,\mu)$ is continuous.
\end{proof}
If $\pi\!: E\to M$ is a $C^r_\K$-vector bundle
and~$U$ an open subset of~$M$, then
$\pi|_{E|_U}^U\!: E|_U\to U$
makes the open submanifold $E|_U:=\pi^{-1}(U)$ of~$E$ a
$C^r_\K$-vector bundle
over the base~$U$.
\begin{defn}\label{boredea}
Let $\K$ be a locally compact topological field,
$\pi\!: E\to M$ be a $C^r_\K$-vector bundle
over a finite-dimensional
$C^r_\K$-manifold~$M$, 
and $K\sub M$ be compact.
We equip $C^r_K(M,E)$ with the topology induced by
$C^r(M,E)$.
\end{defn}
\begin{la}\label{restrcts}
Let
$\pi\!: E\to M$ be a $C^r_\K$-vector bundle
over a finite-dimensional
$C^r_\K$-manifold~$M$,
and $U$ be an open subset of~$M$.
Then the following holds:
\begin{itemize}
\item[\rm (a)]
The restriction map
$C^r(M,E)\to C^r(U,E|_U)$,
$\sigma\mto \sigma|_U$ is continuous.
\item[\rm (b)]
Assume that $\K$ is locally compact, $M$ a finite-dimensional
$C^r_\K$-manifold and $K\sub U$ a compact subset.
Then the restriction map
$\rho_U\!: C^r_K(M,E)\to C^r_K(U,E|_U)$ is
an isomorphism of topological vector spaces.
\end{itemize}
\end{la}
\begin{proof}
(a) Since every local trivialization of~$E|_U$ also
is a local trivialization of~$E$,
Part\,(a) is apparent from the definition of the topologies.\vspace{2mm}

(b) As a consequence of (a), also $\rho_U$ is continuous.
Apparently, it is a linear bijection.
To see that $\rho_U$ is an isomorphism
of topological vector spaces,
we let $\cA_0$ be an atlas of local trivializations
for $E|_{M\,\take\, K}$,
and $\cA_1$ be an atlas for $E|_U$.
Then the topology
on $C^r_K(U,E|_U)$ is initial
with respect to the family of mappings
$\theta_\psi^U\!:C^r_K(U,E|_U)\to C^r(M_\psi,F)$,
$\sigma\mto \sigma_\psi$, for $\psi\in \cA_1$.
Furthermore, by Lemma~\ref{anyatlas},
the topology on $C^r_K(M,E)$ is
initial with respect to
family of mappings $\theta^M_\psi\!:
C^r_K(M,E)\to C^r(M_\psi,F)$
for $\psi\in \cA_0\cup\cA_1$
(defined analogously).
As $\theta^U_\psi\circ \rho_U=\theta^M_\psi$
for all $\psi\in \cA_1$ and $\theta^M_\psi=0$
for all $\psi\in \cA_0$,
the assertion easily follows.
\end{proof}
A variant of Lemma~\ref{amend} (and Lemma~\ref{famgivessec})
will be needed.
\begin{la}\label{amend2}
Let $\pi\!: E\to M$ be a $C^r_\K$-vector
bundle, and $(U_i)_{i\in I}$ be an open cover
of~$M$. For each $i\in I$,
let $\rho_i\!: C^r(M,E)\to C^r(U_i, E|_{U_i})$
be the restriction map. Then
\[
\rho:=(\rho_i)_{i\in I} : \; C^r(M,E)\to \prod_{i\in I} C^r(U_i, E|_{U_i})\,,
\quad
\rho(\sigma)\,:=\, (\sigma|_{U_i})_{i\in I}
\]
is a topological embedding, with closed image.
\end{la}
\begin{proof}
Each $\rho_i$ being continuous
by Lemma~\ref{restrcts}\,(a),
also $\rho$ is continuous.
To see that $\rho$ is an embedding,
consider the set
$\cA$ of all
local trivializations
$\psi\!: E|_{M_\psi}\to M_\psi\times F$ of~$E$ such
$M_\psi\sub U_i$ for some $i\in I$;
here $F$ is the typical fibre of~$E$.
Then $\cA$ is an atlas of local trivializations
for~$E$, and the topology on $C^r(M,E)$ is
initial with respect to the family $(\theta_\psi)_{\psi\in \cA}$
of the maps
$\theta_\psi\!: C^r(M,E)\to C^r(M_\psi,F)$, $\theta_\psi(\sigma):=
\sigma_\psi$ (Lemma~\ref{anyatlas}). If $M_\psi\sub U_i$,
then $\theta_\psi(\sigma)=\sigma_\psi=(\sigma|_{U_i})_\psi
=(\rho_i(\sigma))_\psi$
shows that $\theta_\psi$ is continuous with respect
to the topology induced by $\rho$ on~$C^r(M,E)$.
As a consequence, $\rho$ is a topological embedding.
To complete the proof, note that
\[
H\; :=\; \Big\{ (\sigma_i)_{i\in I}\in
\prod_{i\in I} C^r(U_i,E|_{U_i})\!: (\forall i,j\in I)\;
\sigma_i|_{U_i\cap U_j}=\sigma_j|_{U_i\cap U_j}\Big\}
\]
is closed in $\prod_{i\in I} C^r(U_i,E|_{U_i})$,
the restriction maps $C^r(U_i,E|_{U_i})\to
C^r(U_i\cap U_j,E|_{U_i\cap U_j})$
being continuous (Lemma~\ref{restrcts}\,(a)).
Clearly $\,\im(\rho)\sub H$.
If, conversely, $(\sigma_i)_{i\in I}\in H$,
then $\sigma(x):=\sigma_i(x)$ if $x\in U_i$
gives a well-defined section $\sigma\!: M\to E$.
Clearly $\sigma\in C^r(M,E)$ and $\rho(\sigma)=
(\sigma_i)_{i\in I}$. Thus $\,\im(\rho)=H$
is closed. 
\end{proof}
\begin{defn}\label{defcpsupp2}
Let $\K$ be a locally compact topological field,
$M$ be a
paracompact
$C^r_\K$-manifold
modeled on a finite-dimensional $\K$-vector space~$Z$,
and $\pi\!: E\to M$ be a $C^r_\K$-vector bundle
with fibre an arbitrary
(Hausdorff, not necessarily locally convex)
topological $\K$-vector space~$F$.
Then the set
\[
C^r_c(M,E):=\{\sigma\in C^r(M,E)\!: \;\mbox{$\Supp(\sigma)$
is compact}\,\}
\]
of compactly supported
$C^r_\K$-sections
is a $\K$-vector subspace
of $C^r(M,E)$, and
$C^r_c(M,E)=\bigcup_{K\in \cK(M)}\, C^{\,r}_K(M,E)$,
where $\cK(M)$ denotes the set of all
compact subsets of~$M$.
In the following,
we consider three vector topologies on
$C^r_c(M,E)$:
\begin{itemize}
\item[(a)]
We write
$C^r_c(M,E)_\tvs$
for $C^r_c(M,E)$, equipped
with the finest (a priori not necessarily
Hausdorff) vector topology
making the inclusion maps $\lambda_K\!:
C^r_K(M,E)\to C^r_c(M,E)$
continuous for each compact subset $K\sub M$.
Thus $C^r_c(M,E)_\tvs=\dl\, C^r_K(M,E)$\vspace{-.8mm}
in the category of
not necessarily Hausdorff topological $\K$-vector spaces
and continuous $\K$-linear maps.
\item[(b)]
If $F$ is locally convex, we
write
$C^r_c(M,E)_\lcx$
for $C^r_c(M,E)$, equipped
with the finest (a priori not necessarily
Hausdorff) locally convex vector topology
making the inclusion maps $\lambda_K\!:
C^r_K(M,E)\to C^r_c(M,E)$
continuous for each compact subset $K\sub M$.
Thus $C^r_c(M,E)_\lcx=\dl\, C^r_K(M,E)$\vspace{-.8mm}
in the category of
not necessarily Hausdorff, locally
convex topological $\K$-vector spaces
and continuous $\K$-linear maps.
\item[(c)]
As $M$ is paracompact and locally compact,
there exists a locally finite cover
$\cU=(U_i)_{i\in I}$
of $M$ by relatively compact, open subsets
$U_i\sub M$.
We define
\[
\rho_\cU\!: C^r_c(M,E)\to \bigoplus_{i\in I}\, C^r(U_i,E|_{U_i}),
\quad \rho_\cU(\sigma):=(\rho_i(\sigma))_{i\in I}=(\sigma|_{U_i})_{i\in I}\, ,
\]
where $\rho_i\!: C^r_c(M,E)\to C^r(U_i,E|_{U_i})$
is the restriction map for $i\in I$.
We write $C^r_c(M,E)_\bx$ for $C^r_c(M,E)$,
equipped with the topology $\cO_\cU$ induced by $\rho_\cU$,
where the direct sum is endowed with
the box topology.
\end{itemize}
\end{defn}
\begin{la}\label{boxindep2}
In the situation of Definition~{\rm \ref{defcpsupp2}\,(c)},
assume that both $\cU=(U_i)_{i\in I}$
and $\cV=(V_j)_{j\in J}$
are locally finite covers of~$M$
by relatively compact open sets.
Then $\cO_\cU=\cO_\cV$.
In other words,
the box topology on $C^r_c(M,E)$
is independent of the choice of $\cU$.
\end{la}
\begin{proof}
\!The topologies $\cO_\cU$ and $\cO_\cV$
are induced by
$\rho_\cU\!:
C^r_c(M,E)\to \bigoplus_{i\in I}\,
C^r(U_i,E|_{U_i})$,
$\rho_\cU(\sigma):=(\sigma|_{U_i})_{i\in I}$
and
$\rho_\cV\!:
C^r_c(M,E)\to \bigoplus_{j\in J}\,
C^r(V_j,E|_{V_j})$, $\rho_\cV(\sigma):=(\sigma|_{V_j})_{j\in J}$,
respectively.
Using Lemma~\ref{amend2}
instead of Lemma~\ref{amend},
we can repeat the proof
of Lemma~\ref{boxindep}
verbatim to get the desired result.
\end{proof}
\begin{prop}\label{comparetop2}
Let $M$ be a
paracompact, finite-dimensional $C^r_\K$-manifold
over a locally compact topological field~$\K$.
Let $\pi\!: E\to M$ be a $C^r_\K$-vector bundle
over~$M$, with typical fibre a topological $\K$-vector space~$F$.
Then the following holds:
\begin{itemize}
\item[\rm (a)]
The box topology on $C^r_c(M,E)_\bx$
is Hausdorff.
For every locally finite cover $\cU=(U_i)_{i\in I}$
of~$M$ by relatively compact, open subsets $U_i\sub M$,
the map
\[
\rho_\cU\!: C^r_c(M,E)_\bx \to \bigoplus_{i\in I}\, C^r(U_i,E|_{U_i}),\quad
\rho_\cU(\sigma):= (\sigma|_{U_i})_{i\in I}
\]
has closed image, and {\/\em $\rho_\cU|^{\text{im}\, \rho_\cU}$\/}
is an isomorphism of topological vector spaces.
The inclusion map $C^r_c(M,E)_\bx\to C^r(M,E)$
is continuous.
If $F$ is locally convex,
then $C^r_c(M,E)_\bx$ is locally convex.
\item[\rm (b)]
The inclusion map
$\lambda_K\!: C^r_K(M,E)\to C^r_c(M,E)_\bx$
is continuous and induces the given topology
on $C^r_K(M,E)$, for each compact subset
$K\sub M$.
\item[\rm (c)]
The map $\Phi\!: C^r_c(M,E)_\tvs\to C^r_c(M,E)_\bx$,
$\Phi(\gamma):=\gamma$
is continuous.
Thus $C^r_c(M,E)_\tvs$ is Hausdorff
and induces the given topology on each $C^r_K(M,E)$.
If $\K\not=\C$ and $M$ is $\sigma$-compact,
then $\Phi$ is an isomorphism of topological
$\K$-vector spaces.
\item[\rm (d)]
If $F$ is locally convex, then
$\Psi\!: C^r_c(M,E)_\lcx\to C^r_c(M,E)_\bx$,
$\Psi(\gamma):=\gamma$
is continuous.
Hence $C^r_c(M,E)_\lcx$ is Hausdorff
and induces the given topology on each $C^r_K(M,E)$.
If $\K\not=\C$ and $M$ is $\sigma$-compact,
then $\Psi$ is an isomorphism of topological
$\K$-vector spaces.
\item[\rm (e)]
If $\K$ is a local field
and $\cU=(U_i)_{i\in I}$
is a cover of~$M$ by mutually disjoint,
compact open sets $($cf.\ Lemma~{\rm \ref{onlyopen}\,(b)}$)$,
then
\[
\rho_\cU\!: C^r_c(M,E)_\bx\to \bigoplus_{i\in I}\, C^r(U_i,E|_{U_i}),\quad
\rho_\cU(\sigma):=(\sigma|_{U_i})_{i\in I}
\]
is an isomorphism of topological vector spaces
onto the direct sum, equipped with the box topology.
\item[\rm (f)]
If $\K$ is a local field
and $F$ is locally convex, then $\Psi$
is an isomorphism of topological $\K$-vector spaces,
i.e., $C^r_c(M,E)_\lcx=C^r_c(M,E)_\bx$.
\end{itemize}
In particular,
$C^r_c(M,E)_\bx=
C^r_c(M,E)_\tvs= C^r_c(M,E)_\lcx$
if $\K\not=\C$ and $M$
is $\sigma$-compact.
\end{prop}
\begin{proof}
(a) Using Lemma~\ref{amend2}
instead of Lemma~\ref{amend}
and Remark~\ref{seclcx}
instead of
Proposition~\ref{propprop}\,(b),
the proof of Proposition~\ref{comparetop}\,(a)
carries over.\vspace{2mm}

(b) Using Lemma~\ref{restrcts}\,(b)
and Lemma~\ref{amend2} instead of
Lemma~\ref{restrK} and Lemma~\ref{amend},
respectively,
the proof of Proposition~\ref{comparetop}\,(b)
can be repeated verbatim.\vspace{2mm}

(c) Note that, $C^r(U_i,E|_{U_i})$ being a topological
$C^r(U_i,\K)$-module (Corollary~\ref{aretopmod}),
the multiplication operator
$m_h\!: C^r(U_i,E|_{U_i})\to C^r(M,E|_{U_i})$, $m_h(\sigma)(x):=h(x)\sigma(x)$
is continuous, for each $i\in I$ and $h\in C^r(U_i,\K)$.
The assertion therefore follows
along the lines of
the proof of Proposition~\ref{comparetop}\,(c)
(with $\F:=\K$),
using Lemma~\ref{restrcts}\,(b)
instead of Lemma~\ref{restrK}.\vspace{2mm}

(d) We argue as in the proof of Proposition~\ref{comparetop}\,(d).\vspace{2mm}

(e) We argue as in the proof of Proposition~\ref{comparetop}\,(e),
taking $\F:=\K$.\vspace{2mm}

(f) Using Lemma~\ref{restrcts}\,(b) instead of
Lemma~\ref{restrK}, we can proceed
as in the proof of Proposition~\ref{comparetop}\,(e)
(taking $\F:=\K$).
\end{proof}
Throughout the following,
spaces of compactly supported sections
in vector bundles
will always be equipped
with the box topology.
We abbreviate $C^r_c(M,E):=C^r_c(M,E)_\bx$.
\begin{rem}
If $\K=\R$,
$M$ is $\sigma$-compact and the fibre $F$ is locally convex,
then the box topology on $C^r_c(M,E)$
coincides with the locally convex topology traditionally
considered on this space of compactly supported sections,
as a consequence of
Proposition~\ref{comparetop2}\,(d)
and
Proposition~\ref{propprop}\,(d).
\end{rem}
\begin{rem}\label{patchd3}
Let $M$ be a paracompact,
finite-dimensional $C^r_\K$-manifold over
a locally compact field~$\K$
(where $r\in \N_0\cup\{\infty\}$), and $\pi\!: E\to M$ be a
$C^r_\K$-vector bundle, with fibre an arbitrary
topological $\K$-vector space.
Let $(U_i)_{i\in I}$
be a locally finite cover
of~$M$ by relatively compact, open subsets
$U_i\sub M$
and $\rho_i\!: C^r_c(M,E)\to C^r(U_i,E|_{U_i})$,
$\rho_i(\sigma):=\sigma|_{U_i}$ be the restriction map
for $i\in I$. Then
\[
(C^r_c(M,E),(\rho_i)_{i\in I})
\]
is a patched topological vector space, by
Proposition~\ref{comparetop2}\,(a).\vspace{2mm}
\end{rem}
\begin{center}
{\bf The {\boldmath $\Omega$\/}-Lemma with Parameters}
\end{center}
In the following, we prove generalizations
of the so-called ``$\Omega$-Lemma'' (see \cite[Thm.\,8.7]{Mic}),
formulated in \cite{Mic}
for mappings between subsets of spaces of compactly supported
smooth sections in finite-dimensional real vector bundles.
An essential ingredient of the proof will be a
version of Proposition~\ref{diffpatch} for functions depending on
parameters.
\begin{la}\label{laomop}
Let $\K$ be a locally compact topological field
and
$\pi\!: E\to M$ be a $C^r_\K$-vector bundle
over a paracompact, finite-dimensional
$C^r_\K$-manifold~$M$.
Let $\Omega\sub E$ be an open subset.
Then
\[
C^r_c(M,\Omega)\;  :=\; \{\sigma\in C^r_c(M,E)\!:
\sigma(M)\sub \Omega\}
\]
is an open $($possibly empty$)$ subset of $C^r_c(M,E)$.
\end{la}
\begin{proof}
The present proof is not the shortest one;
it is stated as follows
because in exactly this
form it can be re-used to prove Theorem~\ref{OmegaP}.\\[3mm]
Let $\sigma\in C^r_c(M,\Omega)$.
Using the paracompactness and local compactness
of~$M$,
we find locally finite covers~$(U_i)_{i\in I}$
and $(M_i)_{i\in I}$ of~$M$ by relatively
compact, open sets such that
$K_i:=\wb{U_i}\sub M_i$
and $E|_{M_i}$ is a trivial vector bundle,
for each $i\in I$.\,\footnote{Every $x\in M$ has a relatively
compact open neighbourhood $Q_x$ such that
$E|_{Q_x}$ is trivial; let $P_x$ be an
open neighbourhood of~$x$ such that $\wb{P_x}\sub Q_x$.
Since $M$ is paracompact, there exists a
locally finite open cover $(U_i)_{i\in I}$
subordinate to $(P_x)_{x\in M}$.
Then $U_i\sub P_{x_i}$
for some $x_i\in M$, whence $U_i$ is relatively compact.
By Lemma~\ref{lathicken},
there exists a locally finite cover
$(\wt{U}_i)_{i\in I}$
of~$M$ by relatively compact open sets
such that $\wb{U_i}\sub \wt{U}_i$
for each~$i$. We define $M_i:=\wt{U}_i\cap Q_{x_i}$;
then $(U_i)_{i\in I}$ and $(M_i)_{i\in I}$
have the desired properties.}
Let $\psi_i\!: E|_{M_i}\to M_i\times F$
be a trivialization of $E|_{M_i}$, where
$F$ is the typical fibre of~$E$.
Then
\[
\kappa : \; C^r_c(M,E)\to \bigoplus_{i\in I}\,
C^r(M_i,F)\,,\quad
\kappa (\tau):=(\tau_{\psi_i})_{i\in I}=\big((\tau|_{M_i})_{\psi_i}
\big)_{i\in I}
\]
is a topological embedding,
by definition of the box topology on
$C^r_c(M,E)$ and Lemma~\ref{anyatlas}.
For each $i\in I$,
the set $\Omega_i:=\psi_i(\Omega\cap E|_{M_i})$
is an open neighbourhood of the compact subset
$\{(x,\sigma_{\psi_i}(x))\!: x\in K_i\}$
of $M_i\times F$.
Since $\sigma_{\psi_i}$ is continuous,
a standard compactness argument provides
finite families $(U_{i,j})_{j\in J_i}$
and $(M_{i,j})_{j\in J_i}$
of relatively compact, open subsets $M_{i,j}\sub M_i$
and relatively compact, open subsets $U_{i,j}\sub M_{i,j}$
such that $K_i\sub \bigcup_{j\in J_i}U_{i,j}$,
and open $0$-neighbourhoods $W_{i,j}\sub F$
such that
\begin{equation}\label{cbck}
M_{i,j}
\times (\sigma_{\psi_i}(M_{i,j}) +W_{i,j})\sub \Omega_i
\qquad \mbox{for each $i\in I$ and $j\in J_i\,$.}
\end{equation}
Set $L:=\{(i,j)\!: i\in I,\, j\in J_i\}$.
After replacing $I$ by $L$, $(U_i)_{i\in I}$ by
$(U_{i,j})_{(i,j)\in L}$, $(M_i)_{i\in I}$
by $(M_{i,j})_{(i,j)\in L}$,
and $(\psi_i)_{i\in I}$ by $(\psi_i|_{E|_{M_{i,j}}})_{(i,j)\in L}$,
instead of (\ref{cbck})
we may assume without loss of generality
that there exist open $0$-neighbourhoods $W_i\sub F$
such that
\begin{equation}\label{cbck2}
M_i\times (\sigma_{\psi_i}(M_i)+W_i)\sub \Omega_i
\end{equation}
for each $i\in I$. Then $V_i:=\sigma_{\psi_i}(K_i)+W_i$
is an open neighbourhood of $\sigma_{\psi_i}(K_i)$
in~$F$, for each $i\in I$, and it is a $0$-neighbourhood
in~$F$ for all but finitely many~$i$.
Then $\lfloor K_i, V_i\, \rfloor_r:=\{\gamma\in
C^r(M_i,F)\!: \gamma(K_i)\sub V_i\}$
is an open subset of $C^r(M_i,F)$ and
\[
V\; :=\; \kappa^{-1}\left( \, \bigoplus_{i\in I}\;\lfloor K_i,V_i\rfloor_r
\right)
\]
is an open neighbourhood of $\sigma$ in $C^r_c(M,E)$.
Since $(K_i)_{i\in I}$ is a cover of $M$,
we deduce from (\ref{cbck2}) that $V\sub C^r_c(M,\Omega)$.
Thus $C^r_c(M,\Omega)$ is a neighbourhood of~$\sigma$.
\end{proof}
\begin{thm}[{\boldmath $\Omega$\/}-Lemma with Parameters]\label{OmegaP}
Let $\K$ be a locally compact topological field
and $r,k\in \N_0\cup\{\infty\}$.
Let $\pi_1\!: E_1\to M$ and $\pi_2\!: E_2\to M$
be $C^{r+k}_\K$-vector bundles
over the same paracompact,
finite-dimensional $C^{r+k}_\K$-manifold~$M$,
whose fibres $F_1$, resp., $F_2$
are arbitrary topological $\K$-vector spaces.
Let $P\sub Z$ be an open subset
of a finite-dimensional $\K$-vector space~$Z$,
$\Omega\sub E_1$
be an open subset, and
\[
f\!: \Omega\times P \to E_2
\]
be a mapping of class~$C^{r+k}_\K$
such that $f_p:=f(\sbull,p)\!: \Omega\to E_2$
is a bundle map, for each $p\in P$.
We assume that there exists a compact subset
$K\sub M$ such that $0_x\in \Omega$
for each $x\in M\setminus K$,
and $f(0_x,p)=0_x$ for each
$x\in M\setminus K$ and $p\in P$ $($using the same
symbol $0_\sbull$
for the $0$-section of $E_1$ and $E_2$, resp.$)$
Then the mapping
\[
\phi\!:
C^r_c(M,\Omega)\times P \to C^r_c(M,E_2)\,,\quad
\phi(\sigma,p ):=f_p\circ \sigma
\]
is $C^k_\K$. In particular, if $k=r=\infty$,
then $\phi\!:
C^\infty_c(M,\Omega)\times P \to C^\infty_c(M,E_2)$
is smooth.
\end{thm}
\begin{proof}
Note first that for $\sigma\in C^r_c(M,\Omega)$ and $p\in P$,
we have
\[
\Supp(\phi(\sigma,p))\; \sub \; K \, \cup\, \Supp(\sigma)\,.
\]
Hence $\phi(\sigma,p)$ is indeed compactly supported,
and thus $\, \im(\phi)\sub C^r_c(M,E_2)$.
To see that $\phi$ is of class $C^k_\K$, we now
fix $\sigma\in C^r_c(M,\Omega)$.
We let $(U_i)_{i\in I}$, $(M_i)_{i\in I}$,
$(\psi_i)_{i\in I}$,
$(K_i)_{i\in I}$, 
$(V_i)_{i\in I}$,
$\kappa : \; C^r_c(M,E_1)\to \bigoplus_{i\in I}\,
C^r(M_i,F_1)$, $\kappa(\tau):=(\tau_{\psi_i})_{i\in I}$
and $V:=
\kappa^{-1}\left( \, \bigoplus_{i\in I}\;\lfloor K_i,V_i\rfloor_r
\right)$ be as in the proof of Lemma~\ref{laomop}
(applied with $E:=E_1$
and $F:=F_1$).
Because we can choose each $Q_x$ so small that
also $E_2|_{Q_x}$
is trivial in the proof of Lemma~\ref{laomop},
we may assume without loss of generality
that also $E_2|_{M_i}$ is trivial, for each $i\in I$.
Let $\phi_i\!: E_2|_{M_i}\to M_i\times F_2$
be a trivialization.
Abbreviate $\omega_i:=\phi_i|_{E_2|_{U_i}}^{U_i\times F_2}\!:
E_2|_{U_i}\to U_i\times F_2$.
Then
\[
\kappa_2 : \; C^r_c(M,E_2)\to \bigoplus_{i\in I}\,
C^r(U_i,F_2)\,,\quad
\kappa_2(\tau):=(\tau_{\omega_i})_{i\in I}
=\big((\tau|_{U_i})_{\omega_i}\big)_{i\in I}
\]
is a topological embedding onto a closed
vector subspace, as a consequence of
Proposition~\ref{comparetop2}\,(a) and Lemma~\ref{anyatlas}.
For each $i\in I$,
the map
\[
f_i : \, M_i\times V_i\times P\to F_2\,,\quad
f_i(x,y,p)\, :=  \, \pr_2(\phi_i(f(\psi_i^{-1}(x,y),p)))
\]
is of class $C^{r+k}_\K$, where $\pr_2\!: M_i\times F_2\to
F_2$ is the second coordinate projection.
Hence, by Proposition~\ref{crucial}\,(a),
the map
\[
g_i:\,
\lfloor K_i, V_i\rfloor_r
\times P
\to C^r(U_i,F_2)\,,\quad
g_i(\gamma,p)(x)\, := \, f_i(x,\gamma(x),p)
\]
is of class $C^k_\K$,
where
$\lfloor K_i, V_i\rfloor_r :=\{\gamma\in
C^r (M_i, F_1)\!: \gamma(K_i)\sub V_i\}$.
Furthermore, $0\in \lfloor K_i, V_i\rfloor_r$
and $g_i(0,p)=0$ whenever $M_i\cap (K\cup \Supp(\sigma))=\emptyset$,
which is the case for all but finitely many
$i\in I$. Hence, by Proposition~\ref{sumspara},
the map
\[
g:\, \left( \bigoplus_{i\in I}
\, \lfloor K_i, V_i\rfloor_r \right)
\times P
\,\to \, \bigoplus_{i\in I}\,
C^r(U_i,F_2)\,,\quad
g\Big( \sum_{i\in I} \gamma_i,\, p\Big)\,:=\,
\sum_{i\in I} g_i(\gamma_i,p)
\]
is of class $C^k_\K$.
Since the diagram
\[
\begin{array}{rccc}
C^r_c(M,\Omega)\times P\;\, \supseteq &
V\times P & \stackrel{\kappa|_V\times \text{id}_P}{\longrightarrow}&
\left(\bigoplus_{i\in I}\, \lfloor K_i,V_i\rfloor_r\right)\times P
\vspace{1mm}\\
&\phi|_{V\times P}  \downarrow\;\;\;\;\;\;\;& & \downarrow g\vspace{1mm}\\
&C^r_c(M,E_2) & \stackrel{\kappa_2}{\longrightarrow} &
\bigoplus_{i\in I}\, C^r(U_i,F_2)
\end{array}
\]
commutes and $\kappa_2$ is an embedding
of topological vector spaces with closed image,
we deduce with Lemma~\ref{corestr}
that $\phi|_{V\times P}$ is $C^k_\K$
on the open neighbourhood $V\times P$ of $\{\sigma\}\times P$
in $C^r_c(M,\Omega)\times P$.
As $\sigma$ was arbitrary, $\phi$ is $C^k_\K$.
\end{proof}
Specializing to
a singleton set of parameters,
we obtain:
\begin{cor}[{\boldmath $\Omega$\/}-Lemma]\label{Omeg}
Let $\K$ be a locally compact topological field
and $r,k\in \N_0\cup\{\infty\}$.
Let $\pi_1\!: E_1\to M$ and $\pi_2\!: E_2\to M$
be $C^{r+k}_\K$-vector bundles
over the same paracompact,
finite-dimensional $C^{r+k}_\K$-manifold~$M$,
whose fibres $F_1$, resp., $F_2$
are arbitrary topological $\K$-vector spaces.
Let $\Omega\sub E_1$
be an open neighbourhood
of the image of a section $\sigma_0\in C^r_c(M,E_1)$,
and $f\! : \Omega\to E_2$
be a bundle map of class~$C^{r+k}_\K$,
such that $f\circ\sigma_0$ has compact support. 
Then the map
\[
\phi\!: C^r_c(M,\Omega)\to C^r_c(M,E_2)\,,\quad
\phi(\sigma):=f\circ \sigma
\]
is of class $C^k_\K$.\Punkt
\end{cor}
\begin{rem}\label{remOmeg}
The following special cases of Corollary~\ref{Omeg}
are of particular interest:
\begin{itemize}
\item[(a)]
If $k=r=\infty$,
then $C^\infty_c(M,f)\!:
C^\infty_c(M,\Omega)\to C^\infty_c(M,E_2)$, $\sigma\mto f\circ \sigma$
is smooth.
\item[(b)]
If $r\in \N_0$ and $k=0$,
then $C^r_c(M,f)\!: C^r_c(M,\Omega)\to C^r_c(M,E_2)$,
$\sigma\mto f\circ \sigma$
is continuous.
\end{itemize}
\end{rem}
To illustrate the results,
let us prove that spaces of compactly supported
sections
are topological modules
over the corresponding test function algebras.
First, we observe:
\begin{numba}\label{sumsvble2}
If $\K$ is locally compact and
the manifold $M$ is finite-dimensional
and paracompact
in the situation of~{\bf \ref{sumsvble}},
then apparently also the linear mapping
\[
C^r_c(M,E_1)\times C^r_c(M,E_2)\to C^r_c(M,E_1\oplus E_2),\;\;\;
(\sigma_1,\sigma_2)\mto (x\mto (\sigma_1(x),\sigma_2(x)))
\]
is an isomorphism of topological $\K$-vector spaces.
\end{numba}
As an immediate consequence of Theorem~\ref{pushfbdl}
and Corollary~\ref{Omeg}, we now obtain:
\begin{cor}\label{aretopmod2}
Let $\K$ be a locally compact topological field,
$M$ be a paracompact, finite-dimensional $C^r_\K$-manifold,
and $\pi\!: E\to M$ be a $C^r_\K$-vector bundle,
whose fibre is a topological $\K$-vector
space~$F$.
Then $C^r_c(M,E)$ is a topological $C^r_c(M,\K)$-module.
\end{cor}
\begin{proof}
The proof of Corollary~\ref{aretopmod}
carries over verbatim,
using Corollary~\ref{Omeg}
and {\bf \ref{sumsvble2}}
instead of Theorem~\ref{pushfbdl}\,(b)
and {\bf \ref{sumsvble}}.
\end{proof}
\begin{rem}\label{bundmodu}
If $A$ is an associative topological $\K$-algebra
and $M$ a $C^r_\K$-manifold,
we define a {\em bundle of topological
$A$-modules\/}
as a $C^r_\K$-vector bundle
$\pi\!: E\to M$
whose typical fibre~$F$
is a topological $A$-module,
and equipped with an atlas $\cA$ of local trivializations
such that $\im(g_{\phi,\psi})$
consists of topological $A$-module
automorphisms of~$F$, for all $\phi,\psi\in\cA$.
In this case,
we see as in the proof of Corollary~\ref{aretopmod}
that $C^r(M,E)$ is a topological $C^r(M,A)$-module
(under pointwise
operations). If $\K$ is locally compact and
$M$ is finite-dimensional
and paracompact, then $C^r_c(M,E)$ is
a topological $C^r_c(M,A)$-module (cf.\ proof
of Corollary~\ref{aretopmod2}).\vspace{2mm}
\end{rem}
\begin{center}
{\bf Almost local mappings between spaces of compactly
supported sections}
\end{center}
Our considerations from Section~\ref{secalmloc}
carry over directly to the case of mappings
between spaces of compactly supported
sections in vector bundles (equipped with the box
topology).
Compare the earlier works \cite{SEC} and \cite{DIF}
for a discussion of the real locally convex case,
based on the locally convex direct limit topology
on spaces of compactly supported sections.
\begin{defn}\label{defnal2}
Let
$\K$ be the field of real numbers
or a local field.
Given
$r,s,k \in \N_0\cup\{\infty\}$,
let $\pi_1\!: E_1\to M$ be
a $C^r_\K$-vector bundle
over a paracompact, finite-dimensional
$C^r_\K$-manifold~$M$,
with fibre an arbitrary topological $\K$-vector space~$F_1$.
Let $\pi_2\!: E_2\to N$ be
a $C^s_\K$-vector bundle
over a paracompact, finite-dimensional
$C^s_\K$-manifold~$N$,
with fibre an arbitrary topological $\K$-vector space~$F_2$.
Finally, let
$f\!: P \to C^s_c(N,E_2)$ be a mapping,
defined on an open subset $P\sub C^r_c(M,E_1)$.
\begin{itemize}
\item[\rm (a)]
$f$
is called {\em almost local\/}
if there exist locally finite covers
$(U_i)_{i\in I}$
of $M$ and $(V_i)_{i\in I}$
of~$N$
by relatively
compact, open sets
such that, for all $i\in I$
and $\sigma,\tau\in P$
with $\sigma|_{U_i}=\tau|_{U_i}$,
we have $f(\sigma)|_{V_i}=f(\tau)|_{V_i}$.
\item[\rm (b)]
$f$ is called {\em locally almost
local\/} if every $\sigma \in P$ has an open neighbourhood
$Q\sub P$ such that $f|_Q$ is almost local.
\item[\rm (c)]
In the special case where $M=N$,
we call $f\!: P\to C^s_c(M,E_2)$
a {\em local\/} mapping if, for all $x\in M$
and $\sigma \in P$,
the element $f(\sigma)(x)$ only depends on the
germ of~$\sigma$ at~$x$.\footnote{More precisely,
we require
$f(\sigma)(x)=f(\tau)(x)$
for all $x\in M$ and $\sigma,\tau\in P$ with
the same germ at~$x$.}
\end{itemize}
\end{defn}
It is easy to see that every local mapping
is almost local.
\begin{thm}[Smoothness Theorem]\label{smoothy2}
Let $f\!:
C^r_c(M,E_1)\supseteq P\to C^s_c(N,E_2)$
be a map as described in Definition~{\rm \ref{defnal2}}.
If
$f_K:=f|_{P\cap C^r_K(M,E_1)}$
is of class $C^k_\K$
for every compact subset $K\sub M$
and~$f$ is locally almost local,
then~$f$ is of class $C^k_\K$.
\end{thm}
\begin{proof}
Given $\sigma\in P$,
there exists an open neighbourhood~$Q$
of $\sigma$ in $P$ such that $f|_Q$
is almost local.
As $\sigma$ was arbitrary,
the assertion
will follow if we can show that $f|_W$
is of class $C^k_\K$
for some open neighbourhood $W$ of~$\sigma$
in~$Q$.
To this end, it suffices to show
that the mapping
$g\!: Q-\sigma\to C^s_c(N, E_2)$,
$g(\tau):=f(\sigma+\tau)-f(\sigma)$
is of class $C^k_\K$
on some open
zero-neighbourhood. As $f|_Q$ is
almost local, we find
locally finite covers $(U_i)_{i\in I}$
of $M$
and $(V_i)_{i\in I}$ of $N$,
with each $U_i$ and $V_i$ relatively compact
and open,
such that $f(\tau)|_{V_i}$ only depends
on $\tau|_{U_i}$, for all $\tau\in Q$.
Then apparently also
$g(\tau)|_{V_i}=g(\kappa)|_{V_i}$
for all $\tau,\kappa\in Q-\sigma$ such that
$\tau|_{U_i}=\kappa|_{U_i}$,
showing that also~$g$ is almost local.
Furthermore, given a compact subset
$K\sub M$, the map
$g|_{(Q-\sigma)\cap C^r_K(M,E_1)}$
is of class $C^k_\K$,
since so is the restriction of~$f$
to $Q\cap C^r_{K\cup\sSup(\sigma)}(M,E_1)$.
We abbreviate $R:= Q-\sigma$.\\[3mm]
Next, we pick a locally finite open cover $(\wt{U}_i)_{i\in I}$
of~$M$
such that
$\wb{U_i}\sub\wt{U}_i$
holds for the compact
closures, for all $i\in I$;
such a ``thickening''
exists by Lemma~\ref{lathicken}. 
For each $i\in I$,
we pick
a mapping $h_i\in C^r(\wt{U}_i,\K)$,
with compact support $K_i:=\Supp(h_i)$,
which is constantly~$1$ on $U_i$
(see Lemma~\ref{cutthroat} if $\K$ is a local field;
the real case is standard).\\[3mm]
By Remark~\ref{patchd3}, the family
$(\rho_i)_{i\in I}$ of restriction maps
$\rho_i\!:
C^r_c(M,E_1)\to C^r(\wt{U}_i, E_1|_{\wt{U}_i})$
is a patchwork for $C^r_c(M,E_1)$.
We let $\rho\!: C^r_c(M, E_1)\to \bigoplus_{i\in I}
C^r(\wt{U}_i,E_1|_{\wt{U}_i})=:S$
be the corresponding embedding
taking $\tau$ to $\sum_{i\in I} \rho_i(\tau)$.
Similarly, the family $(\xi_i)_{i\in I}$
of restriction maps
$\xi_i\!: C^s_c(N,E_2)\to C^s(V_i,E_2|_{V_i})$
is a patchwork for $C^s_c(N,E_2)$.\\[3mm]
The mapping $\rho$ being a topological embedding,
we find an open $0$-neighbourhood $H\sub S$
such that $\rho^{-1}(H)\sub R$.
The direct sum being
equipped with the box topology,
after shrinking~$H$ we may assume that $H=\bigoplus_{i\in I}
A_i$ for a family
$(A_i)_{i\in I}$ of open $0$-neighbourhoods
$A_i\sub C^r(\wt{U}_i, E_1|_{\wt{U}_i})$.
As a consequence of Corollary~\ref{aretopmod},
the multiplication operator $\mu_{h_i}\!:
C^r(\wt{U}_i, E_1|_{\wt{U}_i})
\to C^r_{K_i}(\wt{U}_i, E_1|_{\wt{U}_i})$,
$\tau\mto h_i\cdot \tau$
is continuous linear.
Hence, we find an open zero-neighbourhood
$W_i\sub A_i$ such that $h_i\cdot W_i\sub R$,
where we identify
$C^{\, r}_{K_i}(\wt{U}_i, E_1|_{\wt{U}_i})$
with $C^{\, r}_{K_i}(M, E_1)$
as a topological $\K$-vector space
in the natural way, extending sections by~$0$
(cf.\ Lemma~\ref{restrcts}\,(b)).
Then $W:=\rho^{-1}(\bigoplus_{i\in I} W_i)\sub R$
is an open zero-neighbourhood in
$C^r_c(M,E_1)$ such that $\rho_i(W)\sub W_i$
for each $i\in I$.
We define
\[
g_i\!: W_i \to C^s(V_i,E_2|_{V_i}),
\;\;\; g_i:= \xi_i \circ
g|_{R\cap C^r_{K_i}(M,E_1)}
\circ \mu_{h_i}|_{W_i}^R\, .
\]
Then $g_i$ is of class $C^k_\K$,
being a composition of $C^k_\K$-maps.
Note that $\xi_i(g(\tau))
=g(\tau)|_{V_i}=g(h_i \cdot \tau)|_{V_i}
=g_i(\tau|_{\wt{U}_i})$ for each $\tau \in W$ and $i\in I$.
Thus $(g_i)_{i\in I}$ is compatible with $g|_W$
in the sense of Definition~\ref{defnpama}.
We have shown that $g|_W$
is a patched mapping which is of class $C^k_\K$ on the patches.
By Proposition~\ref{diffpatch},
$g|_W$ is of class $C^k_\K$.
\end{proof}
$\;$\vfill\pagebreak

$\;\;\;\;$\\
$\mbox{$\;\;\;$}$\\
\vfill\pagebreak
$\;$\\
\mbox{$\;\;$}\\
\thispagestyle{empty}
$\;$\hfill\pagebreak
\twocolumn
\noindent
{\Large \bf Index}\vspace{3.1 mm}\\
\thispagestyle{plain}
\noindent\label{nowindex}
absolute value\dotfill\pageref{convents}\\
absorbing set\dotfill\pageref{secdirsum}\\
algebra\\
--- continuous inverse algebra\dotfill\pageref{defnci}\\
--- locally finite algebra\dotfill\pageref{refback}\\
almost local map\dotfill\pageref{defnal}, \pageref{defnal2}\\
analytic maps\dotfill\pageref{findimgp}\\
atlas\dotfill\pageref{verybas}\\
Baker-Campbell-Hausdorff Lie group\dotfill\pageref{remmentBCH}\\
balanced set\dotfill\pageref{secdirsum}\\
ball\dotfill\pageref{defnuing}\\
--- metric\dotfill\pageref{defnuing}\\
Boman's theorem\dotfill\pageref{appboman}\\
box neighbourhood\dotfill\pageref{defnbxxx}\\
box topology\\
--- (direct limit properties)\dotfill\pageref{boxisuniv}\\
--- locally convex\dotfill\pageref{lcx2}\\
--- on direct sum\dotfill\pageref{defnbxxx}\\
--- on spaces of sections\dotfill\pageref{defcpsupp2}\\
--- on spaces of test functions\dotfill\pageref{defbxtop}\\
bundle map\dotfill\pageref{defnparvect}\\
bundle of topological modules\dotfill\pageref{bundmodu}\\
canonical vector topology\dotfill\pageref{canonic}\\
chain rule\dotfill\pageref{chainr}\\
chart\dotfill\pageref{verybas}\\
compatible with patched map\dotfill\pageref{defnpama}\\
complete valued field\dotfill\pageref{convents}\\
composition map\dotfill\pageref{compcomp}\\
continuous inverse algebra (CIA)\\
--- (definition)\dotfill\pageref{defnci}\\
--- locally convex, over $\R$ or $\C$\dotfill\pageref{remmentBCH}\\
continuously differentiable maps\\
--- (definition)\dotfill\pageref{defnCr}\\
--- in Michal-Bastiani sense\dotfill\pageref{MichalB}\\
--- into cartesian products\dotfill\pageref{inprod}\\
--- into closed vector subspaces\dotfill\pageref{corestr}\\
--- into projective limits\dotfill\pageref{inpl}\\
--- into real locally convex spaces\dotfill\pageref{MichalB}\\
--- into complex locally convex spaces\dotfill\pageref{MichalB}\\
--- Keller's $C^k_c$-maps\dotfill\pageref{MichalB}\\
--- (local property)\dotfill\pageref{Crlocal}\pagebreak

$\;$\\[0.33cm]
convenient differential calculus\dotfill\pageref{bscsconv}\\
conveniently $\K$-smooth\dotfill\pageref{bscsconv}\\
diffeomorphism group\\
--- of a ball\dotfill\pageref{diffball}\\
--- of a finite-dimensional manifold\dotfill\pageref{thmdiffpara},
\pageref{defnDM}\\
difference quotient map $f^{]1[}$\dotfill\pageref{gargel}\\
differential\dotfill\pageref{differentials}\\
--- (higher)\dotfill\pageref{differentials}\\
direct sums (of topological vector spaces)\\
--- (definition)\dotfill\pageref{basicdsum}\\
--- differentiable maps between\dotfill\pageref{mapsdirsums},
\pageref{sumspara}\\
directional derivative\dotfill\pageref{differentials}\\
evaluation map\dotfill\pageref{evalCk}, \pageref{seclook}\\
exponential law\dotfill\pageref{nownow}, \pageref{expmetriz}\\
extended difference quotient map $f^{[k]}$\\
--- (definition)\dotfill\pageref{defnCr}\\
--- (symmetry properties)\dotfill\pageref{firstsym}\\
extension of scalars\dotfill\pageref{extendscale}\\
finite topology\dotfill\pageref{defnftop}\\
function spaces\\
--- $C_K(X,E)$\dotfill\pageref{defncsp}\\
--- $C^r(U,E)$\dotfill\pageref{top1}\\
--- $C^r(M,E)$\dotfill\pageref{top2}\\
--- $C^r_K(M,E)$\dotfill\pageref{defcrsk}\\
--- $C^r_c(M,E)$\dotfill\pageref{defcpsupp}, \pageref{anoconv}\\
--- locally convex\dotfill\pageref{propprop}, \pageref{lcx1},
\pageref{comparetop}\\
--- properties\dotfill\pageref{propprop}\\
Grothendieck's theorem\dotfill\pageref{thmgroth}\\
group of units\dotfill\pageref{defnci}\\
--- as a Lie group\dotfill\pageref{ciasmooth}\\
$k$-space\dotfill\pageref{defnkspa}\\
Keller's $C^k_c$-maps\dotfill\pageref{MichalB}\\
Lie algebra of a Lie group\dotfill\pageref{verybas}\\
Lie group over~$\K$, alias $\K$-Lie group\dotfill\pageref{verybas}\\
local description of Lie group structures\dotfill\pageref{locchar}\\
local field\dotfill\pageref{convents}\\
locally almost local map\dotfill\pageref{defnal}, \pageref{defnal2}\\
locally convex topological vector space\\
--- (definition)\dotfill\pageref{deflcx}\\
--- (function spaces)\dotfill\pageref{propprop},
\pageref{lcx1}, \pageref{comparetop}\\
--- (direct sums)\dotfill\pageref{lcx2}\\
locally finite algebra\dotfill\pageref{refback}\\
local map\dotfill\pageref{defnal}, \pageref{defnal2}\\
local trivialization\dotfill\pageref{defnbdle}\\
Mackey-Cauchy sequence\dotfill\pageref{defnMC}\\
Mackey complete\dotfill\pageref{defnMC}\\
manifold, $C^r_\K$-manifold\\
--- modeled on topological vector space\dotfill\pageref{verybas}\\
--- paracompact finite-dimensional\dotfill\pageref{paracamen}\\
--- " ", over local fields\dotfill\pageref{onlyopen}\\
--- $\sigma$-compact finite-dimensional\dotfill\pageref{onlyopen}\\
mapping algebras\dotfill\pageref{mapcia}\\
mapping group\dotfill\pageref{defnmgp}\\
metric ball\dotfill\pageref{defnuing}\\
multiplication operator\dotfill\pageref{multop}\\
norm\dotfill\pageref{dfblls}\\
--- operator norm\dotfill\pageref{dfblls}\\
--- maximum norm\dotfill\pageref{dfblls}\\
normalization\dotfill\pageref{dfblls}\\
$\Omega$-Lemma\dotfill\pageref{Omeg}\\
$\Omega$-Lemma with parameters\dotfill\pageref{OmegaP}\\
partition\\
--- into balls\dotfill\pageref{onlyopen}\\
--- into metric balls\dotfill\pageref{disjballs}\\
--- of unity\dotfill\pageref{defnparun}\\
patched mapping

\thispagestyle{plain}\noindent{}--- $C^k_\K$ on the patches\dotfill\pageref{defnpm2}\\
--- (definition)\dotfill\pageref{defnpama}\\
--- differentiability properties\dotfill\pageref{diffpatch}\\
patched topological vector space\\
--- (definition)\dotfill\pageref{defnpatched}\\
--- spaces of sections are such\dotfill\pageref{patchd3}\\
--- test function spaces are such\dotfill\pageref{patchd}\\
patches\dotfill\pageref{pagerpa}\\
patchwork\dotfill\pageref{pagerpa}\\
polynormed topological vector space\dotfill\pageref{findimgp}\\
pullback\\
--- of $C^r$-maps\dotfill\pageref{pullback}, \pageref{pb2}\\
--- of test functions\dotfill\pageref{pbproper}\\
pushforward\\
--- of continuous maps\dotfill\pageref{backbone}, \pageref{pushforw}\\
--- of $C^r$-maps\dotfill\pageref{linearcase}, \pageref{globcruc},
\pageref{pushforw2}, \pageref{crucial}\\ 
restriction map\dotfill\pageref{restrK},
\pageref{restrcts}\vfill\pagebreak

\noindent
sequentially complete\dotfill\pageref{defnseqc}\\
smooth map\dotfill\pageref{defnCr}\\
--- continuous linear maps\dotfill\pageref{linbilin}\\
--- continuous bilinear maps\dotfill\pageref{linbilin}\\
smoothness theorem\dotfill\pageref{smoothy}, \pageref{smoothy2}\\
spaces of compactly supported sections\\
--- comparison of topologies\dotfill\pageref{comparetop2}\\
--- (definition of topology)\dotfill\pageref{defcpsupp2}\\
symmetry properties of $f^{[k]}$\dotfill\pageref{firstsym}\\
tensor product\dotfill\pageref{tensorprod}\\
thickening locally finite covers\dotfill\pageref{lathicken}\\
test function algebras\dotfill\pageref{disctefgp}\\
test function groups\dotfill\pageref{disctefgp}\\
test function spaces\\
--- comparison of topologies\dotfill\pageref{comparetop},
\pageref{traditio}\\
--- (definition)\dotfill\pageref{defcpsupp}\\
--- differentiable maps between these\dotfill\pageref{smoothy}\\
--- direct limit properties\dotfill\pageref{boxdirl1}\\
--- local convexity\dotfill\pageref{comparetop}\\
--- properties\dotfill\pageref{comparetop}\\
topological fields (non-discrete)\dotfill\pageref{convents}\\
topological space (Hausdorff)\dotfill\pageref{convents}\\
topological vector space\\
--- canonical Hausdorff vector topology\hfill\pageref{canonic}\\
--- locally convex\dotfill\pageref{deflcx}\\
--- Mackey complete\dotfill\pageref{defnMC}\\
--- patched\dotfill\pageref{defnpatched}\\
--- polynormed\dotfill\pageref{findimgp}\\
--- sequentially complete\dotfill\pageref{defnseqc}\\
typical fibre\dotfill\pageref{defnbdle}\\
ultrametric field\dotfill\pageref{convents}\\
ultrametric seminorm\dotfill\pageref{deflcx}\\
uniformizing element $\pi$\dotfill\pageref{defnuing}\\
unit group\dotfill\pageref{defnci}\\
--- as a Lie group\dotfill\pageref{ciasmooth}\\
valuation ring $\bO$\dotfill\pageref{convents}\\
vector bundle\dotfill\pageref{defnbdle}\\
weak direct products of Lie groups\dotfill\pageref{propweakdp}\\
weakly $C^k$\dotfill\pageref{thmgroth}\\
weakly smooth\dotfill\pageref{thmgroth}\\
Whitney sum\dotfill\pageref{sumsvble}\pagebreak

\noindent
{\Large\bf List of Symbols}\\[2mm]\label{nowsymb}
\thispagestyle{plain}
$\;\;$\vspace{-4 mm}\\
{\bf Standard Symbols:}\vspace{1 mm}\\
$\C$ the field of complex numbers\\
$\N=\{1,2,3,\ldots\,\}$\\
$\N_0=\{0,1,2,\ldots\,\}$\\
$\R$ the field of real numbers\\
$\Q_p$ the field of $p$-adic numbers\\
$\sph^1=\{z\in\C\!: |z|=1\,\}$\\
$\Z$ the integers\vspace{1 mm}\\
{\bf Operations on Sets and Maps}\vspace{1 mm}\\
$Y^0$ interior of $Y\sub X$ in~$X$\\
$\wb{Y}$ closure of a subset $Y$ of a topological\\
   \spa space (unless re-defined)\\
$f|_Y$ restriction of~$f$ to~$Y$\\
$f|^Y$ corestriction of~$f$ to~$Y\supseteq \im(f)$\vspace{3 mm}\\
{\bf Special Symbols}\vspace{1 mm}\\
$|.|$ (absolute value)\dotfill\pageref{convents}\\
$\bO$ (valuation ring)\dotfill\pageref{convents}\\
$B_\ve^E(x)$, $B_\ve(x)$ (balls)\dotfill\pageref{dfblls}\\
$\|.\|$ (norm, operator norm)\dotfill\pageref{dfblls}\\
$\|.\|_\infty$ (maximum norm)\dotfill\pageref{dfblls}\\
$\cL(E,F)$, $\cL(E)$\dotfill\pageref{dfblls}\\
$U^{[k]}$\dotfill\pageref{domains}\\
$C^k$, $C^k_\K$\dotfill\pageref{defnCr}\\
$f^{[k]}$, $f^{[k]}_\K$\dotfill\pageref{defnCr}\\
$U^{]1[}$\dotfill\pageref{gargel}\\
$f^{]1[}$\dotfill\pageref{gargel}\\
$df(x,v)$, $d^kf(x,v_1,\ldots, v_k)$,
$D_vf$\dotfill\pageref{differentials}\\
$C^k_{MB}$\dotfill\pageref{MichalB}\\
$C^k$-manifold\dotfill\pageref{verybas}\\
$L(G)$\dotfill\pageref{verybas}\\
$A^\times$\dotfill\pageref{defnci}\\
$M_n(A)$\dotfill\pageref{matrices}\\
$F\tensor_\K A$\dotfill\pageref{tensorprod}\\
$A_\bL$\dotfill\pageref{extendscale}\\
$\bigcup_{n\in \N}A_n$\dotfill\pageref{refback}\\
$C_K(X,E)$, $C_K(X,U)$\dotfill\pageref{defncsp}\\
$f(\sbull,p)_*$\dotfill\pageref{backbone},
\pageref{globcruc}, \pageref{crucial}\\
$f_*$\dotfill\pageref{pushforw},
\pageref{pushforw2}, \pageref{pfwcsupp}\\
$C_K(X,f)$\dotfill\pageref{Cf}\pagebreak

$\;$\\[5.5mm]
$C^r(U,E)$\dotfill\pageref{top1}\\
$\lfloor K, W\rfloor$\dotfill\pageref{thesets}\\
$C^r(f,E)$\dotfill\pageref{pullback}, \pageref{pb2}\\
$m_f$ (multiplication operator)\dotfill\pageref{multop}\\
$C^r(M,E)$\dotfill\pageref{top2}\\
$\theta_\kappa$\dotfill\pageref{dag}\\
$C^r(M,\lambda)$\dotfill\pageref{linearcase}\\
$C^r_K(M,E)$\dotfill\pageref{defcrsk}\\
$C^r_K(M,f)$\dotfill\pageref{Cf2}\\
$\lfloor K, U\rfloor_r$\dotfill\pageref{flooropen}\\
$C^r_K(M,G)$ (mapping group)\dotfill\pageref{defnmgp}\\
$C^r_K(M,A)$ (mapping algebra)\dotfill\pageref{mapcia}\\
$\bigoplus_{i\in I}E_i$\dotfill\pageref{basicdsum}\\
$\bigoplus_{i\in I}U_i$\dotfill\pageref{defnbxxx},
\pageref{mapsdirsums}\\
$U^{\{k\}}$, $f^{\{k\}}$\dotfill\pageref{verystrr}\\
$\oplus_{i\in I} f_i$\dotfill\pageref{mapsdirsums}\\
$\prod_{i\in I}^*G_i$ (weak direct product)\dotfill\pageref{propweakdp}\\
$X=\coprod_{i\in I} X_i$\dotfill\pageref{paracamen}\\
$\pi$ (uniformizing element)\dotfill\pageref{defnuing}\\
$B=\bO^d$\dotfill\pageref{defnuing}\\
$C^r_c(M,E)$\dotfill\pageref{defcpsupp}\\
$\cK(M)$\dotfill\pageref{defcpsupp}\\
$C^r_c(M,E)_\tvs$, $C^r_c(M,E)_\lcx$,
$C^r_c(M,E)_\bx$\hfill\pageref{deftvstp}\\
$C^r_c(M,U)$\dotfill\pageref{pfwcsupp}\\
$C^r_c(M,f)$\dotfill\pageref{functcsupp}\\
$C^r_c(M,A)$ (test function algebra)\dotfill\pageref{disctefgp}\\
$C^r_c(M,G)$ (test function group)\dotfill\pageref{disctefgp}\\
$\ve$ (evaluation)\dotfill\pageref{evalCk}, \pageref{seclook}\\
$\Gamma$ (composition map)\dotfill\pageref{compcomp}\\
$f^\vee$\dotfill\pageref{halfcartesian}\\
$\Phi$\dotfill\pageref{halfcartesian}\\
$g^\wedge$\dotfill\pageref{nownow}\\
$c^\infty$-open\dotfill\pageref{bscsconv}\\
$c^\infty_\K$-map\dotfill\pageref{bscsconv}\\
$\Diff^\infty(B)$\dotfill\pageref{diffball}\\
$\End^\infty(B)$\dotfill\pageref{diffball}\\
$\Theta_\psi$\dotfill\pageref{sltly}\\
$\Diff^\infty(M)$\dotfill\pageref{thmdiffpara}\\
$C^\infty_c(M,TM)\wt{\;}$\dotfill\pageref{defspatil}\\
$\End^r_c(U)$, $\cE^r_c(U)$\dotfill\pageref{defendom}\pagebreak

\noindent
$\beta_r$\dotfill\pageref{defbetr}\\
$C^\infty_c(U,\K^d)\wt{\;}$,
$\cE^\infty_c(U)\wt{\;}$,
$\End_c^\infty(U)\wt{\;}$\dotfill\pageref{wildetilde}\\
$m_{r,k}$, $\wt{m}$, $\rho_{r,\eta}$\dotfill\pageref{cpendsmooth}\\
$\Diff^r_c(U)$, $\Diff^\infty_c(U)\wt{\;}$\dotfill\pageref{cpendsmooth}\\
$C^r_c(\phi,E)$\dotfill\pageref{pbproper}\\
$\Diff^r(U)$\dotfill\pageref{DiffU}\\
$\Diff^r(M)$\dotfill\pageref{defnDM}\\
$C^r(M,E)_D$\dotfill\pageref{defnDjs}\\
$D^j$\dotfill\pageref{defnDjs}\\
$T^jM$\dotfill\pageref{defnDjs}\\
\thispagestyle{plain}
$g_{\phi,\psi}$, $G_{\phi,\psi}$\dotfill\pageref{chnge}\\
$\Supp(\sigma)$\dotfill\pageref{bscssec}\\
$C^r_K(M,E)$, $C^r(M,E)\!$
(spaces of sections)\dotfill\pageref{bscssec}\\
$0_\sbull$ (zero-section)\dotfill\pageref{bscssec}\\
$\sigma_\psi$\dotfill\pageref{sigsi}\\
$f_{\phi,\psi}$\dotfill\pageref{defnparvect}\\
$C^r(M,F)$\dotfill\pageref{pushfbdl}\\
$C^r(M,E)$ as a topological module\dotfill\pageref{aretopmod}\\
$E_1\oplus E_2$\dotfill\pageref{sumsvble}\\
$C^r_K(M,E)$\dotfill\pageref{boredea}\\
$C^r_c(M,E)_\tvs$, $C^r_c(M,E)_\lcx$,
$C^r_c(M,E)_\bx$\hfill\dotfill\pageref{defcpsupp2}\\
$C^r_c(M,\Omega)$\dotfill\pageref{laomop}\\
$C^r_c(M,E)$ as a topological module\dotfill\pageref{aretopmod2}
\end{document}